\documentclass[11pt]{amsart}

\usepackage{array}
\usepackage[square]{stmaryrd}
\usepackage[mathscr]{eucal}
\usepackage{amssymb}
\usepackage[centertags]{amsmath}
\usepackage{graphics,graphicx}
\usepackage{bm}
\usepackage{subfig}
\usepackage[small,nohug,heads=vee]{diagrams} 
\usepackage{longtable}

\diagramstyle[labelstyle=\scriptstyle]

\usepackage[margin=2.9cm]{geometry}

\numberwithin{equation}{subsection}
\numberwithin{figure}{subsection}

\newtheorem{theorem2}{Theorem}[subsection]{\bf}{\rm}
\newtheorem{proposition2}[theorem2]{Proposition}{\bf}{\rm}
\newtheorem{lemma2}[theorem2]{Lemma}{\bf}{\rm}
\newtheorem{corollary2}[theorem2]{Corollary}{\bf}{\rm}
\newtheorem{definition2}[theorem2]{Definition}{\bf}{\rm}
\newtheorem{remark2}[theorem2]{Remark}{\it}{\rm}
\newtheorem{example2}[theorem2]{Example}{\bf}{\rm}

\newtheorem{assumption}[theorem2]{Assumption}{\bf}{\rm}
\newtheorem{condition}[theorem2]{Condition}{\bf}{\rm}

\newcommand{\R}{\mathbb{R}}
\newcommand{\Z}{\mathbb{Z}}
\newcommand{\C}{\mathbb{C}}

\newcommand{\CP}[1]{\mathbb{CP}^{#1}}
\newcommand{\RP}[1]{\mathbb{RP}^{#1}}
\newcommand{\To}{\rightarrow}
\newcommand{\del}[2]{\frac{\partial #1}{\partial #2}}

\newcommand{\url}[1]{}

\newarrow{Eq}{=}{=}{=}{=}{=}

\begin{document}

\title{Homological mirror symmetry for Calabi-Yau hypersurfaces in projective space}
\thanks{This research was mostly carried out while the author was a graduate student at MIT. 
This research was partially supported by NSF Grant DMS-0652620. }

\begin{abstract}
We prove Homological Mirror Symmetry for a smooth $d$-dimen\-sional Calabi-Yau hypersurface in projective space, for any $d \ge 3$ (for example, $d=3$ is the quintic three-fold). 
The main techniques involved in the proof are: the construction of an immersed Lagrangian sphere in the `$d$-dimensional pair of pants'; the introduction of the `relative Fukaya category', and an understanding of its grading structure; a description of the behaviour of this category with respect to branched covers (via an `orbifold' Fukaya category); a Morse-Bott model for the relative Fukaya category that allows one to make explicit computations; and the introduction of certain graded categories of matrix factorizations mirror to the relative Fukaya category.
\end{abstract}

\author{Nick Sheridan}

\address{Princeton University and the Institute for Advanced Study, Einstein Drive, Princeton, NJ 08540, United States}             \email{nicks@math.ias.edu} 

\maketitle

\tableofcontents

\section{Introduction}
\label{sec:intro}

\subsection{Mirror symmetry context}

The mirror symmetry phenomenon was discovered by string theorists. 
In its original form, it dealt with Calabi-Yau K\"{a}hler manifolds. 
On such a manifold, one can define symplectic invariants (the `$A$-model') and complex invariants (the `$B$-model').   
Broadly, mirror symmetry says that there exist pairs of manifolds $(M,N)$ such that the $A$-model on $M$ is equivalent to the $B$-model on $N$, and vice versa.

The closed-string version of the $A$-model is encoded in Gromov-Witten invariants (the quantum cohomology ring), and the closed-string version of the $B$-model is encoded in period integrals of the holomorphic volume form. 
Mirror symmetry first came to the attention of mathematicians in 1991, when Candelas, de la Ossa, Green and Parkes \cite{Candelas1991} applied it to make predictions about rational curve counts on the quintic three-fold $M$. 
They constructed a manifold $N$ which ought to be mirror to $M$, and computed the closed-string $B$-model on $N$. 
Assuming mirror symmetry to hold, this allowed them to predict the closed-string $A$-model on $M$, from which they extracted rational curve counts.
This closed-string version of mirror symmetry was proven in 1996 for Calabi-Yau complete intersections in toric varieties, by Givental  \cite{Givental1996} and Lian-Liu-Yau \cite{Lian1997,Lian1999,Lian1999a,Lian2000}.

In the meantime, Kontsevich had introduced the open-string version of mirror symmetry, called the Homological Mirror Symmetry conjecture \cite{Kontsevich1994} in 1994.
He  proposed that the $A$-model should be encoded in the Fukaya category, and the $B$-model should be encoded in the category of coherent sheaves. 
Then if $M$ and $N$ are mirror Calabi-Yaus, there should be an equivalence of the Fukaya category of $M$ with the category of coherent sheaves on $Y$ (on the derived level).
Complete or partial proofs of homological mirror symmetry for closed Calabi-Yau varieties are known for elliptic curves \cite{Polishchuk1998,Polishchuk2000}, abelian varieties \cite{Fukaya2002a} (see \cite{Abouzaid2010d} for the case of the four-torus), Strominger-Yau-Zaslow dual torus fibrations \cite{Kontsevich2001}, and the quartic K3 surface \cite{Seidel2003}. 

In this paper we consider a smooth Calabi-Yau hypersurface $M^n \subset \CP{n-1}$, with symplectic form in the cohomology class $[\omega] = n^2 c_1(\mathcal{O}(1))$, and its (predicted) mirror $N^n$. 

\begin{remark2}
\label{remark2:awkward}
Note that we define $M^n$ to be $(n-2)$-dimensional as a complex manifold, {\bf not} $n$-dimensional. 
We reassure the reader that this notational awkwardness is worth it, in that it makes almost all of our subsequent computations more readable: if we adopted the more intuitive convention that $M^n$ is an $n$-dimensional Calabi-Yau hypersurface, then all of our equations would have an ugly `$n+2$' in them.
\end{remark2}

Our main result is to prove one direction of homological mirror symmetry: we prove that there is a quasi-equivalence of triangulated $A_{\infty}$ categories between the split-closed derived Fukaya category of $M^n$ and (a suitable DG enhancement of) the bounded derived category of coherent sheaves on $N^n$ (see Theorem \ref{theorem2:2main} for the precise statement). 
$M^3$ is the elliptic curve considered by Polishchuk and Zaslow, $M^4$ is the quartic surface considered by Seidel, and $M^5$ is the quintic three-fold. 
We remark that Nohara and Ueda have also considered the case of the quintic three-fold \cite{Nohara2012}, using the results of Seidel \cite{Seidel2003} and our earlier paper \cite{Sheridan2011a}.

In future work, we will extend this result to the case of Fano hypersurfaces in projective space, and we plan also to consider the case of hypersurfaces of general type. 

\subsection{Statement of the main result}
\label{subsec:2maintheorem}

First let us introduce the coefficient field of the categories we will be considering. 

\begin{definition2}
\label{definition2:2novikov}
We define the {\bf universal Novikov field} $\Lambda$, whose elements are formal sums
\[ \psi(r) = \sum_{j=0}^{\infty} c_j r^{\lambda_j},\]
where $c_j \in \C$, and $\lambda_j \in \R$ is an increasing sequence of real numbers such that
\[ \lim_{j \To \infty} \lambda_j = +\infty.\]
It is a non-Archimedean valued field, with norm
\[ \left\| \sum_{j=0}^{\infty} c_j r^{\lambda_j} \right\| := e^{-\lambda_0}\]
(assuming $c_0 \neq 0$). 
It is algebraically closed (see \cite[Appendix 1]{Fukaya2010c}).
\end{definition2}

Observe that there is an inclusion $\C \llbracket r \rrbracket \subset \Lambda$. 
Thus, if $\mathscr{C}$ is a $\C \llbracket r \rrbracket$-linear category, then we can form the $\Lambda$-linear category $\mathscr{C} \otimes_{\C \llbracket r \rrbracket} \Lambda$.

Now we need to consider ($r$-adically) continuous automorphisms of $\Lambda$.
Note that automorphisms of $\C \llbracket r \rrbracket$ are in one-to-one correspondence with elements $\psi \in \C \llbracket r \rrbracket$ such that $\psi(0) \neq 0$: the corresponding automorphism 
\[ \psi^*: \C \llbracket r \rrbracket \To \C \llbracket r \rrbracket\]
is defined by its value on $r$:
\[ \psi^*(r) := r \psi(r).\]

\begin{lemma2}
Any continuous automorphism $\psi^*$ of $\C \llbracket r \rrbracket$ can be lifted (non-uniquely) to a continuous automorphism
\[ \hat{\psi}^*: \Lambda \To \Lambda.\]
\end{lemma2}
\begin{proof}
Consider the abelian group $A \subset (\Lambda,\cdot)$ of elements of norm $1$ (i.e., the elements whose leading-order term is a non-zero multiple of $r^0$).  
It follows easily from the proof of \cite[Lemma 14.1]{Fukaya2010c} that $A$ is divisible. 
Therefore, the homomorphism
\[ \psi: (\Z,+) \to A\]
which sends $1 \mapsto \psi(r)$, extends to a homomorphism
\[ \hat{\psi}: (\R,+) \To A\]
by Baer's criterion. 
We can now define the automorphism $\hat{\psi}^*$ by
\[ \hat{\psi}^*\left(r^\alpha\right) = r^\alpha \hat{\psi}(\alpha).\]
It is clear that $\hat{\psi}^*$ is an $r$-adically continuous automorphism lifting $\psi^*$.
 \end{proof}

\begin{definition2}
\label{definition2:formchang}
Suppose that $R$ is a ring, and $\mathscr{C}$ is an $R$-linear category, and
\[ \varphi: R \To R\]
an automorphism. 
We define a new $R$-linear category $\varphi \cdot  \mathscr{C}$: it is the same category as $\mathscr{C}$, but the $R$-module structure of the morphism spaces is changed via $\varphi$. 
We say that $\varphi \cdot \mathscr{C}$ is obtained from $\mathscr{C}$ via a {\bf change of variables} in $R$.
\end{definition2}

\begin{remark2}
\label{remark2:tenschange}
Suppose that $\mathscr{C}$ is a $\C \llbracket r \rrbracket$-linear category, $\psi^*$ an automorphism of $\C \llbracket r \rrbracket$, and $\hat{\psi}^*$ a lift of $\psi^*$ to an automorphism of $\Lambda$. 
Then 
\[ \left( \psi^* \cdot \mathscr{C}  \right) \otimes_{\C \llbracket r \rrbracket} \Lambda \cong \hat{\psi}^* \cdot \left( \mathscr{C} \otimes_{\C \llbracket r \rrbracket } \Lambda \right).\]
\end{remark2}

Now we introduce the relevant categories. 

\begin{definition2}
\label{definition2:2fermdef}
On the symplectic side, let $M^n \subset \CP{n-1}$ be a smooth hypersurface of degree $n$. 
$M^n$ is an $(n-2)$-dimensional Calabi-Yau.
The Fukaya category $\mathcal{F}(M^n)$ (as defined in \cite{Abouzaid2012}) is a $\Z$-graded $\Lambda$-linear $A_{\infty}$ category.
The split-closed derived Fukaya category $D^{\pi}\mathcal{F}(M^n)$ (see \cite[Section I.4]{Seidel2008} for the definition of the split-closed derived category of an $A_{\infty}$ category) is a split-closed $\Lambda$-linear triangulated $A_{\infty}$ category. 
We remark that the Fukaya category is a symplectic invariant (up to $A_{\infty}$ quasi-isomorphism), so it does not matter which smooth hypersurface $M^n$ we choose (they are all symplectomorphic by Moser's theorem).
\end{definition2}

\begin{definition2}
\label{definition2:2mirror}
On the algebraic side, we define
\[ w_{nov} := u_1 \ldots u_{n} + r \sum_{j=1}^{n} u_j^{n} \in \Lambda[u_1, \ldots, u_n].\]
We set $\widetilde{N}^n_{nov} := \{w_{nov} = 0\} \subset \mathbb{P}_{\Lambda}^{n-1}$. 
We equip $\widetilde{N}^n_{nov}$ with the action of a finite group. 
Observe that the group
\[ \tilde{\Gamma}^*_n :=  (\Z_n)^{n}/(1,1,\ldots,1)\]
(where we denote $\Z_n := \Z/(n)$)
acts on $\mathbb{P}_{\Lambda}^{n-1}$ by multiplying the homogeneous coordinates $u_j$ by $n$th roots of unity (we include the `$*$' for consistency with later notation).
We have a homomorphism
\[ \tilde{\Gamma}^*_n \To \Z_n,\]
given by summing the coordinates. 
We denote the kernel of this homomorphism by $\Gamma^*_n$.
Note that the action of $\Gamma^*_n$  preserves $\widetilde{N}^n_{nov}$, so $\Gamma^*_n$ acts on $\widetilde{N}^n_{nov}$. 
We define
\[ N^n_{nov} := \widetilde{N}^n_{nov}/\Gamma^*_n.\]
We consider a suitable DG enhancement of the bounded derived category of coherent sheaves on $N^n_{nov}$, which is
\[ D^bCoh(N^n_{nov}) \cong D^b\left(Coh^{\Gamma^*_n}\left( \widetilde{N}^n_{nov} \right)\right).\]
It is another triangulated $A_{\infty}$ category over $\Lambda$.
\end{definition2}

\begin{theorem2}
\label{theorem2:2main}
If $n \ge 5$, then there exists a power series $\psi \in \C\llbracket r \rrbracket $, $\psi(0) = \pm 1$, such that, for any lift $\hat{\psi}^*$ of the corresponding automorphism $\psi^*$ of $\C \llbracket r \rrbracket$ to an automorphism of $\Lambda$, there is a quasi-equivalence of $\Lambda$-linear triangulated $A_{\infty}$ categories
\[ D^{\pi} \mathcal{F}(M^n) \cong \hat{\psi}^* \cdot D^bCoh(N^n_{nov}).\]
In fact, we will see that $\psi \in \C\llbracket r^n \rrbracket  \subset \C\llbracket r \rrbracket $.
\end{theorem2}

\begin{remark2}
Observe that, by taking cohomology, Theorem \ref{theorem2:2main} yields an equivalence of ordinary (not $A_{\infty}$) triangulated categories.
\end{remark2}

\begin{remark2}
\label{remark2:neq12}
The requirement that $n \ge 5$ can be removed without difficulty, but would require some ad-hoc arguments in the cases $n=3$ and $n=4$ (see Remark \ref{remark2:nnot4}). 
Since the corresponding result in these cases is already known ($n=3$ is the elliptic curve and $n=4$ is the quartic surface), and it would distract from the main flow of the argument, we leave these cases out. 
\end{remark2}

We remark that the reference \cite{Abouzaid2012}, in which the full Fukaya category $\mathcal{F}(M^n)$ is defined, and one of the results we use in the proof of Theorem \ref{theorem2:2main} is proven (we have stated it as Theorem \ref{theorem2:2afooo}), is still in preparation. 
This is not an ideal situation, but we reassure the reader that this paper is written with minimal reliance on \cite{Abouzaid2012}. 
Let us explain what we mean.

\begin{sloppypar}
We rigorously define a version of the Fukaya category called the `relative' Fukaya category $\mathcal{F}(M^n,D)$. 
It is an $A_{\infty}$ category defined over the coefficient ring $R$, which is the completion of the polynomial ring
\[ \widetilde{R} := \C [r_1, \ldots, r_n]\]
in a certain category of graded algebras. 
In fact, it will turn out that
\[ R \cong \C \llbracket T  \rrbracket [r_1, \ldots, r_n]/(T - r_1 \ldots r_n).\]
We also define a certain category of matrix factorizations, denoted by $\bm{p}_1^*MF^{\bm{G}}(S,w)$, which is a DG category over the same coefficient ring $R$.
We introduce full subcategories
\[ \widetilde{\mathscr{A}} \subset \mathcal{F}(M^n,D),\,\,\,\,\, \widetilde{\mathscr{B}} \subset \bm{p}_1^*MF^{\bm{G}}(S,w),\]
and prove (without reference to \cite{Abouzaid2012}) that there exists a power series $\tilde{\psi} \in \C \llbracket T \rrbracket$, $\tilde{\psi}(0) = \pm 1$, giving rise to an automorphism $\tilde{\psi}^*$ of $R$, with
\[ \tilde{\psi}^*(r_j) = \tilde{\psi}(T) r_j,\]
such that there is an $A_{\infty}$ quasi-isomorphism
\[ \widetilde{\mathscr{A}} \cong \tilde{\psi}^* \cdot \widetilde{\mathscr{B}} .\]
For a more precise statement, see Theorem \ref{theorem2:2main2}. 
\end{sloppypar}

Now there are homomorphisms of $\C$-algebras
\[ R \To \C \llbracket r \rrbracket \hookrightarrow \Lambda,\]
where the left homomorphism maps all $r_j$ to $r$. 
It follows from a theorem of Orlov that there is an $A_\infty$ quasi-equivalence
\[  D^bCoh(N^n_{nov}) \cong  \bm{p}_1^*MF^{\bm{G}}(S,w) \otimes_R \Lambda, \]
and the full subcategory $\widetilde{\mathscr{B}} \otimes_R \Lambda$ generates. 

Observe that $\tilde{\psi}^*$ lifts to the automorphism $\psi^*$ of $C \llbracket r \rrbracket$, where $\psi(r) = \tilde{\psi}(r^n)$, and $\psi^*$ lifts to an automorphism of $\Lambda$.
It follows that, if $\hat{\psi}^*$ is such a lift, then we have a fully faithful $A_{\infty}$ embedding 
\[ \hat{\psi}^* \cdot D^bCoh(N^n_{nov}) \hookrightarrow D^b \left(\mathcal{F}(M^n,D) \otimes_R \Lambda \right)\]
(see Remark \ref{remark2:tenschange}).
We can come to this point without reference to \cite{Abouzaid2012}, but no further.

In Section \ref{sec:2splitgen}, we make the assumption (Assumption \ref{assumption:relfull}) that there is a fully faithful embedding
\[ \mathcal{F}(M^n,D) \otimes_R \Lambda \hookrightarrow \mathcal{F}(M^n),\]
where $\mathcal{F}(M^n)$ is defined as in \cite{Abouzaid2012}.
We justify this assumption `heuristically', and show that the full subcategory
\[ \widetilde{\mathscr{A}} \otimes_R \Lambda \subset  \mathcal{F}(M^n)\]
split-generates, by applying the split-generation criterion of \cite{Abouzaid2012}, which we state as Theorem \ref{theorem2:2afooo}.
This completes the proof of Theorem \ref{theorem2:2main}.

In the rest of this introduction, we give an overview of the main techniques introduced in the rest of the paper, and how they are used to prove Theorem \ref{theorem2:2main}. 
We will make a few  imprecise statements and definitions, in the interest of giving the reader the correct intuitive picture.

\subsection{Affine and relative Fukaya categories}
\label{subsec:2eqaffint}

Let $M$ be a K\"{a}hler manifold, with a smooth normal-crossings divisor $D = D_1 \cup \ldots \cup D_k \subset M$, where $D_j$ are the smooth irreducible components of $D$. 
We require that each $D_j$ be Poincar\'{e} dual to some positive multiple of the K\"{a}hler form $\omega$.
In this paper, we consider three versions of the Fukaya category: the {\bf affine Fukaya category} $\mathcal{F}(M \setminus D)$, the {\bf relative Fukaya category} $\mathcal{F}(M,D)$, and the {\bf full Fukaya category} $\mathcal{F}(M)$.

First let us describe the affine Fukaya category. 
It is closely related to the exact Fukaya category of the exact symplectic manifold $M \setminus D$, as defined by Seidel in \cite{Seidel2008}.

First we explain the objects of the affine Fukaya category.
If $P$ is a symplectic manifold, we denote by $\mathcal{G}P \To P$ the Lagrangian Grassmannian of $P$, whose fibre over $p \in P$ is the space of Lagrangian subspaces of $T_pP$. 
Any smooth Lagrangian immersion $i:L \To P$ comes with an associated lift $i_*: L \To \mathcal{G}P$.
We define an {\bf anchored Lag\-rangi\-an brane} $L^{\#}$ in $P$ to be a Lag\-rang\-i\-an immersion $i:L \To P$, together with a lift
\[ \begin{diagram}
&& \widetilde{\mathcal{G}}P \\
& \ruTo^{i^{\#}}& \dTo \\
L & \rTo_{i_*} & \mathcal{G}P,
\end{diagram}
\]
where $\widetilde{\mathcal{G}}P$ is the universal abelian cover of the total space of $\mathcal{G}P$ (i.e., the one corresponding to the commutator subgroup of $\pi_1(\mathcal{G}P)$), together with a choice of Pin structure on $L$. 
The `anchored' terminology comes from \cite{Fukaya2010a}, where a closely-related concept is studied.

We define the objects of the affine Fukaya category $\mathcal{F}(M \setminus D)$ to be compact, exact, embedded, anchored Lagrangian branes in $M \setminus D$.
The category is $\C$-linear, and morphism spaces and $A_{\infty}$ structure maps are defined in exactly the same way as in the exact Fukaya category of the exact symplectic manifold $M \setminus D$. 
Namely, the morphism spaces are
\[ CF_{\mathcal{F}(M \setminus D)}^*\left(L_0,L_1 \right) := \bigoplus_{p \in L_0 \cap L_1} \C \cdot p\]
(with appropriate modifications to allow for non-transverse intersection of $L_0$ with $L_1$). 
The $A_{\infty}$ structure map
\[ \mu^s: CF^*\left(L_{s-1},L_s\right) \otimes \ldots \otimes CF^* \left( L_0,L_1\right) \To CF^* \left(L_0,L_s \right)\]
is defined by counting rigid holomorphic disks with $s+1$ boundary punctures in $M \setminus D$. 
Namely, the coefficient of $p_0$ in $\mu^s(p_s, \ldots, p_1)$ is given by the signed count of rigid holomorphic disks
\[ u: \mathbb{D} \setminus \{ \zeta_0, \ldots, \zeta_s \} \To M \setminus D,\]
sending the $j$th boundary component to Lagrangian $L_j$, and asymptotic at puncture $\zeta_j$ to intersection point $p_j$.
However, we treat the grading differently. 

To start with, we equip each morphism space with a grading in the abelian group $H_1(M \setminus D)$. 
The $A_{\infty}$ structure maps respect this grading, essentially because if there is a holomorphic disk contributing to some $A_{\infty}$ product in $M \setminus D$, then its boundary is nullhomologous in $M \setminus D$.
Gradings of this type have appeared in, for example, \cite{Fukaya2010a,Seidel2008a,Sheridan2011a}.

\begin{remark2}
Another way of defining this grading (following \cite{Fukaya2010a}) would be to equip $M$ with a basepoint $q \in M \setminus D$, and define an object of the Fukaya category to be a Lagrangian $L \subset M \setminus D$, equipped with a path from $q$ to a point $q_L \in L$, inside $M \setminus D$.
Then, given an intersection point $p \in L_0 \cap L_1$, which by definition is a generator of $CF^*(L_0,L_1)$, we can define a class in $H_1(M \setminus D)$: start from $q$ and follow the path from $q$ to $q_{L_0}$, then follow a path in $L_0$ to the intersection point $p$, then follow a path in $L_1$ to $q_{L_1}$, then follow the path back to $q$. 
The class of this path in $H_1(M \setminus D)$ defines the grading of $p$. 
\end{remark2}

Secondly, if we were given a quadratic complex volume form on $M \setminus D$, we could define a $\Z$-grading of the morphism spaces (as in \cite{Seidel2008}).
This defines a $\Z \oplus H_1(M \setminus D)$-graded category. 
However, this formulation is unsatisfactory: the $\Z$-grading and $H_1$-grading are related. 
For example, changing the quadratic volume form has the effect of changing the $\Z \oplus H_1$ grading by an automorphism preserving the $H_1$ factor. 
We really want a new notion of grading.

In Section \ref{sec:2deftheory}, we define a {\bf grading datum} $\bm{G}$ to be an abelian group $Y$ together with a morphism $f: \Z \To Y$. 
We say that an $A_{\infty}$ category is {\bf $\bm{G}$-graded} if its morphism spaces are $Y$-graded, and the $A_{\infty}$ structure map $\mu^s$ has degree $f(2-s)$. 
If $Y = \Z$ and $f$ is the identity, then this coincides with the usual definition of a $\Z$-graded $A_{\infty}$ category.
We also study the deformation theory of $\bm{G}$-graded $A_{\infty}$ algebras and categories, and prove various classification results about them.

In Section \ref{subsec:2hmindex}, we introduce a grading datum $\bm{G}(M,D)$ associated to $M \setminus D$, as follows: we consider the fibre bundle $\mathcal{G}(M \setminus D)$, with the associated fibration
\[ \begin{diagram}
\mathcal{G}_p(M \setminus D) && \rInto && \mathcal{G}(M \setminus D) \\
&& && \dTo \\
&& && M \setminus D.
\end{diagram}
\]
Taking the abelianization of the associated exact sequence of homotopy groups, we obtain an exact sequence
\[ H_1 \left(\mathcal{G}_p(M \setminus D) \right) \To H_1 \left(\mathcal{G}(M \setminus D)\right)  \To H_1(M \setminus D).\]
We observe that $\mathcal{G}_p(M \setminus D)$ is the Lagrangian Grassmannian of the symplectic vector space $T_pM$, so the Maslov index defines an isomorphism $H_1(\mathcal{G}_p(M \setminus D)) \cong \Z$ (see \cite{Arnold1967}).
Thus we can define the grading datum $\bm{G}(M,D)$ to be given by the first morphism in this exact sequence. 
We show that the affine Fukaya category is naturally $\bm{G}(M,D)$-graded.
Observe that the second map in the exact sequence gives rise to the $H_1(M \setminus D)$-grading mentioned earlier.
 
\begin{remark2}
\label{remark2:fukvers}
For any covering of $\mathcal{G}(M \setminus D)$, we could similarly define a version of $\mathcal{F}(M \setminus D)$ with grading in the covering group. 
For example, in \cite{Seidel2008}, the cover is taken to be a fibrewise universal cover, and the corresponding Fukaya category is $\Z$-graded. 
Depending on which cover we choose, the resulting Fukaya category will have different objects (it can only deal with Lagrangians that lift to that cover), and different grading. 
In this paper, whenever we write $\mathcal{F}(M \setminus D)$ (or $\mathcal{F}(M,D)$), we are referring to the version corresponding to the universal abelian cover. 
\end{remark2}

Next, we introduce the relative Fukaya category, denoted $\mathcal{F}(M , D)$ (following \cite{Seidel2002,Seidel2003}).
Its objects are exactly the same as those of $\mathcal{F}(M \setminus D)$: compact, exact, embedded, anchored Lagrangian branes in $M \setminus D$. 
It is defined over a certain coefficient ring, which we denote by $R(M,D)$. 
To define it, we first introduce the polynomial ring 
\[ \widetilde{R}(M, D) := \C [ r_1, \ldots, r_k ],\]
with one generator for each irreducible component of $D$, and equip it with a $\bm{G}(M,D)$-grading. 
We define $R(M,D)$ to be the completion of $\widetilde{R}(M,D)$ with respect to the degree filtration, in the category of $\bm{G}(M,D)$-graded algebras. 
This will in general be a subring of the completion of $\widetilde{R}(M,D)$ in the category of algebras, which is of course the power series ring $\C \llbracket r_1, \ldots, r_k \rrbracket$. 
We will often write $R$ instead of $R(M,D)$ when no confusion is possible.

We define morphism spaces by 
\[ CF^*_{\mathcal{F}(M,D)}(L_0,L_1) := \bigoplus_{p \in L_0 \cap L_1} R\cdot p.\]
The $A_{\infty}$ structure maps $\mu^s$ count rigid boundary-punctured holomorphic disks in $M$. 
Namely, the coefficient of $p_0$ in $\mu^s(p_s, \ldots, p_1)$ is given by a signed count of rigid holomorphic disks
\[ u: \mathbb{D} \setminus \{ \zeta_0, \ldots, \zeta_s \} \To M\]
with boundary and asymptotic conditions as before. 
Each such disk $u$ contributes a term
\[ r^{u \cdot D} := r_1^{u \cdot D_1} \ldots r_k^{u \cdot D_k} \in R,\]
where $u \cdot D_j$ denotes the topological intersection number of $u$ with $D_j$. 
We observe that $u \cdot D_j \ge 0$ by positivity of intersection, so these coefficients lie in $\widetilde{R}$. 
It remains to show that the sum of these coefficients over all rigid disks is an element of $R$; we will explain this presently.

First we must explain that $\mathcal{F}(M,D)$ is still a $\bm{G}(M,D)$-graded category, but it is important that the coefficient ring $R$ has a non-trivial grading for this to be true. 
For example, consider the $H_1(M\setminus D)$-grading that we mentioned earlier. 
It is no longer true that a holomorphic disk $u$ contributing to an $A_{\infty}$ product $\mu^s$ has nullhomologous boundary in $M \setminus D$, because the disk now maps into $M$, not $M \setminus D$ as it did for the affine Fukaya category. 
However, if we remove small balls surrounding each intersection point of $u$ with a divisor $D_j$, then the resulting surface defines a homology in $M \setminus D$ between the boundary of $u$ and a collection of meridian loops around the divisors. 
Thus, if we define the $H_1(M \setminus D)$-grading of the generator $r_j \in R$ to be the class of a meridian loop around divisor $D_j$, then the $A_{\infty}$ structure respects the $H_1$-grading. 
Of course there remains more to check to show that $\mathcal{F}(M,D)$ is $\bm{G}(M,D)$-graded, but this is the basic idea -- the details can be found in Section \ref{sec:2relfuks}.

This grading is what allows us to define the coefficient ring $R(M,D)$ by taking the completion of the polynomial ring $\widetilde{R}(M,D)$ in the category of $\bm{G}(M,D)$-graded algebras, rather than in the category of algebras. 
The coefficient of $p_0$ in $\mu^s(p_s, \ldots, p_1)$ is a sum over rigid disks $u$ of coefficients $r^{u \cdot D}$.
The number of disks contributing a term $r^{\bm{d}} p_0$ to this product, for fixed $\bm{d}$, is finite by Gromov compactness, but there may be infinitely many monomials $r^{\bm{d}}$ for which the coefficient of $r^{\bm{d}}p_0$ in this product is non-zero. 
That is why we need to take a completion in the definition of our coefficient ring $R(M,D)$. 
However, because the maps $\mu^s$ are $\bm{G}(M,D)$-graded, we only ever need an infinite sum of terms $r^{\bm{d}}$ with the same grading, hence it suffices for the coefficient ring to be the completion in the category of $\bm{G}(M,D)$-graded algebras.

Now consider the map $R \To \C$ which sends all $r_j$ to $0$. 
We have
\[ \mathcal{F}(M,D) \otimes_R \C \cong \mathcal{F}(M \setminus D).\]
That is because the category on the left is the zeroth-order part of $\mathcal{F}(M,D)$, in which a holomorphic disk $u: \mathbb{D} \To M$ contributes to $\mu^s$ if and only if $u \cdot D_j = 0$ for all $j$. 
This corresponds to counting only holomorphic disks which avoid the divisors $D$, i.e., which lie in $M \setminus D$.
By definition, this corresponds to the affine Fukaya category $\mathcal{F}(M \setminus D)$. 
We therefore say that $\mathcal{F}(M , D)$ is a $\bm{G}(M,D)$-{\bf graded deformation} of $\mathcal{F}(M \setminus D)$ over $R(M,D)$.

Finally, we recall the definition of the full Fukaya category $\mathcal{F}(M)$, as in \cite{fooo,Abouzaid2012}, where $M$ is a Calabi-Yau symplectic manifold. 
Its objects are Lagrangians $L \subset M$ equipped with a grading and spin structure. 
Note that, because $M$ is Calabi-Yau, there is a fibrewise universal cover of $\mathcal{G}M$, with covering group $\Z$. 
A grading of $L$ is a lift to this cover: this equips $\mathcal{F}(M)$ with a $\Z$-grading in accordance with Remark \ref{remark2:fukvers}.

$\mathcal{F}(M)$ is $\Lambda$-linear. 
Its morphism spaces are $\Z$-graded $\Lambda$-vector spaces (where $\Lambda$ has degree $0 \in \Z$), defined by
\[ CF^*_{\mathcal{F}(M)}(L_0,L_1) := \bigoplus_{p \in L_0 \cap L_1} \Lambda \cdot p.\]
The $A_{\infty}$ structure maps are defined by signed counts of rigid boundary-punctured holomorphic disks $u: \mathbb{D} \To M$ as before. 
Each such disk $u$ contributes a term $r^{\omega(u)} \in \Lambda$ to the corresponding coefficient of $\mu^s$. 

If we define the ring homomorphism
\begin{eqnarray*} 
R & \To & \Lambda \\
r_j & \mapsto & r \mbox{ for all $j$,}
\end{eqnarray*}
then $\Lambda$ becomes an $R$-algebra. 
We expect that there is a fully faithful embedding of $\Z$-graded $\Lambda$-linear $A_{\infty}$ categories,
\[ \mathcal{F}(M, D) \otimes_{R} \Lambda \hookrightarrow \mathcal{F}(M).\]
That is because the $A_{\infty}$ structure maps are counting the same holomorphic disks, and by Stokes' theorem each disk contributes with the same coefficient in both categories. 
However, the details of the definitions of the two sides are quite different, because the approach to transversality in the two versions of the Fukaya category is different, so we do not rigorously establish that such an embedding exists. 
We state this as Assumption \ref{assumption:relfull}, and give some justification in Remark \ref{remark2:relfull}.

The properties of the three versions of the Fukaya category are summarized in Table \ref{table:threefuk}, and the relationship between the three can be summarized as

\[ \mathcal{F}(M \setminus D) \xrightarrow{\mbox{$\bm{G}(M,D)$-graded deformation}} \mathcal{F}(M, D) \xrightarrow{ \otimes_{R} \Lambda} \mathcal{F}(M). \]

\renewcommand{\arraystretch}{1.5}
\begin{table}
\begin{tabular}{>{\raggedright}p{.145\textwidth}|>{\raggedright}p{.23\textwidth}|>{\raggedright}p{.23\textwidth}
|>{\raggedright\arraybackslash}p{.247\textwidth}}
Notation & Affine: $ \mathcal{F}(M \setminus D)$ & Relative: $\mathcal{F}(M,D)$ & Full:  $\mathcal{F}(M)$ \\ \hline
Coefficients & $\C$ & $R(M,D)$ & $\Lambda$ \\ \hline
Objects & Closed, anchored Lagrangian branes $L^{\#} \subset M \setminus D$ & Closed, anchored Lagrangian branes $L^{\#} \subset M \setminus D$ & Graded, spin Lagrangians $L \subset M$ \\ \hline
Morphisms & $\C \langle L_0 \cap L_1 \rangle$ & $R\langle L_0 \cap L_1 \rangle$  & $\Lambda \langle L_0 \cap L_1 \rangle$ \\ \hline
$A_{\infty}$ maps $\mu^s$ & $\#\{u: \mathbb{D} \To M \setminus D \mbox{ hol.}\}$ & $\#\{u: \mathbb{D} \To M\mbox{ hol.}\}$, with coefficient $r_1^{u \cdot D_1} \ldots r_k^{u \cdot D_k} \in R$ & $\#\{u: \mathbb{D} \To M\mbox{ hol.}\}$, with coefficient $r^{\omega(u)} \in \Lambda$ \\ \hline
Grading & $\bm{G}(M,D)$ & $\bm{G}(M,D)$ & $\Z$ (if $M$ is Calabi-Yau) 
\end{tabular}
\caption{Three versions of the Fukaya category. \label{table:threefuk}}
\end{table}

\subsection{The $B$-model mirror to the affine, relative and full Fukaya categories}
\label{subsec:introb}

We define the smooth Calabi-Yau Fermat hypersurface
\[ M^n := \left\{ \sum_{j=1}^n z_j^n = 0\right\} \subset \CP{n-1},\]
with ample divisors $D_j := \{z_j = 0\}$ for $j = 1, \ldots, n$. 
In this section, we will introduce the $B$-models which ought to be mirror to $\mathcal{F}(M^n \setminus D)$, $\mathcal{F}(M^n,D)$, and $\mathcal{F}(M^n)$.

We define 
\[ R := R(M^n,D) \subset \C \llbracket r_1, \ldots, r_n \rrbracket ,\]
which is the coefficient ring of $\mathcal{F}(M^n,D)$.  
We define the $R$-algebra
\[ S := R[u_1, \ldots, u_n],\]
and equip it with the $\Z$-grading so that $R$ is concentrated in degree $0$, and each $u_j$ has degree $1$.
We define the element 
\[ w = u_1 \ldots u_n + \sum_{j=1}^n r_j u_j^n \in S,\]
of degree $n$.

Now note that $\mathrm{Proj}(S) = \mathbb{P}^{n-1}_R$. 
We consider the variety
\[ \widetilde{N}^n := \{w = 0\} \subset \mathbb{P}^{n-1}_R,\]
and equip it with the action of $\Gamma^*_n$, exactly as we did for $\widetilde{N}^n_{nov}$ in Definition \ref{definition2:2mirror}, then define
\[ N^n := \widetilde{N}^n/ \Gamma^*_n.\]

We note that the algebraic torus 
\[ \mathbb{T} := \left\{ (\lambda_1, \ldots, \lambda_n) \in (\C^*)^n: \lambda_1 \ldots \lambda_n = 1 \right\}\]
acts on $S$, by sending
\begin{eqnarray*}
u_j & \mapsto & \lambda_j u_j,\\
r_j & \mapsto & \lambda_j^{-n} r_j.
\end{eqnarray*}
This action preserves $w$, and commutes with the action of $\Gamma_n^*$, and therefore defines an action of $\mathbb{T}$ on $N^n$. 

We note that there are homomorphisms
\[ \C  \leftarrow R \To \Lambda,\]
given by
\[ 0\,\, \mbox{\reflectbox{$\mapsto$}} \,\, r_j \mapsto r\]
for all $j$.
Hence, by base change, we obtain
\[ \begin{diagram}
N^n_0 & \rTo & N^n & \lTo & N^n_{nov}\\
\dTo && \dTo && \dTo \\
\mathrm{Spec}(\C) & \rTo & \mathrm{Spec}(R) & \lTo & \mathrm{Spec}(\Lambda).
\end{diagram} \]
We call $N^n$ the {\bf total space} of the family over $\mathrm{Spec}(R)$, we call $N^n_0$ the {\bf special fibre}, and $N^n_{nov}$ the {\bf generic fibre}. 
We observe that this $N^n_{nov}$ coincides with the definition given in Definition \ref{definition2:2mirror}.

We expect mirror relationships as follows:

\begin{center}
\begin{tabular}{c|c|c|c}
Coefficients & $\C$ & $R$ &   $\Lambda$ \\ \hline
$B$-model & $\mathrm{Perf}(N^n_0)$ & $\mathrm{Perf}(N^n)$ & $Coh(N^n_{nov})$\\ \hline
$A$-model & $\mathcal{F}(M^n \setminus D)$ & $\mathcal{F}(M^n,D)$ &  $\mathcal{F}(M^n)$ 
\end{tabular}
\end{center}
(however note that we do {\bf not} necessarily claim to prove all of these equivalences -- this table is included to assist the reader in seeing the `big picture').

Note that the map from $R$ to $\C$ is naturally $\mathbb{T}$-equivariant.
Thus, we can define $\mathbb{T}$-equivariant sheaves on $N^n_0$ and $N^n$.
However  it is impossible to equip $\Lambda$ with a $\mathbb{T}$-action so that the map $R \To \Lambda$ is $\mathbb{T}$-equivariant, so we can not talk about $\mathbb{T}$-equivariant sheaves on $N^n_{nov}$.

We expect that $\mathbb{T}$-equivariant sheaves correspond to anchored Lagrangian branes. 
Namely, we have seen that the morphism spaces between anchored Lagrangian branes admit a grading in the group $H_1(M^n \setminus D)$, and hence an action of its character group. 
In this case, there is a natural isomorphism
\[ \mathrm{Hom}\left(H_1(M^n \setminus D),\C^*\right) \cong \mathbb{T},\]
and we expect the mirror correspondences in the above table to be $\mathbb{T}$-equi\-vari\-ant (excluding the last column). 
In fact, we expect something stronger: they should be equivalences of $\bm{G}$-graded categories.

In \cite[Theorem 7.4]{Sheridan2011a}, we proved that there is a fully faithful, $\mathbb{T}$-equivariant embedding
\[ \mathrm{Perf}(N^n_0) \hookrightarrow D^{b}\mathcal{F}(M^n \setminus D).\]
In this paper, we extend this to prove results about the other columns.

Our ultimate aim is to understand the final column, and give a proof of Theorem \ref{theorem2:2main}. 
Thus we need a method of making computations in $D^bCoh(N^n_{nov})$. 
For this purpose, we use the category of {\bf graded matrix factorizations}. 
Let us denote
\[ S_{nov} := S \otimes_R \Lambda = \Lambda[u_1, \ldots, u_n],\]
with
\[ w_{nov} := w \otimes 1 \in S_{nov}.\]
Then $S_{nov}$ is $\Z$-graded, and $w_{nov}$ is homogeneous of degree $n$. 
In \cite{Orlov2009}, Orlov introduced the DG category of graded matrix factorizations of a homogeneous superpotential, $\mathrm{GrMF}(S_{nov},w_{nov})$, and proved that there is an $A_\infty$ quasi-equivalence
\[ \mathrm{GrMF}(S_{nov},w_{nov}) \cong D^bCoh \left( \widetilde{N}^n_{nov} \right),\]
where
\[ \widetilde{N}^n_{nov} := \{ w_{nov} = 0\} \subset \mathbb{P}^{n-1}_{\Lambda}.\]

\begin{remark2}
\label{remark2:lunts}
Actually, Orlov proved \cite[Theorem 3.11]{Orlov2009} that there is a triangulated equivalence on the level of homotopy categories, and it follows from work of Lunts and Orlov \cite[Theorem 2.13]{Lunts2010} on uniqueness of DG enhancements that this extends to a DG quasi-equivalence. 
Recall that two DG categories $\mathscr{A}$ and $\mathscr{B}$ are said to be DG quasi-equivalent if there is a chain of DG functors between DG categories,
\[ \mathscr{A} \leftarrow \mathscr{C}_1 \rightarrow \ldots \leftarrow \mathscr{C}_k \rightarrow \mathscr{B},\]
such that each DG functor is a quasi-equivalence, i.e., its cohomology-level functor is an equivalence.  
Since any DG quasi-equivalence is {\it a fortiori} an $A_\infty$ quasi-equivalence, and any $A_\infty$ quasi-equivalence over a field admits an $A_\infty$ inverse up to homotopy (see, for example, \cite[Corollary 1.14]{Seidel2008}), DG quasi-equivalence is a stronger notion than $A_\infty$ quasi-equivalence. 
\end{remark2}

Similarly, there is an $A_\infty$ quasi-equivalence of the corresponding $\Gamma_n^*$-equi\-vari\-ant categories,
\[ \mathrm{GrMF}(S_{nov},w_{nov})^{\Gamma_n^*} \cong D^bCoh \left( N^n_{nov} \right).\]

Orlov's theorem applies because we work over the field $\Lambda$, and $\widetilde{N}^n_{nov} = \{w_{nov} = 0\}$ is smooth and Calabi-Yau. 
However, recall that by passing from the variety $N^n$, defined over $R$, to the variety $N^n_{nov}$, defined over $\Lambda$, we lose the $\mathbb{T}$-action. 
This is a disadvantage, because $\mathbb{T}$-equivariance constrains the algebraic structures we consider significantly, and makes our classification problems tractable. 

Therefore, we introduce (in Section \ref{sec:2bmodel}) the category $\mathrm{GrMF}(S,w)$ of graded matrix factorizations of $w \in S$, over the coefficient ring $R$.
The coefficient ring $R$ and the graded $R$-algebra $S$ admit an action of $\mathbb{T}$, which preserves $w$. 
Therefore we can talk about $\mathbb{T}$-equivariant objects in $\mathrm{GrMF}(S,w)$, and furthermore there is a fully faithful embedding
\[ \mathrm{GrMF}(S,w) \otimes_R \Lambda \hookrightarrow \mathrm{GrMF}(S_{nov},w_{nov}).\]
It seems reasonable to hope that there is a relationship between the categories $\mathrm{GrMF}(S,w)$ and $\mathrm{Perf}(\widetilde{N}^n)$ in the vein of Orlov's theorem, but we do not pursue this.

In fact, we first introduce a differential $\bm{G}$-graded category $MF^{\bm{G}}(S,w)$ of $\bm{G}$-graded matrix factorizations of $w \in S$ (these combine the $\mathbb{T}$-action with the $\Z$-grading, in the same way that anchored Lagrangian branes combine the $H_1(M^n \setminus D)$-grading with the $\Z$-grading in the Fukaya category). 
We show that $\mathrm{GrMF}(S,w)^{\Gamma_n^*}$ is some `orbifolding' of it.

The starting point for this is the observation  that $\mathrm{GrMF}(S,w)$ is a $\Z_n$-equivariant version of $MF^{\bm{G}}(S,w)$ (compare \cite{Polishchuk2011,Caldararu2010}). 
In fact there is an action of $\tilde{\Gamma}_n^*$ on $MF^{\bm{G}}(S,w)$, and we show that there is a fully faithful embedding
\[MF^{\bm{G}}(S,w)^{\tilde{\Gamma}_n^*} \hookrightarrow \mathrm{GrMF}(S,w)^{\Gamma_n^*}\]
(recall that $\tilde{\Gamma}_n^*$ is an extension of $\Z_n$ by $\Gamma_n^*$; the $\Z_n$ got eaten up turning $MF$ into $\mathrm{GrMF}$). 
Actually, our notation for the $\tilde{\Gamma}_n^*$-equivariant category in the main body of the text is different -- we will write it as 
\[ \bm{p}_1^* MF^{\bm{G}}(S,w) \equiv MF^{\bm{G}}(S,w)^{\tilde{\Gamma}_n^*}.\]
We will not give the precise meaning of `$\bm{p}_1^*$' in this introduction, but will continue to use it for consistency with our later notation.
 
We consider the object $\mathcal{O}_0$ of $MF^{\bm{G}}(S,w)$, corresponding to the ideal $(u_1, \ldots, u_n) \subset S$. 
We denote its endomorphism algebra by $\mathscr{B}$. 
It is a deformation of the exterior algebra
\[ \mathrm{Ext}^*_{Coh(\C^n)}(\mathcal{O}_0,\mathcal{O}_0) \cong \Lambda^* \C^n\]
over the power series ring $R$.
Deformations of $\Lambda^* \C^n$ are governed by Hochschild cohomology, which is given by polyvector fields, by the Hochschild-Kostant-Rosenberg isomorphism:
\[ HH^*(\Lambda^* \C^n) \cong \C\llbracket u_1, \ldots, u_n \rrbracket [\theta_1, \ldots, \theta_n],\]
where the variables $u_j$ commute and the variables $\theta_j$ anti-commute.
We construct a minimal model for this $A_{\infty}$ deformation of $\Lambda^* \C^n$, and prove that its deformation classes are given exactly by the coefficients of $w$ (following \cite{Efimov2009}). 
We prove a classification theorem (Theorem \ref{theorem2:2typea}), which shows that these deformation classes, together with the $\bm{G}$-grading, are enough to determine the deformation up to $A_{\infty}$ quasi-isomorphism and formal change of variables.

We then consider the full subcategory of  $\bm{p}_1^*MF^{\bm{G}}(S,w)$ whose objects are the equivariant twists of $\mathcal{O}_0$. 
We denote it by $\widetilde{\mathscr{B}}$. 
It can be determined completely from $\mathscr{B}$. 
We also introduce a full subcategory 
\[ \widetilde{\mathscr{A}} \subset \mathcal{F}(M^n,D)\]
(the remaining sections of this introduction will consist of an explanation of how to compute $\widetilde{\mathscr{A}}$).
We prove the following generalization of \cite[Theorem 7.4]{Sheridan2011a}:

\begin{theorem2}
\label{theorem2:2main2}
There exists a $\tilde{\psi} \in \C\llbracket T \rrbracket  \subset R$, where $T = r_1  \ldots r_n$, with $\tilde{\psi}(0) = \pm 1$, and a quasi-isomorphism of $\bm{G}$-graded $R$-linear $A_{\infty}$ categories
\[ \tilde{\psi}^* \cdot \widetilde{\mathscr{B}} \cong \widetilde{\mathscr{A}}.\]
\end{theorem2}

By our previous discussion, there is a fully faithful $A_{\infty}$ embedding
\[ \bm{p}_1^* MF^{\bm{G}}(S,w) \otimes_R \Lambda \hookrightarrow D^bCoh(N^n_{nov})\cong D^bCoh^{\Gamma_n^*}\left(\widetilde{N}^n_{nov} \right).\]
The images of the equivariant twists of $\mathcal{O}_0$ correspond to equivariant twists of the restrictions of the Beilinson exceptional collection $\Omega^j(j)$ restricted to $\widetilde{N}^n_{nov}$ (for $j = 0, 1, \ldots, n-1$) by characters of $\Gamma_n^*$.

$D^bCoh(N^n_{nov})$ is split-closed, and the equivariant twists of the restrictions of the Beilinson exceptional collection split-generate it. 
It follows immediately from Theorem \ref{theorem2:2main2}, by tensoring with $\Lambda$, that there is a quasi-equivalence of triangulated $A_{\infty}$ categories
\[ \hat{\psi}^* \cdot D^bCoh(N^n_{nov}) \cong D^\pi \left( \widetilde{\mathscr{A}} \otimes_R \Lambda \right) \subset D^{\pi} \mathcal{F}(M^n)\]
(under our assumption (Assumption \ref{assumption:relfull}) that the second embedding above exists).

Finally, to complete the proof, we wish to show that $\widetilde{\mathscr{A}} \otimes_R \Lambda$ split-generates the Fukaya category. 
We do this by applying the split-generation result of \cite{Abouzaid2012}, which says that, if the closed-open string map
\[ \mathcal{CO}: QH^*\left(M^n\right) \To HH^*\left(\widetilde{\mathscr{A}}\otimes_R \Lambda \right)\]
is non-zero in the top degree $2(n-2)$, then $\widetilde{\mathscr{A}} \otimes_R \Lambda$ split-generates $D^{\pi} \mathcal{F}(M^n)$. 

\begin{sloppypar}
We observe that 
\[ \mathcal{CO}([\omega]) = r \del{\mu^*}{r} \in HH^2 \left(\widetilde{\mathscr{A}} \otimes_R \Lambda \right)\]
(in words, the image of the class of the symplectic form under the closed-open string map is the class in $HH^2$ corresponding to deforming the Fukaya category by scaling the symplectic form). 
We now observe that $\mathcal{CO}$ is a $\Lambda$-algebra homomorphism, so
\[ \mathcal{CO}([\omega]^{n-2}) = \left(r \del{\mu^*}{r}\right)^{n-2}.\]
We then compute that this class is non-zero in the Hochschild cohomology. 
It follows from the split-generation criterion that $\widetilde{\mathscr{A}} \otimes_R \Lambda$ split-generates $D^{\pi} \mathcal{F}(M^n)$. 
This completes the proof of Theorem \ref{theorem2:2main}.
\end{sloppypar}

In the rest of this introduction, we explain how we make computations in the relative Fukaya category $\mathcal{F}(M^n,D)$, which are sufficient to prove Theorem \ref{theorem2:2main2}. 

\subsection{Behaviour of the Fukaya category under branched covers}
\label{subsec:2branchbehav}

Suppose that $N$ and $M$ are compact K\"{a}hler manifolds with smooth normal-crossings divisors $E \subset N$ and $D \subset M$ as before, and that 
\[\phi: (N,E) \To (M,D)\]
is a branched cover ramified about the divisors $E$, sending divisor $E_j$ to divisor $D_j$, and with ramification of degree $a_j$ about divisor $E_j$. 
We aim to understand how the affine and relative Fukaya categories of $(N,E)$ and $(M,D)$ are related.

First, we observe that the map
\[ \phi: N \setminus E \To M \setminus D\]
is an unbranched cover. 
Therefore, any holomorphic disk in $M \setminus D$ lifts to $N \setminus E$, because it is contractible. 
It follows that the problem of relating $\mathcal{F}(N \setminus E)$ to $\mathcal{F}(M \setminus D)$ is essentially one of algebraic bookkeeping: we need to keep track of how the holomorphic disks lift, but do not need to compute any new moduli spaces of disks. 
This leads one to the statement that $\mathcal{F}(N \setminus E)$ is a `semi-direct product of $\mathcal{F}(M \setminus D)$ with the character group of the covering group of $\phi$' (see \cite[Section 8b]{Seidel2003} and \cite[Section 9]{Seidel2008a}).
We rephrase this in Section \ref{subsec:2affcov} using the language of $\bm{G}$-graded categories, in which we write
\[ \mathcal{F}(N \setminus D) \cong  \bm{p}^*\mathcal{F}(M \setminus D) \]
(we won't explain this notation in the introduction).
 
Now we try to understand the behaviour of the relative Fukaya category with respect to branched covers.
This is not as simple as the unramified cover case, because holomorphic disks may pass through the branching locus, and then they do not lift to the cover. 
In order to relate $\mathcal{F}(N, E)$ to $\mathcal{F}(M,D)$, we introduce a `smooth orbifold relative Fukaya category'  $\mathcal{F}(M,D,\bm{a})$, where $\bm{a} = (a_1, \ldots,a_k)$ denotes the degrees of ramification of the cover about the divisors (but we could define the smooth orbifold relative Fukaya category for any tuple $\bm{a}$ of $k$ positive integers).

The objects and generators of the morphism spaces of $\mathcal{F}(M, D,\bm{a})$ are the same as for $\mathcal{F}(M,D)$. 
The coefficient ring $R(M,D,\bm{a})$ is slightly different: it is still a completion of $\widetilde{R}(M,D)$ in the category of $\bm{G}(M,D)$-graded algebras, but $\widetilde{R}(M,D)$ is equipped with a different $\bm{G}(M,D)$-grading from that used in the definition of $R(M,D)$. 
For example, the $H_1$-grading of $r_j \in \widetilde{R}(M,D)$ is defined to be $a_j$ times the class of a meridian loop about divisor $D_j$.
 
The $A_{\infty}$ structure maps $\mu^s$ count holomorphic disks $u: \mathbb{D} \To M$ that have ramification of degree $a_j$ about divisor $D_j$ wherever they intersect it. 
Each such disk contributes 
\[r_1^{\#(u \cap D_1)} \ldots r_k^{\#(u \cap D_k)} \in R(M,D,\bm{a})\]
(an intersection point of $u$ with $D_j$ contributes $1$ to $u \cap D_j$, although it contributes $a_j$ to the topological intersection number $u \cdot D_j$).
In particular, if $\bm{a} = (1,1, \ldots, 1)$ then we recover the relative Fukaya category.
The category is still $\bm{G}(M,D)$-graded.

The holomorphic disks $u: \mathbb{D} \To M$ contributing to the orbifold relative Fukaya category, now do lift to holomorphic disks $u: \mathbb{D} \To N$ (by the homotopy lifting criterion).
Thus, the relative Fukaya category $\mathcal{F}(N,E)$ is related to $\mathcal{F}(M,D,\bm{a})$ in exactly the same way that the affine Fukaya category $\mathcal{F}(N \setminus E)$ is related to $\mathcal{F}(M \setminus D)$: in the language of $\bm{G}$-graded categories,
\[ \mathcal{F}(N,D) \cong \bm{p}^* \mathcal{F}(M,D,\bm{a}).\]

It now remains to relate $\mathcal{F}(M,D,\bm{a})$ to $\mathcal{F}(M,D)$. 
In fact, we are only able to relate the `first-order' parts of the categories (but this turns out to be enough for our purposes).
The first-order relative Fukaya category is defined to be
\[ \mathcal{F}(M,D)/\mathfrak{m}^2 := \mathcal{F}(M,D) \otimes_R R/\mathfrak{m}^2,\]
where $\mathfrak{m} \subset R$ is the maximal ideal. 
It is linear over $R/ \mathfrak{m}^2$.
It retains only the information about rigid holomorphic disks $u: \mathbb{D} \To M$ passing through a {\bf single} divisor $D_j$ (with multiplicity $1$). 

Let us write
\[ \mu^* = \mu_0^* + \mu_1^*\]
for the $A_{\infty}$ structure maps $\mu^*$  in $\mathcal{F}(M,D)/\mathfrak{m}^2$, where $\mu_0^*$ gives the affine Fukaya category and $\mu_1^*$ gives the first-order terms.
Then the $A_{\infty}$ relations tell us that $\mu_1^*$ is a Hochschild cocycle, hence defines an element
\[ \sum_{j=1}^k r_j \alpha_j \in HH^*(\mathcal{F}(M \setminus D)) \otimes \mathfrak{m} / \mathfrak{m}^2.\]
We call $\alpha_j$ the {\bf first-order deformation classes} of $\mathcal{F}(M,D)$.

We prove (Theorem \ref{theorem2:2defclassram}) that, if $\mathcal{F}(M,D)$ has first-order deformation classes $\alpha_j$, then $\mathcal{F}(M,D,\bm{a})$ has first-order deformation classes $\alpha_j^{a_j}$, where the power  is taken with respect to the Yoneda product on Hochschild cohomology. 
The proof looks very similar to the proof that the map $QH^*(M) \To HH^*(\mathcal{F}(M))$ is a ring homomorphism.

\begin{remark2}
At first sight, this may seem a strange result: first-order deformation classes of a category live in $HH^2$, and the Yoneda product respects the $\Z$-grading, so one would expect the class $\alpha_j^{a_j}$ to no longer live in $HH^2$ and therefore not be an appropriate first-order deformation class. 
The solution lies in the fact that the coefficient rings $R$ have non-trivial gradings, and in fact the coefficient rings for $\mathcal{F}(M,D)$ and $\mathcal{F}(M,D,\bm{a})$ have {\bf different} gradings: thus, both $r_j \alpha_j$ and $r_j \alpha_j^{a_j}$ have degree $2$ in the respective Hochschild cohomology groups in which they live.
\end{remark2}

Combining this result with our previous observations, if we have a branched cover $\phi:(N,E) \To (M,D)$, then we can compute $\mathcal{F}(N,E)$ to first order if we know $\mathcal{F}(M,D)$ to first order. 

\subsection{The Fukaya category of $M^n$}
\label{subsec:2onedim}

We will now explain how to compute the Fukaya category of $M^n$. 
We will keep the one-dimensional case ($n=3$) as a running example throughout, despite the fact that we do not prove this case of Theorem \ref{theorem2:2main} completely in this paper. 
We do this because one can see all of the holomorphic disks in the Fukaya category in this case, and gain intuition for the various versions of the Fukaya category that we introduce, and results that we prove about it.

We consider the Fermat hypersurfaces
\[ M^n_a := \left\{ \sum_{j=1}^n z_j^a = 0 \right\} \subset \CP{n-1}\]
with the smooth ample normal-crossings divisors $D_j = \{z_j = 0\}$ for $j = 1, \ldots, n$.
There is a branched cover
\begin{eqnarray*}
\phi: (M^n_a, D) & \To & (M^n_1,D), \\
{[}z_1: \ldots : z_{n}{]} & \mapsto & {[}z_1^a: \ldots : z_{n}^a{]}.
\end{eqnarray*} 
Thus, applying the results described in the previous section, if we can compute $\mathcal{F}(M^n_1,D)$ to first order, then we can compute $\mathcal{F}(M^n_a,D)$ to first order. 

We observe that $M^n_1 \cong \CP{n-2}$, and $D$ consists of $n$ hyperplanes in general position. 
$M^n_1 \setminus D$ is called the $(n-2)$-dimensional {\bf pair of pants}. 
In Section \ref{sec:2mb}, we construct a Lagrangian immersion
\[ L:S^{n-2} \To M^n_1 \setminus D\]
which has an anchored Lagrangian brane structure. 
This Lagrangian was introduced in \cite{Sheridan2011a}.

\begin{figure}
\centering
\includegraphics[width=0.7\textwidth]{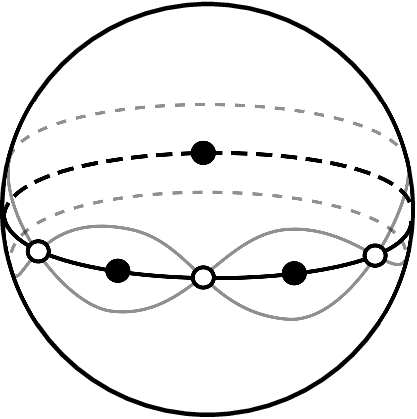}
\caption{$\CP{1}$ with its real locus $\RP{1}$ is shown in black, with the divisors $D$ indicated by solid dots. 
The Lagrangian immersion $L:S^1 \To \CP{1} \setminus D$ is shown in grey, and its self-intersection points are marked with open dots. 
\label{fig:2l1col}}
\end{figure}

\begin{example2}
$M^3_1 = \CP{1}$, and $D$ consists of three divisors (points). 
The Lagrangian immersion $L:S^1 \To \CP{1} \setminus D$ is shown in Figure \ref{fig:2l1col}. 
The $A_{\infty}$ algebra $CF^*(L,L)$ was described in \cite{Seidel2008a}. 
It was introduced as a $\Z \oplus H_1(M^3_1 \setminus D)$-graded category, but it is not hard to see the underlying $\bm{G}(M^3_1,D)$-graded structure.
\end{example2}

\begin{figure}
\centering
\includegraphics[width=0.95\textwidth]{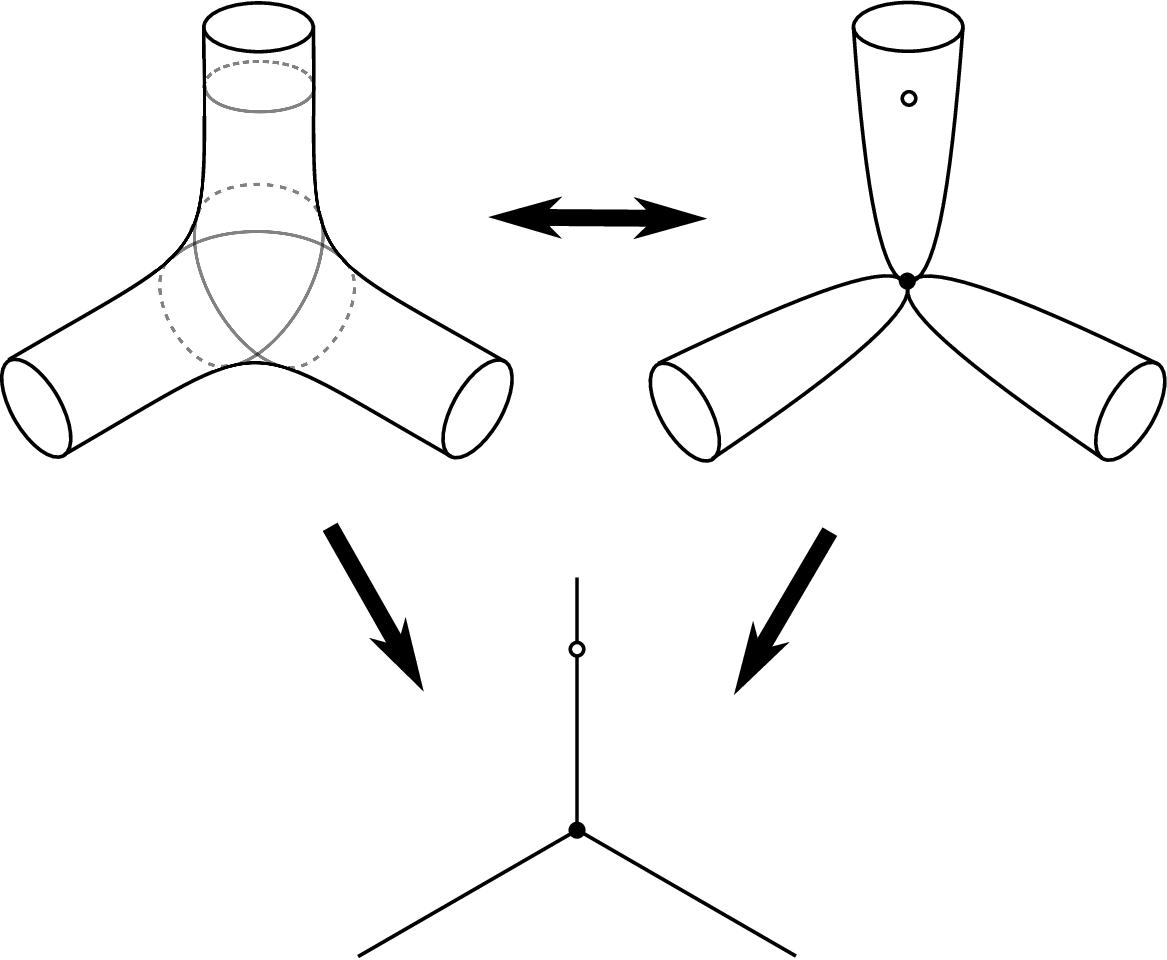}
\caption{The Lagrangian $L:S^1 \To \CP{1} \setminus D$ can also be drawn as a `trefoil', shown here in grey in the centre of the pair of pants (upper left). 
It lies over the central point (shown as a solid dot) in the tropical pair of pants (bottom). 
The mirror is the union of coordinate axes in $\C^3$ (upper right), and the object corresponding to $L$ is the structure sheaf of the origin (shown as a solid black dot). 
Another SYZ torus fibre is drawn in grey on the upper leg of the pair of pants. 
It lies over a point (shown as an open dot) in the tropical pair of pants.
It is mirror to the structure sheaf of a point (shown as an open dot).
\label{fig:2syzpic}}
\end{figure}

The reason this Lagrangian is important is because it can be regarded as a `fibre' in a Strominger-Yau-Zaslow fibration. 
See Figure \ref{fig:2syzpic} for the picture in the one-dimensional case. 
More generally, as shown in \cite{Mikhalkin2004}, the pair of pants $\CP{n-2} \setminus D$ is a singular torus fibration over the `tropical amoeba of the pair of pants', which is some space stratified by affine manifolds.
The torus fibration is non-singular over the top-dimensional faces of the tropical pair of pants. 

We suggest that one should think, not of an SYZ fibration of the pair of pants over the tropical pair of pants, with some singular fibres, but rather of an SYZ {\bf family} of objects of the Fukaya category, parametrized by the tropical pair of pants. 
The immersed Lagrangian sphere $L$ is the object corresponding to the central point in the tropical pair of pants in this picture. 
The objects corresponding to points of the top-dimensional strata are Lagrangian torus fibres (recall that the fibration is non-singular there).
The objects corresponding to points on the in-between strata are lower-dimensional incarnations of $L$, crossed with tori. 
We provided some evidence for this philosophy in \cite{Sheridan2011a}, where we showed that the endomorphism algebra of $L$ in $\mathcal{F}(\CP{n-2} \setminus D)$ is quasi-isomorphic to the endomorphism algebra of the structure sheaf of the origin in the mirror category of matrix factorizations. 

\begin{figure}
\centering
\subfloat[Disks contributing to $CF^*(L,L)$ in the affine Fukaya category, $\mathcal{F}(\CP{1} \setminus D)$.]{
\includegraphics[width=0.45\textwidth]{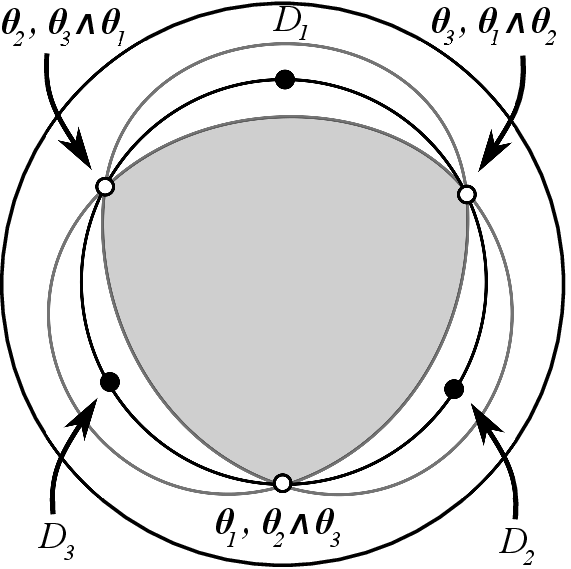} 
\label{subfig:2l1triang}}
\hfill
\subfloat[Disks contributing to the first-order deformation class of $CF^*(L,L)$ in the relative Fukaya category, $\mathcal{F}(\CP{1},D)$.]{ 
\includegraphics[width=0.45\textwidth]{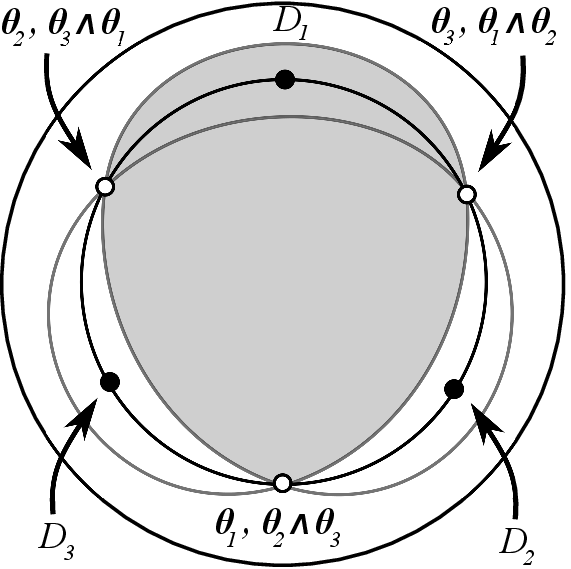} 
\label{subfig:2l1tear}}
\caption{Holomorphic disks contributing to $CF^*(L,L)$.
\label{fig:2l1disks}}
\end{figure}

In Section \ref{sec:2mb}, we compute $CF^*(L,L)$ to first order in $\mathcal{F}(\CP{n-2},D)$, using a Morse-Bott model for the relative Fukaya category, based on the `cluster homology' of \cite{Cornea2006}. 
We compute that the underlying vector space is an exterior algebra:
\[ CF^*(L,L) \cong \Lambda^* \C^n \cong \C[\theta_1, \ldots, \theta_n],\]
where the variables $\theta_j$ anti-commute. 
For example, when $n=3$, $CF^*(L,L)$ is generated by $H^*(S^1)$ (whose two generators we identify as the bottom and top classes $1$ and $\theta_1 \wedge \theta_2 \wedge \theta_3$), together with two generators for each self-intersection point, which we label as in Figure \ref{fig:2l1disks}. 

We next show that the zeroth-order algebra structure, $\mu^2_0$, coincides with the exterior algebra. 
In the case $n=3$, the corresponding holomorphic triangles are shown in Figure \ref{subfig:2l1triang}. 
The shaded triangle can be viewed as having inputs $\theta_1$ and $\theta_2$ and output $\theta_1 \wedge \theta_2$, while the corresponding triangle on the back of the figure can be viewed as having inputs $\theta_2$ and $\theta_1$ and output $-\theta_1 \wedge \theta_2$. 
The other products follow similarly.

$A_{\infty}$ structures with underlying cohomology algebra $\Lambda^* \C^n$ are classified by the Hochschild cohomology, which is given by polyvector fields, by the Hochschild-Kostant-Rosenberg isomorphism:
\[ HH(\Lambda ^*\C^n) \cong \C\llbracket u_1, \ldots, u_n \rrbracket [\theta_1, \ldots, \theta_n],\]
where variables $u_i$ commute and $\theta_i$ anti-commute. 
We show that the endomorphism algebra $CF^*(L,L)$ in the affine Fukaya category $\mathcal{F}(\CP{n-2} \setminus D)$ is completely determined, up to $A_{\infty}$ quasi-isomorphism, by a single higher-order product, having the form
\[ \mu^n(\theta_1, \ldots, \theta_n) = 1,\]
corresponding to the Hochschild cohomology class
\[ u_1 \ldots u_n \in \C\llbracket u_1, \ldots, u_n \rrbracket [\theta_1, \ldots, \theta_n].\]
In the case $n=3$, we can see the corresponding holomorphic disk in Figure \ref{subfig:2l1triang}. 
It is the shaded triangle, which we view as a degenerate $4$-gon having inputs $\theta_1$, $\theta_2$, $\theta_3$, and output a degenerate vertex on one of the sides of the triangle, corresponding to $1$. 

We next compute that the endomorphism algebra of $CF^*(L,L)$ in the first-order relative Fukaya category $\mathcal{F}(\CP{n-2},D)/\mathfrak{m}^2$ is determined by structure maps of the form
\[ \mu^1(\theta_j) = r_j \cdot 1,\]
corresponding to first-order deformation classes 
\[ r_j u_j \in \C\llbracket u_1, \ldots, u_n \rrbracket [\theta_1, \ldots, \theta_n] \otimes \mathfrak{m}/\mathfrak{m}^2.\]
When $n=3$, we can see the corresponding holomorphic disks in Figure \ref{subfig:2l1tear}. 
The shaded `teardrop' shape has one input $\theta_1$, and a degenerate output vertex corresponding to $1$. 
It intersects divisor $D_1$ exactly once, and does not intersect the other divisors, hence contributes with a coefficient $r_1$. 
Thus it gives rise to the term $\mu^1(\theta_1) = r_1 \cdot 1$.

It follows from the result described in Section \ref{subsec:2branchbehav} that the first-order deformation classes of the $A_{\infty}$ algebra
\[ \mathscr{A} :=  CF^*_{\mathcal{F}(\CP{n-2},D,(n,\ldots,n))}(L,L)\]
in the orbifold Fukaya category, are $r_j u_j^n$. 
Thus, the full deformation class of $\mathscr{A}$ is
\[u_1 \ldots u_n + \sum_{j=1}^{n} r_j u_j^n + \mathcal{O}(r^2),\]
which we observe coincides with the defining polynomial $w$ of $N^n$, to first order.

We prove a classification theorem (Theorem \ref{theorem2:2typea}) which shows that this is enough information to determine the full $\bm{G}$-graded deformation, up to $A_{\infty}$ quasi-isomorphism and formal change of variables. 

We show that the $A_{\infty}$ algebra
\[\mathscr{B} := hom^*_{MF^{\bm{G}}(S,w)}(\mathcal{O}_0,\mathcal{O}_0)\]
also has the same underlying algebra, $\bm{G} \cong \bm{G}(\CP{n-2},D)$-grading, and deformation classes (recall the deformation classes were given exactly by $w$ itself). 
Therefore, by the above-mentioned classification theorem, we have a formal change of variables $\tilde{\psi}^*$, and an $A_{\infty}$ quasi-isomorphism
\[ \tilde{\psi}^* \cdot \mathscr{B} \cong  \mathscr{A}.\]

\begin{figure}
\centering
\includegraphics[width=0.9\textwidth]{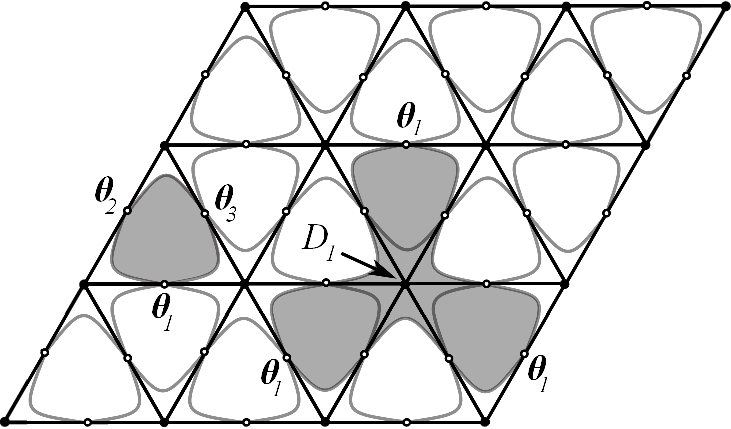} 
\caption{A fundamental domain of the elliptic curve $M^3$, with the divisor $D$ consisting of nine points, indicated by black dots.
The $9$-fold cover $(M^3,D) \To (\CP{1},D)$ has ramification of order $3$ about each divisor. 
The black lines are the pullback of $\RP{1}$. 
The lifts of $L$ are the curves shown in grey, with their intersection points shown as open dots. 
We have shaded two of the holomorphic disks contributing to the $A_{\infty}$ products between lifts of $L$. 
Note that the triangle of Figure \ref{subfig:2l1triang}, which contributes to the affine Fukaya category, lifts directly to $M^3$, whereas the `teardrop' illustrated in Figure \ref{subfig:2l1tear}, which contributes the term $r_1 u_1$ to the deformation class of the relative Fukaya category, does not lift, but rather gives rise to the `clover-leaf' shape, which contributes the deformation class $r_1 u_1^3$.
\label{fig:2l1torus}}
\end{figure}

We now define the full subcategories 
\[\widetilde{\mathscr{A}} := \bm{p}_1^* \mathscr{A} \subset \mathcal{F}(M^n,D)\]
and
\[ \widetilde{\mathscr{B}} := \bm{p}_1^* \mathscr{B} \subset \bm{p}_1^* MF^{\bm{G}}(S,w).\] 
It follows from the preceding argument that we have
\[ \tilde{\psi}^* \cdot \widetilde{\mathscr{B}} \cong  \tilde{\psi}^* \cdot\bm{p}_1^* \mathscr{B} \cong \bm{p}_1^* \mathscr{A} \cong \widetilde{\mathscr{A}}.\]
This completes the proof of Theorem \ref{theorem2:2main2}. See Figure \ref{fig:2l1torus} for a picture in the case $n=3$.

{\bf Acknowledgements: } 
I would like to thank my advisor, Paul Seidel, for all of his help with this work. These results owe a great deal to him, both because his previous work \cite{Seidel2002,Seidel2003,Seidel2008,Seidel2008a} laid the foundations for them, and because he gave me a lot of guidance and made many useful suggestions along the way. 
I would also like to thank Mohammed Abouzaid for many helpful discussions and suggestions, and for showing me a preliminary version of \cite{Abouzaid2012}. 
I would also like to thank Grisha Mikhalkin for helping me to find the construction of the immersed Lagrangian sphere in the pair of pants, on which this work is based.

\section{Graded and equivariant categories}
\label{sec:2deftheory}

The main purpose of this section is to introduce the relevant notions of graded and equivariant algebraic objects, and modify the results in \cite[Section 3]{Seidel2003} to classify such objects.

\subsection{Grading data}
\label{subsec:2graddat}

For the purposes of this section, we fix an integer $n \ge 3$.

\begin{definition2}
An {\bf unsigned grading datum} $\bm{G}$ is an abelian group $Y$ together with a morphism $f: \Z \To Y$. 
We will use the shorthand $\bm{G} = \{ \Z \overset{f}{\To} Y\}$. 
We will often write $\bm{G}$ as
\[ \begin{diagram}
\Z & \rTo^{f} & Y &\rTo^{g} & X & \rTo & 0,
\end{diagram}
\]
where $X$ is the cokernel of $f$.
We say that $\bm{G}$ is {\bf exact} if the map $\Z \To Y$ is injective. 
\end{definition2}

\begin{definition2}
A {\bf morphism} of unsigned grading data, $\bm{p}: \bm{G}_1 \To \bm{G}_2$, is a morphism $p: Y_1 \To Y_2$ that makes the following diagram commute:
\[\begin{diagram}
 \Z & \rTo^{f_1} & Y_1  \\
\dEq & & \dTo<{p}  \\
\Z &  \rTo^{f_2} & Y_2,
\end{diagram}
\]
Composition of morphisms is defined in the obvious way, and this defines a category of unsigned grading data. 
We say that a morphism of unsigned grading data is {\bf injective} (respectively {\bf surjective}) if the map $p$ is injective (respectively surjective). 
We will sometimes write $\bm{p}$ as 
\[\begin{diagram}
 \Z & \rTo^{f_1} & Y_1 & \rTo^{g_1} & X_1 & \rTo & 0  \\
\dEq & & \dTo<{p} & & \dTo<{p_X} \\
\Z &  \rTo^{f_2} & Y_2 & \rTo^{g_2} & X_2 & \rTo & 0,
\end{diagram}
\]
where $p_X$ is the map induced by $p$.
\end{definition2}

\begin{definition2}
We define the {\bf sign grading datum}, $\bm{G}_{\sigma} := \{ \Z \To \Z_2\}$.
\end{definition2}

\begin{definition2}
\label{definition2:gradingdat}
We define a {\bf grading datum}  $(\bm{G},\bm{\sigma})$ to be an unsigned grading datum $\bm{G}$, together with a {\bf sign morphism}, which is a morphism of unsigned grading data,
\[ \bm{\sigma}: \bm{G} \To \bm{G}_{\sigma}.\]
We define a morphism of grading data to be a morphism of unsigned grading data that is compatible with sign morphisms. 
Henceforth, we will often omit the sign morphism $\bm{\sigma}$ from the notation to avoid clutter.
\end{definition2}

The sign morphism is important because it allows us to define certain signs in our algebraic objects, which allows us to work over coefficient rings in which $2 \neq 0$.

\begin{example2}
We define the grading datum $\bm{G}_{\Z} := \{ \Z \overset{\mathrm{id}}{\To} \Z\}$, with the obvious morphism $\bm{\sigma}$.
It is an initial object in the category of grading data.
\end{example2}

In practice, when doing computations we will work with objects called pseudo-grading data.

\begin{definition2}
\label{definition2:2pseudograd}
A {\bf pseudo-grading datum} $\bm{H}$ is a morphism of abelian groups $f: Z \To Y$, together with an element $c \in \mathrm{Hom}(Z,\Z)$, whose image lies inside $2\Z \subset \Z$.
\end{definition2}

\begin{definition2}
\label{definition2:2morps}
A {\bf morphism} of pseudo-grading data, 
\[\bm{p}: \bm{H}_1 \To \bm{H}_2,\]
consists of maps $p_Z$ and $p_Y$ that make the following diagram commute:
\[\begin{diagram}
Z_1 & \rTo^{f_1} & Y_1 \\
\dTo<{p_Z} & & \dTo<{p_Y} \\
Z_2 &  \rTo^{f_2} & Y_2,
\end{diagram}
\]
together with an element $d \in \mathrm{Hom}(Y_1,\Z)$, whose image lies inside $2\Z \subset \Z$, such that
\[ c_1 = p_Z^*(c_2) + f_1^*(d).\]
\end{definition2}

\begin{definition2}
Given morphisms of pseudo-grading data
\[ \bm{H}_0 \overset{\bm{p}_0}{\longrightarrow} \bm{H}_1 \overset{\bm{p}_1}{\longrightarrow} \bm{H}_2,\]
we define their composition $\bm{p} = \bm{p}_1 \circ \bm{p}_0$ by composing maps in the obvious way, and setting
\[ d := d_0 + p_{0,Y}^* d_1.\]
\end{definition2}

\begin{definition2}
\label{definition2:gradfromps}
Given a pseudo-grading datum $\bm{H}$:
\[ \begin{diagram}
Z & \rTo^{f} & Y,
\end{diagram}\]
together with $c$, we define a grading datum $\bm{G}(\bm{H})$ by
\[ \begin{diagram}
\Z & \rTo & (\Z \oplus Y)/Z,
\end{diagram}\]
where we define the map $Z \To \Z \oplus Y$ by 
\[ z \mapsto c(-z) \oplus f(z)\]
and the other maps in the obvious way. 
We define the sign morphism $\bm{\sigma}: \bm{G}(\bm{H}) \To \bm{G}_{\sigma}$ by
\begin{eqnarray*}
\sigma: (\Z \oplus Y)/Z & \To & \Z_2,\\
\sigma(j \oplus y) & := & j.
\end{eqnarray*}
Observe that the condition that the image of $c$ lies in $2\Z$ ensures that $\sigma$ is well-defined.
\end{definition2}

\begin{definition2}
\label{definition2:2morpsmor}
Given a morphism of pseudo-grading data $\bm{p}: \bm{H}_1 \To \bm{H}_2$ as in Definition \ref{definition2:2morps}, we define a corresponding morphism of grading data
\[ \bm{G}(\bm{p}): \bm{G}(\bm{H}_1) \To \bm{G}(\bm{H}_2),\]
where
\[ \bm{G}(\bm{p})(j \oplus y_1) := (j + d(y_1)) \oplus p(y_1).\]
It is not hard to check that this defines a functor from the category of pseudo-grading data to the category of grading data.
\end{definition2}

\begin{example2}
If we denote by $\bm{0} := \{0 \To 0\}$ the zero morphism, then $\bm{G}(\bm{0}) = \bm{G}_{\Z}$.
\end{example2}

\begin{example2}
\label{example2:grmfd}
Given $n \ge 0$, we define the pseudo-grading datum 
\[\bm{H}_{MF(n)} := \{ \Z \overset{\times n}{\To} \Z\},\] 
with $c = 2$.
We denote the corresponding grading datum by $\bm{G}_{MF(n)}$. 
It is exact, and has corresponding short exact sequence
\[ \begin{diagram}
0 & \rTo & \Z & \rTo^{f} &  (\Z \oplus \Z)/(-2,n)  &\rTo^{g} & \Z_n & \rTo & 0,
\end{diagram}\]
This grading datum is important because it controls Orlov's `category of graded matrix factorizations' of a superpotential of degree $n$ (hence our terminology). 
\end{example2}

Now we will introduce another important grading datum, although first we introduce a bit of convenient notation.

\begin{definition2}
We denote $[k] := \{1, \ldots, k\}$ for any positive integer $k$ and, for any $K \subset [k]$, 
\[ y_K := \sum_{j \in K} y_j \in \Z \langle y_1, \ldots, y_{k} \rangle \cong \Z^k.\]
\end{definition2}

\begin{example2}
\label{example2:ses1}
Given $n \ge 1$, we denote by $\bm{H}^n_a$ the pseudo-grading datum 
\[\begin{diagram}
\Z & \rTo^{y_{[n]}} & \Z \langle y_1, \ldots, y_{n} \rangle,
\end{diagram}\]
together with $c =2(n-a)$. 
We denote by $\bm{G}^n_a$ the corresponding grading datum. 
This grading datum is important because it controls both the relative Fukaya category and the category of equivariant matrix factorizations that we will consider.
\end{example2}

Now we prove a Lemma relating some of the grading data that we have introduced. 
We will use it to relate the category of matrix factorizations to the category of coherent sheaves.

\begin{lemma2}
\label{lemma2:squaregrad}
For any $n \ge 1$, there is a commutative square of grading data:
\[ \begin{diagram}
\bm{G}^n_n & \rOnto^{\bm{q}_1} & \bm{G}_{\Z} \\
\dInto^{\bm{p}_1} & & \dInto^{\bm{p}_2} \\
\bm{G}^n_1 & \rOnto^{\bm{q}_2} & \bm{G}_{MF(n)},
\end{diagram}\]
such that $\bm{p}_1$ and $\bm{p}_2$ are injective.
\end{lemma2}
\begin{proof}
The morphism $\bm{p}_2$ is uniquely defined because $\bm{G}_{\Z}$ is an initial object. 
It is clearly injective.
The morphism $\bm{q}_1$ comes from the zero morphism of pseudo-grading data, with $d = 0$ (note this is a morphism of pseudo-grading data because $c=0$ in $\bm{H}^n_n$). 

The morphism $\bm{p}_1$ comes from the morphism of pseudo-grading data,
\[ \begin{diagram}
 \Z & \rTo^{y_{[n]}} & \Z^{n} \\
\dTo^{\times n} &  & \dTo^{ \times n} \\
 \Z & \rTo^{y_{[n]}} & \Z^{n} , \\
\end{diagram}\]
with $d = 2(1-n)y_{[n]}$ (where here we denote by $y_{[n]}$ the element of the dual space $(\Z^{n})^{\vee} \cong \Z^{n}$). 
It is clearly injective.

The morphism $\bm{q}_2$ comes from the morphism of pseudo-grading data,
\[ \begin{diagram}
 \Z & \rTo^{y_{[n]}} & \Z^{n}\\
\dTo^{\times (-1)} &  & \dTo^{-y_{[n]}} \\
 \Z & \rTo^{\times n} & \Z , \\
\end{diagram}\]
with $d = 2y_{[n]}$.

It is a simple exercise, applying Definition \ref{definition2:2morpsmor}, to check that the diagram commutes in the category of grading data (although it does not commute in the category of pseudo-grading data).
 \end{proof}

\subsection{Graded vector spaces}
\label{subsec:2ggrad}

We recall that, if $Y$ is an abelian group, then a {\bf $Y$-graded vector space} is a vector space $V$, together with a collection of vector spaces $V_y$, indexed by $y \in Y$, and an isomorphism
\[ V \cong \bigoplus_{y \in Y} V_y.\]

\begin{definition2}
\label{definition2:2pushgrad}
Let $p: Y_1 \To Y_2$ be a morphism of abelian groups. 
Given a $Y_1$-graded vector space $V$, we define the $Y_2$-graded vector space $p_* V$, where
\[ (p_* V)_b := \bigoplus_{p(a) = b} V_a.\]
In particular, the underlying vector space $V$ does not change.
\end{definition2}

\begin{definition2}
\label{definition2:2pullgrad}
Let $p: Y_1 \To Y_2$ be a morphism of abelian groups. 
Given a $Y_2$-graded vector space $V$, we define the $Y_1$-graded vector space $p^* V$, where
\[ (p^* V)_a := V_{p(a)}.\]
\end{definition2}

\begin{remark2}
\label{remark2:injpull}
If $p: Y_1 \To Y_2$ is injective, then $p^*V$ is just the part of $V$ whose $Y_2$-degree lies in $\mathrm{im}(p)$.
\end{remark2}

\begin{definition2}
Let $\bm{G} = \{ \Z \To Y\}$ be a grading datum.
A $\bm{G}$-graded vector space $V$ is the same thing as a $Y$-graded vector space $V$.
\end{definition2}

In fact, for the purposes of this section, `$\bm{G}$-graded' is virtually identical to `$Y$-graded'. 
Things will get more complicated in the subsequent sections.

\begin{definition2}
Given an element $y \in Y$, and a $\bm{G}$-graded vector space $V$, we define $V[y]$ to be $V$ with grading shifted by $y$. 
\end{definition2}

\begin{definition2}
Given a morphism $\bm{p}$ of grading data, we define operations $\bm{p}_*$ and $\bm{p}^*$ on $\bm{G}$-graded vector spaces to be identical to $p_{*}$ and $p^*$.
\end{definition2}

\begin{definition2}
\label{definition2:2z2grad}
If $V$ is a $\bm{G}$-graded vector space, then it automatically becomes a $\Z_2$-graded vector space: 
\[ V \cong \bm{\sigma}_* V,\]
where $\bm{\sigma}$ is the sign morphism of $\bm{G}$. 
Given $v \in V$ of pure degree with respect to this $\Z_2$-grading, we denote its $\Z_2$-degree by $\sigma(v)$.
\end{definition2}

\begin{definition2}
A $\bm{G}$-graded algebra is a $Y$-graded algebra, i.e., one whose multiplication respects the $Y$-grading.
\end{definition2}

\begin{example2}
\label{example2:vecu}
Recall the grading datum $\bm{G}^n_1$ of Example \ref{example2:ses1}.
We define the $n$-dimensional $\bm{G}^n_1$-graded vector space
\[ U_n := \C \langle u_1, \ldots, u_n \rangle, \]
where we equip $u_j$ with degree $(-1, y_j) \in (\Z \oplus Y_n) /( 2(1-n) \oplus y_{[n]})$.
\end{example2}

\begin{example2}
\label{example2:vecd}
Let $a$ be an integer.
We define the $\bm{G}^n_1$-graded vector space
\[ V^n_a := \C \langle r_1, \ldots, r_n \rangle,\]
where we equip $r_j$ with degree $(2-2a, ay_j) \in (\Z \oplus Y_n)/(2(1-n) \oplus y_{[n]})$.
\end{example2}

\begin{remark2}
\label{remark2:q1p1r}
We observe that, if $\bm{p}_1$ and $\bm{q}_1$ are the morphisms of grading data defined in Lemma \ref{lemma2:squaregrad}, then 
\[ \bm{q}_{1*} \bm{p}_1^* V^n_n \cong V^n_n\]
is a $\Z$-graded vector space concentrated in degree $0$. 
\end{remark2}

We remark that, if $V$ is a $\bm{G}$-graded vector space, then the exterior algebra $\Lambda(V)$ and symmetric algebra $\mathrm{Sym}(V)$ have natural $\bm{G}$-graded algebra structures. 

\begin{definition2}
\label{definition2:2alga}
We define the $\bm{G}^n_1$-graded exterior algebra 
\[ A_n := \Lambda(U_n) \cong \C [ \theta_1, \ldots, \theta_n ],\]
where the variable $\theta_j$ anti-commute.
\end{definition2}

\begin{definition2}
\label{definition2:2ringra}
We define the $\bm{G}^n_1$-graded polynomial ring $\C [V^{n,\vee}_a]$. 
It has a natural filtration by the degree of the polynomial (we'll call this the {\bf order filtration}). 
We define the ring $R^n_a$ to be the completion of $\C[V^{n,\vee}_a]$ with respect to the order filtration, in the category of $\bm{G}^n_1$-graded $\C$-algebras. 
If we took the completion in the category of $\C$-algebras, the result would be the power series ring 
\[ \C \llbracket V_a^{n,\vee} \rrbracket \cong \C \llbracket r_1, \ldots, r_n \rrbracket.\]
Taking the completion in the category of $\bm{G}^n_1$-graded algebras has a different significance: it means that we take the completion of each graded part separately, then take the direct sum of these.
\end{definition2}

The next Lemma looks abstruse, but will be used in Section \ref{subsec:2coh} to relate equivariant matrix factorizations to equivariant coherent sheaves:

\begin{lemma2}
\label{lemma2:gradsquare}
Suppose we are given a commutative diagram of exact morphisms of grading data:
\[\begin{diagram}
\bm{G}_1 & \rTo^{\bm{q}_1} & \bm{G}_{\Z} \\
\dInto^{\bm{p}_1} & & \dInto^{\bm{p}_2} \\ 
\bm{G} & \rTo^{\bm{q}_2} & \bm{G}_2,
\end{diagram}
\]
where $\bm{p}_1$ and $\bm{p}_2$ are injective, and a $\bm{G}$-graded vector space $V$. 
Taking the $X$-part of the grading data and morphisms, we find there are morphisms 
\[ \begin{diagram}
X_1 & \rTo^{p_{1,X}} & X & \rTo^{q_{2,X}} & X_2
\end{diagram}
\]
whose composition is $0$ (by commutativity of the diagram). 
We define the group $\Gamma$ to be the homology of this sequence:
\[ \Gamma := \mathrm{ker}(q_{2,X})/\mathrm{im}(p_{1,X}).\]
Then  $\bm{p}_2^* \bm{q}_{2*} V$ admits a $\Gamma$-grading, and hence an action of the character group $\Gamma^*$, and there is an isomorphism
\[ (\bm{p}_2^* \bm{q}_{2*} V)^{\Gamma^*} \cong \bm{q}_{1*} \bm{p}_1^* V\]
as $\Z$-graded vector spaces, where the superscript $\Gamma^*$ denotes the $\Gamma^*$-invariant part (or equivalently the part of degree $0 \in \Gamma$).
\end{lemma2}
\begin{proof}
We have, by definition, the degree-$j$ parts
\[ \left(\bm{p}_2^* \bm{q}_{2*} V\right)^j := \bigoplus_{q_{2}(y) = p_{2}(j)} V_y,\]
and
\[ \left(\bm{q}_1^*\bm{p}_{1*} V\right)^j := \bigoplus_{q_{1}(y) = j} V_{p_{1}(y)}.\]
Note that the first of these is equal to the part of $V$ whose $Y$-grading lies in $q_2^{-1}(\mathrm{im}(p_2))$, while the second is equal to the part of $V$ whose $Y$-grading lies in $\mathrm{im}(p_1)$ (using the exactness of $\bm{p}_1$ and $\bm{p}_2$). 
By commutativity of the diagram, 
\begin{eqnarray*}
y &=& p_1(y_1) \\
\Rightarrow q_2(y) &=& p_2(q_1(y_1)) \\
\Rightarrow y & \in & q_2^{-1}(\mathrm{im}(q_1)),
\end{eqnarray*}
so $\mathrm{im}(p_1) \subset q_2^{-1}(\mathrm{im}(p_2))$.
Therefore, we have
\[ \bm{q}_1^*\bm{p}_{1*} V \subset \bm{p}_2^*\bm{q}_{2*} V.\]
Furthermore, the left-hand side is exactly equal to the part of the right-hand side whose $Y$-grading lies in
\[ \mathrm{im}(p_1) \subset q_2^{-1}(\mathrm{im}(p_2)),\]
so we can equip the right-hand side with a grading in
\[ q_2^{-1}(\mathrm{im}(p_2))/ \mathrm{im}(p_1) \cong \Gamma,\]
and the left-hand side is equal to the part of degree $0 \in \Gamma$. 

The fact that the $\Z$-gradings match up also follows from commutativity of the diagram:
\[ q_1(y) = j \Rightarrow p_2(j) = q_2(p_1(y)).\]
This completes the proof.
 \end{proof}

\subsection{$\bm{G}$-graded $A_{\infty}$ algebras and Hochschild cohomology}

We now define appropriate notions of $\bm{G}$-graded $A_{\infty}$ algebras and Hochschild cohomology. 
For the purposes of this section, let $\bm{G}$ be an exact grading datum.

\begin{definition2}
\label{definition2:2cc}
Let $R$ be a $\bm{G}$-graded algebra, and let $A,B$ be $\bm{G}$-graded $R$-bimodules.
For each $s \ge 0$, we define a $Y \oplus \Z$-graded $R$-bimodule, whose degree-$(y,s)$ part is
\[ CC^y_c(A,B|R)^s := \mathrm{Hom}_{R-bimod}(A^{\otimes_R s},B)_{y-f(s)},\]
called {\bf compactly supported Hochschild cochains of length $s$ and degree $y$}. 
Note that the $Y$-grading is not quite the obvious one: if $\phi^s$ changes $Y$-degree by $y'$, then we define the $Y$-grading of $\phi^s$ to be $y := f(s) + y'$. 
We will omit the `$|R$' from the notation unless it is necessary to avoid confusion. 
There is a natural filtration of $CC^*_c(A,B)$, called the {\bf length filtration}, given by
\[ F^s CC^*_c(A,B) := \bigoplus_{s' \ge s} CC_c^*(A,B)^{s'}.\]
We define the $\bm{G}$-graded Hochschild cochain complex $CC^*_{\bm{G}}(A,B)$ to be the completion of $CC^*_c(A,B)$ with respect to the length filtration, in the category of $\bm{G}$-graded $R$-bimodules.
Explicitly, the degree-$y$ part is 
\[ CC^y_{\bm{G}}(A,B) := \prod_{s \ge 0} CC_c^y(A,B)^s,\]
and
\[ CC^*_{\bm{G}}(A,B) := \bigoplus_{y \in Y} CC^y_{\bm{G}}(A,B).\]
If $B = A$, we denote
\[ CC^*_{\bm{G}}(A) := CC^*_{\bm{G}}(A,A).\]
Given $\phi \in CC^*_{\bm{G}}(A,B)$, we write  $\phi^s$ for the length-$s$ component of $\phi$. 
We also define a $\Z$-graded version $f^* CC^*_{\bm{G}}(A,B)$ of the Hochschild cohomology (where $f: \Z \To Y$ is the morphism coming from the grading datum $\bm{G}$). 
For $p \in \Z$, we will write $CC^p_{\bm{G}}(A,B)$ for $CC^{f(p)}_{\bm{G}}(A,B)$.
\end{definition2}

\begin{definition2}
We also define a `truncated' version of the $\Z$-graded version of the Hochschild cochain complex: for $p \in \Z$,
\[ TrCC^p_{\bm{G}}(A,B) := \prod_{s \ge p} CC_c^{f(p)}(A,B)^{s} \subset CC^p_{\bm{G}}(A,B).\] 
\end{definition2}

\begin{definition2}
Suppose that $R$ is a $\bm{G}$-graded algebra, and $A$ and $B$ are $\bm{G}$-graded $R$-bimodules. 
The {\bf Gerstenhaber product} is a map of degree $(f(-1),-1) \in Y \oplus \Z$,
\begin{eqnarray*}
CC^*_c(A,B) \otimes_{R} CC^*_c(A) & \To & CC^*_c(A,B), \mbox{ which we denote by} \\
\phi \otimes \psi & \mapsto & \phi \circ \psi \mbox{, and which is defined by} \\
\phi \circ \psi (a_{n},\ldots,a_1) &:=& 
\end{eqnarray*}
\[\sum_{i+j+k = n} (-1)^\dagger \phi^{i+k+1}(a_{i+j+k},\ldots,a_{i+j+1},\psi^j(a_{i+j},\ldots,a_{i+1}),a_i,\ldots,a_1),\]
where
\[\dagger = (\sigma(\psi) + 1)(\sigma(a_1) + \ldots + \sigma(a_i) - i)\]
(recalling Definition \ref{definition2:2z2grad}). 
If the left and right actions of $R$ on the $R$-bimodules $A$ and $B$ coincide, then the Gerstenhaber product is $R$-bilinear; otherwise it is only $R$-linear in $\phi$.
Because the Gerstenhaber product respects the length filtration and the $\bm{G}$-grading, it defines a product
\[ CC^*_{\bm{G}}(A,B) \otimes_R CC^*_{\bm{G}}(A) \To CC^*_{\bm{G}}(A,B)\]
of degree $f(-1)$, also called the Gerstenhaber product.
\end{definition2}

\begin{definition2}
If $R$ is a $\bm{G}$-graded algebra and $A$ is a $\bm{G}$-graded $R$-bimodule, then we define the {\bf Gerstenhaber bracket}, which is a Lie bracket of degree $f(-1)$ on $CC^*_{\bm{G}}(A)$, by
\[ [\phi, \psi] := \phi \circ \psi - (-1)^{(\sigma(\phi) + 1)(\sigma(\psi) + 1)}\psi \circ \phi.\]
\end{definition2}

\begin{definition2}
If $R$ is a $\bm{G}$-graded algebra, then a {\bf $\bm{G}$-graded associative $R$-algebra} is a $\bm{G}$-graded $R$-bimodule $A$, together with an element 
\[ \mu^2 \in CC_{c}^2(A)^2\]
satisfying the associativity relation
\[ \mu^2 \circ \mu^2 = 0.\]
\end{definition2}

\begin{remark2}
If $A$ is a $\bm{G}$-graded associative $R$-algebra, then the product
\[a_2 \cdot a_1 := (-1)^{\sigma(a_1)} \mu^2(a_2,a_1)\]
is associative and respects the $Y$-grading, and makes $A$ into a $Y$-graded associative $R$-algebra in the usual sense. 
\end{remark2}

\begin{definition2}
If $(A,\mu^2)$ is a $\bm{G}$-graded associative $R$-algebra, then we define the {\bf Hochschild differential}
\begin{eqnarray*}
\delta: CC^*_{\bm{G}}(A) &\To& CC^*_{\bm{G}}(A) \\
\delta (\tau) &:=&\left[\mu^2, \tau \right].
\end{eqnarray*}
It has degree $f(1) \in Y$ and increases length by $1$ (i.e., it is induced by a similar differential of degree $(f(1),1)$ on $CC^*_c(A)$).
It follows from the fact that $\mu^2 \circ \mu^2 = 0$ that $\delta$ is a differential, i.e., $\delta^2 = 0$. 
We define the {\bf Hochschild cohomology} of $A$ to be its cohomology,
\[HH^*_{\bm{G}}(A) := H^*(CC^*_{\bm{G}}(A),\delta),\]
which is a $\bm{G}$-graded $R$-bimodule (as $\delta$ has pure degree in $Y$).
Furthermore, because $\delta$ is pure of degree $(f(1),1)$ on $CC^*_c(A)$, we can also define the compactly-supported version $HH^*_c(A)$, which is $Y \oplus \Z$-graded.
\end{definition2}

\begin{definition2}
\label{definition2:2rainfalg}
If $R$ is a $\bm{G}$-graded algebra, then a {\bf $\bm{G}$-graded $A_{\infty}$ algebra} over $R$ is a $\bm{G}$-graded $R$-bimodule $A$, together with an element
\[ \mu \in CC_{\bm{G}}^2(A),\]
satisfying $\mu^0 = 0$, and such that the {\bf $A_{\infty}$ associativity relation}
\[ \mu \circ \mu = 0\]
is satisfied. 
We denote $A_{\infty}$ algebras by $\mathcal{A} := (A,\mu^*)$. 
If $\mu^1 =0$ (or equivalently, if $\mu$ sits inside $TrCC_{\bm{G}}^2(A)$), we say that $\mathcal{A}$ is {\bf minimal}.
If we have $\mu^*$ such that $\mu^0 \neq 0$, but still $\mu \circ \mu = 0$, then $(A,\mu^*)$ is called a {\bf curved} $A_{\infty}$ algebra.
\end{definition2}

\begin{remark2}
\label{remark2:ainfgrad}
In other words, $A$ is a $Y$-graded $R$-bimodule, equipped with $R$-multilinear maps
\[ \mu^s: A^{\otimes_R s} \To A\]
of degree $f(2-s) \in Y$, for all $s \ge 1$, satisfying the $A_{\infty}$ associativity relations.
\end{remark2}

\begin{definition2}
If $A$ and $B$ are $\bm{G}$-graded $R$-bimodules, we define a new `product' 
\begin{eqnarray*}
CC^*_{\bm{G}}(B) \otimes_{R} TrCC^1_{\bm{G}}(A,B) &\To & CC^*_{\bm{G}}(A,B) \mbox{, denoted by}\\
\phi \otimes \psi &\mapsto& \phi \diamond \psi \mbox{, and which we define by} \\
 (\phi \diamond \psi)^n(a_n,\ldots,a_1) &:=&
\end{eqnarray*} 
\[\sum_{i_1+\ldots+i_j = n} \phi^j(\psi^{i_1}(a_{i_1+\ldots+i_j},\ldots),\psi^{i_2}(a_{i_2+\ldots+i_j},\ldots),\ldots,\psi^{i_j}(a_{i_j},\ldots,a_1)).\]
It is $R$-linear only in the first variable $\phi$. 
Note that $\diamond$ has degree $0 \in Y$, and is clearly associative: $(F \diamond G) \diamond H = F \diamond (G \diamond H)$.
\end{definition2}

\begin{definition2}
\label{definition2:2amorph}
If $\mathcal{A} = (A,\mu)$ and $\mathcal{B} = (B,\eta)$ are $\bm{G}$-graded $A_{\infty}$ algebras over $R$, then an $A_{\infty}$ morphism from $\mathcal{A}$ to $\mathcal{B}$ is an element $F \in TrCC^1_{\bm{G}}(A,B)$ such that
\[ F \circ \mu - \eta \diamond F = 0 \in CC_{\bm{G}}^2(A,B).\]
Composition of two $A_{\infty}$ morphisms is defined using the product $\diamond$. 
$F$ is called {\bf strict} if $F^j = 0$ for all $j \ge 2$. 
An {\bf $A_\infty$ morphism} $F$ from $\mathcal{A}$ to $\mathcal{B}$ induces a homomorphism of graded associative algebras on the level of cohomology, which we denote by
\[ H^*(F): H^*(\mathcal{A}) \To H^*(\mathcal{B}).\]
\end{definition2}

\begin{definition2}
\label{definition2:quasiisomor}
If an $A_{\infty}$ morphism induces an isomorphism on the level of cohomology, then it is said to be a {\bf quasi-isomorphism}. 
\end{definition2}

In fact, when our $A_{\infty}$ algebras are minimal, there is an easier notion, that of {\bf formal diffeomorphism} (see \cite[Section 1c]{Seidel2008}).

\begin{definition2}
\label{definition2:2formdiff}
If $R$ is a $\bm{G}$-graded $\C$-algebra, and $A$ and $B$ are $\bm{G}$-graded $R$-modules, then a {\bf $\bm{G}$-graded formal diffeomorphism} from $A$ to $B$ is an element
\[ F \in TrCC^1_{\bm{G}}(A,B)\]
such that
\[ F^1: A \To B \]
is an isomorphism of $R$-modules. 
Formal diffeomorphisms can be composed using $\diamond$: If $F$ is a $\bm{G}$-graded formal diffeomorphism from $A$ to $B$ and $G$ is a $\bm{G}$-graded formal diffeomorphism from $B$ to $C$, then $G \diamond F$ is a $\bm{G}$-graded formal diffeomorphism from $A$ to $C$. 
\end{definition2}

\begin{lemma2}
\label{lemma2:formdiffinv}
$\bm{G}$-graded formal diffeomorphisms can be strictly inverted: if $F$ is a $\bm{G}$-graded formal diffeomorphism from $A $ to $B $, then there exists a unique $\bm{G}$-graded formal diffeomorphism $F^{-1}$ from $B$ to $A$ such that 
\[F^{-1} \diamond F = \mathrm{id} \mbox{ and } F \diamond F^{-1} = \mathrm{id},\]
where `$\mathrm{id}$' denotes the formal diffeomorphism from $A$ to itself, given by
\[ (\mathrm{id})^s = \left\{ \begin{array}{ll}
					\mathrm{id}:A \To A & \mbox{ if $s = 1$,} \\
					0 & \mbox{otherwise}.
				\end{array} \right.\]
\end{lemma2}
\begin{proof}
We construct a left inverse $G$ for $F$, inductively in the length filtration: $G^1 = (F^1)^{-1}$, and if $G$ is determined to length $\le s-1$, then at order $s$ we have
\[ G^s(F^1(a_s), F^1(a_{s-1}) , \ldots, F^1(a_1)) = -\sum_{s' \le s-1}G^{s'}(F(\ldots), \ldots, F(\ldots)), \]
and since $F^1$ is an isomorphism, this determines $G^s$ uniquely. 
One can similarly prove that $F$ has a unique strict right inverse. 
Because $\diamond$ is associative, the left and right inverses coincide. 
 \end{proof}

\begin{definition2}
\label{definition2:fstar}
If $F \in TrCC^1_{\bm{G}}(A,B)$ is a formal diffeomorphism, we define a $\bm{G}$-graded, $R$-linear map
\begin{eqnarray*}
 F_*: CC_{\bm{G}}(A) & \To & CC_{\bm{G}}(B), \\
F_* \eta &:=& (F \circ \eta) \diamond F^{-1}.
\end{eqnarray*}
\end{definition2}

\begin{lemma2}
\label{lemma2:presgerst}
$F_*$ preserves the Gerstenhaber bracket: 
\[ \left[ F_* \alpha, F_* \beta \right] = F_* [ \alpha, \beta].\]
\end{lemma2}
\begin{proof}
First, if $F$ is a formal diffeomorphism, then some simple algebraic manipulation yields 
\[ (G \diamond F) \circ H = (G \circ F_*H) \diamond F.\]
It follows that
\begin{eqnarray*}
F_* \alpha \circ F_* \beta &=& \left( (F \circ \alpha) \diamond F^{-1} \right) \circ F_* \beta \\
&=& \left( (F \circ \alpha) \circ F^{-1}_* F_* \beta \right) \diamond F^{-1} \\
&=& \left((F \circ \alpha) \circ \beta \right) \diamond F^{-1}.
\end{eqnarray*}
Hence,
\begin{eqnarray*}
\left[ F_* \alpha, F_* \beta \right] &=& \left( (F \circ \alpha) \circ \beta - (-1)^{\dagger} (F \circ \beta) \circ \alpha \right) \diamond F^{-1} \\
&=& \left(F \circ \left( \alpha \circ \beta - (-1)^{\dagger} \beta \circ \alpha \right)\right) \diamond F^{-1} \\
&=& F_* [\alpha,\beta],
\end{eqnarray*}
where $\dagger = (\sigma(\alpha)+1)(\sigma(\beta)+1)$.
 \end{proof}

\begin{corollary2}
\label{corollary2:formdiffpush}
If $(A,\mu)$ is a minimal $\bm{G}$-graded $A_{\infty}$ algebra over $R$, and $F \in TrCC^1_{\bm{G}}(A,B)$ is a formal diffeomorphism from $A$ to $B$, then $F_* \mu$ is a $\bm{G}$-graded minimal $A_{\infty}$ structure on $B$, and $F$ defines an $A_{\infty}$ quasi-isomorphism from $(A,\mu)$ to $(B,F_* \mu)$.
\end{corollary2}
\begin{proof}
It follows from Lemma \ref{lemma2:presgerst} that $[F_*\mu,F_*\mu] = F_*[\mu,\mu] = 0$, so $F_*\mu$ is an $A_\infty$ structure on $B$.
By construction,
\[ F \circ \mu = (F_* \mu) \diamond F,\]
so $F$ defines an $A_{\infty}$ morphism from $(A,\mu)$ to $(B,F_*\mu)$; because $F^1$ is an isomorphism, $F$ is a quasi-isomorphism.
 \end{proof}

\begin{remark2}
\label{remark2:formdiffgp}
We observe that
\[ (G \diamond F)_* \mu = G_*(F_* \mu).\]
\end{remark2}

\begin{definition2}
\label{definition2:2gact}
If $A$ is a $\bm{G}$-graded vector space over $\C$, then there is a group
\[ \mathfrak{G}(A) := \{F \in TrCC^1_{\bm{G}}(A): F^1 = \mathrm{id}\},\]
the {$\bm{G}$-graded formal diffeomorphisms} from $A$ to itself whose leading term is the identity. 
It follows from Lemma \ref{lemma2:formdiffinv}, Corollary \ref{corollary2:formdiffpush} and Remark \ref{remark2:formdiffgp}  that $\mathfrak{G}(A)$ is a group, and that it acts on the space of minimal $\bm{G}$-graded $A_{\infty}$ structures on $A$. 
Note that this action preserves the underlying algebra $(A,\mu^2)$, and that $F$ defines an $A_{\infty}$ quasi-isomorphism from $\mu$ to $F_* \mu$. 
\end{definition2}

\begin{definition2}
If we are given a $\bm{G}$-graded associative algebra $A$ over $\C$, we define $\mathfrak{A}(A)$, the set of $\bm{G}$-graded minimal $A_{\infty}$ algebras $\mathcal{A} = (A,\mu)$ over $\C$, with $\mu^2$ coinciding with the product on $A$.
\end{definition2}

\begin{definition2}
Suppose that $\mathcal{A} = (A,\mu) \in \mathfrak{A}(A)$, and $\mu^s = 0$ for $2<s<d$. 
Then the $A_{\infty}$ associativity relations $\mu \circ \mu = 0$ imply that $\mu^d$ is a Hochschild cocycle for $(A,\mu^2)$. 
The class
\[ \left[ \mu^d \right] \in HH_{c,\bm{G}}^2(A)^{d}\]
is called the {\bf order-$d$ deforming class} of $\mathcal{A}$.
\end{definition2}

\begin{remark2}
In \cite{Seidel2003} and \cite{Sheridan2011a}, the class $[\mu^d]$ is called an order-$d$ {\bf deformation class}. 
However a large part of this paper is devoted to studying different objects (see Definition \ref{definition2:2defclass}), also elements of Hochschild cohomology, and also called deformation classes in \cite{Seidel2003}. 
That is why, to avoid confusing the reader, and only for the purposes of the current paper, we use the terminology `deforming class' to distinguish this object.
\end{remark2}

We recall a versality result from \cite{Seidel2008}, appropriately modified to take into account $\bm{G}$-grading:

\begin{proposition2}
\label{proposition2:2vers1}
Suppose that $A$ is a $\bm{G}$-graded associative algebra over $\C$, and there exists $d>2$ such that
\[HH_{c,\bm{G}}^2(A)^{s} \cong \left\{\begin{array}{ll}
                                                    \C & \mbox{for $s=d$}\\
                                                    0 & \mbox{ for all $s>2$, $s \neq d.$}
                                              \end{array}\right.\]
Suppose that $\mathcal{A}_1 = (A,\mu_1)$ and $\mathcal{A}_2 = (A,\mu_2)$ both lie in $\mathfrak{A}(A)$, satisfy $\mu_1^s = \mu_2^s = 0$ for all $2<s<d$, and have non-trivial order-$d$ deforming class in $HH_{c}^2(A)^{d}$.
Then $\mathcal{A}_1$ and $\mathcal{A}_2$ are related by a formal diffeomorphism. 
\end{proposition2}
\begin{proof}
The proof is by a straightforward order-by-order construction of a formal diffeomorphism $F$ such that $F_* \mu_1 = \mu_2$, showing that all obstructions to the existence of $F$ vanish (see \cite[Lemma 3.2]{Seidel2003}).
 \end{proof}

Now we define Hochschild cohomology of an $A_{\infty}$ algebra:

\begin{definition2}
\label{definition2:2hhainf}
Suppose that $\mathcal{A}=(A,\mu)$ is a $\bm{G}$-graded $A_{\infty}$ algebra.
We define the Hochschild differential
\begin{eqnarray*}
\delta: CC^*_{\bm{G}}(A) &\To &CC^*_{\bm{G}}(A),\\
\delta(\tau) &:=& \left[ \mu, \tau \right].
\end{eqnarray*}
It follows from the fact that $\mu \circ \mu = 0$, and that the Gerstenhaber bracket satisfies (a version of) the Jacobi relation, that $\delta$ is a differential, i.e., $\delta^2 = 0$. 
We define the Hochschild cohomology of $\mathcal{A}$ to be the cohomology,
\[ HH^*_{\bm{G}}(\mathcal{A}) := H^*(CC^*_{\bm{G}}(A),\delta).\]
The Hochschild differential has pure degree $f(1)$, so $HH^*_{\bm{G}}(\mathcal{A})$ is $\bm{G}$-graded. 
However, note that the Hochschild differential is no longer pure with respect to the length, so we can not define the $Y \oplus \Z$-graded compactly-supported version, as we could for an associative algebra. 
However, $\delta$ does always increase (or preserve) the length; therefore, we do still have the length filtration on the Hochschild cochain complex.
\end{definition2}

\begin{definition2}
If $\mathcal{A} = (A,\mu)$ is a $\bm{G}$-graded {\bf minimal} $A_{\infty}$ algebra (i.e., $\mu \in TrCC_{\bm{G}}^2(A)$), then the Hochschild differential preserves the truncated Hochschild cochains. 
Thus it makes sense to define the {\bf truncated Hoch\-schild complex} $(TrCC_{\bm{G}}^*(\mathcal{A}),\delta)$, and call its cohomology the {\bf truncated Hoch\-schild cohomology} $TrHH^*_{\bm{G}}(\mathcal{A})$. 
Recall that it is $\Z$-graded.
\end{definition2}

\begin{remark2}
\label{remark2:ss}
Suppose that $A$ is a $\bm{G}$-graded associative algebra and $\mathcal{A} \in \mathfrak{A}(A)$.
The length filtration on the Hochschild cochain complex $CC^*_{\bm{G}}(\mathcal{A})$ yields a `$\bm{G}$-graded spectral sequence' $(E_d^{*,*},\delta_d^{*,*})$. 
The $E_d$ page is $Y \oplus \Z$ graded, and $\delta_d$ has degree $(f(1), d)$. 
The spectral sequence starts on page $E_1$, where it is given by
\[ E_1^* = CC^*_{c,\bm{G}}(A),\]
with differential $\delta_1 = [\mu^2, -]$. 
The cohomology of this differential is, by definition, the Hochschild cohomology of the algebra $A$:
\[E_2^* = HH_{c,\bm{G}}^*(A).\]
When we need to prove that the spectral sequence converges to $HH^*_{\bm{G}}(\mathcal{A})$, we will apply the `complete convergence theorem' \cite[Theorem 5.5.10]{Weibel1994}, in the abelian category of $\bm{G}$-graded modules.
The length filtration is clearly bounded above, because all Hochschild cochains have length $s \ge 0$, and hence it is also exhaustive. 
It is also complete in the category of $\bm{G}$-graded $R$-bimodules, because the Hochschild cochain complex is defined as a completion with respect to the length filtration. 
Therefore, to prove that the spectral sequence converges to $HH^*_{\bm{G}}(\mathcal{A})$, we must show that it is {\bf regular}:  for each $(y,s) \in Y \oplus \Z$, the differentials
\[ \delta_d: E_d^{(y,s)} \To E_d^{(y + f(1), s + d)}\]
vanish for sufficiently large $d$.
When we need to prove that the spectral sequence converges, we will in fact show that the differentials $\delta_d$ vanish whenever $d$ is sufficiently large (independent of $y,s$).
\end{remark2}

\begin{remark2}
\label{remark2:ssfirst}
Suppose that $A$ is a $\bm{G}$-graded associative algebra over $\C$, and $\mathcal{A} \in \mathfrak{A}(A)$ has $\mu^s = 0$ for $2 < s < d$, and order-$d$ deforming class $[\mu^d]$. 
Then the first non-zero differential in the spectral sequence is $\delta_{d-1}$, and it is given by
\begin{eqnarray*}
\delta_{d-1}: HH_{c,\bm{G}}^{(y,s)}(A) & \To & HH_{c,\bm{G}}^{(y+f(1),s+d-1)}(A), \\
\delta_{d-1}(\phi) & = & \left[\left[\mu^d \right], \phi \right],
\end{eqnarray*}
where we observe that the Gerstenhaber bracket $[-,-]$ descends to the cohomology. 
\end{remark2}

Now we will define the appropriate notions of $\bm{G}$-graded $A_{\infty}$ categories and their Hochschild cohomology.

\begin{definition2}
\label{definition2:gprecat}
Let $\bm{G} = \{\Z \To Y\}$ be a grading datum, and $R$ a $\bm{G}$-graded algebra. 
A {\bf $\bm{G}$-graded pre-category} $\mathcal{C}$ over $R$ is a set of objects $Ob(\mathcal{C})$, together with morphism spaces
\[ hom^*_{\mathcal{C}}(K,L)\]
which are $\bm{G}$-graded $R$-bimodules, and an action of $Y$ on $Ob(\mathcal{C})$ by `shifts' $[y]$, together with isomorphisms
\[ S_{y_1,y_2}:   hom^*_{\mathcal{C}}(K,L)  \overset{\cong}{\To}  hom^*_{\mathcal{C}}(K[y_1],L[y_2])[y_1-y_2]\]
which are compatible in the obvious way. 
We can think of this as equipping $\mathcal{C}$ with an $R[Y]$-bimodule structure.
We define the $\bm{G}$-graded group $CC^*_{\bm{G}}(\mathcal{C})$, by analogy with Definition \ref{definition2:2cc}, restricting to the parts that respect the above isomorphisms. 
We can think of this strict equivariance requirement as taking $CC^*_{\bm{G}}(\mathcal{C}|R[Y])$, where $\mathcal{C}$ acquires an $R[Y]$-linear structure from the $Y$-action. 
Explicitly, it means that for $\phi^s \in CC^s_{\bm{G}} (\mathcal{C})$, and any $y_0, \ldots, y_s \in Y$,
\[ \phi^s \circ ( S_{y_{s-1},y_s} \otimes \ldots \otimes S_{y_0,y_1}) = S_{y_0,y_s} \circ \phi^s,\]
as maps
\[ hom^*(L_{s-1},L_s) \otimes \ldots \otimes hom^*(L_0,L_1) \To hom^*(L_0[y_0],L_s[y_s])[y_0-y_s].\]
\end{definition2}

\begin{definition2}
We define a {\bf $\bm{G}$-graded $A_{\infty}$ category} $\mathcal{C}$ over $R$ to be a $\bm{G}$-graded pre-category together with
\[ \mu \in CC^2_{\bm{G}}(\mathcal{C})\]
satisfying $\mu^0 = 0$ and $\mu \circ \mu = 0$. 
An $A_{\infty}$ category is said to be {\bf cohomologically unital} if its cohomological category is unital.
We define the Hochschild cohomology $HH^*_{\bm{G}}(\mathcal{C})$ by analogy with Definition \ref{definition2:2hhainf}. 
We also consider the case where $\mu^0 \neq 0$; in this case we say that $\mathcal{C}$ is {\bf curved}. 
\end{definition2}

\begin{definition2}
We define $\bm{G}$-graded $A_{\infty}$ functors by analogy with Definition \ref{definition2:2amorph}. 
An $A_\infty$ functor induces an ordinary functor on the level of cohomology; if this functor is an equivalence, we call the functor a {\bf quasi-equivalence}. 
\end{definition2}

\begin{definition2}
\label{definition2:2minfun}
A $\bm{G}$-graded $A_{\infty}$ category is said to be {\bf minimal} if $\mu$ lies in $TrCC^2_{\bm{G}}(\mathcal{C})$. 
\end{definition2}

\begin{remark2}
\label{remark2:minimality}
The notions of unitality and equivalence for minimal $A_{\infty}$ categories are simpler than for non-minimal categories. 
Because there is no differential $\mu^1$ on the morphisms spaces, $\mu^2$ is strictly associative and therefore defines a category. 
Thus, if $(\mathcal{C},\mu)$ is minimal and cohomologically unital, then $(\mathcal{C},\mu^2)$ is a category, in particular is unital. 
We say that two objects of $\mathcal{C}$ are {\bf quasi-isomorphic} if they are quasi-isomorphic as objects of $(\mathcal{C},\mu^2)$. 
An $A_{\infty}$ functor $\mathscr{F}:(\mathcal{C},\mu) \To (\mathcal{D},\eta)$ between minimal $A_{\infty}$ categories is a quasi-equivalence if and only if the functor $\mathscr{F}^1: (\mathcal{C},\mu^2) \To (\mathcal{D},\eta^2)$ is an equivalence. 
\end{remark2}

\begin{lemma2}
\label{lemma2:mininv}
If $\mathscr{F}:(\mathcal{C},\mu) \To (\mathcal{D},\eta)$ is an $A_{\infty}$ quasi-equivalence between minimal $A_{\infty}$ categories, then there exists an $A_{\infty}$ functor $\mathscr{G}:(\mathcal{D},\eta) \To (\mathcal{C},\mu)$ such that $\mathscr{F}^1$ and $\mathscr{G}^1$ are mutually (strictly) inverse quasi-equivalences.
\end{lemma2}
\begin{proof}
Follows from Lemma \ref{lemma2:formdiffinv} (see \cite[Theorem 2.9]{Seidel2008}).
 \end{proof}

\begin{remark2}
\label{remark2:invertquasi}
If the coefficient ring $R$ is a field, then Lemma \ref{lemma2:mininv} is true even if the categories $\mathcal{C}$ and $\mathcal{D}$ are not minimal (see \cite[Theorem 2.9]{Seidel2008}). 
However, if the ring $R$ is not a field, then the hypothesis of minimality in Lemma \ref{lemma2:mininv} is crucial. 
\end{remark2}

\begin{definition2}
\label{definition2:2algtocat}
If $\mathcal{A}$ is a $\bm{G}$-graded $A_{\infty}$ algebra, then we denote by $\underline{\mathcal{A}}$ the smallest $\bm{G}$-graded $A_{\infty}$ category with an object whose endomorphism algebra is $\mathcal{A}$. 
Namely, $\underline{\mathcal{A}}$ has objects $K[y]$ indexed by $Y$, and morphism spaces
\[ hom^*(K[y_1],K[y_2]) := \mathcal{A}[y_2 - y_1].\]
By definition, there is an isomorphism
\[CC^*_{\bm{G}}\left(\mathcal{A}\right) \cong CC^*_{\bm{G}}\left(\underline{\mathcal{A}}\right).\]
We define the $A_{\infty}$ structure maps $\mu^*$ on $\underline{\mathcal{A}}$ to be the image of those on $\mathcal{A}$ under this isomorphism.
\end{definition2}

Now we explain how $\bm{G}$-graded $A_{\infty}$ categories can be `pulled back' along injective morphisms of grading data, and how this operation affects the Hoch\-schild cohomology.

\begin{definition2}
\label{definition2:2prepullainf}
Let $\bm{p}:\bm{G}_1 \hookrightarrow \bm{G}_2$ be an injective morphism of grading data, $R$ a $\bm{G}_1$-graded algebra, and $\mathcal{C}$ a $\bm{G}_2$-graded pre-category over $\bm{p}_*R$. 
We define $\bm{p}^*\mathcal{C}$, a $\bm{G}_1$-graded pre-category over $R$, to have the same objects as $\mathcal{C}$, but $\bm{G}_1$-graded morphism spaces
\[ hom^*_{\bm{p}^* \mathcal{C}}(K,L) := \bm{p}^* hom^*_{\mathcal{C}}(K,L)\]
(so it is not necessarily a full sub-pre-category).
We note that $\bm{p}^*\mathcal{C}$ still has an action of $Y_2$ on objects by shifts, so the subgroup $Y_1 \subset Y_2$ acts, equipping $\bm{p}^*\mathcal{C}$ with the structure of a $\bm{G}_1$-graded pre-category over $R$.
\end{definition2}

\begin{remark2}
\label{remark2:pullcc}
Because $Y_2$ acts on $\bm{p}^*\mathcal{C}$ by shifts (shifting all objects simultaneously by the same $y \in Y_2$), it acts on $CC^*_{\bm{G}_1}(\bm{p}^*\mathcal{C})$. 
However the action of the subgroup $Y_1$ is trivial by definition, since we restrict to Hochschild cochains that respect the shifts by $Y_1$. 
So there is an action of the group $\Gamma := Y_2/Y_1$ on $CC^*_{\bm{G}_1}(\bm{p}^*\mathcal{C})$. 
It is not hard to see that the $\Gamma$-fixed part is isomorphic to $\bm{p}^*CC^*_{\bm{G}_2}(\mathcal{C})$. 
Thus, we have 
\[ \bm{p}^*CC^*_{\bm{G}_2}(\mathcal{C}) \cong CC^*_{\bm{G}_1}(\bm{p}^* \mathcal{C})^{\Gamma},\]
and it follows that, for any integer $j$,
\[ CC^j_{\bm{G}_2}(\mathcal{C}) \cong CC^j_{\bm{G}_1}(\bm{p}^*\mathcal{C})^{\Gamma} \subset CC^j_{\bm{G}_1}(\bm{p}^*\mathcal{C}).\]
\end{remark2}

Thus we can make the following:

\begin{definition2}
\label{definition2:2pullainf}
If $\mathcal{C}$ is a $\bm{G}_2$-graded $A_{\infty}$ category over $\bm{p}_* R$, with structure maps $\mu^*$, then we define $\bm{p}^*\mathcal{C}$, a $\bm{G}_1$-graded $A_{\infty}$ category over $R$, whose $A_{\infty}$ structure maps are given by the image of $\mu^*$ under the inclusion
\[ CC^2_{\bm{G}_2}(\mathcal{C}) \cong CC^2_{\bm{G}_1}(\bm{p}^*\mathcal{C})^{\Gamma} \subset CC^2_{\bm{G}_1}(\bm{p}^*\mathcal{C}).\]
\end{definition2}

\begin{remark2}
\label{remark2:pullhh}
It follows that there is an isomorphism
\[ \bm{p}^* HH^*_{\bm{G}_2}(\mathcal{C}) \cong HH^*_{\bm{G}_1}(\bm{p}^*\mathcal{C})^{\Gamma} \subset  HH^*_{\bm{G}_1}(\bm{p}^*\mathcal{C}).\]
\end{remark2}

Now we explain how a $\bm{G}$-graded $A_{\infty}$ category can be `pushed forward' along a surjective morphism of grading data.

\begin{definition2}
Let $\bm{p}: \bm{G}_1 \To \bm{G}_2$ be a surjective morphism of grading data, with kernel $\Gamma \subset Y_1$, and $\mathcal{C}$ a $\bm{G}_1$-graded $A_{\infty}$ category over a $\bm{G}_1$-graded algebra $R$. 
We now define a $\bm{G}_2$-graded pre-category $\bm{p}_* \mathcal{C}$ over $\bm{p}_* R$, as follows: First, observe that there are canonical isomorphisms of $\bm{G}_2$-graded vector spaces
\[ \bm{p}_* hom^*_{\mathcal{C}}(K[y],L) \cong \bm{p}_* hom^*_{\mathcal{C}}(K,L)\]
for any $y \in \Gamma$. 
We define the set of objects of $\bm{p}_* \mathcal{C}$ to be the quotient of the set of objects of $\mathcal{C}$ by the action of $\Gamma \subset Y_1$. 
We define the $\bm{G}_2$-graded morphism spaces to be
\[ hom^*_{\bm{p}_* \mathcal{C}}(K,L) := \bm{p}_* hom^*_{\mathcal{C}}(K,L).\]
This is well-defined by our previous remark. 
Furthermore, because $\bm{p}$ is surjective, there is an obvious action of $Y_2$ on $\bm{p}_* \mathcal{C}$ by shifts, so $\bm{p}_* \mathcal{C}$ is a $\bm{G}_2$-graded pre-category. 
In some sense we have
\[ \bm{p}_* \mathcal{C} = \mathcal{C} \otimes_{R[Y_1]} R[Y_2].\]
\end{definition2}

\begin{remark2}
\label{remark2:pushcat}
Note that
\[CC^y_{c,\bm{G}_2}(\bm{p}_* \mathcal{C}) \cong \left(\bm{p}_* CC^*_{c,\bm{G}_1}(\mathcal{C})\right)^y.\]
It follows that $CC^y_{\bm{G}_2}(\bm{p}_* \mathcal{C})$ is just the completion of $\left(\bm{p}_* CC^*_{\bm{G}_1}(\mathcal{C})\right)^y$ with respect to the length filtration, for all $y \in Y_2$ (observe that the completion is only needed when $\bm{p}$ has an infinite kernel).
It follows that there is an inclusion
\[ \bm{p}_*CC^*_{\bm{G}_1}(\mathcal{C}) \hookrightarrow CC^*_{\bm{G}_2}(\bm{p}_* \mathcal{C}).\]
\end{remark2}

\begin{definition2}
\label{definition2:2pushainfcat}
If $\mathcal{C}$ is a $\bm{G}_1$-graded $A_{\infty}$ category over $R$, with structure maps $\mu^*$, and $\bm{p}: \bm{G}_1 \To \bm{G}_2$ a surjective morphism of grading data, then we define $\bm{p}_*\mathcal{C}$, a $\bm{G}_2$-graded $A_{\infty}$ category over $\bm{p}_*R$, whose $A_{\infty}$ structure maps are given by the image of $\mu^*$ under the inclusion
\[ CC^2_{\bm{G}_1}(\mathcal{C}) \hookrightarrow CC^2_{\bm{G}_2}(\bm{p}_*\mathcal{C}).\]
\end{definition2}

\subsection{Deformations of $A_{\infty}$ algebras}
\label{subsec:2defainf}

For the purposes of this section, let us fix a grading datum $\bm{G}$ and a $\bm{G}$-graded vector space $V$. 
Let $R$ be the $\bm{G}$-graded ring which is the completion of the $\bm{G}$-graded polynomial ring $\C[V]$ with respect to the order filtration, in the category of $\bm{G}$-graded rings. 
We denote the order-$j$ part of $R$ by $R^j$.
There is a natural projection $R \To R^0 \cong \C$, given by setting all $r_j = 0$, and we denote the kernel of this projection by $\mathfrak{m} \subset R$. 

We also denote by $R_0$ the part of $R$ of degree $0 \in Y$. 

\begin{definition2}
Given $\psi \in R_0$, there is a $\bm{G}$-graded algebra homomorphism
\begin{eqnarray*} 
\psi^*: R & \To & R, \\
\psi^*(r^j) & := & \psi^j \cdot r^j, \mbox{ where $r^j \in R^j$ has order $j$.}
\end{eqnarray*}
We now define a group, which by abuse of notation we call $\mathrm{Aut}(R)$, by setting
\begin{eqnarray*}
\mathrm{Aut}(R) &:=& \{ \psi \in R_0: \psi(0) \neq 0 \}, \\
\phi \cdot \psi & := & \phi^*(\psi) \phi,
\end{eqnarray*}
so that
\[ (\phi \cdot \psi)^* = \phi^* \circ \psi^*.\]
Thus, we have an action of $\mathrm{Aut}(R)$ on $R$ by $\bm{G}$-graded algebra isomorphisms.
Note that the condition $\psi(0) \neq 0$ ensures that $\psi$ has a unique inverse in $\mathrm{Aut}(R)$, which can be constructed order-by-order.
\end{definition2}

\begin{definition2}
If $A$ is a $\bm{G}$-graded vector space, then $A \otimes R$ is a $\bm{G}$-graded $R$-bimodule, and we have an isomorphism
\[ CC^*_{\bm{G}}(A \otimes R, A \otimes R|R) \cong CC^*_{\bm{G}}(A,A \otimes R|\C).\]
If we have
\begin{eqnarray*}
\phi &\in& CC^*_{\bm{G}}(A,A \otimes R),\\
\phi &=& \sum_{j \ge 0} \phi_j, \mbox{ where } \\
\phi_j & \in & CC^*(A, A \otimes R^j),
\end{eqnarray*}
then we call $\phi_j$ the {\bf order-$j$ component} of $\phi$.
\end{definition2}

\begin{definition2}
\label{definition2:2graddef}
Let $\mathcal{A} = (A,\mu_0)$ be a $\bm{G}$-graded minimal $A_{\infty}$ algebra over $\C$. 
A {\bf $\bm{G}$-graded deformation of $\mathcal{A}$ over $R$} is an element
\[ \mu \in CC_{\bm{G}}^2(A \otimes R|R) \cong CC_{\bm{G}}^2(A,A \otimes R|\C)\]
that makes $A \otimes R$ into an $A_{\infty}$ algebra over $R$ (i.e., $\mu \circ \mu = 0$), and whose order-$0$ component is $\mu_0$. 
If $\mu \in TrCC_{\bm{G}}^2(A,A \otimes R)$, then the deformation is said to be {\bf minimal}. 
\end{definition2}

\begin{remark2}
\label{remark2:r00mod}
Observe that $CC_{\bm{G}}^2(A,A \otimes R)$ and $TrCC_{\bm{G}}^2(A,A \otimes R)$ are $R_0$-modules in the obvious way. 
\end{remark2}

\begin{definition2}
The group $\mathrm{Aut}(R)$ acts on $A \otimes R$, where $\psi$ acts by the automorphism
\[ \psi^* := \mathrm{id} \otimes \psi^*: A \otimes R \To A \otimes R.\]
Hence it acts on $CC^*_{\bm{G}}(A \otimes R|R)$, where $\psi$ acts by 
\[ (\psi^* \cdot \phi)(a_s, \ldots, a_1) := \psi^* \phi \left( (\psi^*)^{-1}a_s, \ldots, (\psi^*)^{-1} a_1 \right).\]
It is not difficult to see that this action preserves the Gerstenhaber product, so if $\mu$ is one $\bm{G}$-graded deformation of $\mathcal{A}$ over $R$, then $\psi^* \cdot \mu$ is another. 
We say that $\psi^* \cdot \mu$ is obtained from $\mu$ by the {\bf formal change of variables} $\psi^*$ (compare Definition \ref{definition2:formchang}). 
\end{definition2}

Now suppose that $(A,\mu_0)$ and $(B,\eta_0)$ are $\bm{G}$-graded $A_{\infty}$ algebras over $\C$, and $(A,\mu)$ and $(B,\eta)$ are $\bm{G}$-graded $A_{\infty}$ deformations of these over $R$. 
We recall (from Definition \ref{definition2:2amorph}) that a $\bm{G}$-graded $A_{\infty}$ morphism from $(A,\mu)$ to $(B,\eta)$ over $R$ is an element 
\[ F \in TrCC_{\bm{G}}^1(A \otimes R, B \otimes R|R) \cong TrCC_{\bm{G}}^1(A,B \otimes R|\C)\]
such that 
\[ F \circ \mu - \eta \diamond F = 0.\]
Once again, we write
\[ F = \sum_{j \ge 0} F_j,\]
where $F_j$ is the order-$j$ component of $F$.
Observe that $F_0$ is a $\bm{G}$-graded $A_{\infty}$ morphism from $(A,\mu_0)$ to $(B,\eta_0)$. 

The notion of $A_{\infty}$ morphisms over a general ring $R$ is not as well-behaved as over a field. 
For example, it is not clear that $A_{\infty}$ quasi-isomorphisms can be inverted over $R$. 
It turns out that we will only need invertibility of quasi-isomorphisms in two situations: when the coefficient ring is a field, and when our $A_\infty$ algebra is {\bf minimal}. 
For minimal $A_\infty$ algebras, an $A_\infty$ quasi-isomorphism is necessarily a formal diffeomorphism. 
Here we recall the notion of $\bm{G}$-graded formal diffeomorphisms from Definition \ref{definition2:2formdiff}, and make appropriate modifications for the case of $A_{\infty}$ algebras defined over $R$:

\begin{definition2}
\label{definition2:2minformdiffr}
If $A$ and $B$ are $\bm{G}$-graded vector spaces over $\C$, and $R$ a $\bm{G}$-graded power series ring, then a {\bf $\bm{G}$-graded formal diffeomorphism} from $A \otimes R$ to $B \otimes R$ is an element
\[ F \in TrCC^1_{\bm{G}}(A,B\otimes R)\]
such that
\[ F^1: A \otimes R \To B \otimes R\]
is an isomorphism of free $R$-modules.
\end{definition2}

As before (see Lemma \ref{lemma2:formdiffinv} and Corollary \ref{corollary2:formdiffpush}), formal diffeomorphisms form a group with multiplication $\diamond$, and they can be used to push forward minimal $\bm{G}$-graded $A_{\infty}$ structures over $R$, so that $F$ defines a quasi-isomorphism from $(A,\mu)$ to $(B,F_*\mu)$, with strict inverse $F^{-1}$.
This last point is particularly important, because (as we stated above), there is no reason for an arbitrary $A_{\infty}$ quasi-isomorphism over $R$ to be invertible.

Now we introduce an analogue of Definition \ref{definition2:2gact} for minimal deformations of $A_{\infty}$ algebras.

\begin{definition2}
If $\mathcal{A} = (A,\mu_0)$ is a $\bm{G}$-graded minimal $A_{\infty}$ algebra over $\C$, we consider the group of $\bm{G}$-graded formal diffeomorphisms from $A$ to itself, whose leading-order term is the identity:
\[ \mathfrak{G}_R(A) := \{ F \in TrCC^1_{\bm{G}}(A,A \otimes R): F^1 = \mathrm{id} + \mathfrak{m}\}.\]
The group $\mathfrak{G}_R(A)$ acts on the set of $\bm{G}$-graded minimal $A_{\infty}$ structures on $A \otimes R$.
As before, the action of $F$ on $\mu \in TrCC^2_{\bm{G}}(A,A \otimes R)$ is denoted by $F_* \mu \in TrCC^2_{\bm{G}}(A,A \otimes R)$, and $F$ defines an $A_{\infty}$ morphism from $\mu$ to $F_* \mu$. 
\end{definition2}

\begin{lemma2}
\label{lemma2:quasiorder}
Let $\mathcal{C}$ be a minimal $A_\infty$ category over $\C$, and let $\mathscr{C}$ be a minimal deformation of $\mathcal{C}$ over some power series ring $R$. 
Let $X$ and $Y$ be objects of $\mathscr{C}$. 
If $X$ and $Y$ are quasi-isomorphic as objects of $\mathcal{C}$ (i.e., at $0$th order), then they are quasi-isomorphic as objects of $\mathscr{C}$.
\end{lemma2}
\begin{proof}
Let $f \in hom^*_\mathcal{C}(X,Y)$ be an isomorphism in $H^*(\mathcal{C})$. 
Now regard $f$ (more precisely, $f \otimes 1$) as an element of $hom^*_\mathscr{C}(X,Y)$. 
It is closed because $\mathscr{C}$ is minimal. 
Furthermore, for any object $Z$, the map of free $R$-modules
\[ \mu^2(f,\cdot): hom_\mathscr{C}^*(Z,X) \To hom_\mathscr{C}^*(Z,Y)\]
is an isomorphism to $0$th order, because $f$ is an isomorphism in $\mathcal{C}$; it follows that it is an isomorphism to all orders. 
Similarly, $\mu^2(\cdot,f)$ is an isomorphism, and it follows that $f$ defines an isomorphism in $H^*(\mathscr{C})$.

\end{proof}

\begin{definition2}
\label{definition2:2defclass}
Suppose that $\mu$ is a $\bm{G}$-graded deformation of the $A_{\infty}$ algebra $\mathcal{A} = (A,\mu_0)$ over $R$.
The first-order component of the $A_{\infty}$ relation $\mu \circ \mu = 0$ tells us that
\[ \mu_1 \in CC_{\bm{G}}^2\left(A,A \otimes R^1\right) \]
is a Hochschild cochain. 
Thus, we obtain an element
\[ [\mu_1] \in HH_{\bm{G}}^2(\mathcal{A}, \mathcal{A} \otimes R^1),\]
which we call the {\bf first-order deformation class} of $\mu$.
\end{definition2}

\begin{definition2}
If $\mu$ is a $\bm{G}$-graded {\bf minimal} deformation of the {\bf minimal} $A_{\infty}$ algebra $\mathcal{A}$ over $R$,
then the first-order component of $\mu$ defines an element in the {\bf truncated} Hochschild cohomology,
\[ [\mu_1] \in TrHH_{\bm{G}}^2(\mathcal{A}, \mathcal{A} \otimes R^1),\]
which we also call the {\bf first-order deformation class} of the minimal deformation $\mu$.
\end{definition2}

We are now almost ready to prove our main classification result for deformations of $A_{\infty}$ algebras. 
It turns out that in our particular situation, we need to incorporate a finite group action into the picture, so we now briefly explain how to do that.

\begin{definition2}
Let $H$ be a finite group. 
An {\bf action} of $H$ on a grading datum $\bm{G}$ is an action $\alpha$ of $H$ on $Y$ by group homomorphisms,
\[ \alpha: H \To \mathrm{Hom}_{Ab}(Y,Y),\]
which preserves $\Z$.
We will denote $\alpha(h)$ by $\alpha_h$.
\end{definition2}

\begin{example2}
\label{example2:hgradact}
In the case of Example \ref{example2:ses1}, there is an action of the symmetric group $H = S_n$ on $\bm{G}_n$, by permuting the generators of $\Z^{n}$.
\end{example2}

\begin{definition2}
Suppose that $H$ acts on the grading datum $\bm{G}$, and $V$ is a $\bm{G}$-graded vector space, and we have an action
\[ \rho: H \To \mathrm{Hom}_{Vect}(V,V).\]
We say that the action $\rho$ is $\bm{G}$-{\bf graded} if $\rho(h)$ maps $V_y$ to $V_{\alpha_h(y)}$.
\end{definition2}

\begin{definition2}
Suppose that we have compatible $\bm{G}$-graded actions of $H$ on the $\bm{G}$-graded algebra $R$, and on the $\bm{G}$-graded $R$-bimodules $A$ and $B$. 
Then there is a $\bm{G}$-graded action of $H$ on $CC^*_{\bm{G}}(A,B)$, via
\[ (h \cdot \phi)^s (a_s,\ldots,a_1) := h \cdot \phi^s (h^{-1} \cdot a_s,\ldots,h^{-1} \cdot a_1).\]
For $j$ an integer, we denote the $H$-invariant part of $CC^j_{\bm{G}}(A,B)$ by 
\[ CC^j_{\bm{G}}(A,B)^H := \{\phi \in CC^j_{\bm{G}}(A,B): h \cdot \phi = \phi \mbox{ for all $h \in H$}\}.\] 
\end{definition2}

\begin{definition2}
\label{definition2:2stricteq}
We say that a $\bm{G}$-graded $A_{\infty}$ algebra $\mathcal{A} = (A,\mu)$ over $R$ is {\bf strictly $H$-equivariant} if $\mu$ lies in $CC_{\bm{G}}^2(A)^H \subset CC_{\bm{G}}^2(A)$. 
Equivalently, we have
\[ \mu^k(h \cdot a_k, \ldots, h \cdot a_1) = h \cdot \mu^k(a_k, \ldots,a_1)\]
for all $k$ and all $h \in H$.
\end{definition2}

\begin{remark2}
We remark that $TrCC_{\bm{G}}^2(A,A \otimes R)^H$ is naturally a $R^H_{0}$-module (compare Remark \ref{remark2:r00mod}).
\end{remark2}

Now we prove our main classification result for $\bm{G}$-graded deformations of $A_{\infty}$ algebras over $R$:

\begin{proposition2}
\label{proposition2:2vershard}
Suppose that $\mathcal{A}=(A,\mu_0)$ is a $\bm{G}$-graded minimal $A_{\infty}$ algebra over $\C$, and furthermore is strictly $H$-equivariant with respect to the action of some finite group $H$ on $A$. 
Suppose that 
\[ TrHH_{\bm{G}}^2(\mathcal{A}, \mathcal{A} \otimes R)^H\]
is generated, as an $R^H_{0}$-module, by its first-order part
\[ TrHH_{\bm{G}}^2(\mathcal{A},\mathcal{A} \otimes R^1)^H,\]
and this first-order part is one-dimensional as a $\C$-vector space.
Then any two strictly $H$-equivariant $\bm{G}$-graded minimal deformations of $\mathcal{A}$, whose first-order deformation classes are non-zero in 
\[TrHH_{\bm{G}}^2(\mathcal{A},\mathcal{A} \otimes R^1)^H,\]
 are related by an element of $\mathrm{Aut}(R)^H$ composed with a $\bm{G}$-graded formal diffeomorphism. 
\end{proposition2}
\begin{proof}
Suppose that $(A,\mu)$ and $(A,\eta)$ are two such deformations.
We will construct, order-by-order, elements $\psi \in \mathrm{Aut}(R)^H$ and $F \in TrCC_{\bm{G}}^1(A,A \otimes R)^H$ so that $\psi \cdot \mu = F_* \eta$. 
The equation that $\psi$ and $F$ must satisfy is
\[ \psi \cdot \mu = F_* \eta.\]
We call this the $A_{\infty}$ relation for the purposes of this proof. 

We denote 
\begin{eqnarray*}
\psi &=& \sum_{j \ge 0} \psi_j \mbox{, where } \psi_j \in (R^j)^H_0\mbox{, and} \\
F &=& \sum_{j \ge 0} F_j \mbox{, where } F_j \in TrCC_{\bm{G}}^1(A,A \otimes R^j)^H.
\end{eqnarray*}
We start with $F_0 = \mathrm{id}$.
The order-zero component of the $A_{\infty}$ equation says that $\mu_0 = \eta_0$, which is true by assumption.
 
Now suppose, inductively, that we have determined $F_j$ and $\psi_{j-1}$ for all $j \le k-1$, $H$-invariant and $\bm{G}$-graded, and that 
\[ (\psi \cdot \mu - F_* \eta)_j = 0 \mbox{ for all $j \le k-1$.}\] 
We show that it is possible to choose $F_k$ and $\psi_{k-1}$ so that 
\[ (\psi \cdot \mu - F_* \eta)_k = 0.\]
The left hand side lies in $TrCC_{\bm{G}}^2(A,A \otimes R^k)^H$. 

First, we observe that
\[ [ \psi \cdot \mu - F_* \eta, \psi \cdot \mu + F_* \eta] = 0,\]
by expanding out the brackets: the cross-terms vanish by symmetry ($\psi \cdot \mu$ and $F_* \eta$ both have degree $2$), and the other terms vanish because $\psi \cdot \mu$ and $F_* \eta$ are $A_{\infty}$ structures. 
Now note that 
\[ (\psi \cdot \mu + F_* \eta)_0 = \mu_0 + \eta_0 = 2 \mu_0,\]
so taking the order-$k$ component of the previous equation gives us
\begin{eqnarray*}
 [ \psi \cdot \mu - F_* \eta, \mu_0]_k & = & 0 \\
\Rightarrow \delta \left((\psi \cdot \mu - F_* \eta)_k\right) = 0,
\end{eqnarray*}
regardless of what $F_k$ and $\psi_{k-1}$ are (here $\delta$ is the Hochschild differential).

Now we pick out the terms in $(\psi \cdot \mu - F_* \eta)_k$ that involve $F_k$ and $\psi_{k-1}$. 
First, we have
\[ (\psi \cdot \mu)_k = \psi_{k-1} \cdot \mu_1 + \mbox{ (terms involving $\psi_j$ for $j \le k-1$)}.\]
Next, following the proof of Lemma \ref{lemma2:formdiffinv}, we see that $F^{-1}_0 = \mathrm{id}$, and
\[ F^{-1}_k = - F_k + \mbox{(terms involving $F_j$ for $j \le k-1$).}\]
Thus, 
\begin{eqnarray*}
(F_* \eta)_k &=& \left((F \circ \eta) \diamond F^{-1} \right)_k  \\
&=& (F_k \circ \eta_0) \diamond \mathrm{id} + (\mathrm{id} \circ \eta_0) \circ (-F_k) + \mbox{(terms involving $F_j$, $j \le k-1$)}\\
&=& F_k \circ \eta_0 - \eta_0 \circ F_k \\
&=& \delta(F_k)
\end{eqnarray*}   
(the signs work out because $F$ has degree $1$). 
Therefore,
\[ (\psi \cdot \mu - F_* \eta)_k = \psi_{k-1} \cdot \mu_1 + \delta(F_k) + D_k,\]
where $D_k$ contains all the terms that do not involve $\psi_{k-1}$ or $F_k$.
Note that if we set $F_k =0$ and $\psi_{k-1} = 0$, our previous argument shows that $\delta(D_k) = 0$, so $D_k$ defines a class
\[ [D_k] \in TrHH_{\bm{G}}^2(\mathcal{A},\mathcal{A} \otimes R^k)^H.\]

Thus, we need to choose $\psi_{k-1}$ so that $\psi_{k-1} \cdot [\mu_1] = -[D_k]$ in the truncated Hochschild cohomology $TrHH_{\bm{G}}^2(\mathcal{A}, \mathcal{A} \otimes R^k)^H$. 
We can do this by our assumption that $[\mu_1]$ is non-zero, hence generates the one-dimensional first-order component $TrHH_{\bm{G}}^2(\mathcal{A},\mathcal{A} \otimes R^1)^H$, which generates $TrHH_{\bm{G}}^2(\mathcal{A},\mathcal{A} \otimes R)^H$ as a $R^H_0$-module. 
We then choose $F_k$ to effect the Hochschild coboundary between $\psi_{k-1} \cdot \mu_1$ and $-D_k$. 
We can make $F_k$ $H$-invariant by averaging over $H$.

Finally, note that at first order, we have 
\[\psi_0 \cdot [\mu_1] = [\eta_1],\]
from which it follows that $\psi_0 \neq 0$, because $[\mu_1]$ and $[\eta_1]$ are both non-zero, so indeed $\psi \in \mathrm{Aut}(R)^H$.
 \end{proof}

\subsection{Computations}
\label{subsec:2comp}

In this section, we introduce specific $\bm{G}$-graded $A_{\infty}$ algebras and deformations, and use the results of the previous sections to prove classification theorems for them.

Let us introduce some notation. We fix an integer $n \ge 4$ (we will be considering hypersurfaces in $\CP{n-1}$) and $a \ge 1$ (this will be the degree of the hypersurface in $\CP{n-1}$ that we will consider). 
In our intended applications in the current paper, $a$ will be either $1$ or $n$.

Throughout this section, we will be using the grading datum $\bm{G} := \bm{G}^n_1$ from Example \ref{example2:ses1}. 
We denote
\[ Y := \Z \langle y_1, \ldots, y_n \rangle,\]
so that $\bm{G}$ is given by
\[ \Z \To \widetilde{Y} \cong (\Z \oplus Y)/(2(1-n) \oplus y_{[n]}).\]
For an element $y \in Y$, we will denote
\[ |y| := y_{[n]} \cdot y.\]
We will denote by $H$ the symmetric group on $n$ elements, and recall that it acts on $\bm{G}$ (see Example \ref{example2:hgradact}).

We recall the $\bm{G}$-graded exterior algebra 
\[A := A_n := \Lambda(U_n)\]
of Definition \ref{definition2:2alga}.
For each subset $K \subset [n]$, we denote the corresponding element of $A$ by
\[ \theta^K := \bigwedge _{j \in K} \theta^j.\]
We equip the vector space $U_n$ with an $H$-action, which up to sign is the obvious action by permuting basis elements. 
In other words,
\[ h(u_j) = \pm u_{h(j)}.\]
We will not need to specify the actual signs. 
There is an induced action of $H$ on $A$.

We recall (from Definition \ref{definition2:2ringra}) the $\bm{G}$-graded ring
\[R_a := R^n_a,\]
which is the completion of $\C [ r_1,\ldots, r_{n} ]$ with respect to the length filtration in the category of $\bm{G}$-graded algebras.
We give a name to one important element of $R_a$: we denote
\[ T := r^{y_{[n]}} = r_1 \ldots r_{n}.\]
We equip $R_a$ with an $H$ action, which up to sign is the obvious action by permuting basis elements. 
We furthermore denote
\[ R:= R_{n}\]
because $a=n$ is the most important case we will consider. 
Finally, we will denote by $R^j_a$ the order-$j$ part of $R_a$. 
Note the change of notation from Definition \ref{definition2:2ringra}, where the superscript denoted the number of generators.
We hope this does not cause confusion.

\begin{definition2}
\label{definition2:2phihkr}
Suppose that $A \cong \Lambda(U)$ is an exterior algebra over $\C$. 
We define the {\bf Hochschild-Kostant-Rosenberg (HKR) map}
\begin{eqnarray*} 
\Phi: CC^*_c(A)&\To& \C[ U ]  \otimes A, \\
\Phi(\alpha) & := & \sum_{s=0}^{\infty} \alpha^s(\bm{u},\ldots,\bm{u}),
\end{eqnarray*}
where
\[ \bm{u} := \sum_{j=1}^{n} u_j \theta^j.\]
This defines an isomorphism on cohomology (the {\bf HKR isomorphism} \cite{Hochschild1962})
\[ \Phi: HH^*_c(A) \overset{\cong}{\To} \C[U] \otimes A.\]
We will slightly abuse notation and also write
\[ \Phi: HH^*(A) \overset{\cong}{\To} \C \llbracket U \rrbracket \otimes A\]
for the non-compactly supported version. 

Now suppose that $R$ is an arbitrary commutative $\C$-algebra. 
We define the $R$-algebra $A \otimes R$, and observe that $\Phi$ also induces an isomorphism
\[ \Phi: HH^*_c(A \otimes R|R)  \overset{\cong}{\To} R[U] \otimes A\]
(because tensoring with $R$ is exact), and hence also an isomorphism
\[ \Phi: HH^*(A \otimes R|R)  \overset{\cong}{\To} R \llbracket U \rrbracket \otimes A.\]
We also call these HKR maps.

We observe that HKR maps are `natural', in the sense that there is a commutative diagram
\[\begin{diagram} 
CC^*(A \otimes R|R) \otimes_R S & \rTo^{\Phi_R \otimes \bm{1}} & \left(R\llbracket U \rrbracket  \otimes A\right) \otimes_R S \\
\dTo && \dTo \\
CC^*(A \otimes S|S) & \rTo^{\Phi_S} & S\llbracket U \rrbracket  \otimes A.
\end{diagram}\]
\end{definition2}

\begin{definition2}
\label{definition2:2typea}
Let us fix $a = n$.
We say that a $\bm{G}$-graded $A_{\infty}$ algebra $\mathscr{A}$ over $R$ has {\bf type A} if it satisfies the following properties:
\begin{itemize}
\item Its underlying $R$-module and order-$0$ cohomology algebra is 
\[ (\mathscr{A},\mu^2_0) \cong A \otimes R;\]
\item It is strictly $H$-equivariant;
\item It satisfies
\[ \Phi(\mu^{\ge 3}) = \pm u_1 \ldots u_{n} + \sum_{j=1}^{n} \pm r_j u_j^{n} + \mathcal{O}(r^2).\]
\end{itemize}
The last expression may seem a little strange because $\mu^{\ge 3}$ is not a Hochschild cochain; here `$\Phi$' refers to the map 
\[ \Phi:CC^*_{\bm{G}}(A \otimes R|R) \To R \llbracket U \rrbracket \otimes A.\]
\end{definition2}

Now we state the main result we will prove in this section:

\begin{theorem2}
\label{theorem2:2typea}
Suppose that $\mathscr{A}_1 = (A \otimes R, \mu)$ and $\mathscr{A}_2 = (A \otimes R,\eta)$ are two $\bm{G}$-graded $A_{\infty}$ algebras over $R$ of type A. 
Then there exists $\psi \in R_0^H$, of the form
\[ \psi(r_1, \ldots, r_n) = \pm 1 + \sum_{j=1}^{\infty} c_j T^j,\]
and a $\bm{G}$-graded formal diffeomorphism $F \in \mathfrak{G}_R(A)$, such that
\[ \psi \cdot \mathscr{A}_1 = F_* \mathscr{A}_2.\]
\end{theorem2}

To begin the proof of Theorem \ref{theorem2:2typea}, we first show that the unspecified signs in $\Phi(\mu^{\ge 3})$ actually amount only to a single unspecified sign:

\begin{lemma2}
\label{lemma2:signsrid}
Any $A_\infty$ algebra of type A is strictly isomorphic to one with
\[ \Phi(\mu^{\ge 3}) = u_1 \ldots u_n \pm \sum_{j=1}^n r_j u_j^n + \mathcal{O}(r^2);\]
henceforth we shall assume that all our $A_\infty$ algebras of type A have this property.
\end{lemma2}
\begin{proof}
Note that the terms $r_ju_j^n$ all have the same sign, by $H$-equivariance. 
If the leading term has a `$+$' sign we are done; if it has a `$-$' sign, let $\eta$ be an $n$th root of $-1$, and change the basis for $A$ by multiplying all generators $\theta^j$ by $\eta$.
   
\end{proof}

We will now give a brief outline of the remainder of the proof of Theorem \ref{theorem2:2typea}. 
The first step is to prove (in Corollary \ref{corollary2:aunique}, using Proposition \ref{proposition2:2vers1}) that the order-$0$ parts $(A,\mu_0)$ and $(A,\eta_0)$ are related by a formal diffeomorphism.
The classification of the order-$0$ part is governed by the Hochschild cohomology $HH^*_{\bm{G}}(A)$, which is determined via the HKR isomorphism. 
We thereafter denote by $\mathcal{A} \cong (A,\mu_0)$ this (unique up to formal diffeomorphism) order-$0$ part. 

We then study deformations of $\mathcal{A}$ over $R$. 
The classification of such deformations is governed by the Hochschild cohomology $HH^*_{\bm{G}}(\mathcal{A}, \mathcal{A} \otimes R)$, which is also determined from the HKR isomorphism, via a spectral sequence. 
We apply Proposition \ref{proposition2:2vershard} to show that such deformations are unique up to $A_{\infty}$ quasi-isomorphism and the action of $\mathrm{Aut}(R)^H$.

Now let us begin. 
We first explain how to use the HKR isomorphism to calculate $HH^*_{c,\bm{G}}(A)$, and more generally $HH^*_{\bm{G}}(A,A \otimes R_a)$. 
Taking gradings into account, the HKR isomorphism tell us that
\[  HH^*_{c,\bm{G}}(A) \cong \C[U] \otimes \Lambda(U)\]
as $\widetilde{Y} \oplus \Z$-graded vector spaces. 
We will make the $\bm{G}$-grading of the right-hand side explicit in Lemma \ref{lemma2:hhrest}.

More generally, it follows that there is an isomorphism of $\widetilde{Y} \oplus \Z$-graded vector spaces,
\[ HH_{c,\bm{G}}^*(A,A \otimes R_a) \cong \C[U] \otimes \Lambda(U)  \otimes R_a. \]
$\C[U] \otimes \Lambda(U) \otimes R_a$ is generated by terms $u^{\bm{b}} \theta^K r^{\bm{c}}$, where $\bm{b}, \bm{c} \in Y_{\ge 0}$ (recall $Y \cong \Z^n$, so we define $Y_{\ge 0} := \Z_{\ge 0}^n$) and $K \subset [n]$.
If $\bm{b} = \sum_j b_j y_j$, then this is the image under the HKR map of a Hochschild cochain which sends 
\[ \bigotimes_j (\theta^j)^{\otimes b_j} \mapsto r^{\bm{c}} \theta^K.\]

We start by examining what the various gradings on a Hochschild cochain tell us.

\begin{lemma2}
\label{lemma2:ccrest}
If a generator $\tau \in CC_{c,\bm{G}}^{s+t}\left(A,A \otimes R_a^j\right)^s$ sends
\[ \bigotimes_{i=1}^{s} \theta^{K_i} \mapsto r^{\bm{c}} \theta^{K_0},\]
then we have
\begin{eqnarray}
j &=& |\bm{c}| \label{eqn:ccr1} \\
a \bm{c} + y_{K_0} - \sum_{i=1}^s y_{K_i} &=& q y_{[n]} \mbox{ in $Y$, for some $q$ } \label{eqn:ccr2} \\
t &=& (n-2)q + (2-a)j. \label{eqn:ccr3}
\end{eqnarray}
\end{lemma2}
\begin{proof}
Equation (\ref{eqn:ccr1}) follows by definition. 
To prove Equations (\ref{eqn:ccr2}) and (\ref{eqn:ccr3}), we recall that the grading of $\theta_j$ is $(-1,y_j)$ (Example \ref{example2:vecu}) and the grading of $r_j$ is $(2-2a,ay_j)$ (Example \ref{example2:vecd}). 
We alter the grading datum $\bm{G}$ by an automorphism sending
\[ (j,y) \mapsto (j+y_{[n]} \cdot y, y).\]
This is equivalent to considering the pseudo-grading datum with the same exact sequence as $\bm{H}_n$, but 
\[c = 2(n-1) - y_{[n]} \cdot y_{[n]} = n-2.\]
Then $\theta_j$ has grading $(0,y_j)$ and $r_j$ has grading $(2-a,a y_j)$. 

If the grading of $\tau$ is $f(t) = (t,0) \in  (\Z \oplus Y)/Z = \widetilde{Y}$, then
\[ (0, y_{K_0}) + ((2-a)|\bm{c}|, a \bm{c}) - \sum_{j=1}^s (0,y_{K_j}) \cong (t,0) \mbox{ modulo $Z$.}\]
Recalling that the image of $Z$ in $\Z \oplus Y$ is given by
\[ ((2-n)q, q y_{[n]})\]
(in the altered grading datum), we have
\[(0, y_{K_0}) + ((2-a)|\bm{c}|, a \bm{c}) - \sum_{j=1}^s (0,y_{K_j}) = (t,0) + q(2-n,y_{[n]}),\]
from which the result follows.
 \end{proof}

Now recalling the HKR isomorphism, a generator of $HH_{c,\bm{G}}^{s+t}(A,A \otimes R_a)^s$ has the form $u^{\bm{b}} \theta^K r^{\bm{c}}$, where $\bm{b}$ and $\bm{c}$ are elements of $Y_{\ge 0}$, and $K \subset [n]$. 
We examine what the gradings tell us about such a generator of the Hochschild cohomology.

\begin{lemma2}
\label{lemma2:hhrest}
If $u^{\bm{b}} \theta^K r^{\bm{c}}$ is a generator of $HH_{c,\bm{G}}^{s+t}\left(A,A \otimes R^j_a\right)^s$, then the following equations hold:
\begin{eqnarray}
j &=& |\bm{c}| \label{eqn:hhr1} \\
s &=& |\bm{b}| \label{eqn:hhr2} \\
y_K + a\bm{c} - \bm{b} &=& q y_{[n]} \mbox{ in $Y$, for some $q$} \label{eqn:hhr3} \\
t &=& (n-2)q + (2-a)j \label{eqn:hhr4} \\
|K| &=& s+t + 2(q-j). \label{eqn:hhr5}\\
(n-2)|K| &=& (n-2)s + nt + 2(a-n)j \label{eqn:hhr6}
\end{eqnarray}
\end{lemma2}
\begin{proof}
Equations (\ref{eqn:hhr1}) and (\ref{eqn:hhr2}) hold by definition.
Equation (\ref{eqn:hhr3}) follows from Equation (\ref{eqn:ccr2}).
Equation (\ref{eqn:hhr4}) follows from Equation (\ref{eqn:ccr3}).
The first step in proving Equations (\ref{eqn:hhr5}) and (\ref{eqn:hhr6}) is to take the dot product of Equation (\ref{eqn:hhr3}) with $y_{[n]}$:
 \begin{eqnarray*}
y_{[n]} \cdot (y_K + a\bm{c} - \bm{b}) &=& qy_{[n]} \cdot y_{[n]} \\
 \Rightarrow |K| &=& |\bm{b}| - a |\bm{c}| + nq \\
&=& s - a j + nq.
 \end{eqnarray*}
To prove equation (\ref{eqn:hhr5}), we use equation (\ref{eqn:hhr4}) to substitute for $nq$:
\begin{eqnarray*}
|K|&=& s - a j + t + (a-2) j + 2q\\
&=& s+t + 2(q-j).
 \end{eqnarray*}
We use the same equation to prove equation (\ref{eqn:hhr6}), but this time we first multiply by $n$ then substitute in equation (\ref{eqn:hhr4}):
\begin{eqnarray*}
(n-2)|K| &=& (n-2)(s - aj + nq) \\
 &=& (n-2)s - (n-2)aj + n(t+(a-2)j) \\
&=& (n-2)s + nt + 2(a-n)j.
\end{eqnarray*}
 \end{proof}

We now set about determining the order-$0$ part of an $A_{\infty}$ algebra of type A. 
We identify the possible $\bm{G}$-graded $A_{\infty}$ algebras $\mathcal{A} = (A,\mu_0)$ over $\C$ with underlying vector space and $\mu^2_0$ given by $A$.

\begin{lemma2}
\label{lemma2:mun2}
We have
\[ CC_{c,\bm{G}}^{s+t}(A)^s \cong 0\]
unless $t$ is divisible by $(n-2)$.
\end{lemma2}
\begin{proof}
Follows from Equation (\ref{eqn:hhr4}) with $j = 0$.
 \end{proof}

\begin{lemma2}
\label{lemma2:hhw}
The Hochschild cohomology of $A$ satisfies, for $s>2$,
\[ HH_{c,\bm{G}}^2(A)^{s} \cong \left\{ \begin{array}{ll}
                                        \C \cdot u^{y_{[n]}} & \mbox{if $s = n$,} \\
                                        0 & \mbox{otherwise,}
                                 \end{array} \right. \]
where $u^{y_{[n]}} = u_1 \ldots u_{n}$.
\end{lemma2}
\begin{proof}
Let $u^{\bm{b}} \theta^K$ be a generator of $HH_{c,\bm{G}}^2(A)^{s}$. 
Equation (\ref{eqn:hhr4}), with $j=0$, yields
\[ s = 2-(n-2)q.\]
We want $s>2$, so we have $q < 0$.
Now equation (\ref{eqn:hhr5}), with $j = 0$, yields
\[ 2 + 2q = |K| \ge 0. \]
Thus, $-1 \le q < 0$, so $q = -1$ and $|K| = 0$, so $K = \emptyset$.

Equation (\ref{eqn:hhr3}) now yields $\bm{b} = y_{[n]}$. 
Therefore $u^{\bm{b}} \theta^K = u^{y_{[n]}}$, as required.
 \end{proof}

\begin{corollary2}
\label{corollary2:aunique}
There is a formal diffeomorphism $F \in TrCC^1_{\bm{G}}(A)$ between the $0$th-order parts of $\mathscr{A}_1$ and $\mathscr{A}_2$:
\[\eta_0 = F_* \mu_0.\]
\end{corollary2}
\begin{proof}
First, we observe that $CC_{c,\bm{G}}^2(A)^1 \cong 0$ by Equation (\ref{eqn:ccr3}) (with $j=0$, and assuming $n \ge 4$), so the $A_{\infty}$ algebras are minimal. 
Furthermore, $\mu^d_0 = 0$ for $2<d<n$ because $CC_{c,\bm{G}}^2(A)^{d} \cong 0$ by Lemma \ref{lemma2:mun2}.
The result now follows from Proposition \ref{proposition2:2vers1}, because the deforming class of an algebra of type A is
\[ \Phi_0(\mu_0) = u_1 \ldots u_n =  u^{y_{[n]}}\]
by definition.
 \end{proof}

Henceforth, we will denote by $\mathcal{A}$ any $\bm{G}$-graded minimal $A_{\infty}$ algebra over $\C$ with cohomology algebra given by $A$, and whose order-$n$ deforming class is a non-zero multiple of $u^{y_{[n]}}$. 
The previous lemma says that $\mathcal{A}$ is well-defined up to formal diffeomorphism.

We will consider deformations of $\mathcal{A}$ over $R_a$, which are controlled by the Hochschild cohomology with coefficients in $R_a$,
\[ HH_{\bm{G}}^*(\mathcal{A}, \mathcal{A} \otimes R_a).\]

Recall from Remark \ref{remark2:ss} that the filtration by length on the Hochschild complex $CC^*_{\bm{G}}(A)$ yields a spectral sequence for the Hochschild cohomology, with 
\[E_2^{*,*} = HH_{c,\bm{G}}^*(A,A).\]

\begin{lemma2}
\label{lemma2:ssconv}
The spectral sequence induced by the length filtration on the Hochschild cochain complex $CC^*_{\bm{G}}(\mathcal{A})$ converges to the Hochschild cohomology $HH^*_{\bm{G}}(\mathcal{A})$.
\end{lemma2}
\begin{proof}
By Remark \ref{remark2:ss}, it suffices to prove that the spectral sequence is regular.

We first work out the $\widetilde{Y} \oplus \Z$ grading of a generator $u^{\bm{b}} \theta^K $ of $HH_{c,\bm{G}}^y(A)^s$, following Lemma \ref{lemma2:hhrest}. 
Firstly, by definition, $s = |\bm{b}|$. 
Next, this generator changes the $\widetilde{Y}$-degree by 
\[ (-|K|+|\bm{b}|,-\bm{b} + y_K) \in \Z \oplus Y/(2(1-n),y_{[n]}).\]
Therefore, its $\widetilde{Y}$-grading is
\[ (-|K|+2|\bm{b}|,-\bm{b}+y_K)\]
(recalling the conventions for grading Hochschild cochains). 
Altering the grading datum by an automorphism 
\[ (j, y) \mapsto (j + 2y_{[n]} \cdot y,y),\]
this becomes equivalent to a grading
\[(|K|,-\bm{b}+y_K)\]
(note: $c=2$ in the correspondingly altered pseudo-grading datum).

Now recall that the differential on page $d$ of the spectral sequence maps
\[ \delta_d: E_d^{(y,s)} \To E_d^{(y+f(1),s+d)}.\]
If both domain and codomain of the differential are to be non-zero, then we must have generators $u^{\bm{b}_1}\theta^{K_1}$ and $u^{\bm{b}_2} \theta^{K_2}$ such that
\[ (|K_1| + 1,-\bm{b}_1+y_{K_1}) = (|K_2|,-\bm{b}_2+y_{K_2}) + q(2,y_{[n]})\]
and
\[ |\bm{b}_2| = |\bm{b}_1| + d.\]
Because $0 \le |K_1|,|K_2| \le n$, and $|K_1| +1 = |K_2| + 2q$, we must have $q \le (n+1)/2$. 
We also have
\[ \bm{b}_2 = \bm{b}_1 + y_{K_2} - y_{K_1} + q y_{[n]},\]
and hence that 
\[ |\bm{b}_2| = |\bm{b}_1| + |K_2| - |K_1| + nq = |\bm{b}_1| + 1 + (n-2)q.\]
Therefore, for the differential $\delta_d$ to be non-zero, we must have
\[ d =|\bm{b}_2| - |\bm{b}_1| = (n-2)q +1 \le \frac{(n+1)(n-2)}{2}+1.\] 
Hence, for $d$ sufficiently large, the differential vanishes, so the spectral sequence is regular. 
 \end{proof}

Now, if we are to apply the deformation theory of Section \ref{subsec:2defainf}, we need to know that our deformations are minimal. 

\begin{lemma2}
\label{lemma2:mindef}
If $a=n$, and $n \ge 5$, then any $\bm{G}$-graded deformation of $\mathcal{A}$ over $R$ is minimal.
\end{lemma2}
\begin{proof}
Equation (\ref{eqn:ccr3}) with $a=n$ shows that
\[CC_{c,\bm{G}}^2(A,A \otimes R)^s \cong 0\]
unless $t$ is divisible by $(n-2)$. 
The $\mu^0$ and $\mu^1$ terms live in the spaces with $t=2$ and $t=1$, respectively. 
So when $n \ge 5$, any deformation of $\mathcal{A}$ over $R$ satisfies $\mu^0 = 0$ and $\mu^1 = 0$, hence is minimal.
 \end{proof}

\begin{remark2}
\label{remark2:nnot4}
Lemma \ref{lemma2:mindef} is the only place where we need the assumption $n \ge 5$, as opposed to $n \ge 4$. 
If $n=4$, then $CC_{c,\bm{G}}^2(A,A \otimes R)^0$ is one-dimensional, generated by the element $T\theta^{[4]}$. 
This would correspond to some curvature in our category. 
In the Fukaya category of the quartic surface, this type of deformation can be ruled out by choosing an almost-complex structure with no holomorphic disks with boundary on our Lagrangians, but we will avoid this complication as that case is dealt with completely in \cite{Seidel2003}.
\end{remark2}

The next step is to determine the first-order deformation space.

\begin{lemma2}
\label{lemma2:defclasshh}
If $n \ge 4$, then the vector space
\[ HH_{\bm{G}}^2(A, A \otimes R^1_a)\]
is generated by the elements
\[ r_j u_j^a,\]
for $j = 1,\ldots,n$.
\end{lemma2}
\begin{proof}
We apply Lemma \ref{lemma2:hhrest} with $s+t = 2$ and $j = 1$. 
 Equation (\ref{eqn:hhr5}) yields
\[|K| = 2q \ge 0.\]
 Furthermore, Equation (\ref{eqn:hhr4}) gives
 \begin{eqnarray*}
 2-s &=& (n-2)q + 2-a \\
 \Rightarrow a-s &=& (n-2)q.
 \end{eqnarray*}
  
 Furthermore, we have $\bm{c} = y_k$ for some $k$.  
 Taking the dot product of $y_k$ with Equation (\ref{eqn:hhr3}) gives us
 \begin{eqnarray*}
  y_k \cdot (y_K + a\bm{c} - \bm{b}) &=& q y_k \cdot y_{[n]} \\
  \Rightarrow \mbox{($0$ or $1$)} &=& -a + y_k \cdot \bm{b} + q \\
  & \le & -a + s + q \\
  &=& (3-n)q.
  \end{eqnarray*}
  If $n \ge 4$, then $3-n < 0$ and $q \ge 0$, so we must have $q = 0$. 
  Thus, we have $|K| = 2q = 0$, so $K=\emptyset$. 
  From Equation (\ref{eqn:hhr3}), we obtain $\bm{b} = a y_j$, so the generator has the form
  \[ r_j u_j^a.\]
 \end{proof}

Having determined the first-order deformation space, we now seek to understand the higher-order parts of the deformation. 
The following lemma will be useful:

\begin{lemma2}
\label{lemma2:j2}
If  $r^{\bm{c}} u^{\bm{b}} \theta^K$ is a generator of  $TrHH_{\bm{G}}^2(A,A \otimes R_a)$, and $j=|\bm{c}| \ge 2$, then $\bm{c} \ge y_{[n]}$.
\end{lemma2}
\begin{proof}
We apply Lemma \ref{lemma2:hhrest}, with $s+t = 2$. 
Equation (\ref{eqn:hhr5}) yields
\[ 0 \le |K| = 2(1+q-j).\]
It follows that $q \ge j-1$. 
Now we split into two cases, depending on whether $q = j-1$ or $q > j-1$.

{\bf Case 1:} $q=j-1$. In this case, we have $|K| = 0$ so $K = \emptyset$. 
Thus Equation (\ref{eqn:hhr3}) gives
\[ a \bm{c} = q y_{[n]} + \bm{b}.\]
Taking the dot product with $y_k$ yields
\[ a c_k = q+ b_k \ge q = j-1 > 0,\]
since $j \ge 2$. 
Hence $c_k \ge 1$ for all $k$, and the result is proved.

{\bf Case 2:} $q>j-1$. In this case, we have $q > j-1 \ge 1$, so $q \ge 2$. 
Thus, taking the dot product of Equation (\ref{eqn:hhr3}) with $y_k$ gives
\[ a  c_k = q + b_k - y_k \cdot y_K \ge q - 1 > 0.\]
Hence $c_k \ge 1$ for all $k$, and the result is proved.
 \end{proof}

Now if $a=n$, then the grading of $T = r^{y_{[n]}} \in R$ is $0$ in $\widetilde{Y}$, and furthermore $T$ is clearly $H$-invariant. 
Thus $\C\llbracket T \rrbracket  \subset R^H_0$ (in fact one can show they are equal).
It follows that $TrHH_{\bm{G}}^2(\mathcal{A},\mathcal{A} \otimes R)^H$ is a $\C\llbracket T \rrbracket $-module. 

\begin{lemma2}
\label{lemma2:thhversal}
Suppose that $a=n$, and $n \ge 4$. 
Then the $H$-invariant part of the truncated Hochschild cohomology group
\[ TrHH_{\bm{G}}^2(\mathcal{A},\mathcal{A} \otimes R)^H\]
is generated, as a $\C\llbracket T \rrbracket $-module, by its first-order part
\[ TrHH_{\bm{G}}^2(\mathcal{A},\mathcal{A} \otimes R^1)^H.\]
\end{lemma2}
\begin{proof}
We prove the result by induction on the order $j$, and using the spectral sequence induced by the length filtration.

When $j=0$, we must show that $TrHH_{\bm{G}}^2(\mathcal{A}) \cong 0$. 
Let $u^{\bm{b}} \theta^K$ be a generator of $TrHH_{\bm{G}}^2(A)$. 
By equation (\ref{eqn:hhr4}) with $j=0$, we have
\[ t = (n-2)q \le 0\]
(since we are considering the truncated Hochschild cohomology). 
By equation (\ref{eqn:hhr5}) with $j=0$, we have
\[ |K| = 2(1+q) \ge 0.\]
Thus we have $-1 \le q \le 0$. 
We have two cases: 

{\bf Case 1:} $q = 0$. In this case, $|K| = 2$, $t = 0$, so $s=2$. 
Equation (\ref{eqn:hhr3}) yields
\[ y_K = \bm{a},\]
so our generator has the form $u_i u_k \theta^i \wedge \theta^k$. 
But now, if $\sigma_{ik} \in H$ denotes the transposition of elements $i$ and $k$, we have 
\[ \sigma_{ik}( u_i u_k \theta^i \wedge \theta^k) =u_k u_i \theta^k \wedge \theta^i = -u_i u_k \theta^i \wedge \theta^k.\]
Hence, the coefficient of this term in any $H$-invariant element of the truncated Hochschild cohomology must be $0$.

\begin{remark2}
Recall that $H$ may differ from the obvious action by permutation of basis elements, by some signs. 
However, each sign appears twice (once for $u_i$, once for $\theta^i$), so they do not affect this computation. 
\end{remark2}

{\bf Case 2:} $q = -1$. In this case, $|K| = 0$, so $K = \emptyset$. 
Equation (\ref{eqn:hhr3}) yields
\[ y_{[n]} = \bm{b}.\]
Thus, the generator is $u^{y_{[n]}}$, the deforming class identified in Lemma \ref{lemma2:hhw}. 
However, this generator gets killed by the first non-trivial differential of the spectral sequence. 
To see this, observe that, because $\mu^{d} = 0$ for $2<d<n$ by Lemma \ref{lemma2:mun2}, the first non-trivial differential is $\delta_{n-1}$, and is given by 
\[ \delta_{n-1}(\phi) = [[\mu^{n}],\phi] = \left[u^{y_{[n]}},\phi\right]\]
(see Remark \ref{remark2:ssfirst}). 
The Gerstenhaber bracket on $CC^*(A)$ gets carried to the Schouten bracket on polyvector fields $\C\llbracket U \rrbracket  \otimes \Lambda(U)$. 
It follows quickly from the explicit form of the Schouten bracket (see \cite[Equation (3.7)]{Seidel2008a}) that for any $W \in \C\llbracket U \rrbracket $,
\[ [W,\eta] = \iota_{dW}( \eta),\]
and therefore that the cohomology of $[W,-]$ is the cohomology of the Koszul complex associated to $dW$. 
In particular, the class $W$ itself gets killed by taking the cohomology of this differential.
In our case, this means that $u^{y_{[n]}}$ is killed by $\delta_{n-1}$.

This completes the proof that $TrHH_{\bm{G}}^2(\mathcal{A}) \cong 0$.

When $j=1$, the statement is simply that the first-order part is generated by the first-order part. 
So it remains to prove that, for $j \ge 2$, the order-$j$ part is generated by the first-order part.

Suppose, inductively, that 
\[TrHH_{\bm{G}}^2(\mathcal{A},\mathcal{A} \otimes R^k) \subset \C\llbracket T \rrbracket  \cdot TrHH_{\bm{G}}^2(\mathcal{A},\mathcal{A}\otimes R^1)\]
for all $k \le j-1$, for some $j \ge 2$. 
We prove that the same holds for $k=j$. 
Let 
\[r^{\bm{c}}u^{\bm{b}}\theta^K \in TrHH_{\bm{G}}^2(A,A \otimes R^j)\]
be a generator.
Because $j \ge 2$, it follows from Lemma \ref{lemma2:j2} that $r^{\bm{c}} = T \cdot r^{\bm{c}'}$, where
\[ \bm{c}' = \bm{c} - y_{[n]} \ge 0.\]
Because $T\in R_0$, we have
\[ r^{\bm{c}'} u^{\bm{b}} \theta^K \in TrHH_{\bm{G}}^2(A,A \otimes R^{j-n}).\] 
It follows that
\[ TrHH_{\bm{G}}^2(A,A \otimes R^j) = T \cdot TrHH_{\bm{G}}^2(A,A \otimes R^{j-n}),\]
and hence, from the spectral sequence induced by the length filtration, that
\[ TrHH_{\bm{G}}^2(\mathcal{A},\mathcal{A} \otimes R^j) = T \cdot TrHH_{\bm{G}}^2(\mathcal{A},\mathcal{A} \otimes R^{j-n}).\]
The result now follows by taking $H$-invariant parts of this equation.

This completes the inductive step, and hence the proof.
 \end{proof}

\begin{corollary2}
There exists a $\psi \in \C\llbracket T \rrbracket $, with $\psi(0) = \pm 1$, and a $\bm{G}$-graded formal diffeomorphism $F$, such that
\[ \psi \cdot \mathscr{A}_1 = F_* \mathscr{A}_2.\]
\end{corollary2}
\begin{proof}
By Corollary \ref{corollary2:aunique}, there is a formal diffeomorphism $F_0$ from $(A,\eta_0)$ to $(A,\mu_0)$. 
Therefore we can push forward $\eta$ by $F_0$, and reduce to the case where $\mu_0 = \eta_0$. 
Now note that the deformations are minimal (Lemma \ref{lemma2:mindef}), that
\[ TrHH^2_{\bm{G}}(\mathcal{A},\mathcal{A} \otimes R^1)^H \cong \C \langle r_1 u_1^n + \ldots + r_n u_n^n \rangle\]
is one-dimensional (by Lemma \ref{lemma2:defclasshh}), and that $TrHH^2_{\bm{G}}(\mathcal{A},\mathcal{A} \otimes R)^H$ is generated, as $R^H_0$-module, by its first-order part (Lemma \ref{lemma2:thhversal}).
Thus, the result follows from Proposition \ref{proposition2:2vershard}. 
The fact that $\psi(0) = \pm 1$ follows from the proof of Proposition \ref{proposition2:2vershard}: if the signs of $r_1u_1^n + \ldots + r_n u_n^n$ are the same for both $\mathscr{A}_1$ and $\mathscr{A}_2$, then $\psi(0) = 1$; if they are different then $\psi(0) = -1$.
 \end{proof}

This completes the proof of Theorem \ref{theorem2:2typea}.

We now make one final computation of Hochschild cohomology. 
For this computation, we will be interested in an extra structure on Hochschild cohomology: the {\bf Yoneda product}, which makes it into an associative algebra. 
In fact, together with the Gerstenhaber bracket, this makes the Hochschild cohomology  into a {\bf Gerstenhaber algebra} (see \cite{Gerstenhaber1963}).

\begin{definition2}
\label{definition2:2yon}
Let $\mathcal{A}$  be a $\bm{G}$-graded $A_{\infty}$ algebra.  
We define the {\bf Yoneda product}, which is a map
\[ CC^*_{\bm{G}}(\mathcal{A}) \otimes  CC^*_{\bm{G}}(\mathcal{A}) \To  CC^*_{\bm{G}}(\mathcal{A}),\]
which we denote by
\[ \phi \otimes \psi \mapsto \phi \smile \psi,\]
and which is defined by
\[ (\phi \smile \psi)^n (a_n, \ldots , a_1) := \]
\[ \sum_{0 \le j \le k \le l \le m \le n} (-1)^{\dagger} \mu^{\alpha}(a_{n}, \ldots, \phi^{\beta}(a_{m}, \ldots), a_{l},\ldots, \psi^{\gamma}(a_{k},\ldots),a_j,\ldots, a_1), \]
where 
\[ \alpha = n+2+l+j-m-k, \,\, \beta = m-l, \,\, \gamma = k-j,\]
and
\[ \dagger = (\sigma(\phi) + 1)(\sigma(a_1) + \ldots + \sigma(a_l) - l) + (\sigma(\psi) + 1)(\sigma(a_1) + \ldots + \sigma(a_j) - j).\]
\end{definition2}

\begin{remark2}
\label{remark2:yon}
We record the following useful information about the Yoneda product:
\begin{itemize}
\item The Hochschild differential $\delta$ and the Yoneda product $\smile$ extend to an $A_{\infty}$ algebra structure on $CC^*_{\bm{G}}(\mathcal{A})$ (see \cite[Proposition 1.7]{Getzler1993}); in fact one can check that this $A_{\infty}$ algebra is $\bm{G}$-graded. 
\item In particular, $\smile$ defines a $\bm{G}$-graded, associative product on $HH^*_{\bm{G}}(\mathcal{A})$.
\item If $\bm{k}$ is a field, and $U$ a $\bm{k}$-vector space, then the HKR isomorphism
\[ \Phi: HH^*(\Lambda(U)) \To \bm{k}\llbracket U \rrbracket  \otimes \Lambda(U),\]
is an isomorphism of $\bm{k}$-algebras, where the product on $HH^*(\Lambda(U))$ is the Yoneda product, and the product on polyvector fields is the wedge product (see \cite[Section 8]{Kontsevich2003}).
\item The spectral sequence $E_d^{s,t} \Rightarrow HH^*(\mathcal{A})$ induced by the length filtration (Remark \ref{remark2:ss}) respects the multiplicative structure of the Yoneda product.
\end{itemize}
\end{remark2}

We now carry out our final Hochschild cohomology computation. 
Suppose that $\mathscr{A}$ is an $A_{\infty}$ algebra of type A, and $\underline{\mathscr{A}}$ its extension to a $\bm{G}^n_1$-graded category. 
Let  $\bm{q}_1$ and $\bm{p}_1$ be the morphisms of grading data defined in Lemma \ref{lemma2:squaregrad}. 
We denote 
\[ \widetilde{\mathscr{A}} :=  \bm{p}_1^* \underline{\mathscr{A}}.\]
Note that $\bm{q}_{1*} \widetilde{\mathscr{A}}$ is a $\bm{G}_{\Z}$-graded  (or equivalently, $\Z$-graded) $A_{\infty}$ category, over the $\Z$-graded coefficient ring $\bm{q}_{1*} \bm{p}_1^* R \cong R$, which has degree $0 \in \Z$, by Remark \ref{remark2:q1p1r}. 
Therefore, the homomorphism 
\begin{eqnarray*}
R & \To & \Lambda,\\
r_j & \mapsto & r \mbox{ for all $j$}
\end{eqnarray*} 
(where $\Lambda$ is the universal Novikov field of Definition \ref{definition2:2novikov})
respects the $\Z$-grading, and we can form the $\bm{G}_{\Z}$-graded $\Lambda$-linear $A_{\infty}$ category 
\[\widetilde{\mathscr{A}}_{nov} := \bm{q}_{1*}\widetilde{\mathscr{A}} \otimes_R \Lambda.\]
Our aim now is to compute a certain part of
\[ HH^*_{\bm{G}_{\Z}}\left(\left. \widetilde{\mathscr{A}}_{nov} \right| \Lambda  \right).\]

We observe that $Y^n_1$ acts on $\underline{\mathscr{A}}$ by shifts, hence also acts on $\widetilde{\mathscr{A}}_{nov}$, hence also acts on $CC^* _{\bm{G}_{\Z}}(\widetilde{\mathscr{A}}_{nov})$. 
The action of $\mathrm{ker}(Y^n_n \To \Z) \subset Y^n_1$ is trivial, because shifts in this subgroup become trivial when we take $\bm{q}_{1*}$ of our category. 
The action of $\Z$ is also trivial, because we consider only $\Z$-equivariant cochains, by definition. 
It follows that the action of the image of $Y^n_n$ in $Y^n_1$ is trivial, so there is an action of
\[ \tilde{\Gamma}^*_n \cong Y^n_1/Y^n_n\]
on $CC^*_{\bm{G}_{\Z}}(\widetilde{\mathscr{A}}_{nov})$. 
Note that this is the same $\tilde{\Gamma}_n^*$ as in Definition \ref{definition2:2mirror}.

\begin{lemma2}
\label{lemma2:hhtypea}
There is an isomorphism of $\Z$-graded $\Lambda$-algebras,
\[ HH^*_{\bm{G}_{\Z}}\left(\widetilde{\mathscr{A}}_{nov}\right)^{\tilde{\Gamma}_n^*} \cong \Lambda [ \alpha]/\alpha^{n-1},\]
where $\alpha$ has degree $2$.
\end{lemma2}
\begin{proof}
We consider the spectral sequence induced by the length filtration on 
\[ CC^*_{\bm{G}_{\Z}}\left(\widetilde{\mathscr{A}}_{nov}\right)^{\tilde{\Gamma}_n^*}.\]
As we saw in Remark \ref{remark2:ss}, the filtration is bounded above, hence exhaustive, and is also complete, so if we can prove that it is regular, then it must converge to the Hochschild cohomology.

First, note that $\mu^s = 0$ in $\mathscr{A}$ unless $s-2$ is divisible by $(n-2)$ (by Equation (\ref{eqn:ccr3}) with $a=n$ and $s+t=2$). 
In particular, $\mathscr{A}$ and hence $\widetilde{\mathscr{A}}_{nov}$ are minimal. 
So the spectral sequence has $E_1$ page
\[ CC^*_{c,\bm{G}_{\Z}}\left(\widetilde{\mathscr{A}}_{nov}\right)^{\tilde{\Gamma}_n^*}.\]

We observe that there is an obvious morphism of chain complexes,
\[ CC^*_{\bm{G}_{\Z}} \left(\bm{q}_{1*} \bm{p}_1^* \underline{\mathscr{A}} \right) \otimes_R \Lambda \To CC^*_{\bm{G}_{\Z}}\left(\widetilde{\mathscr{A}}_{nov}\right),\]
and indeed the $A_{\infty}$ structure maps on $\widetilde{\mathscr{A}}_{nov}$ are the image of the $A_{\infty}$ structure maps of $\mathscr{A}$ under this morphism. 
However, this is not necessarily a quasi-isomorphism, because the Hochschild cochain complex is defined as a direct product over cochains of all lengths $s \ge 0$, and arbitrary direct products do not commute with $\otimes_R \Lambda$.

On the other hand, finite products do commute with $\otimes_R \Lambda$, so the above morphism of chain complexes does induce an isomorphism on compactly-supported Hochschild cochains:
\[ CC^*_{c,\bm{G}_{\Z}} \left(\bm{q}_{1*} \bm{p}_1^* \underline{\mathscr{A}} \right) \otimes_R \Lambda \cong CC^*_{c,\bm{G}_{\Z}}\left(\widetilde{\mathscr{A}}_{nov}\right).\] 

Now, applying Remarks \ref{remark2:pushcat} and \ref{remark2:pullcc}, we have isomorphisms
\begin{eqnarray*}
CC^*_{c,\bm{G}_{\Z}} \left( \bm{q}_{1*} \bm{p}_1^* \underline{\mathscr{A}} \right)^{\tilde{\Gamma}_n^*} & \cong &  \bm{q}_{1*} CC^*_{c,\bm{G}^n_n}\left( \bm{p}_1^* \underline{\mathscr{A}} \right)^{\tilde{\Gamma}_n^*}  \\
& \cong & \bm{q}_{1*} \bm{p}_1^* CC^*_{c,\bm{G}}\left( \underline{\mathscr{A}} \right) \\
& \cong & \bm{q}_{1*} \bm{p}_1^* CC^*_{c,\bm{G}} \left( \mathscr{A} \right).
\end{eqnarray*}

It follows that the spectral sequence has $E_1$ page
\[  \bm{q}_{1*} \bm{p}_1^* CC^*_{c,\bm{G}} \left( \mathscr{A}\right) \otimes_R \Lambda.\]

The $E_2$ page is
\[  \bm{q}_{1*} \bm{p}_1^* HH^*_{c,\bm{G}} \left( \mathscr{A}, \mu^2\right) \otimes_R \Lambda.\]
We show that there is an algebra isomorphism $(\mathscr{A}, \mu^2_0) \cong (\mathscr{A},\mu^2)$. 
I.e., the higher-order terms in the product $\mu^2$ can be absorbed into the order-$0$ product $\mu^2_0$, which we know to be the exterior product.
This follows immediately from the calculation
\[ TrHH^2_{\bm{G}}(A, A \otimes R)^{0,H} \cong 0,\]
which was carried out in the proof of Lemma \ref{lemma2:thhversal} (as `Case 1'), together with a deformation theory argument. 
Thus, we can replace $\mathscr{A}$ by a quasi-isomorphic $A_{\infty}$ algebra $\mathscr{A}'$, whose underlying algebra is the exterior algebra
\[ \left(\mathscr{A}',\mu^2\right) \cong A \otimes R.\]
Furthermore, $\mathscr{A}'$ and $\mathscr{A}$ coincide to order $n$, because the higher-order terms in $\mu^2$ were functions of $T = r_1 \ldots r_n$. 
So $\mathscr{A}'$ is still of type A. 
Henceforth we'll simply write $\mathscr{A}$ for this replacement.

It now follows by the HKR isomorphism that the $E_2$ page of the spectral sequence is
\begin{eqnarray*}
 \bm{q}_{1*} \bm{p}_1^* HH^*_{c,\bm{G}}(A \otimes R) \otimes_R \Lambda &\cong & \bm{q}_{1*} \bm{p}_1^* (R[U] \otimes A) \otimes_R \Lambda \\
&\cong & \left( R[U] \otimes A \right)^{\tilde{\Gamma}_n^*} \otimes_R \Lambda \\
&\cong & \left(\Lambda[U] \otimes A \right)^{\tilde{\Gamma}_n^*}.
\end{eqnarray*}
The last step follows because the action of $\Gamma_n^*$ on $R$ is trivial

Because $\mu^s = 0$ for $2<s<n$, the first non-trivial differential is $\delta_{n-1}$, which is given by Gerstenhaber bracket with the order-$n$ deformation class $[\mu^n_{nov}]$. 
The Gerstenhaber bracket gets carried to the Schouten bracket under the HKR isomorphism, and  $\mu^n_{nov}$ gets carried to 
\[ w_{nov} := \Phi([\mu^n_{nov}]) \in   \left( \Lambda [U] \otimes A \right)^{\tilde{\Gamma}_n^*}.\]
By the `naturality' property of the HKR map (see Definition \ref{definition2:2phihkr}), $w_{nov} = w \otimes 1$, where 
\[w \in \left(R[U] \otimes A\right)^{\tilde{\Gamma}_n^*}\]
is the image of $[\mu^n]$ under the HKR map.

Now, by the grading computations in the proof of Lemma \ref{lemma2:thhversal}, $w$ has the form 
\[ w = f_1(T)u_1 \ldots u_n + f_2(T)\sum_{j=1}^n r_j u_j^n \in R[ u_1, \ldots, u_n ] ,\]
where $f_1, f_2 \in \C\llbracket T \rrbracket $, and we recall that $T := r_1 \ldots r_n$. 
By the definition of `type A', we have $f_1(0) = 1$ and $f_2(0) = 1$. 
It follows that 
\[ w_{nov} = f_1(r^n) u_1 \ldots u_n + f_2(r^n) \sum_{j=1}^n r u_j^n \in \Lambda [ u_1, \ldots, u_n ]^{\tilde{\Gamma}_n^*}.\]

In particular, $w_{nov}$ lies in $\Lambda[ U ]$, i.e., there are no non-trivial polyvector field terms appearing.
Therefore, as we saw in the proof of Lemma \ref{lemma2:thhversal}, 
\[ \delta_{n-1}(-) = \left[ \Phi(\mu^n), - \right] = \iota_{dw_{nov}}(-)\]
gives the Koszul complex for the sequence
\[ \del{w_{nov}}{u_1}, \ldots, \del{w_{nov}}{u_n} \in \Lambda [ U ].\]

Now we show that this sequence is regular. 
This follows because $w_{nov}$ has an isolated singularity at the origin. 
To see this, observe that we have relations
\begin{eqnarray*}
 f_1(r^n)\frac{u_1 \ldots u_n}{u_j} & \equiv & -n r f_2(r^n) u_j^{n-1}\\
 \Rightarrow \frac{u_1 \ldots u_n}{u_j} & \equiv & rf(r) u_j^{n-1},
\end{eqnarray*}
in the ring
\[ \Lambda [u_1, \ldots, u_n ]/(\partial_1 w_{nov}, \ldots, \partial_n w_{nov}), \]
where
\[ f(r) = -n \frac{f_2(r^n)}{f_1(r^n)} \in \C\llbracket r^n \rrbracket , \mbox{ $f(0) \neq 0$.}\]
Taking the product of these relations gives
\[(u_1 \ldots u_n)^{n-1} \equiv r^nf(r)^n (u_1 \ldots u_n)^{n-1},\]
and hence that
\[ (u_1 \ldots u_n)^{n-1} \equiv 0,\]
because $1 + \mathcal{O}(r)$ is invertible in $\Lambda$.
Returning to the original relation, we have
\[ (rf(r) u_j^n)^{n-1} = (u_1 \ldots u_n)^{n-1} = 0,\]
and hence sufficiently high powers of each generator $u_j$ vanish (recalling $f \neq 0$, so $f \in \Lambda^*$ is invertible). 
Therefore $w_{nov}$ has an isolated singularity at $0$, so the sequence is regular, so the cohomology of the Koszul complex is the Jacobian ring.

Now 
\[ \Lambda [ U ] ^{\tilde{\Gamma}_n^*} \]
is generated, as a $\Lambda$-algebra, by the elements
\[ u_1 \ldots u_n, u_1^n, \ldots, u_n^n,\]
and each of these has degree $2 \in \Z$.
To see this, note that for $u^{\bm{b}}$ to be invariant under the action of $\tilde{\Gamma}_n^*$, its degree must be in the image of $Y^n_n$ in $Y^n_1$. 
From the definition, this means that
\[ (|\bm{b}|, - \bm{b}) = (t + 2(n-1)|\bm{m}|, -n \bm{m}) - q(2(1-n), y_{[n]})\]
for some $t, q \in \Z$ and $\bm{m} \in \Z^n$, or in other words, setting $s = |\bm{b}|$ for the length,
\begin{eqnarray*}
s & = & t + 2(n-1)|\bm{m}| + 2(n-1)q,\\
\bm{b} &=& n \bm{m} + q y_{[n]}.
\end{eqnarray*}
It is not hard to show that, since $\bm{b} \ge \bm{0}$, we can arrange that $\bm{m} \ge \bm{0}$ and $q \ge 0$, so $u^{\bm{b}}$ is a product of the generators $u_1 \ldots u_n, u_1^n, \ldots, u_n^n$, as claimed. 
Furthermore, one can show that $s+t = 2q + 2|\bm{m}|$, where $s = |\bm{b}|$ is the length and $t$ is the $\Z$-grading. 
It follows that the degree of each generator is $2$, when regarded as an element of Hochschild cohomology.

It follows that the $E_n$ page of our spectral sequence is given by
\[ \Lambda [ u_1, \ldots,u_n ] /(\partial_1 w_{nov}, \ldots, \partial_n w_{nov}) \cap \Lambda [ u_1 \ldots u_n, u_1^n, \ldots, u_n^n ] .\]
From the relations in the Jacobian ring, we have
\[ u_1 \ldots u_n = r f(r) u_j^n \mbox{ where $f(r) \neq 0$,}\]
so we only need a single generator $\alpha := u_1 \ldots u_n$, and furthermore this generator satisfies 
\[ \alpha^{n-1} = 0,\]
as we showed above. 
It is easy to check using Gr\"{o}bner bases that $\alpha^{n-2}$ does not lie in the ideal generated by the partial derivatives $\partial_j w_{nov}$.
It follows that the $E_n$ page of the spectral sequence is given by
\[ \Lambda[\alpha]/\alpha^{n-1}.\]
Because $\alpha$ has degree $2 \in \Z$, $E_n$ is graded in even degrees so all subsequent differentials in the spectral sequence vanish, so the spectral sequence is regular. 
Since the length filtration is bounded above and complete, it follows by \cite[Theorem 5.5.10]{Weibel1994} that the spectral sequence converges to the Hochschild cohomology. 
Thus the $E_n$ page is isomorphic to the associated graded algebra of the Hochschild cohomology; but since it is one-dimensional in each degree, it is in fact isomorphic to the Hochschild cohomology. 
This completes the proof.
 \end{proof}

Now we observe that the element 
\[ r_j \del{\mu^*}{r_j} \in CC^2_{\bm{G}^n_n} \left(\widetilde{\mathscr{A}}\right)\]
is a Hochschild cochain (this follows by applying $r_j \partial / \partial r_j$ to the $A_{\infty}$ relation $\mu \circ \mu = 0$). 
Hence it defines an element in $HH^*(\widetilde{\mathscr{A}})$. 
We denote by $\beta$ the element of $HH^*(\widetilde{\mathscr{A}}_{nov})$ that is the image of
\[ \left(  r_j \del{\mu^*}{r_j} \right) \otimes 1 \in HH^* \left(\widetilde{\mathscr{A}}\right) \otimes_R \Lambda\]
under the obvious map.

\begin{lemma2}
\label{lemma2:defom}
The element $\beta$ lies in the $\tilde{\Gamma}_n^*$-equivariant part of $HH^*(\widetilde{\mathscr{A}}_{nov})$, and corresponds to $g \cdot \alpha$ for some invertible $g \in \Lambda^*$ under the isomorphism of Lemma \ref{lemma2:hhtypea}.
\end{lemma2}
\begin{proof}
The fact that $\beta$ is $\tilde{\Gamma}_n^*$-equivariant follows immediately from the fact that $\mu^*$ is. 
We recall from the proof of Lemma \ref{lemma2:hhtypea} that, firstly, we arrange that $\mu^2$ is independent of $r_j$ and, secondly, the image of $\mu^{\ge 3}$ under the HKR map to $R\llbracket U \rrbracket  \otimes A$ has the form
\[ w = f_1(T) u_1 \ldots u_n + f_2(T) \sum_{j=1}^n r_j u_j^n.\]
It follows that the image of $r_j \partial \mu^*/ \partial r_j$ under the HKR map has the form 
\[ r_j \del{w}{r_j} = Tf_1'(T) u_1 \ldots u_n + f_2(T) r_j u_j^n + Tf_2'(T) \sum_{j=1}^n r_j u_j^n.\]
It follows from naturality of the HKR map that the image of $ \beta = (r_j \partial \mu^* / \partial r_j) \otimes 1$ on the $E_2$ page of the spectral sequence has the form
\[ r^n f_1'(r^n) u_1 \ldots u_n + f_2(r^n)  ru_j^n  + r^n f_2'(r^n) \sum_{j=1}^n r u_j^n.\]

We now recall from the proof of Lemma \ref{lemma2:hhtypea} that the $E_n$ page of the spectral sequence is the Jacobian ring \[\Lambda [ u_1, \ldots, u_n ] /(\partial_1 w_{nov}, \ldots, \partial_n w_{nov}) \cap \Lambda [u_1 \ldots u_n, u_1^n, \ldots, u_n^n ].\]
We recall (from the proof of Lemma \ref{lemma2:hhtypea}) that, in the Jacobian ring, we have relations
\[ f_1(r^n) u_1 \ldots u_n \equiv -n  f_2(r^n) r u_j^n,\]
and we set $\alpha := u_1 \ldots u_n$, so the image of $\beta$ on the $E_n$ page of the spectral sequence is equivalent to
\[ \left(r^n f_1'(r^n) - \frac{1}{n}f_1(r^n) - r^n \frac{f_2'(r^n)f_1(r^n)}{f_2(r^n)} \right) \alpha = g(r^n) \alpha,\]
where $g(r^n) \in \C\llbracket r^n \rrbracket $ and $g(0) \neq 0$.

Recalling that the spectral sequence degenerates at the $E_n$ page, this completes the proof.
 \end{proof}

\subsection{First-order deformations}
\label{subsec:2firstod}

In this section, we will consider a very specific situation. 
Let $\mathcal{A} = (A,\mu_0)$ be a $\bm{G}$-graded minimal $A_{\infty}$ algebra over $\C$, and $V$ a $\bm{G}$-graded vector space. 
We consider the $\bm{G}$-graded $\C$-algebra
\[ R = \C[V]/ \mathfrak{m}^2,\]
where $\mathfrak{m} \subset \C[V]$ is the maximal ideal corresponding to $0$.
A $\bm{G}$-graded {\bf first-order deformation} of $\mathcal{A}$ over $V$ is an $R/\mathfrak{m}^2$-linear $A_\infty$ algebra $\mathscr{A}$ such that $\mathscr{A}/\mathfrak{m}$ is strictly isomorphic to $\mathcal{A}$. 

Equivalently, it is an element
\[ \mu = \mu_0 + \mu_1 \in CC_{\bm{G}}^2(A \otimes R,A \otimes R),\]
whose order-$0$ component agrees with $\mu_0$, and such that
\[ \mu_1 \in CC_{\bm{G}}^2(\mathcal{A},\mathcal{A} \otimes V^{\vee})\]
is a Hochschild cocycle. 
The class $[\mu_1] \in HH_{\bm{G}}^2(\mathcal{A},\mathcal{A} \otimes V^{\vee})$ is called the {\bf deformation class} of the deformation.

\begin{definition2}
\label{definition2:similar}
Suppose we are given two $\bm{G}$-graded, $R/\mathfrak{m}^2$-linear $A_\infty$ algebras $\mathscr{A}$ and $\mathscr{B}$, which are first-order deformations of $\mathcal{A} = \mathscr{A}/\mathfrak{m}$ and $\mathcal{B} = \mathscr{B}/\mathfrak{m}$, respectively, over $V$. 
$\mathscr{A}$ and $\mathscr{B}$ are said to be {\bf first-order quasi-equivalent} if there is a quasi-isomorphism of $\C$-linear $A_\infty$ algebras,
\[ \mathcal{A} \cong \mathcal{B},\]
and the resulting isomorphism of Hochschild cohomology groups,
\[ HH^2_{\bm{G}}(\mathcal{A},\mathcal{A} \otimes V^\vee) \cong HH^2_{\bm{G}}(\mathcal{B},\mathcal{B} \otimes V^\vee),\]
equates the first-order deformation classes of $\mathscr{A}$ and $\mathscr{B}$. 
\end{definition2}

\begin{lemma2}
\label{lemma2:similar}
Let $\mathscr{C}$ be a $\bm{G}$-graded, $R/\mathfrak{m}^2$-linear $A_\infty$ category, and let $\mathcal{C} := \mathscr{C}/\mathfrak{m}$. 
Suppose that objects $X$ and $Y$ are quasi-isomorphic in the $\C$-linear $A_\infty$ category $\mathcal{C}$; then the $\bm{G}$-graded, $R/\mathfrak{m}^2$-linear $A_\infty$ algebras
\[ \mathscr{A}:= hom^*_{\mathscr{C}}(X,X), \mathscr{B} := hom^*_\mathscr{C}(Y,Y)\]
 are first-order quasi-equivalent (in the sense of Definition \ref{definition2:similar}).
\end{lemma2}
\begin{proof}
Let $\mathcal{C}_X$ be the full subcategory of $\mathcal{C}$ with the single object $X$; define $\mathcal{C}_Y$ similarly. 
It is standard that the inclusions
\[ \mathcal{C}_X \hookrightarrow \mathcal{C} \hookleftarrow \mathcal{C}_Y\]
are quasi-equivalences; it follows that there is a quasi-isomorphism
\[ hom^*_\mathcal{C}(X,X) \cong hom^*_\mathcal{C}(Y,Y).\]

Furthermore, it follows by Morita invariance (see, for example, \cite[Lemma 2.6]{Seidel2008}) that the restriction maps
\[ HH_{\bm{G}}^*(\mathcal{C}_X, \mathcal{C}_X \otimes V^\vee) \leftarrow HH_{\bm{G}}^*(\mathcal{C},\mathcal{C}\otimes V^\vee) \rightarrow HH_{\bm{G}}^*(\mathcal{C}_Y,\mathcal{C}_Y\otimes V^\vee)\]
are isomorphisms. 
The deformation class of $\mathscr{C}$ lies in the middle Hochschild cohomology group, and restricts to the deformation classes of $\mathscr{C}_X, \mathscr{C}_Y$ on the left and right, respectively. 
It follows that the deformation classes of $\mathscr{A} = hom^*_\mathscr{C}(X,X)$ and $\mathscr{B} = hom^*_\mathscr{C}(Y,Y)$ match up under Morita invariance; therefore $\mathscr{A}$ and $\mathscr{B}$ are first-order quasi-equivalent.

\end{proof}

\section{The affine Fukaya category}
\label{sec:2affeqi}

In this section, we introduce the affine Fukaya category $\mathcal{F}(M)$ of an exact symplectic manifold, and explain its relationship with the exact Fukaya category as defined in \cite{Seidel2008}, which we will denote by $\mathcal{F}'(M)$. 
The difference is essentially that the affine Fukaya category has fewer objects, but a richer grading structure.

\subsection{Grading data from the Lagrangian Grassmannian}
\label{subsec:2hmindex}

We recall some notions from \cite[Chapters 11, 12]{Seidel2008}. 
Let $M$ be a symplectic manifold. 
We denote by $\mathcal{G}M$ the bundle of Lagrangian subspaces of $TM$. 
Observe that, because we have a fibration
\[\begin{diagram}
\mathcal{G}_pM& \rInto & \mathcal{G}M \\
& & \dTo \\ 
&  & M,
\end{diagram}
\]
there is an associated exact sequence 
\[ \ldots \To \pi_2(M) \To \pi_1( \mathcal{G}_pM) \To \pi_1( \mathcal{G}M) \To \pi_1(M) \To *.\]
Now because abelianization is right-exact in an appropriate sense, and 
\[ \pi_1(\mathcal{G}_pM) \cong H_1(\mathcal{G}_pM) \cong \Z\]
is abelian (where the isomorphism is given by the Maslov index), the sequence
\[ \pi_2(M) \To  H_1(\mathcal{G}_pM) \To H_1(\mathcal{G}M) \To H_1(M) \To 0\]
is exact.

\begin{definition2}
\label{definition2:2hm}
We define a grading datum $\bm{G}(M)$:
\[\begin{diagram}
 \Z & \rTo^{f} & Y(M) & \rTo^{g} & X(M) & \rTo & 0 \\
\dEq & & \dEq & & \dEq && \\
H_1(\mathcal{G}_pM) &  \rTo &  H_1( \mathcal{G}M) & \rTo &  H_1(M) & \rTo & 0.
\end{diagram}
\]
To define the sign morphism $\bm{\sigma}$, we must specify a map
\[ \sigma: H_1(\mathcal{G}M) \To \Z_2.\]
To do this, we consider the real vector bundle $\mathscr{L} \To \mathcal{G}M$, whose fibre over a point is identified with the Lagrangian subspace at that point. 
The first Stiefel-Whitney class defines an element
\[ w_1(\mathscr{L}) \in H^1\left(\mathcal{G}M; \Z_2\right),\]
and $\sigma$ is defined by pairing with this element. 
\end{definition2}

We will see that it is natural to define the Fukaya category as a $\bm{G}(M)$-graded category. 
That is because of the relationship between $\bm{G}(M)$ and index theory of Cauchy-Riemann operators in $M$, which we now explain.

\begin{remark2}
One could define a non-abelian grading datum to be a morphism $\Z \To Y$, where now $Y$ is allowed to be non-abelian, and rework the theory of grading data (Section \ref{sec:2deftheory}) appropriately. 
This could be applied to study Fukaya categories of symplectic manifolds with non-abelian fundamental group.
However, apart from being more complicated, when one studies the the relative Fukaya category (Section \ref{sec:2relfuks}), it becomes absolutely necessary to consider abelian grading data. 
At any rate, all of the manifolds we consider in this paper have abelian fundamental group.
\end{remark2}

We recall (from \cite[Appendix C.3]{mcduffsalamon}) the definition of a {\bf bundle pair} $(\Sigma,E,F)$ over a Riemann surface with boundary $\Sigma$. 
It is a complex vector bundle $E \To \Sigma$ together with a totally real subbundle $F \subset E|_{\partial \Sigma}$. 
We recall that a bundle pair defines a Cauchy-Riemann operator
\[ D: \Omega^0_F(\Sigma,E) \To \Omega^{0,1}_F (\Sigma,E),\]
whose zeroes are holomorphic sections of $E$ whose boundary values lie in $F$.
We recall also the {\bf boundary Maslov index} of a bundle pair, which is an integer $\mu(\Sigma,E,F)$ such that the index of (an appropriate Sobolev-space version of) the Cauchy-Riemann operator $D$ associated to the bundle pair $(\Sigma,E,F)$ is
\[ \mathrm{ind}(D) = n \chi(\Sigma) + \mu(\Sigma,E,F),\]
where $n$ is the complex dimension of a fibre of $E$. 

Now suppose we are given a map $u: \Sigma \To M$, together with a Lagrangian subbundle $F \subset u^*TM|_{\partial \Sigma}$. 
$F$ defines a lift 
\[\begin{diagram}
&  & \mathcal{G}M \\
& \ruTo^{\rho} & \dTo \\ 
\partial \Sigma & \rTo_{u|_{\partial \Sigma}}  & M,
\end{diagram}
\]
which defines a class $[\rho] \in H_1(\mathcal{G}M) = Y(M)$, which is a lift of
\[ u|_{\partial \Sigma} = \partial u = 0 \in H_1(M) = X(M).\] 
So by exactness of the sequence $\bm{G}(M)$, $[\rho]$ lies in the image of $\Z$. 

\begin{lemma2}
\label{lemma2:maslovgrad}
We have
\[ [\rho] = f(\mu(\Sigma,u^*TM,F)),\]
in $Y(M)$, where $f: \Z \To Y(M)$ comes from the grading datum.
\end{lemma2}
\begin{proof}
Define a decomposition of bundle pairs,
\[(\Sigma, u^*TM,F) = (\Sigma_1,u^*TM,F) \cup (\Sigma_2,u^*TM,F),\]
 where $\Sigma_1$ is a small ball in the interior of $\Sigma$, $\Sigma_2$ is its complement, and the totally real subbundle of $u^*TM|_{\partial \Sigma_1}$ is defined by a lift
\[ \rho': \partial \Sigma_1 \To \mathcal{G}M,\]
chosen in such a way that there is a trivialization
\[ (\Sigma_2,u^*TM,F) \cong (\Sigma_2, \Sigma_2 \times \C^n, \Sigma_2 \times \R^n).\]
It follows quickly from the properties of the boundary Maslov index (see \cite[Theorem C.3.5]{mcduffsalamon}) that 
\[ \mu(\Sigma_2,u^*TM,F) = 0.\]
Then, by the composition property of the boundary Maslov index,
\begin{eqnarray*}
 \mu(\Sigma,u^*TM,F) &=& \mu(\Sigma_1, u^*TM,F) + \mu(\Sigma_2,u^*TM,F) \\
&=& \mu(\Sigma_1,u^*TM,F).
\end{eqnarray*}
By our definition of $(\Sigma_2,u^*TM,F)$, there is a lift
\[\begin{diagram}
&  & \mathcal{G}M \\
& \ruTo^{P} & \dTo \\ 
 \Sigma_2 & \rTo_{u|_{\Sigma_2}}  & M,
\end{diagram}
\] 
and $P$ defines a homology between the class $[\rho]$ and the class $[\rho']$.
If the ball $\Sigma_1$ is centred on a point $p$, then there is an obvious isomorphism
\[ u^*TM|_{\Sigma_1} \cong \Sigma_1 \times T_pM,\]
so $\rho'$ defines a class in $H_1(\mathcal{G}_pM) \cong \Z$, which is exactly equal to $\mu(\Sigma_1,u^*TM,F)$ (essentially from the definition of the Maslov class). 
The result follows.
 \end{proof}

\begin{remark2}
\label{remark2:masexact}
Consider the abelianization of the long exact sequence of homotopy groups of the fibration $\mathcal{G}M$:
\[ \ldots \To \pi_2(M) \To H_1(\mathcal{G}_pM) \To H_1(\mathcal{G}M) \To H_1(M) \To 0.\]
It follows from Lemma \ref{lemma2:maslovgrad} and the definition of the long exact sequence of homotopy groups that the map 
\[ \pi_2(M) \To H_1(\mathcal{G}_p M) \cong \Z\]
is given by evaluation of $2c_1(TM)$.
It follows that, if this map vanishes, then the grading datum $\bm{G}(M)$ is exact (recall that this means the morphism $\Z \To Y(M)$ is injective).
\end{remark2}

\subsection{Anchored branes}
\label{subsec:2eqbrane}

In this section, we will define the notion of an {\bf anchored Lagrangian brane} in $M$. 
These will be the objects of the affine Fukaya category $\mathcal{F}(M)$ (when $M$ is exact).
First we recall the notion of a (non-anchored) graded Lagrangian brane from \cite[Chapter 11]{Seidel2008}, which is an object of $\mathcal{F}'(M)$.

If $2c_1(M) = 0$, and $M$ is equipped with a quadratic complex volume form $\eta$, then we can construct a {\bf phase map}
\[ \alpha_M: \mathcal{G}M \To S^1.\]
Now if $i:L \To M$ is a Lagrangian immersion, there is a canonical lift of 
\[\begin{diagram}
& & \mathcal{G}M \\
& \ruTo^{i_*} & \dTo \\
L & \rTo_{i} & M.
\end{diagram}
\]
Hence, given $\eta$, there is a map 
\begin{eqnarray*}
\alpha_L: L &\To& S^1,\\
\alpha_L & := & \alpha_M \circ i_*.
\end{eqnarray*}
A (non-anchored) {\bf graded Lagrangian brane} in $M$ is a compact embedded Lagrangian $L \subset M$, together with a lift $\alpha^{\#}_L$ of $\alpha_L$ to $\R$, and a Pin structure on $L$.

Now we introduce a new notion. 
Let 
\[\pi: \widetilde{\mathcal{G}}M \To \mathcal{G}M\] 
denote the universal abelian cover of the manifold $\mathcal{G}M$ (i.e., the one associated to the commutator subgroup of $\pi_1(\mathcal{G}M)$). 

\begin{definition2}
An {\bf anchored Lagrangian brane} $L^{\#}$ in $M$ is a Lag\-rang\-i\-an immersion $i:L \To M$ of a compact manifold $L$ into $M$, together with a lift $i^{\#}$ as follows:
\[\begin{diagram}
& & \widetilde{\mathcal{G}}M \\
& \ruTo^{i^{\#}} & \dTo^{\pi} \\
L & \rTo_{i_*} & \mathcal{G}M,
\end{diagram}
\]
and a Pin structure. 
Note that we do {\bf not} need a quadratic complex volume form to define the notion of an anchored Lagrangian brane.
\end{definition2}

We observe that there is a natural action of the covering group $H_1(\mathcal{G}M) \cong Y(M)$ on $\widetilde{\mathcal{G}}M$, and hence on anchored Lagrangian branes. 
We denote the action of $y \in Y(M)$ on $L^{\#}$ by $y \cdot L^{\#}$.

Now we explain the relationship between anchored Lagrangian branes and (non-anchored) graded Lagrangian branes.
We observe that, given a quadratic complex volume form $\eta_M$, there exists a lift of the squared phase map
\[\begin{diagram}
\widetilde{\mathcal{G}}M & \rTo^{\tilde{\alpha}_M} & \R \\
\dTo & & \dTo\\
\mathcal{G}M & \rTo{\alpha_M} &  S^1.
\end{diagram}
\]

\begin{definition2}
\label{definition2:2forgetbrane}
Given a quadratic volume form $\eta_M$, a lift $\tilde{\alpha}_M$ as above, and an embedded anchored Lagrangian brane $L^{\#} = (L,i^{\#})$, we define a Lagrangian brane $\mathfrak{f}(L^{\#})$ in $M$ ($\mathfrak{f}$ for `forgetting' the anchored structure), with the same underlying Lagrangian $L$ and Pin structure, and 
\[ \alpha_L^{\#} = \tilde{\alpha}_M \circ i^{\#}.\]
\end{definition2}

\subsection{The affine Fukaya category}
\label{subsec:2eqaffuk}

Let $M$ be an exact symplectic manifold with convex boundary.
In this section, we define the {\bf affine Fukaya category}, $\mathcal{F}(M)$. 
It will be a $\bm{G}(M)$-graded $A_{\infty}$ category. 
The definition is very closely related to the definition of the exact Fukaya category $\mathcal{F}'(M)$, given in \cite[Section 12]{Seidel2008}, to which the reader is referred for all technical details. 
We will see that $\mathcal{F}(M)$ is essentially a full subcategory of $\mathcal{F}'(M)$ with a richer grading structure -- in particular, the analytic details of defining moduli spaces of pseudoholomorphic disks to define the $A_{\infty}$ structure maps are completely analogous.

Objects of $\mathcal{F}(M)$ are embedded anchored Lagrangian branes. 
For each pair of objects, we choose a Floer datum on $M$ (in the sense of \cite[Section 8e]{Seidel2008}), and for all moduli spaces of boundary-punctured holomorphic disks with boundary components labelled by anchored Lagrangian branes, we make a consistent universal choice of perturbation data on $M$ (in the sense of \cite[Section 9h]{Seidel2008}). 
We assume that the action of $Y(M)$ on anchored Lagrangian branes lifts to an action on Floer and perturbation data (i.e., if we change some of the boundary labels of a boundary-punctured holomorphic disk by the action of $Y(M)$, then the perturbation datum does not change).

\begin{sloppypar}
We now define the $\bm{G}(M)$-graded morphism spaces in $\mathcal{F}(M)$.
Now, given anchored Lagrangian branes $L_0^{\#}, L_1^{\#}$, we define the morphism space $CF^*(L_0^{\#},L_1^{\#})$ to be generated by paths $p: [0,1] \To M$ satisfying $p(0) \in L_0$, $p(1) \in L(1)$, which are flowlines of the Hamiltonian vector field associated with the corresponding Floer datum. 
\end{sloppypar}

Given such a $p$, we define its grading $y \in Y(M)$ to be the unique element such that $p$ lifts to a path from $L_0^{\#}$ to $y \cdot L_1^{\#}$ in $\widetilde{\mathcal{G}}M$, which has Maslov index $0 \in \Z$. 
To explain what this means, we observe that there is a commutative diagram
\[\begin{diagram}
\widetilde{\mathcal{G}}M & \rTo & \widetilde{M} \\
\dTo && \dTo\\
\mathcal{G}M & \rTo & M,
\end{diagram}
\]
where $\widetilde{M} \To M$ is the universal abelian cover of $M$. 
Thus, associated with any anchored Lagrangian brane $L^{\#}$ is a lift, $\widetilde{L}$, of $L$ to $\widetilde{M}$. 
The fact that $p$ must lift to $\widetilde{\mathcal{G}}M$ implies that it must lift to $\widetilde{M}$; this already defines $y \in Y(M)$ up to addition of an element in the image of $\Z \To Y(M)$.
Now we observe that the fibres of the bundle 
\[\widetilde{\mathcal{G}}M \To \widetilde{M}\]
are the universal covers of the fibres of the Lagrangian Grassmannian $\mathcal{G}\widetilde{M}$. 
Thus, the anchored brane structures $L_0^{\#}, L_1^{\#}$ equip $\widetilde{L}_0, \widetilde{L}_1$ with the structure of  `abstract Lagrangian branes' (see \cite[Section 12a]{Seidel2008}). 
Therefore, if the path $p$ lifts to a path $\tilde{p}$ from $\widetilde{L}_0$ to $y \cdot \widetilde{L}_1$ in $\widetilde{M}$, then we can define the Maslov index $i$ of any lift of $\tilde{p}$ to $\widetilde{\mathcal{G}}M$, and it is equal to the relative Maslov index of the abstract linear Lagrangian branes at either end of $\tilde{p}$. 
It is this index that we require to be $0$. 
Given $p$, it is clear that we can find $y' \in Y(M)$ such that the path $p$ lifts to a path $\tilde{p}$ from $\widetilde{L}_0$ to $y' \cdot \widetilde{L}_1$, but the Maslov index $i$ may not be zero.
However, we then necessarily have
\[ y = y' - f(i),\]
so the $Y(M)$-grading of $p$ is well-defined.

We define the $A_{\infty}$ structure maps in $\mathcal{F}(M)$ by counting rigid pseudo-holomorphic disks in $M$. 
That is, given objects 
\[L_0^{\#}, \ldots, L_s^{\#}\]
and morphisms
\[ p_j \in CF^*(L_{j-1}^{\#},L_j^{\#}) \mbox{ for $j =1,\ldots,s$}\]
and
\[ p_0 \in CF^*(L_0^{\#},L_s^{\#}),\]
we define the coefficient of $p_0$ in $\mu^s(p_s,\ldots,p_1)$ to be the count of rigid pseudolomorphic disks in $M$ with boundary conditions on $L_j$, asymptotic to the generators $p_j$.
 
We now explain why these structure maps are $\bm{G}(M)$-graded. 
Firstly, observe that the structure maps respect the action of $Y(M)$ on objects, because we chose the perturbation data to do so (compare Definition \ref{definition2:gprecat}).
From \cite[Section 11l]{Seidel2008}, we recall the definition of an orientation operator $D_{p}$ corresponding to a generator $p$ of $CF^*(L_0^{\#},L_1^{\#})$. 
We lift $p$ to a path 
\[ \rho: [0,1] \To \widetilde{\mathcal{G}}M\]
connecting $L_0^{\#}(p(0))$ to $y \cdot L_1^{\#}(p(1))$, where $y \in Y(M)$ is the degree of $p$.
Now define a smooth, non-decreasing function $\psi: \R \to [0,1]$ such that $\psi(s) = 0$ for $s \ll 0$ and $\psi(s) = 1$ for $s \gg 0$. 
We consider the Hermitian vector bundle over the upper half plane $\R \times \R_{\ge 0}$, with fibre over $(s,t)$ given by $T_{p(\psi(s))}M$.
We introduce Lagrangian boundary conditions along the real axis, given by $\rho(\psi(s))$.
These Lagrangian boundary conditions define a Cauchy-Riemann operator, which we denote by $D_{p}$.
$D_{p}$ is Fredholm, its index is equal to the relative Maslov index of the abstract Lagrangian branes at either end of $p$, which is $0$ by the definition of $y$, and its determinant line is canonically isomorphic to the orientation line $o_p$ of $p$.

Given a holomorphic disk $u$ contributing to an $A_{\infty}$ product $\mu^s$, we denote the linearized operator of the pseudoholomorphic holomorphic curve equation at $u$ (with fixed domain $S$) by $D_{S,u}$. 
It is a Cauchy-Riemann operator on the trivial Hermitian vector bundle $u^* TM$ over $S$.
We can glue the orientation operators $D_{p_1}, \ldots, D_{p_s}$ and $D_{p_0}^{\vee}$ to $D_{S,u}$ along the strip-like ends to obtain a new Cauchy-Riemann operator over the closed disk.
We denote this operator by $\overline{D}$. 
The gluing formula then implies that
\begin{eqnarray*}
i(\overline{D}) &=& i(D_{S,u}) + i(p_1) + \ldots + i(p_s) +(n- i(p_0))\\
& = & i(D_{S,u}) + n,
\end{eqnarray*}
and there is a canonical isomorphism
\[\mathrm{det}(\overline{D}) \cong \mathrm{det}(D_u) \otimes o_{p_1} \otimes \ldots \otimes o_{p_s} \otimes o_{p_0}^{\vee}.\]
Now the Cauchy-Riemann operator $\overline{D}$ is given by a bundle pair $(D^2,E,F)$, which is equivalent to a bundle pair $(D^2,u^*TM,F)$, where $u:D^2 \To M$ is obtained from the original disk $u$ (which had strip-like ends converging to the generators $p_j$) by gluing the orientation operators onto the ends.
The boundary conditions for $\overline{D}$ define a map $\rho: \partial D^2 \To \mathcal{G}M$ which lifts the boundary map $\partial u$.

If we think of $\rho$ as a map
\[ \rho: [0,1] \To \mathcal{G}M\]
such that $\rho(0) = \rho(1)$ lies on $L_0$, then we have a lift
\[ \begin{diagram}
&& \widetilde{\mathcal{G}}M \\
& \ruTo^{\tilde{\rho}} & \dTo \\
[0,1] & \rTo^{\rho} & \mathcal{G}M.
\end{diagram}\]
The boundary conditions $L_j$ lift to $(y_1 + \ldots + y_j) \cdot L_j^{\#}$, and finally $\rho(1)$ lands on $(y_1 + \ldots + y_s - y_0) \cdot L_0^{\#}$. 
It follows that
\[ [\rho] = -y_0 + \sum_{j=1}^s y_j.\]
Lemma \ref{lemma2:maslovgrad} now implies that
\begin{eqnarray*}
 f\left(i(\overline{D}) - n \chi \left( D^2 \right)\right) &=&  -y_0 + \sum_{j=1}^s y_j \\
\Rightarrow f(i(D_{S,u})) &= & -y_0 + \sum_{j=1}^s y_j \mbox{ in $Y(M)$.}
\end{eqnarray*}
We now recall that, for the disk to be rigid, the {\bf extended} linearized operator $D_u$ (in which the modulus of the domain is allowed to vary, as well as the map) should have index zero. 
The dimension of the moduli space of disks with $s+1$ marked boundary points is $s-2$, so this means that for a rigid disk,
\[ y_0 =  f(2-s) + \sum_{j=1}^s y_j\]
in $Y(M)$.
It follows that the affine Fukaya category is a $\bm{G}(M)$-graded $A_{\infty}$ category (see Remark \ref{remark2:ainfgrad}).
We observe that the $A_{\infty}$ associativity equations are satisfied, by the same argument as for $\mathcal{F}'(M)$ (\cite[Proposition 12.3]{Seidel2008}).

We will now explain how $\mathcal{F}(M)$ is related to the $\Z$-graded exact Fukaya category $\mathcal{F}'(M)$, as defined in \cite{Seidel2008}.
We recall that, to define $\Z$-gradings on $\mathcal{F}'(M)$, we require that $2c_1(M) = 0$ and equip $\eta$ with a quadratic volume form $\eta_M$.

Recall from Definition \ref{definition2:2forgetbrane} that, if we equip $M$ with a quadratic volume form $\eta_M$, then we obtain a squared phase map
\[ \alpha_M: \mathcal{G}M \To S^1,\]
and if we define a lift 
\[ \tilde{\alpha}_M: \widetilde{\mathcal{G}}M \To \R,\]
then we obtain a forgetful map $\mathfrak{f}$ from anchored Lagrangian branes to (non-anchored) Lagrangian branes. 

Now on the level of $H_1$, $\alpha_M$ induces a map
\[ Y(M) \To \Z.\]
This defines a morphism of grading data, $\bm{p}^{\eta}: \bm{G}(M) \To \bm{G}_{\Z}$. 
It follows that $\bm{p}^{\eta}_* \mathcal{F}(M)$ is a $\Z$-graded $A_{\infty}$ category. 
We have:

\begin{lemma2}
\label{lemma2:eqaffemb}
The forgetful map $\mathfrak{f}$ on objects extends to a fully faithful embedding of $\Z$-graded $A_{\infty}$ categories,
\[ \mathfrak{f}: \bm{p}^{\eta}_* \mathcal{F}(M) \To \mathcal{F}'(M).\]
\end{lemma2}

\begin{remark2}
The image of this embedding consists of all (non-anchored) Lagrangian branes $L$ such that the image of $H_1(L)$ in $H_1(M)$ is zero.
\end{remark2}

\subsection{Covers}
\label{subsec:2affcov}

We explain how the affine Fukaya category behaves with respect to finite covers (essentially following \cite[Section 8b]{Seidel2003}). 
Suppose that $M,N$ are exact symplectic manifolds with convex boundary, with assumptions as in Section \ref{subsec:2eqaffuk}, and $\phi: N \To M$ is an exact symplectic covering (i.e., a covering such that the Liouville one-form on $N$ is pulled back from that on $M$ via $\phi$). 

\begin{condition}
\label{condition:cov}
We assume that the covering group of $\phi$ is abelian, and the induced map
\[ \phi_*: H_1(N) \To H_1(M)\]
is injective.
\end{condition}

Then we have an induced covering
\[\begin{diagram}
\mathcal{G}N & \rTo & N \\
\dTo<{\phi_*} && \dTo>{\phi} \\
\mathcal{G}M & \rTo &M,
\end{diagram}
\]
and hence an injective morphism of grading data, $\bm{p}: \bm{G}(N) \To \bm{G}(M)$, given by
\[\begin{diagram}
H_1(\mathcal{G}_pN) & \rTo & H_1(\mathcal{G}N) & \rTo & H_1(N) & \rTo & 0 \\
\dEq^{\phi_*} && \dTo^{\phi_*} && \dTo^{\phi_*} && \\
H_1(\mathcal{G}_{\phi(p)}M) & \rTo & H_1(\mathcal{G}M) & \rTo & H_1(M) & \rTo & 0.
\end{diagram}
\]

\begin{proposition2}
\label{proposition2:2affcovers}
If Condition \ref{condition:cov} holds, then there is a fully faithful embedding of $\bm{G}(N)$-graded categories,
\[ \bm{p}^*(\mathcal{F}(M)) \hookrightarrow \mathcal{F}(N).\]
\end{proposition2}
\begin{proof}
Note that $\phi$ induces an inclusion of commutator subgroups,
\[ [ \pi_1(N),\pi_1(N)] \hookrightarrow [ \pi_1(M),\pi_1(M)].\]
Because the covering group is abelian by Condition \ref{condition:cov}, $[ \pi_1(M),\pi_1(M)]$ lies in the image of $\pi_1(N)$. 
Because the map
\[ \pi_1(N)/[\pi_1(N) , \pi_1(N)] \To \pi_1(M)/[\pi_1(M), \pi_1(M)]\]
is injective by Condition \ref{condition:cov}, this means that $\phi$ induces an isomorphism of commutator subgroups,
\[ [ \pi_1(N),\pi_1(N)] \cong [ \pi_1(M),\pi_1(M)].\]
It follows that $\phi_*$ induces an isomorphism of commutator subgroups of $\mathcal{G}N$ and $\mathcal{G}M$.
It follows that the universal abelian covers of $\mathcal{G}N$ and $\mathcal{G}M$ are isomorphic, so after choosing such an isomorphism, there is a bijective correspondence between anchored Lagrangian branes in $N$ and in $M$.

One can similarly set up a correspondence between morphism spaces and moduli spaces of pseudoholomorphic disks defining the $A_{\infty}$ structure maps, and show that the gradings correspond, so the categories are strictly equivalent. 
The only difference is that some Lagrangians which are embedded in $N$ may not be embedded when projected to $M$.
 \end{proof}

\subsection{The relative case}
\label{subsec:2relcase}

Now we specialize to a particular type of exact symplectic manifold with convex boundary.

\begin{definition2}
\label{definition2:2kahlp}
A {\bf K\"{a}hler pair} $(M,D)$ consists of:
\begin{itemize}
\item A compact complex manifold $M$.
\item A smooth normal-crossings divisor
\[ D = \bigcup_{j=1}^k D_j,\]
with smooth irreducible components $D_j \subset M$. 
\item A cohomology class $c \in H^2(M)$, and integers $d_j$, such that
\[ P.D.([D_j]) = d_j c \mbox{ for all $j$}\]
(where $P.D.([D_j])$ is the cohomology class Poincar\'{e} dual to $D_j$). We call $d_j$ the {\bf multiplicity} of $D_j$.
\item A K\"{a}hler form $\omega$ on $M$.
\item On $M \setminus D$ (which we call the {\bf affine part} of $M$), a K\"{a}hler potential $h$, so that
\[ \omega = dd^c h,\]
and $h$ is proper.
\item It follows, in particular, that $\alpha := d^ch$ is a primitive for $\omega|_{M \setminus D}$. 
We call it the {\bf Liouville one-form}. 
\end{itemize} 

We define the {\bf linking number} of $\alpha$ with $D_j$: for any disk $y: D^2 \To M$ embedded so that it intersects the divisor $D_j$ positively and transversely at the origin and nowhere else, and doesn't intersect any of the other divisors, the linking number of $\alpha$ with $D_j$ is
\[ l_j := - \lim_{\delta \To 0} \alpha(y(\gamma_\delta)),\]
where $\gamma_\delta$ is the (positively-oriented) circle of radius $\delta$ in $D^2$. 
One easily checks that this definition is independent of choices. 
It follows from the assumption of properness of $h$ that each  $l_j$ is positive. 
Note that the linking number is called the `wrapping number' in \cite[Lemma 5.17]{McLean2012} (which uses the opposite sign convention); we prefer `linking' to `wrapping' to avoid confusion with the wrapped Fukaya category.

We introduce some standing notation: for any subset $K \subset \{1, \ldots, k\}$ we denote
\[ D_K := \bigcap_{j \in K} D_j.\]
\end{definition2}

\begin{lemma2}
\label{lemma2:stokes}
Let $\Sigma$ be a two-dimensional surface with boundary $\partial \Sigma$. 
For any smooth map $u:(\Sigma,\partial \Sigma) \To (M,M \setminus D)$, we have
\[ \omega(u) = \alpha(\partial u) + \sum_{j=1}^k l_j (u \cdot D_j).\]
Applying this in the case $\partial \Sigma = \emptyset$, it follows that the cohomology class of the K\"{a}hler form $\omega$ is
\[ [\omega] = \left(\sum_{j=1}^k d_j l_j \right) c.\]
\end{lemma2}
\begin{proof}
Suppose $\Sigma$ intersects the divisors $D_j$ transversely. 
Let $B \subset \Sigma$ be a union of small balls surrounding each intersection point with the divisor $D$. 
Applying  Stokes' theorem to $\Sigma \setminus B$ shows that
\begin{eqnarray*}
\omega(u) & \approx & \omega\left(u|_{ \Sigma \setminus B}\right) \\
& =& \alpha(\partial u) -\alpha(\partial B) \\
& \approx & \alpha(\partial u) + \sum_{j=1}^k l_j (u \cdot D_j).
\end{eqnarray*}
Now take the limit that the size of $B$ goes to zero.

\end{proof}

\begin{example2}
\label{example2:kahlerline}
We describe an important source of examples of K\"{a}hler pairs, following \cite[Chapter I]{Griffiths1978}.
Suppose we have a complex manifold $M$, equipped with a positive line bundle $\mathcal{L}$, and a smooth normal-crossings divisor $D = \bigcup D_j$, such that each component $D_j$ is the (transverse) vanishing locus of a section $s_j$ of $\mathcal{L}^{\otimes d_j}$. 
Define $c := c_1(\mathcal{L})$. 
Then we have $P.D.([D_j]) = d_j c$, i.e., the multiplicity of $D_j$ is $d_j$.
If the line bundle $\mathcal{L}^{\otimes d_j}$ is equipped with a Hermitian metric $\|\cdot\|$, then we can define a K\"{a}hler form $\omega_j$ which is the curvature of the corresponding connection, with cohomology class 
\[ [\omega_j] = c_1\left(\mathcal{L}^{\otimes d_j}\right) = d_j c.\] 
We can define a K\"{a}hler potential for this K\"{a}hler form on $M \setminus D_j$, by
\[ h_j := -\log \left(\| s_j\|^2\right),\]
so that $\alpha_j := d^ch_j$ is a primitive for $\omega_j|_{M \setminus D_j}$. 
One easily checks that $\alpha_j$ has linking number $1$ with $D_j$, because $s_j$ vanishes transversely along $D_j$.
Hence, if all of the line bundles $\mathcal{L}^{\otimes d_j}$ are equipped with Hermitian metrics, we can define the K\"{a}hler form 
\[ \omega := \sum_{j=1}^k l_j \omega_j,\]
for any choice of linking numbers $l_j>0$. 
This K\"{a}hler form has K\"{a}hler potential
\[ h := \sum_{j=1}^k l_j h_j.\] 
These have the required properties to make $(M,D)$ into a K\"{a}hler pair.
\end{example2}

\begin{example2}
\label{example2:fermat}
We consider the Fermat hypersurface,
\[ M^n_a := \left\{ \sum_{j=1}^n z_j^a = 0 \right\} \subset \CP{n-1},\]
with line bundle $\mathcal{L} = \mathcal{O}(1)$, and the ample divisors
\[ D_j := \{z_j = 0\}\]
for $j = 1, \ldots, n$. 
We choose linking numbers $l_j = a$ for all $j$ (for reasons which will become apparent in Example \ref{example2:fermatcover}). 
Thus we can construct a K\"{a}hler pair $(M^n_a,D)$ by Example \ref{example2:kahlerline}, whose K\"{a}hler form has cohomology class
\[ [\omega] = an \cdot c_1(\mathcal{O}(1)).\]
\end{example2}

\begin{definition2}
\label{definition2:gmd}
Now recall that there is a grading datum $\bm{G}(M \setminus D)$ associated to $M \setminus D$, with exact sequence
\[  \Z \To H_1(\mathcal{G}(M \setminus D)) \To H_1(M \setminus D) \To 0.\]
We will denote
\[ \bm{G}(M,D) := \bm{G}(M \setminus D)\]
in the relative case, and write the exact sequence as 
\[ \Z \To \widetilde{Y}(M,D) \To \widetilde{X}(M,D) \To 0.\]
\end{definition2}

\begin{lemma2}
\label{lemma2:totreal}
The grading datum $\bm{G}(M,D)$ is independent of the choice of K\"{a}hler form $\omega$ on $M$; it depends only on the complex structure of $M \setminus D$.
\end{lemma2}
\begin{proof}
For any almost-complex manifold $X$, let $\mathscr{R}(X)$ denote the bundle of totally real subspaces of the tangent bundle of $X$. 
For any symplectic form compatible with the almost-complex structure, the bundle of Lagrangian subspaces $\mathcal{G}(X)$ embeds inside $\mathscr{R}(X)$. This embedding is a homotopy equivalence. 
This follows from the well-known that, for a fixed vector space equipped with compatible symplectic and complex structures, the Lagrangian Grassmannian is a deformation retract of the totally real Grassmannian; applying the five-lemma to the long exact sequences of homotopy groups for the two fibre bundles shows that the embedding induces an isomorphism of homotopy groups, and hence is a homotopy equivalence by Whitehead's theorem. 

Therefore, the grading datum $\bm{G}(X) := \{H_1(\mathcal{G}_x(X)) \To H_1(\mathcal{G}(X))\}$  is isomorphic to the grading datum $\{H_1(\mathscr{R}_x(X)) \To H_1(\mathscr{R}(X))\}$, which manifestly depends only on the almost-complex structure. 
Applying these arguments to $X = M \setminus D$ yields the result. 

\end{proof}

We will now introduce a pseudo-grading datum $\bm{H}(M,D)$, together with a morphism of grading data
\[ \bm{p}:\bm{G}(\bm{H}(M,D)) \To \bm{G}(M,D).\]
Compare \cite[Section 5]{Ionel2003}.

The homology long exact sequence for the pair $M, M \setminus D$ says that
\[ H_2(M) \To H_2(M,M \setminus D) \To H_1(M \setminus D) \To H_1(M)\]
is exact. 
Note that
\begin{eqnarray*}
 H_2(M,M \setminus D) & \cong &H^{2n-2}(D) \\
 &\cong & \bigoplus_{j=1}^k H^{2n-2}(D_j) \\
&\cong & \Z \langle y_1, \ldots, y_k \rangle.
\end{eqnarray*}
The first step follows by Poincar\'{e} duality, and the second follows by inductively applying Mayer-Vietoris and using the normal crossings condition. 
Explicitly, $y_j$ represents the class of a disk intersecting $D_j$ once positively, with boundary a meridian loop around $D_j$. 

\begin{definition2}
\label{definition2:2gradm}
We define the pseudo-grading datum $\bm{H}(M,D)$ by taking the first two terms of the exact sequence above:
\[\begin{diagram}
 Z(M,D) & \rTo^{f} & Y(M,D) \\
\dEq & & \dEq \\
H_2(M;\Z) &  \rTo^{f} & \Z \langle y_1, \ldots, y_k \rangle,
\end{diagram}\]
where 
\begin{eqnarray*}
f(u) &=& \sum_{i=1}^k (u \cdot D_i) y_i.
\end{eqnarray*}
We define the element $c \in \mathrm{Hom}(Z,\Z)$ to be given by $2 c_1(TM) \in H^2(M)$.
\end{definition2}

\begin{definition2}
\label{definition2:mapq}
We define a map
\[q:Y(M,D) \To H_1(\mathcal{G}(M \setminus D))\]
as follows:
Recall that $y_i$ corresponds to a disk
\[ y_i : (D^2,\partial D^2) \To (M,M \setminus D).\]
We choose a lift
\[\begin{diagram}
& & \mathcal{G}M \\
& \ruTo^{\bar{y}_i} & \dTo \\
D^2 & \rTo_{y_i} & M,
\end{diagram}
\]
then define
\[ q(y_i) := [\partial \bar{y}_i] \in H_1(\mathcal{G}(M \setminus D)).\]
We observe that any two such lifts $\bar{y}_i$ are homotopic, because $\mathcal{G}M \To M$ is a fibration with connected fibres and $D^2$ is contractible, so the homology class $q(y_i)$ is independent of the choice of lift.
Note that the boundary Maslov index vanishes:
\[ \mu(D^2,y_i^* TM, q(y_i)) = 0.\]
\end{definition2}

\begin{lemma2}
\label{lemma2:gradcomm}
The diagram 
\[\begin{diagram}
H_2(M) & \rTo & \Z \langle y_1, \ldots, y_k \rangle \\
\dTo<{2c_1} && \dTo<{q} \\
\Z & \rTo & H_1(\mathcal{G}(M \setminus D))  
\end{diagram}
\]
commutes.
\end{lemma2}
\begin{proof}
Suppose that $u: \Sigma \To M$ is a surface representing a homology class in $H_2(M)$, and intersecting the divisors $D_j$ transversely. 
One side of the square maps 
\[ u \mapsto \sum_{i=1}^k (u \cdot D_i) \partial \bar{y}_i,\]
while the other maps 
\[ u \mapsto f(2c_1(u)).\]
We consider the bundle pair $(\Sigma,E,\emptyset)$ with empty boundary, simply given by the complex vector bundle $E := u^*(TM)$. 
Its Maslov index is $\mu(\Sigma,E,\emptyset) = 2c_1(u)$. 
We now define a decomposition of this bundle pair: $\Sigma = \Sigma_1 \cup \Sigma_2$, where $\Sigma_1$ is a union of small balls around each of the intersection points of $u$ with divisors $D_i$, and $\Sigma_2$ is the rest of $\Sigma$. 
We define the Lagrangian boundary conditions $F$ along $\partial \Sigma_1$ by requiring that the corresponding lift of the boundary $\partial \Sigma_1 \To \mathcal{G}M$ extends to a lift $\Sigma_1 \To \mathcal{G}M$. 
Then $\mu(\Sigma_1,E,F) = 0$, and the composition property for bundle pairs (see \cite[Appendix C.3]{mcduffsalamon}) says that
\[2c_1(u) = \mu(\Sigma,E,\emptyset) = \mu(\Sigma_1,E,F) + \mu(\Sigma_2,E,F) = \mu(\Sigma_2,E,F).\]
We note that the boundary conditions we have associated to a small ball around a positive intersection point of $u$ with the divisor $D_i$ define a map
\[ S^1 \To \mathcal{G}(M \setminus D)\]
representing the class $q(y_i)$, by definition. 
Because $u$ maps $\Sigma_2$ into $M \setminus D$, it now follows from Lemma \ref{lemma2:maslovgrad} that
\[ f(2c_1(u)) = f(\mu(\Sigma_2,E,F)) = \sum_{i=1}^k (u \cdot D_i) q(y_i),\]
so the diagram commutes.
 \end{proof}

\begin{definition2}
\label{definition2:pmorph}
We define a morphism of grading data, 
\[ \bm{p}: \bm{G}(\bm{H}(M,D)) \To \bm{G}(M, D),\]
by
 \[\begin{diagram}
\Z & \rTo & \left( \Z \oplus Y(M,D)\right)/Z(M,D) \\
\dEq && \dTo>{j \oplus y \mapsto f(j) + q(y)} \\
\Z & \rTo & H_1(\mathcal{G}(M \setminus D))
\end{diagram}
\]
(recalling Definition \ref{definition2:gradfromps}). 
It follows from Lemma \ref{lemma2:gradcomm} that this is well-defined.
\end{definition2}

\begin{lemma2}
\label{lemma2:psgrel}
If $H_1(M) = 0$ and if the map
\[ 2c_1(T(M \setminus D)): \pi_2(M \setminus D) \To \Z\]
vanishes, then the morphism of grading data
\[ \bm{p}: \bm{G}(\bm{H}(M,D)) \To \bm{G}(M,D)\]
is an isomorphism.
\end{lemma2}
\begin{proof}
The morphism $\bm{p}$ extends to a commutative diagram
\[\begin{diagram}
0& \rTo & \Z & \rTo & \left(\Z  \oplus Y(M,D)\right)/Z(M,D) & \rTo & H_1(M \setminus D) & \rTo & 0 \\
&&\dEq && \dTo && \dEq &&  \\
0 & \rTo &\Z & \rTo & H_1(\mathcal{G}(M \setminus D)) & \rTo & H_1(M \setminus D) & \rTo & 0.
\end{diagram}
\]
The bottom row is exact by the exact sequence for the homotopy groups of the fibration $\mathcal{G}(M \setminus D)$ and the condition on $\pi_2$ (see Definition \ref{definition2:2hm} and Remark \ref{remark2:masexact}). 
The homology long exact sequence for the pair $M,M \setminus D$ tells us that
\[ H_2(M) \To H_2(M,M \setminus D) \To H_1(M \setminus D) \To 0\]
is exact, because $H_1(M) = 0$. 
It quickly follows that the top row is exact. 
It follows that the middle map is an isomorphism, by the five-lemma.
 \end{proof}

Now suppose that we equip $M$ with a meromorphic $n$-form $\eta$ (i.e., an $(n,0)$-form), whose zeroes and poles lie along the divisors $D_j$.
Then we obtain a quadratic complex volume form $\eta^2$ on $M \setminus D$, and recall that this defines a morphism $\bm{p}^{\eta}: \bm{G}(M, D) \To \bm{G}_{\Z}$, allowing us to equip our category with a $\Z$-grading.

\begin{lemma2}
\label{lemma2:pseudograd}
Consider the morphism of pseudo-grading data
\[ \bm{q}: \bm{H}(M,D) \To \bm{0},\]
defined by setting
\begin{eqnarray*}
 d: \Z \langle y_1, \ldots, y_k \rangle & \To & \Z, \\
d(y_j) & = & 2p_j,
\end{eqnarray*}
where $p_j$ is the order of the pole of $\eta$ along divisor $D_j$. 
Then the diagram
\[ \begin{diagram}
\bm{G}(\bm{H}(M,D)) & \rTo^{\bm{p}}& \bm{G}(M,D) \\
& \rdTo_{\bm{G}(\bm{q})} & \dTo>{\bm{p}^{\eta}} \\
&& \bm{G}_{\Z}
\end{diagram}\]
commutes.
\end{lemma2}
\begin{proof}
Follows essentially from the definition of the boundary Maslov index, see \cite[Theorem C.3.5, `Normalization' property]{mcduffsalamon}.
 \end{proof}

\begin{remark2}
\label{remark2:cygrad}
If $M$ is Calabi-Yau, then it admits a nowhere-vanishing holomorphic volume form $\eta$, so there is a canonical morphism of grading data
\[ \bm{G}(M, D) \To \bm{G}_{\Z},\]
which is induced by the zero morphism of pseudo-grading data, in accordance with Lemma \ref{lemma2:pseudograd}.
\end{remark2}

\begin{lemma2}
\label{lemma2:fermatgrad}
The grading datum associated to the Fermat hypersurfaces with coordinate divisors, $(M^n_a,D)$ (see Example \ref{example2:fermat}) is
\[ \bm{G}(M^n_a,D) \cong \bm{G}(\bm{H}^n_a) \cong \bm{G}^n_a,\]
where $\bm{H}^n_a$ is the pseudo-grading datum of Example \ref{example2:ses1}. 
\end{lemma2}
\begin{proof}
First we observe that 
\[\bm{H}(M^n_a,D) \cong \bm{H}^n_a.\]
This follows from the fact that $H_2(M^n_a) \cong \Z$, generated by the class of a line $[P]$, that
\[ [P] \cdot D_j = 1\]
for all $j$, and that
\[ c_1([P]) = n-a.\]
Next we observe that $\bm{G}(M^n_a,D) \cong \bm{G}(\bm{H}(M^n_a,D))$. 
This follows from Lemma \ref{lemma2:psgrel}: we have $H_1(M^n_a) \cong 0$, and the map
\[ 2c_1(TM^n_a): \pi_2(M^n_a \setminus D) \To \Z\]
vanishes because $c_1(TM^n_a)$ is Poincar\'{e} dual to a multiple of the hyperplane class $[D_j]$, and hence vanishes on $H_2(M^n_a \setminus D)$.
 \end{proof}

We would now like to consider branched covers of K\"{a}hler pairs. 
We run into a problem: the pullback of a K\"{a}hler form by a branched cover is degenerate along the branch locus, and therefore is not a K\"{a}hler form.
We will see in Lemma \ref{lemma2:degen} that it is possible to remedy this situation by adding a small exact two-form supported in a neighbourhood of the branch locus. 
To keep track of these regions where the K\"{a}hler form must be perturbed, we make the following:

\begin{definition2}
\label{definition2:kahlerplus}
A {\bf K\"{a}hler$^+$ pair} $(M,D^+)$ is a K\"{a}hler pair $(M,D)$, together with a choice of open neighbourhood $D^+$ of $D$.
\end{definition2}

\begin{definition2}
\label{definition2:2brcov}
Suppose that $(N,E^+)$ and $(M,D^+)$ are K\"{a}hler$^+$ pairs of the same dimension (each with $k$ divisors), and $\bm{a} = (a_1, \ldots, a_k)$ is a tuple of positive integers.
An {\bf $\bm{a}$-branched cover} of K\"{a}hler$^+$ pairs,
\[ \phi: (N,E^+) \To (M,D^+), \]
is a holomorphic branched cover $\phi: N \To M$ which near any point $p \in E_K$ (recall $E_K$ is the intersection of all divisors $E_j$ such that $j \in K$, for $K \subset [k]$) has the local form
\[ \phi(z_1,\ldots,z_n) = \left(z_1^{a_{j_1}},\ldots, z_m^{a_{j_m}},z_{m+1},\ldots,z_n\right),\]
where $K = \{j_1,\ldots,j_m\}$, and for all $1 \le i \le m$, $\{z_i = 0\}$ corresponds to the divisor $E_{j_i}$ in the domain, and $D_{j_i}$ in the codomain. 
Thus, $\phi$ maps divisor $E_j$ to $D_j$, and has branching of order $a_j$ along divisor $E_j$ (and no branching anywhere else). 

We require that $\phi$ respects the cohomology classes $c$:
\[  \phi^* c_M = c_N.\]
We also require that $\phi^{-1}(D^+) = E^+$, and that $\phi$ respects the K\"{a}hler potentials (hence also symplectic forms) away from $D^+$, in the following sense: there is a function $\rho$ on $N$, supported in $E^+$, such that
\[ \phi^* h_M = h_N + \rho|_{N \setminus E}.\]
It follows that
\[ \omega_N = \phi^*\omega_M + dd^c \rho,\]
and in particular, that $[\omega_N] = \phi^*[\omega_M]$.
It also follows that
\[ a_j d^N_j = d^M_j \mbox{ and } w^N_j = a_j w^M_j.\]

Note that $\phi|_{N \setminus E}: N \setminus E \To M \setminus D$ is an unbranched cover of manifolds. 
\end{definition2}

\begin{lemma2}
\label{lemma2:degen}
(Compare \cite[Proposition 10]{Auroux2000}). Let $(M,D^+)$ be a K\"{a}hler$^+$ pair, and $N$ be a complex manifold of the same dimension, equipped with a normal-crossings divisor $E = \cup E_j$ as in the definition of a K\"{a}hler pair, together with a holomorphic branched cover $\phi: (N,E) \To (M,D)$ of the local form required in the definition of an $\bm{a}$-branched cover of K\"{a}hler pairs.
Then there exists a K\"{a}hler$^+$ structure on $(N,E)$, so that $\phi: (N,E^+) \To (M,D^+)$ is a branched cover of K\"{a}hler$^+$ pairs in the sense of Definition \ref{definition2:2brcov}.
\end{lemma2}
\begin{proof}
We define 
\begin{eqnarray*}
E^+ &:=& \phi^{-1}(D^+) \\
c_N &:=& \phi^* c_M \\
d^N_j &:=& \frac{d^M_j}{a_j}\\
w^N_j &:=& a_j w^M_j \\
\tilde{h}_N &:=& \phi^* h_M\\
\tilde{\alpha}_N &:=& d^c \tilde{h}_N = \phi^* \alpha_M\\
\tilde{\omega}_N &:=& \phi^* \omega_M.
\end{eqnarray*}
It is automatic that $[E_j] = d^N_j c_N$, and the linking of $\tilde{\alpha}_N$ with $E_j$ is $w^N_j$. 

We could set $\rho = 0$ and be done, except for one thing: $\tilde{\omega}_N$ is not a K\"{a}hler form. 
It is  closed and has the right cohomology class, and it is non-negative, in the sense that
\[ \tilde{\omega}_N(v,Jv) = \omega_M(\phi_*v, J\phi_*v) \ge 0\]
(using holomorphicity of $\phi$). 
However, it is not strictly positive: in fact, 
\[ \tilde{\omega}_N(v,Jv) = 0 \iff v \in ker(\phi_*),\]
 and $ker(\phi_*)$ is non-trivial along the branching locus. 
To prove the Lemma, it suffices for us to construct a function $\rho:N \To \R$, supported in $E^+$, so that
\[ \omega_N := \tilde{\omega}_N + dd^c\rho\]
is a K\"{a}hler form, i.e., $\omega_N(v,Jv)>0$ for all non-zero $v$.

To construct $\rho$, we first show that, for any $p \in N$, there exists a function $\rho_p: N \To \R$ such that:
\begin{itemize}
\item $\rho_p$ is supported in $E^+$;
\item $d d^c \rho_p(v, Jv) \ge 0$ for all $v \in ker(\phi_*)$;
\item in a neighbourhood of $p$, $dd^c \rho_p(v,Jv) > 0$ for all non-zero $v \in ker(\phi_*)$.
\end{itemize}
To construct $\rho_p$, suppose that we have local holomorphic coordinates around $p$ in which
\[ \phi(z_1, \ldots, z_n) = \left(z_1^{a_{j_1}}, \ldots, z_m^{a_{j_m}}, z_{m+1}, \ldots, z_n\right),\]
where all $a_{j_i} \neq 0$.
So the kernel of $\phi_*$ at a point $(z_1,\ldots,z_n)$ is the subspace
\[ \{(v_1,\ldots,v_n) \in \C^n: v_jz_j = 0 \mbox{ for all }j \le m \}.\]
We define 
\[ \rho_p(z_1, \ldots, z_n) := \chi(|z_1|) \ldots \chi(|z_n|)(|z_1|^2 + \ldots + |z_m|^2), \]
where $\chi$ is a cutoff function identically equal to $1$ in a neighbourhood of $0$. 
One easily checks that $\rho_p$ has the required properties. 

By adding together a finite number of such functions, we obtain a function $\tilde{\rho}$, supported in $E^+$, such that 
\[ dd^c \tilde{\rho}(v,Jv) > 0\]
for all non-zero $v \in ker(\phi_*)$. 
We now show that, for all sufficiently small $\varepsilon > 0$, 
\[ \omega_N^{\varepsilon} := \tilde{\omega}_N + \varepsilon dd^c \tilde\rho\]
is a K\"{a}hler form. 
To check this, we must show that $\omega_N^{\varepsilon}(v,Jv) > 0$ for all $v \in S(TN)$, the unit sphere bundle of the tangent bundle of $N$. 
We denote by $Ker \subset S(TN)$ the kernel of $\phi_*$; it is a closed subset of the compact manifold $S(TN)$, so compact.  
We know that $dd^c \tilde{\rho}(v,Jv) > 0$ on $Ker$, so it is strictly positive on an open neighbourhood $Ker^+$ of $Ker$. 
Hence, since $\tilde{\omega}_N(v,Jv) \ge 0$, it follows that $\omega_N^\varepsilon(v,Jv) > 0$ on $Ker^+$, for any $\varepsilon>0$.

On the complement of $Ker^+$, which is compact, $\tilde{\omega}_N(v,Jv)$ is strictly positive, hence bounded below by a strictly positive constant, and $|dd^c \tilde{\rho}(v,Jv)|$ is bounded above; so for sufficiently small $\varepsilon$, $\omega_N^\varepsilon(v,Jv)>0$ on this set also. 
Therefore, for sufficiently small $\varepsilon>0$, the function $\rho := \varepsilon \tilde{\rho}$ has the required properties: it is supported in $E^+$, and $\omega_N := \tilde{\omega}_N + dd^c \rho$ is a K\"{a}hler form.

\end{proof}

\begin{example2}
\label{example2:fermatcover}
We consider the branched cover
\begin{eqnarray*}
\phi_a: (M^n_a,D) & \To &( M^n_1,D), \\
\phi_a \left( [z_1: \ldots : z_n] \right) & = & [z_1^a : \ldots : z_n^a],
\end{eqnarray*}
which has branching of degree $a$ about each divisor $D_j$. 
It follows from Lemma \ref{lemma2:degen} that there exists a K\"{a}hler$^+$ structure on $(M^n_a,D)$, so that $\phi_a$  is an $(a, \ldots, a)$-branched cover of K\"{a}hler$^+$ pairs. 
Note that this need not be the same K\"{a}hler$^+$ structure on $(M^n_a,D)$ as that defined in Example \ref{example2:fermat}; however, the K\"{a}hler form has the same cohomology class, and the Liouville one-form has the same linking numbers.
\end{example2}

\begin{lemma2}
\label{lemma2:branchps}
Let $\phi: (N , E^+) \To (M , D^+)$ be an $\bm{a}$-branched cover. 
The unbranched cover
\[ \phi: N \setminus E \To M \setminus D\]
induces a morphism of grading data
\[ \bm{p}: \bm{G}(N,E) \To \bm{G}(M,D),\]
as in Section \ref{subsec:2affcov} (of course this branched cover does not respect symplectic forms on $E^+$, but it does respect the complex structure, so it still induces a map of grading data by Lemma \ref{lemma2:totreal}). 
We define a morphism of pseudo-grading data
\[ \tilde{\bm{p}}: \bm{H}(N,E) \To \bm{H}(M,D),\]
where
\begin{eqnarray*}
\tilde{p}_Y(y_j) & = &  a_j y_j, \mbox{ and} \\
d(y_j) &=& 2-2a_j.
\end{eqnarray*}
Then the diagram
\[\begin{diagram}
\bm{G}(\bm{H}(N,E)) & \rTo^{\bm{G}(\tilde{\bm{p}})} & \bm{G}(\bm{H}(M,D)) \\
\dTo<{\bm{p}_N} && \dTo>{\bm{p}_M} \\
\bm{G}(N,E) & \rTo_{\bm{p}} & \bm{G}(M,D)
\end{diagram}\]
commutes.
\end{lemma2}
\begin{proof}
It suffices to prove that 
\[p_Y(\partial \bar{y}_j) = f(2-2a_j) + a_j \partial \bar{y}_j.\]
This follows easily from the definition, and the local form
\[ (z_1, z_2,\ldots, z_n) \mapsto (z_1^{a_j}, z_2, \ldots, z_n)\]
of $\phi$ near divisor $E_j$.
 \end{proof}

\begin{corollary2}
If $\phi_n: (M^n_n ,D) \To (M^n_1,D)$ is the $(n, \ldots, n)$-branched cover of Fermat hypersurfaces introduced in Example \ref{example2:fermatcover}, then the induced morphism of grading data
\[ \bm{p}: \bm{G}(M^n_n,D) \To \bm{G}(M^n_1,D)\]
coincides with the morphism
\[ \bm{p}_1: \bm{G}^n_n \To \bm{G}^n_1\]
of Lemma \ref{lemma2:squaregrad} (recalling Lemma \ref{lemma2:fermatgrad}).  
\end{corollary2}

\section{Moduli Spaces of Disks}
\label{sec:2moduli}

In this section we introduce the various moduli spaces of pseudo-holomorphic disks that we will need to define the versions of the Fukaya category that we will consider.

\subsection{Moduli spaces of holomorphic spheres and disks}
\label{subsec:2dmcomp}

\begin{definition2}
Given an unordered set $\bm{E}$, with $|\bm{E}| \ge 3$, we define $\mathcal{R}_0(\bm{E})$, the moduli space of holomorphic spheres with distinct marked points $q_e$ indexed by $e \in \bm{E}$, up to biholomorphism preserving marked points.
\end{definition2}

\begin{definition2}
Given an ordered tuple $\bm{L} = (L_0,\ldots,L_d)$, a {\bf disk with boundary labels $\bm{L}$} is a disk with $d+1$ distinct boundary marked points, $\zeta_0,\zeta_1,\ldots,\zeta_d$, with the boundary component between $\zeta_i$ and $\zeta_{i+1}$ labelled $L_i$ (understood modulo $d+1$).
\end{definition2}

We define three types of moduli spaces of disks:

\begin{definition2}
\label{definition2:2disks1}
Given a tuple $\bm{L}$, and a set $\bm{E}$, with  $|\bm{L}| + 2|\bm{E}| \ge 3$, we define $\mathcal{R}(\bm{L},\bm{E})$ to be the moduli space of holomorphic disks $S$ with boundary labels $\bm{L}$, together with distinct internal marked points $q_e$ indexed by $e \in \bm{E}$.
We consider these objects up to biholomorphism preserving all marked points.
\end{definition2}

\begin{definition2}
Given a tuple $\bm{L}$, we define $\mathcal{R}_1(\bm{L}) := \mathcal{R}(\bm{L},\{1\})$, the moduli space of holomorphic disks $S$ with boundary labels $\bm{L}$ and a single interior marked point $q$.
\end{definition2}

\begin{figure}
\centering
\includegraphics[width = 0.7\textwidth]{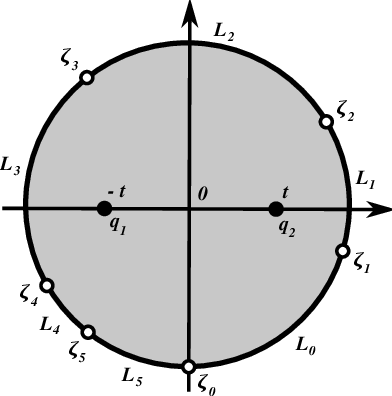}
\caption{The moduli space $\mathcal{R}_2(\bm{L})$, where $|\bm{L}| = 6$.
\label{fig:2m2}}
\end{figure}

\begin{definition2}
Given a tuple $\bm{L}$, we define $\mathcal{R}_2(\bm{L}) \subset \mathcal{R}(\bm{L},\{1,2\})$ to be the  moduli space of holomorphic disks $S$ with Lagrangian labels $\bm{L}$, together with interior marked points $q_1, q_2$, such that there is a biholomorphism of $S$ with the unit disk $\{|w| \le 1\}\subset \C$ sending 
\begin{eqnarray*}
\zeta_0 &\mapsto& -i\\
q_1 & \mapsto &-t \\
q_2 &\mapsto & t 
\end{eqnarray*}
for some $t \in (0,1) \subset \R$ (see Figure \ref{fig:2m2}).
\end{definition2}

Given a point $r$ in one of these moduli spaces, we denote by $S_r$ the corresponding (marked) Riemann surface, with all boundary marked points removed.

\subsection{Deligne-Mumford compactifications}
\label{subsec:2delmum}

We make a universal choice of strip-like and cylindrical ends for each of these moduli spaces.
We denote by $\overline{\mathcal{R}}_0(\bm{E})$, $\overline{\mathcal{R}}(\bm{L},\bm{E})$, $\overline{\mathcal{R}}_1(\bm{L})$, $\overline{\mathcal{R}}_2(\bm{L})$ the Deligne-Mumford compactifications of these moduli spaces by stable spheres and disks.
We now describe these compactifications.

The Deligne-Mumford compactification of $\mathcal{R}_0(\bm{E})$ consists of stable trees of spheres. 
Boundary strata are indexed by stable trees $T$, with semi-infinite edges indexed by $\bm{E}$. 
We denote by $V(T)$ the set of vertices of $T$, and $E(T)$ the set of edges of $T$.
A tree is called {\bf stable} if each vertex has valence $\ge 3$. 
For each vertex $v$ of $T$, we denote by $\bm{E}_v$ the set of edges of $T$ incident to $v$. 
The boundary stratum indexed by $T$ is
\[ \mathcal{R}^T_0(\bm{E}):= \prod_{v} \mathcal{R}_0(\bm{E}_v).\]
Points in this stratum correspond to trees of sphere bubbles, with semi-infinite edges corresponding to marked points, and finite edges corresponding to nodes (see \cite[Section D.3, Figure 4]{mcduffsalamon}).
The Deligne-Mumford (or Grothendieck-Knudsen) compactification, as a set, is the union of all such strata. 
It is a smooth manifold (see \cite[Section D.5]{mcduffsalamon}).
The codimension of the stratum indexed by $T$ is $2(|V(T)|-1)$.

\begin{figure}
\centering
\includegraphics[width=0.7\textwidth]{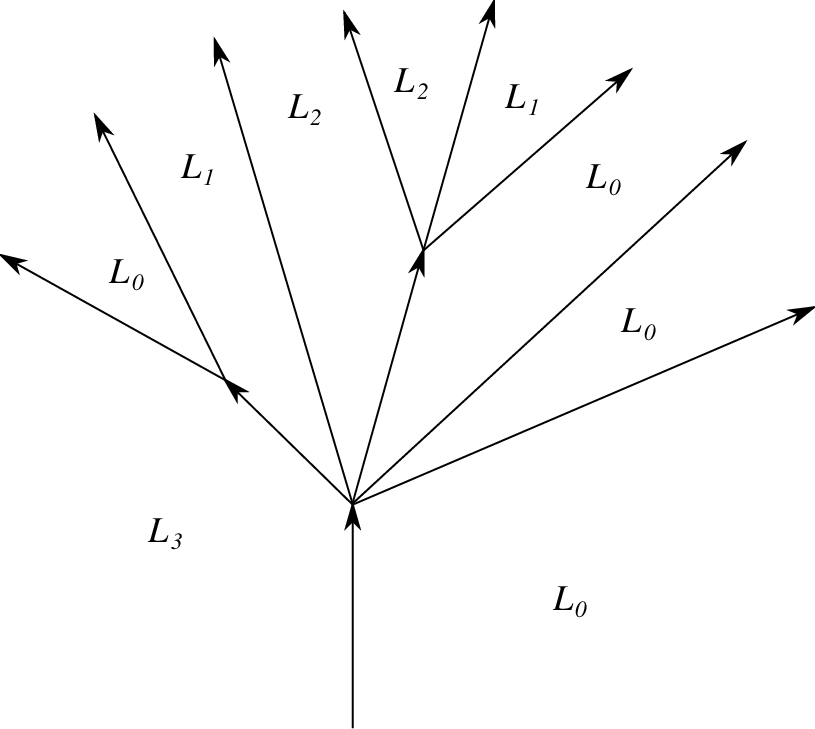}
\caption{If $\bm{L}$ is some tuple, then a $k$-leafed stable tree $T$ is said to have labels $\bm{L}$ if the connected components of $\R^2 \setminus T$ are labeled by the elements of $\bm{L}$, in order. In this figure, $\bm{L} = (L_0,L_0,L_0,L_1,L_2,L_2,L_1,L_0,L_3)$.
A labeling $\bm{L}$ of $T$ induces a labeling $\bm{L}_v$ of the regions surrounding each vertex $v$. 
In this figure, the induced labeling of the regions surrounding the uppermost vertex $v$ is $\bm{L}_v = (L_0,L_1,L_2,L_2)$.
\label{fig:2treelabels}}
\end{figure}

\begin{definition2}
\label{definition2:2plantree}
A {\bf directed $d$-leafed planar tree} is a directed $d$-leafed tree $T$ embedded in $\R^2$.
It consists of the following data:
\begin{itemize}
\item a finite set of vertices $V(T)$;
\item a set of $d$ semi-infinite {\bf outgoing edges};
\item a single semi-infinite {\bf incoming edge}, connected to a vertex $v \in V(T)$ called the {\bf root} of $T$;
\item a set $E(T)$ of internal edges.
\end{itemize}
A vertex is allowed to have zero outgoing edges, but must always have exactly one incoming edge.
We say that a vertex $v \in V(T)$ is {\bf stable} if it has $\ge 2$ outgoing edges, and {\bf semi-stable} if it has $\ge 1$ outgoing edges. 
If all vertices of $T$ are stable, we call $T$ stable; if all vertices are semi-stable, we call $T$ semi-stable.
Given a tuple $\bm{L}$, we say that $T$ {\bf has labels $\bm{L}$} if the connected components of $\R^2 \setminus T$ are labeled by the elements of $\bm{L}$, in order. 
A labeling of $T$ induces a labeling $\bm{L}_v$ of the regions surrounding each vertex $v\in V(T)$ (see Figure \ref{fig:2treelabels}).
\end{definition2}

The Deligne-Mumford compactification of $\overline{\mathcal{R}}(\bm{L},\bm{E})$ consists of stable trees of disk and sphere bubbles with appropriate markings (see \cite[Section 2.3]{Frauenfelder2008}). 
It is a smooth manifold with corners.
Each boundary stratum is indexed by a directed tree $T$, together with a directed planar subtree $T_L$ with labels $\bm{L}$. 
We denote $T_E := T\setminus T_L$.
We require that the semi-infinite edges of $T_E$ are indexed by $\bm{E}$. 
For each vertex $v \in T_L$, we have a labeling $\bm{L}_v$ as above, and denote by $\bm{E}_v$ the set of edges incident to $v$ in $T_E$. 
For each vertex $v \in T_E$, we denote by $\bm{E}_v$ the set of edges incident to $v$ in $T_E$.
We require that the tree is stable, in the sense that for each vertex $v \in T_L$,
\[ |\bm{L}_v| + 2|\bm{E}_v| \ge 3,\]
while for each vertex $v \in T_E$,
\[ |\bm{E}_v| \ge 3.\]

The tree $T$ corresponds to the stratum
\[ \mathcal{R}^T(\bm{L},\bm{E}) \cong \prod_{v \in V(T_L)} \mathcal{R}(\bm{L}_v,\bm{E}_v) \times \prod_{v \in V(T_E)} \mathcal{R}_0(\bm{E}_v).\]
Points in this stratum correspond to nodal disks, with semi-infinite edges of $T_L$ corresponding to boundary marked points,  finite edges of $T_L$ corresponding to boundary nodes, semi-infinite edges of $T_E$ corresponding to internal marked points, and finite edges of $T_E$ corresponding to internal nodes. 
The codimension of this stratum is 
\[ |V(T_L)| + 2|V(T_E)| - 1.\]

\begin{figure}
\centering
\subfloat[The codimension-$1$ part of $\mathcal{R}^{1,T}_2(\bm{L})$.]{
\includegraphics[scale=1]{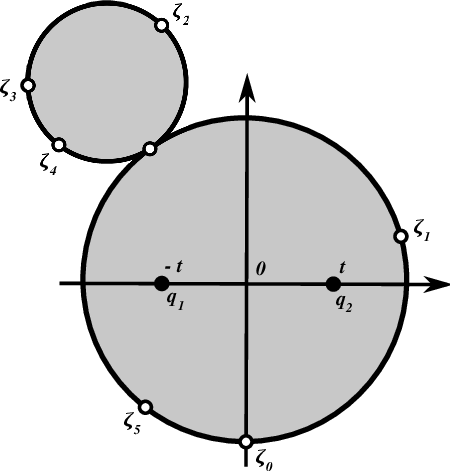}
\label{subfig:2m11}}\\
\subfloat[The codimension-$1$ part of $\mathcal{R}^{2,T}_2(\bm{L})$.]{
\includegraphics[scale=1]{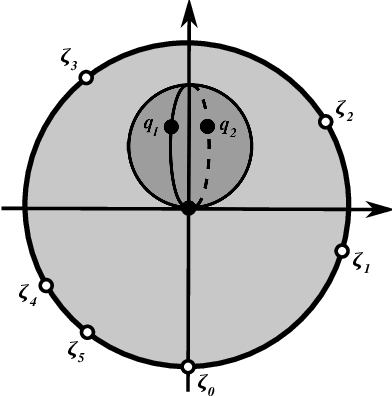}
\label{subfig:2m13}}\\
\subfloat[The codimension-$1$ part of $\mathcal{R}^{3,T}_2(\bm{L})$.]{
\includegraphics[scale=.9]{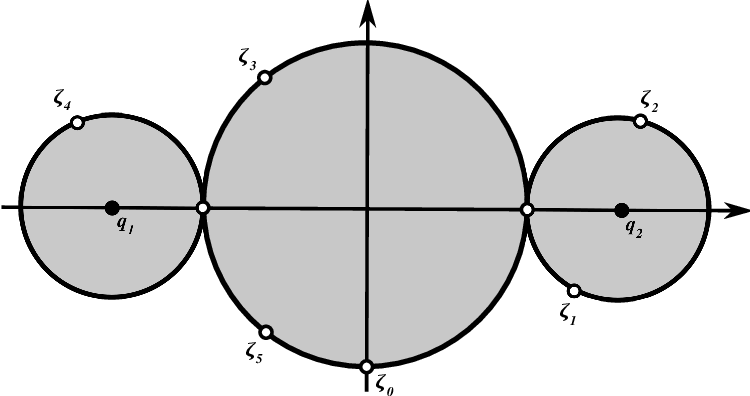}
\label{subfig:2m12}}
\caption{The codimension-$1$ boundary components of $\mathcal{R}_2(\bm{L})$, where $|\bm{L}| = 6$.
\label{fig:2m2bound}}
\end{figure}

The boundary strata of $\overline{\mathcal{R}}_2(\bm{L})$ fall into three types (we have illustrated the codimension-$1$ part of each stratum in Figure \ref{fig:2m2bound}):
\begin{itemize}
\item strata indexed by directed planar trees $T$ with boundary labels $\bm{L}$, together with a distinguished vertex $v_1$, so that all vertices other than possibly $v_1$ have valence $\ge 3$; these correspond to codimension-$(|V(T)|-1)$ strata 
\[ \mathcal{R}^{1,T}_2(\bm{L}) \cong \mathcal{R}_2(\bm{L}_{v_1}) \times \prod_{v \in V(T)\setminus \{v_1\}} \mathcal{R}(\bm{L}_v)\] 
which consist of nodal disks glued together in the obvious way (see Figure \ref{subfig:2m11});
\item another set of strata indexed by directed planar trees $T$ with boundary labels $\bm{L}$, together with a distinguished vertex $v_1$, so that all vertices other than possibly $v_1$ have valence $\ge 3$; these correspond to codimension-$|V(T)|$ strata
\[ \mathcal{R}^{2,T}_2(\bm{L}) \cong \mathcal{R}_1(\bm{L}_{v_1}) \times \prod_{v \in V(T)\setminus \{v_1\}} \mathcal{R}(\bm{L}_v)\]
which consist of nodal disks glued in the obvious way, together with a sphere with three marked points, two of which are the marked points $q_1,q_2$ and one of which is a node, identified to the internal marked point $q$ in the disk coming from the factor $\mathcal{R}_1(\bm{L}_{v_1})$ (see Figure \ref{subfig:2m13});
\item strata indexed by directed planar trees $T$ with boundary labels $\bm{L}$ and two (different) distinguished vertices $v_1, v_2$, so that all vertices other than possibly $v_1$ and $v_2$ have valence $\ge 3$, and the branch of $T$ containing $v_1$ lies strictly to the left of the branch containing $v_2$; these correspond to codimension-$(|V(T)|-2)$ strata
\[ \mathcal{R}^{3,T}_2(\bm{L}) \cong \mathcal{R}_1(\bm{L}_{v_1}) \times \mathcal{R}_1(\bm{L}_{v_2}) \times \prod_{v \in V(T)\setminus \{v_1,v_2\}} \mathcal{R}(\bm{L}_v)\] 
consisting of nodal disks glued together in the obvious way, where the internal marked point in the disk coming from the factor $\mathcal{R}_1(\bm{L}_{v_j})$ corresponds to the marked point $q_j$, for $j = 1,2$  (see Figure \ref{subfig:2m12}).
\end{itemize}

\subsection{Moduli spaces of pseudoholomorphic disks}
\label{subsec:2holcurv}

Let $(M,D)$ be a K\"{a}hler pair (see Definition \ref{definition2:2kahlp}). 
Thus, $M$ is a K\"{a}hler manifold, and $D = D_1 \cup \ldots \cup D_k$ is a union of smooth ample divisors with normal crossings. 
In this section we will define moduli spaces of pseudoholomorphic disks mapping into $M$.

\begin{definition2}
\label{definition2:2labelling}
\begin{sloppypar}
Let $\bm{F}$ be a finite set, and let
\[ \ell: \bm{F} \To [k]\]
be a function from $\bm{F}$ to the set $[k] := \{1,\ldots,k\}$ indexing the divisors $D_1,\ldots,D_k$. 
We call such a function a {\bf labelling} of $\bm{F}$.
Recalling the definition of the pseudo-grading datum $\bm{H}(M,D)$ from Definition \ref{definition2:2gradm}, we denote 
\[ \bm{d}(\ell) := \sum_{j=1}^k |\ell^{-1}(j)| y_j \in Y_{\ge 0}\]
(where $Y_{\ge 0} := \Z_{\ge 0}\langle y_1, \ldots, y_k \rangle$). 
If $\bm{a} = (a_1, \ldots, a_k)$ is a tuple of positive integers, then we define
\[ \bm{d}_{\bm{a}}(\ell) :=  \sum_{j=1}^k a_j |\ell^{-1}(j)| y_j \in Y_{\ge 0}.\] 
We denote $\bm{1} = (1, \ldots, 1)$ ($k$ copies), so that $\bm{d} = \bm{d}_{\bm{1}}$.
\end{sloppypar}
\end{definition2}

\begin{definition2}
\label{definition2:2r0ell}
Let $\bm{E}$ and $\bm{F}$ be finite sets, $|\bm{E}|+|\bm{F}| \ge 3$. 
Let $\ell: \bm{F} \To [k]$ be a labelling of $\bm{F}$. 
We define 
\[\mathcal{R}_0(\bm{E},\ell) := \mathcal{R}_0(\bm{E} \sqcup \bm{F}).\]
We define
\[ \mathcal{R}_0(\ell) := \mathcal{R}_0(\emptyset,\ell).\]
\end{definition2}

\begin{definition2}
Given a tuple of objects $\bm{L}$, finite sets $\bm{E},\bm{F}$ such that $|\bm{L}| + 2|\bm{E}| +2|\bm{F}| \ge 3$, and a labelling $\ell:\bm{F} \To [k]$, we define the moduli space
\[ \mathcal{R}(\bm{L},\bm{E},\ell) := \mathcal{R}(\bm{L},\bm{E} \sqcup \bm{F}).\]
We define
\[ \mathcal{R}(\bm{L},\ell) := \mathcal{R}(\bm{L},\emptyset,\ell).\]
\end{definition2}

\begin{definition2}
\label{definition2:hams}
We denote by $\mathcal{H}(M,D)$ the space of smooth functions $H:M \To \R$ which vanish, along with their first derivatives, along $D$. 
This is the class of Hamiltonian functions we will use to perturb our pseudoholomorphic curve equations. 
\end{definition2}

\begin{definition2}
\label{definition2:acs}
We denote by $\mathcal{J}(M,D)$ the space of smooth almost-com\-plex structures on $TM$ which preserve $TD_j$ for all $j$, and are $\omega$-tame. 
\end{definition2}

Sometimes we will omit the `$(M,D)$' from the notation and simply write $\mathcal{H}$ and $\mathcal{J}$.

For each pair of objects in the affine Fukaya category $\mathcal{F}(M \setminus D)$, we choose a Floer datum on $M$, i.e., maps
\[ H: [0,1] \To \mathcal{H}(M,D), J: [0,1] \To \mathcal{J}(M,D)\]
 (see \cite[Section 8e]{Seidel2008}). 
We make a universal choice of perturbation data on $M$ on each of the moduli spaces $\mathcal{R}_0(\ell)$, $\mathcal{R}(\bm{L},\ell)$, $\mathcal{R}_1(\bm{L})$ and $\mathcal{R}_2(\bm{L})$. 
This consists of a choice of 
\[ K \in \Omega^1(S_r,\mathcal{H}), J \in C^\infty(S_r,\mathcal{J}(M,D)).\]
for all $r$ in the moduli space, varying smoothly with $r$ (see \cite[Section 9h]{Seidel2008}).
Note that the choice of perturbation data on the moduli spaces $\mathcal{R}_0$, $\mathcal{R}$ may be different for different labellings $\ell$, even though they have the same number of boundary components and internal marked points. 

We define a universal choice of perturbation data on the moduli spaces $\mathcal{R}_0(\bm{E},\ell)$ and $\mathcal{R}(\bm{L},\bm{E},\ell)$ by pulling back the perturbation data from $\mathcal{R}_0(\ell)$ and $\mathcal{R}(\bm{L},\ell)$ respectively, via the map forgetting the marked points indexed by $\bm{E}$.

\begin{definition2}
\label{definition2:perthol}
Given a choice of perturbation datum $(K,J)$ on a Riemann surface $S$, we say that a map $u:S \To M$ is a {\bf pseudoholomorphic curve} if $(Du-Y)^{0,1} = 0$. 
This is an equation in $\mathrm{Hom}(TS,u^*TM)$, where for each $\xi \in TS$, $Y(\xi) \in \Omega^1(S,C^\infty(TM))$ is the Hamiltonian vector field associated to $K(\xi)$, and $u^*TM$ is equipped with the complex structure $J$ (see \cite[Equation (8.9)]{Seidel2008}). 
\end{definition2}

We require that
\begin{itemize}
\item the Hamiltonian part of each perturbation datum is $0$ on the moduli spaces $\mathcal{R}_0(\ell)$;
\item on the strip-like end associated to each boundary puncture, the perturbation datum agrees with the associated Floer datum;
\item the choices of perturbation data are consistent with respect to the Deligne-Mumford compactifications outlined in Section \ref{subsec:2dmcomp}, in the sense of \cite[Section 9i]{Seidel2008}, see also \cite[Section 3]{Cieliebak2007};
\item the choices of perturbation data are invariant under shifts of the anchored Lagrangian branes by the covering group action
\[\pi_1(\mathcal{G}(M \setminus D)) \cong \widetilde{Y}(M, D),\] 
as was the case for the affine Fukaya category (see Section \ref{subsec:2eqaffuk}).
\end{itemize}

Our definitions of $\mathcal{H}(M,D)$ and $\mathcal{J}(M,D)$ ensure that any pseudoholomorphic curve $u$ which is not contained in $D_j$ must intersect it positively. 
To see why, recall that solutions of the pseudoholomorphic curve equation $u: S \To M$ are in one-to-one correspondence with $J'$-holomorphic sections of the trivial fibration $S \times M \To S$, where the almost-complex structure $J'$ on the total space is defined using the perturbation data $(K,J)$, see for example \cite[Remark 17.1]{Seidel2008}. 
The conditions we have imposed on our perturbation data ensure that $J'$ makes $S \times D_j \subset S \times M$ into an almost-complex submanifold. 
Then positivity of intersection for pseudoholomorphic curves follows from the usual positivity of intersection of $J$-holomorphic curves with almost-complex submanifolds (see for example \cite[Proposition 7.1]{Cieliebak2007}).

We will use the shorthand 
\[u \cdot D := \sum_{j=1}^k (u \cdot D_j) y_j \in Y(M,D)\]
(see Definition \ref{definition2:2gradm} for the definition of $Y(M,D)$), where $u \cdot D_j$ denotes the topological intersection number, for any class $u \in H_2(M)$ or $H_2(M,M \setminus D)$. 
By positivity of intersection,
\[u \cdot D \in Y(M,D)_{\ge 0}\]
if $u$ is a pseudo-holomorphic disk or sphere that is not contained inside any of the divisors $D_j$. 

\begin{definition2}
Given a tuple $\bm{L} = (L_0^{\#},\ldots,L_s^{\#})$ of anchored Lagrangian branes, an {\bf associated set of generators} is a tuple $\bm{p} = (p_0,\ldots,p_s)$  where $p_j$ is a generator of $CF^*(L_j^{\#},L_{j-1}^{\#})$ for $j \ge 1$, and $p_0$ is a generator of $CF^*(L_0^{\#},L_s^{\#})$.
\end{definition2}

\begin{definition2}
Given a tuple of objects $\bm{L}$ with associated generators $\bm{p}$, finite sets $\bm{E},\bm{F}$ such that $|\bm{L}| + 2|\bm{E}| +2|\bm{F}| \ge 3$, and a labelling $\ell:\bm{F} \To [k]$ as in Definition \ref{definition2:2labelling}, we define an element of the moduli space $\mathcal{M}(\bm{p},\bm{E},\ell)$ to be a pair $(r,u)$, where $r$ is an element of $\mathcal{R}(\bm{L},\bm{E},\ell)$ and $u:S_r \To M$ is a smooth map, such that:
\begin{itemize}
\item $u$ is a pseudo-holomorphic curve, with Lagrangian boundary conditions given by the labels $\bm{L}$;
\item $u$ is asymptotic to the generators $\bm{p}$ along the corresponding strip-like ends;
\item $u(q_f) \in D_{\ell(f)}$ for each $f \in \bm{F}$;
\item $u \cdot D = \bm{d}(\ell)$.
\end{itemize}
(see \cite[Equation (8.10)]{Seidel2008} for the definition of `asymptotic to the generators $\bm{p}$').
Note that, because the perturbation data on $\mathcal{R}(\bm{L},\bm{E},\ell)$ are pulled back from $\mathcal{R}(\bm{L},\ell)$ via the forgetful maps, for any $\bm{E}' \subset \bm{E}$ there is a forgetful map
\[ \mathcal{M}(\bm{p},\bm{E},\ell) \To \mathcal{M}(\bm{p},\bm{E}',\ell).\]
There are evaluation maps 
\[\bm{ev}: \mathcal{M}(\bm{p},\bm{E},\ell) \To M^{\bm{E}}\]
which respect the forgetful maps.
We will sometimes write $\mathcal{M}(\bm{p},\bm{E},\bm{d})$ instead of $\mathcal{M}(\bm{p},\bm{E},\ell)$, where $\bm{d} = \bm{d}(\ell)$, and we denote
\[ \mathcal{M}(\bm{p},\bm{d}):=\mathcal{M}(\bm{p},\emptyset,\bm{d}).\]
\end{definition2}

\begin{remark2}
\label{remark2:diskmarkpts}
Our assumptions on the perturbation data ensure that $u$ intersects $D_j$ in isolated points with positive multiplicity. 
Since each marked point $q_f$ with $\ell(f) = j$ contributes at least $1$ to $u \cdot D_j$, our requirement that $u \cdot D = \bm{d}(\ell)$ ensures that $u$ intersects $D_j$ only at the marked points $q_f$ with $\ell(f)=j$, and each intersection has multiplicity $1$. 
\end{remark2}

\begin{definition2}
\label{definition2:2strips}
For a tuple of $2$ Lagrangian labels $\bm{L}$ with associated generators $\bm{p}$, we define $\mathcal{M}(\bm{p},0)$, the set of holomorphic strips with boundary conditions on $\bm{L}$, intersection number $0$ with the divisors $D$, translation-invariant perturbation coming from the corresponding Floer datum, asymptotic to the generators $\bm{p}$, modulo translation by $\R$ (see \cite[Equation (8.8)]{Seidel2008}).
\end{definition2}

Given a pseudo-holomorphic curve with an internal marked point, we define the notion of `tangency to a divisor to order $k$' at the marked point, in accordance with \cite{Cieliebak2007}: 

\begin{definition2}
Suppose we are given:
\begin{itemize}
\item a Riemann surface $S$ with an internal marked point $q \in S$;
\item a pseudo-holomorphic curve $u:S \To M$;
\item a choice of divisor $D_j \subset M$;
\item an integer $k \ge 1$.
\end{itemize}
We say that $u$ is {\bf tangent to $D_j$ at $q$ to order $k$} if 
\begin{itemize}
\item $u(q) \in D_j$;
\item all partial derivatives of $u$ at $q$ of order $\le k$ lie inside the tangent space $T D_j$.
\end{itemize}
We remark that this does not depend on the choice of coordinates.
\end{definition2}

When $k=0$, this is the same thing as a point constraint $u(q) \in D_j$. 
For $k \ge 1$, one should think of the curve $u$ having `ramification of order $k+1$' about the divisor $D_j$. 
If $u$ is not contained in $D_j$, then the point $q$ contributes at least $k+1$ to the intersection number of $u$ with $D_j$ (see \cite[Proposition 7.1]{Cieliebak2007}).

\begin{definition2}
Suppose that $\bm{a} = (a_1, \ldots, a_k)$ is a tuple of $k$ positive integers. 
We define the moduli space $\mathcal{M}(\bm{p},\bm{E},\ell, \bm{a})$ in exactly the same way as we did $\mathcal{M}(\bm{p},\bm{E},\ell)$, with the following exceptions:
\begin{itemize}
\item For each $f \in \bm{F}$, we require $u$ to be tangent to $D_{\ell(f)}$ at $q_f$ to order $a_{\ell(f)}-1$;
\item We require that $u \cdot D = \bm{d}_{\bm{a}}(\ell)$.
\end{itemize} 
In particular, we have an isomorphism
\[ \mathcal{M}(\bm{p},\bm{E},\ell) \cong \mathcal{M}(\bm{p},\bm{E},\ell,\bm{1}).\]
\end{definition2}

\begin{remark2}
\label{remark2:multinter}
Because each marked point $q_f$ contributes $\ge a_{\ell(f)}$ to $u \cdot D_{\ell(f)}$, elements $u \in \mathcal{M}(\bm{p},\bm{E},\ell, \bm{a})$ do not intersect any of the divisors $D_j$ anywhere other than at marked points $q_f$, where they intersect with multiplicity $a_{\ell(f)}$.
\end{remark2}

\begin{definition2}
Let $\bm{L}$ be a tuple of Lagrangians with associated generators $\bm{p}$, $a$ a positive integer, and $j \in [k]$. We define $\bm{F}$ to be a set with a single element, and $\ell$ a labelling which assigns $j$ to this element. 
We let $\bm{a}$ be any tuple such that $a_j = a$.
 Then we define
\[\mathcal{M}_1(\bm{p},j,a):= \mathcal{M}(\bm{p},\ell,\bm{a}),\]
 the moduli space of pseudo-holomorphic disks with a single internal marked point, which is tangent to divisor $D_j$ to order $a-1$, and no other intersections with the divisors. 
Note that we're not defining anything new; this is just convenient notation for us to have.
\end{definition2}

\begin{definition2}
Given a tuple of Lagrangians $\bm{L}$ with associated generators $\bm{p}$, a positive integer $a$, together with a choice of divisor $D_j$, we define an element of the moduli space $\mathcal{M}_2(\bm{p},j,a)$ to consist of the following data:
\begin{itemize}
\item a point $r \in \mathcal{R}_2(\bm{L})$;
\item a smooth map $u:S_r \To M$,
\end{itemize}
such that:
\begin{itemize}
\item $u$ is a pseudo-holomorphic curve;
\item $u$ is asymptotic to the generators $\bm{p}$ along the strip-like ends;
\item $u \cdot D = (a+1) y_j$;
\item $u$ is tangent to $D_j$ at $q_1$ to order $a-1$, and to $D_j$ at $q_2$ to order $0$.
\end{itemize}
\end{definition2}

\begin{remark2}
As before, elements $u \in \mathcal{M}_2(\bm{p},j,a)$ do not intersect the divisors $D_i$ for $i \neq j$, and intersect $D_j$ only at $q_1$ (with multiplicity $a$) and $q_2$ (with multiplicity $1$).
\end{remark2}

\begin{sloppypar}
Finally, we define the relevant moduli spaces of pseudo-holomorphic spheres, which will appear in the Gromov compactifications of the moduli spaces $\mathcal{M}(\bm{p},\ell)$.
\end{sloppypar}

\begin{definition2}
\label{definition2:mund0k}
Let $\bm{E}, \bm{F}$ be finite sets, $\ell: \bm{F} \To [k]$ a labelling, and $\beta \in H_2(M)$ a homology class. 
We define an element of the moduli space $\mathcal{M}_0(\bm{E},\ell,\beta)$ to be a pair $(r,u)$, where $r \in \mathcal{R}_0(\bm{E},\ell)$ and $u:S_r \To M$ is a smooth map  such that
\begin{itemize}
\item $[u] = \beta$;
\item $u$ is a pseudo-holomorphic curve;
\item $u(q_f) \in D_{\ell(f)}$ for each $f \in \bm{F}$.
\end{itemize}
Note that, because our perturbation data on $\mathcal{R}_0(\bm{E},\ell)$ are pulled back from $\mathcal{R}_0(\ell)$ by the forgetful map, for any $\bm{E}' \subset \bm{E}$ there is a forgetful map
\[ \mathcal{M}_0(\bm{E},\ell) \To \mathcal{M}_0(\bm{E}',\ell),\]
and there are evaluation maps 
\[ \bm{ev}: \mathcal{M}_0(\bm{E},\ell) \To M^{\bm{E}},\]
which respect the forgetful maps in the obvious sense. 
\end{definition2}

\begin{remark2}
Our requirement that the Hamiltonian component of the perturbation data vanishes on $\mathcal{R}_0(\ell)$ implies that any pseudo-holomorphic sphere with $[\beta] = 0$ is constant.
\end{remark2}

\begin{remark2}
To define Gromov-Witten invariants as in \cite{Cieliebak2007}, one should consider the moduli space of pseudo-holomorphic spheres $u$ which are not contained inside $D$, and such that $u \cdot D = \bm{d}(\ell)$. 
That is not our aim here: $\mathcal{M}_0(\bm{E},\ell, \beta)$ is the moduli space of spheres that may appear in Gromov compactifications, and includes all constant spheres and spheres contained inside $D$.
\end{remark2}

\subsection{Indices and orientations}
\label{subsec:2indor}

Each of the moduli spaces of pseudo-holomorphic curves we have defined can be defined as the set of zeroes of a smooth section of a Banach bundle over a Banach manifold of maps (more precisely, the moduli space can be covered by such).
We follow \cite[Chapters 8 and 9]{Seidel2008} in defining the functional analytic framework, with modifications following \cite[Section 6]{Cieliebak2007} to take into account the `orders of tangency' restrictions. 
This just means we have to use $W^{k+1,p}$ maps rather than $W^{1,p}$, so that we can make sense of `derivatives of order $\le k$'.

The linearization of this smooth section defines a Fredholm operator. 
The moduli spaces are said to be regular when the linearization is surjective everywhere. 
When they are regular, the moduli spaces are smooth manifolds, with dimension given by the index of the Fredholm operator.
In this section, we will outline the calculation of this dimension.

\begin{lemma2}
\label{lemma2:gradm}
Let $\bm{L}$ be a tuple of anchored Lagrangian branes with associated generators $\bm{p}$, and $u$ an element of $\mathcal{M}(\bm{p},\bm{E},\ell)$. 
Let $\tilde{y}_j \in \widetilde{Y}(M, D)$ be the degree of $p_j$.
Then the index of the extended linearized operator $D_u$ at $u$ satisfies
\[f( i(D_u) - 2|\bm{E}| + 2-s)=  -\tilde{y}_0 + \sum_{j=1}^s \tilde{y}_j + q(\bm{d}(\ell)),\]
where $q$ is the map from Definition \ref{definition2:mapq}. 
Furthermore, there is a canonical identification of orientation lines
\[ o_{p_0} \cong \mathrm{det}(D_u) \otimes o_{p_1} \otimes \ldots \otimes o_{p_s}.\]
\end{lemma2}
\begin{proof}
Each path $p_j$ lifts to a path from $L_{j-1}^{\#}$ to $\tilde{y}_j \cdot L_j^{\#}$ in $\widetilde{\mathcal{G}}(M \setminus D)$, which we use to define orientation operator $D_{p_j}$ whose index is $0$.

We denote the linearized operator of the pseudo-holomorphic curve equation at $u$ (with fixed domain $S$) by $D_{S,u}$. 
We glue the orientation operators $D_{p_j}$ to $D_{S,u}$ along the strip-like ends to obtain a bundle pair over the closed disk, and hence a Cauchy-Riemann operator $\overline{D}$.
The gluing formula implies that
\begin{eqnarray*}
i(\overline{D}) &=& i(D_{S,u}) + i(D_{p_1}) + \ldots + i(D_{p_s}) +(n- i(D_{p_0})) \\
&=& i(D_{S,u}) + n,
\end{eqnarray*}
and there is a canonical isomorphism
\[ \mathrm{det}(\overline{D}) \cong \mathrm{det}(D_{S,u}) \otimes o_{p_1} \otimes \ldots \otimes o_{p_s} \otimes o_{p_0}^{\vee}.\]

As in Section \ref{subsec:2eqaffuk}, this bundle pair defining $\overline{D}$ is equivalent to a bundle pair $(D^2, u^*TM, \rho)$, where $\rho: \partial D^2 \To \mathcal{G}(M \setminus D)$ which lifts the boundary map $\partial u$, and hence has index 
\[ i(\overline{D}) = n + \mu(D^2, u^*TM, \rho).\]
As in Section \ref{subsec:2eqaffuk}, we can compute the homology class
\[ [\rho] = -\tilde{y}_0 + \sum_{j=1}^s \tilde{y}_j.\]
As in the proof of Lemma \ref{lemma2:gradcomm}, we define a decomposition of the bundle pair $(D^2,u^*TM,\rho)$ into two bundle pairs: $(\Sigma_1,u^*TM,\rho')$, consisting of a union of small balls surrounding each of the points $q_f$ for $f \in \bm{F}$, with boundary conditions given by $\partial \bar{y}_{\ell(f)}$, and $(\Sigma_2,u^*TM,\rho \cup \rho')$ the complement of $\Sigma_1$. 
Then the decomposition property of the boundary Maslov index, together with Lemma \ref{lemma2:maslovgrad}, say that
\begin{eqnarray*}
 \mu(D^2,u^*TM,\rho) &=& \mu(\Sigma_1,u^*TM,\rho) + \mu(\Sigma_2,u^*TM,\rho \cup \rho')\\
 &=& \mu(\Sigma_2,u^*TM,\rho \cup \rho') \\
\Rightarrow f(\mu(D^2,u^*TM,\rho)) &=& [\rho] + [\rho'] \\
&=& -\tilde{y}_0 + \sum_{j=1}^s \tilde{y}_j + q(\bm{d}(\ell)).
\end{eqnarray*}
The result now follows, as
\begin{eqnarray*}
 i(D_u) &=& i(D_{S,u}) + \mathrm{dim}\left( \mathcal{R}(\bm{L},\bm{E}) \right)\\
&=& i(\overline{D}) -n  + s-2 + 2|\bm{E}| \\
&=& \mu(D^2,u^*TM,\rho) + s-2+2|\bm{E}|.
\end{eqnarray*}
The isomorphism of orientation lines follows after we fix an orientation for $\mathcal{R}(\bm{L},\bm{E})$.
 \end{proof}

\begin{lemma2}
\label{lemma2:gradmorb}
Suppose we have objects $\bm{L}$ with associated generators $\bm{p}$, a $k$-tuple $\bm{a} = (a_1, \ldots, a_k)$ of positive integers, and an element $u \in \mathcal{M}(\bm{p},\bm{E},\ell,\bm{a})$. 
Let $\tilde{y}_j \in \widetilde{Y}(M, D)$ be the degree of $p_j$.
Then the index of the extended linearized operator $D_u$ at $u$ satisfies
\[f \left( i(D_u)  +2-s-2|\bm{E}| - 2\left|\bm{d}(\ell) - \bm{d}_{\bm{a}}(\ell)\right|\right) =  -\tilde{y}_0 + \sum_{j=1}^s \tilde{y}_j + q(\bm{d}_{\bm{a}}(\ell)).\]
\end{lemma2}
\begin{proof}
Follows from Lemma \ref{lemma2:gradm}, where we observe that being tangent to $D_{\ell(f)}$ at $q_f$ to order $a_{\ell(f)}$ imposes an additional $2(a_{\ell(f)} - 1)$-dimensional constraint on the disk, leading to the final term on the left-hand side, which is equal to
\[ 2\sum_{f \in \bm{F}} 1- a_{\ell(f)}.\] 
 \end{proof}

We can perform similar calculations for $\mathcal{M}_0$ and $\mathcal{M}_2$.
We obtain:

\begin{lemma2}
\label{lemma2:indor0}
If $D_u$ is the extended linearized operator at $u \in \mathcal{M}_0(\bm{E},\ell,\bm{a})$, then
\[f( i(D_u) -2n + 6 - 2|\bm{E}| - 2 \left| \bm{d}(\ell) - \bm{d}_{\bm{a}}(\ell)\right|) =  q(\bm{d}_{\bm{a}}(\ell))) .\]
\end{lemma2}

\begin{lemma2}
\label{lemma2:indor2}
If $D_u$ is the extended linearized operator at $u \in \mathcal{M}_2(\bm{p},j,a)$, then
\[f( i(D_u) +2a+1-s) =   -\tilde{y}_0 + \sum_{l=1}^s \tilde{y}_l + q((a+1) y_j),\]
and there is a canonical isomorphism of orientation lines as before. 
To clarify, recall that $\tilde{y}_l \in \widetilde{Y}(M,D)$ is the degree of $p_l$, but $y_j$ is the $j$th generator of $Y(M,D)$.
\end{lemma2}

\begin{remark2}
\label{remark2:regularity}
One can prove that the moduli spaces $\mathcal{M}(\bm{p},\bm{E},\ell)$ are regular for generic choices of perturbation data, by essentially the same argument as in \cite[Section 3.2]{mcduffsalamon}, with modifications as in \cite[Section 9k]{Seidel2008}. 
Namely, for each map $u$ in the Banach manifold of maps, one can choose the perturbation datum essentially arbitrarily on an open subset of the domain, and this is enough to achieve transversality. 
There is one exception: the moduli space of holomorphic strips that do not intersect the boundary divisors is defined using a translation-invariant perturbation datum (see Definition \ref{definition2:2strips}). 
It is shown in \cite{Floer1995,oh97} that these moduli spaces are regular for generic choice of Floer data.
\end{remark2}

\begin{remark2}
\label{remark2:indiv}
Regularity is not as straightforward for the moduli spaces of stable pseudo-holomorphic spheres, because our restriction that the almost-complex structure preserves $TD_j$ is non-generic, and allows for the possibility of non-regular sphere bubbles contained inside $D_j$.
However sphere bubbles inside $D$ may appear in the Gromov compactification of our moduli spaces of pseudo-holomorphic disks, so we will need to address their transversality properties in Section \ref{subsec:2gromov}.
\end{remark2}

\subsection{Gromov compactness}
\label{subsec:2gromov}

We now describe Gromov compactifications of the moduli spaces we have defined. 
We observe that each moduli space has bounded energy by Lemma \ref{lemma2:stokes}, and standard Gromov compactness shows that any sequence in one of these moduli spaces has a subsequence which converges, up to sphere and disk bubbling (see \cite{mcduffsalamon} for spheres and \cite{Frauenfelder2008} for disks). 

\begin{sloppypar}
Now we describe the Gromov compactification of the moduli spaces $\mathcal{M}(\bm{p},\bm{E},\ell)$.
\end{sloppypar}

\begin{definition2}
\label{definition2:2mcomptree}
Let $\bm{L}$ be a tuple of elements of the affine Fukaya category, $\bm{p}$ an associated set of generators, $\bm{E},\bm{F}$ finite sets, and $\ell: \bm{F} \To [k]$ a labelling, such that $|\bm{L}| + 2|\bm{E}| + 2|\bm{F}| \ge 2$.
A stratum of $\overline{\mathcal{M}} (\bm{p},\bm{E},\ell)$ is indexed by an object $\bm{T}$, where $\bm{T}$ consists of the following data:
\begin{itemize}
\item A tree $T$, together with a directed planar subtree $T_L$ with labels $\bm{L}$;
\item An indexing of the semi-infinite edges of $T_E := T\setminus T_L$ by $\bm{E} \sqcup \bm{F}$;
\item For each vertex $v \in T_E$, a homology class $\beta_v \in H_2(M)$;
\item For each edge $e$ of $T_L$, a choice of generator $p_e \in CF^*(L_{r(e)},L_{l(e)})$, where $L_{r(e)}, L_{l(e)}$ are the Lagrangian labels to the right and left of $e$ respectively, such that the generators are given by $\bm{p}$ for the external edges.
\end{itemize}
We require that
\[ \sum_{v \in T_E} \beta_v \cdot D + \sum_{v \in T_L} \bm{d}(\ell_v) = \bm{d}(\ell).\]
For each vertex $v \in V(T)$, we denote by $\bm{F}_v$ the set of semi-infinite edges in $T_E$ that are incident to $v$ and have index in $\bm{F}$, and  by $\ell_v: \bm{F}_v \To [k]$ the labelling induced by $\ell$. 
We denote by $\bm{E}_v$ the remaining edges (finite or semi-infinite) in $T_E$ that are incident to $v$. 
For each vertex $v \in T_L$, we denote by $\bm{L}_v$ the tuple of Lagrangians labelling the regions surrounding $v$, and by $\bm{p}_v$ the set of chosen generators for the edges adjacent to $v$. 
The tree $T_L$ is required to be {\bf semi-stable}, in the sense that for each vertex $v \in V(T_L)$ we have
\[ |\bm{L}_v| + 2|\bm{E}_v| + 2|\bm{F}_v| \ge 2.\]
\end{definition2}

\begin{definition2}
\label{definition2:2mcompstrat}
Given such an object $\bm{T}$, we define the corresponding stratum of the Gromov compactification. 
For stable vertices $v \in T_E$, we define
\[ \mathcal{M}^{\bm{T}}(v) := \mathcal{M}_0(\bm{E}_v,\ell_v,\beta_v),\]
where the perturbation data on $\mathcal{R}_0(\bm{E}_v,\ell_v)$ is induced by our choice of perturbation data on the Deligne-Mumford compactification of $\mathcal{R}(\bm{E},\ell)$. 
For unstable vertices $v \in T_E$, we define $\mathcal{M}^{\bm{T}}(v)$ to be the moduli space of $J$-holomorphic spheres with marked points indexed by $\bm{E}_v \sqcup \bm{F}_v$, modulo automorphisms (the constant almost-complex structure $J$ is induced by our choice of perturbation data).  
For $v \in T_L$, we define
\[ \mathcal{M}^{\bm{T}}(v) := \mathcal{M}(\bm{p}_v,\bm{E}_v,\ell_v).\]
Now, letting $\bm{E}_{\mathrm{int}}$ denote the internal (finite) edges of $T_E$, we have an obvious evaluation map
\[\bm{ev}^{\bm{T}}: \prod_{v \in V(T)} \mathcal{M}^{\bm{T}}(v) \To M^{\bm{E}_{\mathrm{int}}} \times M^{\bm{E}_{\mathrm{int}}}.\]
We define
\[ \mathcal{M}^{\bm{T}}(\bm{p},\bm{E},\bm{\ell}) := \left(\bm{ev}^{\bm{T}}\right)^{-1}\left(\triangle^{\bm{T}}\right),\]
where
\[ \triangle^{\bm{T}} \subset  M^{\bm{E}_{\mathrm{int}}} \times M^{\bm{E}_{\mathrm{int}}}\]
denotes the diagonal. 
\end{definition2}

\begin{definition2}
We define the {\bf Gromov compactification}  $\overline{\mathcal{M}}(\bm{p},\bm{E},\ell)$ to be the closure of the image of the obvious map
\[\mathcal{M}(\bm{p},\bm{E},\ell) \hookrightarrow \coprod_{\bm{T}} \mathcal{M}^{\bm{T}}(\bm{p},\bm{E},\ell),\]
where the right-hand side is equipped with the Gromov topology. 
\end{definition2}

\begin{lemma2}
$\overline{\mathcal{M}}(\bm{p},\bm{E},\ell)$ is compact.
\end{lemma2}
\begin{proof}
Standard Gromov compactness \cite{mcduffsalamon,Frauenfelder2008}  shows that $\mathcal{M}(\bm{p},\bm{E},\ell)$ is compact up to sphere and disk bubbling. 
$\overline{\mathcal{M}}(\bm{p},\bm{E},\ell)$ includes all possible sphere bubbling, and all possible {\bf semi-stable} disk bubbling. 
So it remains to rule out unstable disk bubbling, i.e., non-constant disk bubbles with no interior marked points. 
Any non-constant holomorphic disk $u$ with boundary on an anchored Lagrangian brane must intersect $D$, by exactness; and any holomorphic disk which is a Gromov limit of elements of $\mathcal{M}(\bm{p},\bm{E},\ell)$ must have a marked point wherever it intersects $D$ by Remark \ref{remark2:diskmarkpts}, because the intersection point persists in a neighbourhood of the disk in the Gromov topology, by positivity of intersection. 
It follows that $\overline{\mathcal{M}}(\bm{p},\bm{E},\ell)$ is compact. 

\end{proof}

\begin{definition2}
We say that $\overline{\mathcal{M}}(\bm{p},\bm{E},\ell)$ is {\bf regular} if all of its strata are regular and do not contain any unstable sphere bubbles. 
\end{definition2}

\begin{remark2}
Standard gluing theorems imply that, if $\overline{\mathcal{M}}(\bm{p},\bm{E},\ell)$ is regular, then it has the structure of a topological manifold with corners. 
\end{remark2}

For our purpose, which is to define the relative Fukaya category, we need to show that all Gromov compactifications $\overline{\mathcal{M}}(\bm{p},\ell)$ of moduli spaces of pseudo-holomorphic disks of virtual dimension $\le 1$ are regular and compact (the zero-dimensional moduli spaces are used to define the structure maps of the Fukaya category, and the one-dimensional moduli spaces are used to prove the $A_\infty$ relations hold).
For this to be true, we must place some conditions on the K\"{a}hler pair $(M,D)$. 

\begin{definition2}
\label{definition2:posmon}
$M$ is {\bf positively  monotone} if
\[ \omega= \tau c_1(TM)\]
in $H^2(M)$, for some $\tau > 0$.
\end{definition2}

\begin{example2}
If $1 \le a \le n-1$, then the symplectic manifold $M^n_a$ of Example \ref{example2:fermat} is positively monotone.
\end{example2}

\begin{proposition2}
\label{proposition2:posmon}
If $M$ is positively  monotone, then all Gromov compactifications $\overline{\mathcal{M}}(\bm{p},\ell)$ of moduli spaces of pseudo-holomorphic disks of virtual dimension $\le 1$ are regular, for generic choice of perturbation data.
\end{proposition2}
\begin{proof}
Non-constant sphere bubbling is ruled out by the standard argument \cite{Oh1993}: any non-constant sphere $u$ which bubbles off has positive energy, hence carries index $2c_1(u) \ge 2$ with it, leaving behind a disk of virtual dimension $<0$; since moduli spaces of disks are generically regular by  Remark \ref{remark2:regularity}, this does not happen for generic perturbation data. 

Constant sphere bubbles correspond to disks with some marked points coinciding. 
These configurations have virtual codimension $\ge 2$, and their moduli space is generically regular, hence they do not appear in moduli spaces of virtual dimension $\le 1$.
 
\end{proof}

\begin{definition2}
\label{definition2:cy}
$M$ is {\bf  Calabi-Yau} if $c_1(TM) = 0$ in $H^2(M)$.
\end{definition2}

\begin{example2}
If $a=n$, then the symplectic manifold $M^n_a$ of Example \ref{example2:fermat} is Calabi-Yau.
\end{example2}

\begin{proposition2}
\label{proposition2:cyrel}
If $M$ is  Calabi-Yau, then all Gromov compactifications $\overline{\mathcal{M}}(\bm{p},\ell)$ of moduli spaces of pseudo-holomorphic disks of virtual dimension $\le 1$ are regular for generic choice of perturbation data.
\end{proposition2}
\begin{proof}
Non-constant sphere bubbling is ruled out for generic choice of perturbation data by a dimension-counting argument as in \cite{Hofer1995}. 
The basic idea is that non-constant $J$-holomorphic spheres in a Calabi-Yau sweep out a subset of codimension $\ge 4$, whereas a moduli space of pseudo-holomorphic disks of dimension $\le 1$ sweeps out a subset of dimension $\le 3$, so generically the two do not meet and sphere bubbling can not occur (here, `a subset of dimension $d$' means a subset contained in the image of a map from a possibly non-compact manifold of dimension $d$).

For this argument to work, we must ensure that all moduli spaces of holomorphic spheres are sufficiently regular that they do in fact sweep out a subset of codimension $\ge 4$. 
First let us deal with stable holomorphic sphere bubbles. 
Moduli spaces of stable holomorphic spheres $u$ which are not contained in $D$ are generically regular by the argument of Remark \ref{remark2:regularity}, and hence sweep out a subset of dimension
\[ 2n+2c_1(u) - 4 = 2n-4,\]
as required.

However, we can not guarantee regularity of moduli spaces of stable holomorphic spheres $u$ which are contained in some divisor $D_j$, because of our non-generic constraint that the almost-complex structure should preserve $TD_j$. 
Nevertheless, we can perturb the almost-complex structure along $D_j$ essentially arbitrarily away from the other divisors, and extend this perturbation to the rest of $M$. 
This allows us to guarantee regularity of the moduli space of stable $J$-holomorphic spheres contained in $D_j$, but not contained in any of the other divisors $D_k$ for $k \neq j$.
This moduli space sweeps out a subset of $D_j$ of dimension
\[ 2(n-1) + 2c_1(TD_j)(u)-4.\]
Note that this is different from the virtual dimension of the subset swept out by the corresponding moduli space of spheres in $M$, which would be
\[ 2n + 2c_1(TM)(u) - 4.\]
We observe that, because $TM \cong TD_j \oplus ND_j$ along $D_j$, we have
\[ c_1(TM) = c_1(TD_j) + c_1(ND_j),\]
and $c_1(ND_j) = P.D.([D_j])|_{D_j} = d_j c|_{D_j}$. 
For any $J$-holomorphic sphere $u$, we have $\omega(u)\ge 0 \Rightarrow c(u)\ge0$. 
It follows that $0 \ge c_1(TM)(u) \ge c_1(TD_j) (u)$, and hence that the space swept out by stable holomorphic spheres in $D_j$ has codimension $\ge 4$. 
By a similar argument, the moduli space of pseudo-holomorphic disks contained in $D_K$, but not contained in $D_j$ for $j \notin K$, sweeps out a subset of codimension $\ge 4$.

Now observe that, when we remove the sphere bubbles from our nodal disk, the index does not change (because all the sphere bubbles have Chern number zero), so the disk that our sphere bubbles must meet up with still lies in a moduli space of virtual dimension $\le 1$. 
Moduli spaces of pseudo-holomorphic disks are generically regular by Remark \ref{remark2:regularity}, so this moduli space of disks sweeps out a subset of dimension $\le 3$, hence generically avoids the subset swept out by the pseudo-holomorphic spheres, because it is of codimension $\ge 4$. 

Now we need to deal with the case of unstable sphere bubbling. 
Unstable sphere bubbles are $J$-holomorphic for some fixed almost-complex structure $J$. 
Any $J$-holomorphic sphere is a multiple cover of a simple one; so it suffices to rule out bubbling of a simple $J$-holomorphic sphere. 
It is no longer true that `$J$-holomorphic spheres sweep out a submanifold of codimension $\ge 4$', because the almost-complex structures vary in high-dimensional families in our choice of perturbation data. 
However, the moduli space of pseudo-holomorphic disks, together with a simple $J$-holomorphic sphere bubble attached at a point $p$ (where $J$ is the almost-complex structure associated to $p$ by our perturbation data), and not contained in $D$, is generically regular and has negative virtual dimension, hence is empty. 
One can then deal with simple $J$-holomorphic spheres contained inside $D$ as before. 

Thus we have ruled out all non-constant sphere bubbling; constant sphere bubbling can be ruled out as in the proof of Proposition \ref{proposition2:posmon}.

\end{proof}

\begin{proposition2}
\label{proposition2:nodivintcomp}
If $\ell:\bm{F} \To [k]$ is a labelling such that $\ell^{-1}(j) = \emptyset$ for some $j$, then the Gromov compactifications $\overline{\mathcal{M}}(\bm{p},\ell)$ of moduli spaces of pseudo-holomorphic disks are regular for generic choice of perturbation data.
\end{proposition2}
\begin{proof}
As in the proofs of Propositions \ref{proposition2:posmon} and \ref{proposition2:cyrel}, it suffices to rule out non-constant sphere bubbling, because all disks are generically regular. 
Any non-constant sphere bubble has positive intersection number with all divisors $D_j$, hence the pseudo-holomorphic disks $u$ which Gromov converge to the nodal disk containing the sphere bubble must also have positive intersection number with all divisors $D_j$. 
However, by definition, $u \cdot D = \bm{d}(\ell)$, and this contradicts the hypothesis on $\ell$. 

\end{proof}

Now let $\bm{L}$ be a set of Lagrangian labels, $\bm{p}$ an associated set of generators, $a$ be a positive integer, and $D_j$ one of the divisors. 
We define the three types of boundary strata in the Gromov compactification of $\mathcal{M}_2(\bm{p},j,a)$ (compare Figure \ref{fig:2m2bound}). 
We observe that any spheres bubbling off from a sequence in $\mathcal{M}_2(\bm{p},j,a)$ are necessarily constant, because they do not intersect the divisors $D_i$ for $i \neq j$; thus we need only consider strata consisting of disk bubbles.

\begin{definition2}
Let $\bm{T}$ consist of the following data:
\begin{itemize}
\item A directed planar tree $T$ with Lagrangian labels $\bm{L}$, and a distinguished vertex $v_1$;
\item For each edge $e$ of $T$, a generator $p_e \in CF^*(L_{r(e)},L_{l(e)})$,
\end{itemize}
such that all vertices are semi-stable with the possible exception of $v_1$.
We define 
\[\mathcal{M}^{1,\bm{T}}_2(\bm{p},j,k) := \mathcal{M}_2(\bm{p}_{v_1},j,a) \times \prod_{v \neq v_1} \mathcal{M}(\bm{p}_v,\ell_v)\]
(note that for $v \neq v_1$, $\bm{F}_v = \emptyset$, so the $\ell_v$ is irrelevant but we include it in the notation for consistency).
\end{definition2}

The second stratum corresponds to $t \To 0$, so the marked points $z_1$ and $z_2$ come together and bubble off a pseudo-holomorphic sphere. 
This sphere has intersection number $0$ with all the divisors other than $D_j$, hence it must be constant.
Thus, the holomorphic disk attached to the sphere has intersection number $(a+1)$ with the divisor $D_j$, and $0$ with the other divisors, and only intersects $D_j$ at the nodal point $z$ where it is attached to the constant sphere. 
It follow that the disk is tangent to $D_j$ at $z$ to order $a+1$. 
Therefore, for appropriate choice of perturbation data, the disk is an element of $\mathcal{M}_1(\bm{p},j,a+1)$.

\begin{definition2}
Let $\bm{T}$ consist of the following data:
\begin{itemize}
\item A directed planar tree $T$ with Lagrangian labels $\bm{L}$, and a distinguished vertex $v_1$;
\item For each edge $e$ of $T$, a generator $p_e \in CF^*(L_{r(e)},L_{l(e)})$,
\end{itemize}
such that all vertices are semi-stable with the possible exception of $v_1$.
We define 
\[\mathcal{M}^{2,\bm{T}}_2(\bm{p},j,a) := \mathcal{M}_1(\bm{p}_{v_1},j,a+1) \times \prod_{v \neq v_1} \mathcal{M}(\bm{p}_v,\ell_v).\]
\end{definition2}

The third stratum corresponds to $t \To 1$, so the marked points $z_1$ and $z_2$ move to the boundary and bubble off disks at the boundary.

\begin{definition2}
Let $\bm{T}$ consist of the following data:
\begin{itemize}
\item A directed planar tree $T$ with Lagrangian labels $\bm{L}$, and two distinguished vertices $v_1$ and $v_2$;
\item For each edge $e$ of $T$, a generator $p_e \in CF^*(L_{r(e)},L_{l(e)})$,
\end{itemize} 
such that all vertices are semi-stable with the possible exception of $v_1$ and $v_2$, and the branch of $T$ containing $v_1$ lies strictly to the left of the branch containing $v_2$. 
We define 
\[ \mathcal{M}^{3,\bm{T}}_2(\bm{p},j,a) := \mathcal{M}_1(\bm{p}_{v_1},j,a) \times \mathcal{M}_1(\bm{p}_{v_2},j,1) \times \prod_{v \neq v_1,v_2} \mathcal{M}(\bm{p}_v,\ell_v).\]
\end{definition2}

\begin{definition2}
We define the moduli space $\overline{\mathcal{M}}_2 (\bm{p},j,a)$ to be
\[\left(\coprod_{\bm{T}}\mathcal{M}^{1,\bm{T}}_2(\bm{p},j,a)\right) \coprod \left(\coprod_{\bm{T}} \mathcal{M}^{2,\bm{T}}_2(\bm{p},j,a)\right) \coprod \left(\coprod_{\bm{T}} \mathcal{M}^{3,\bm{T}}_2(\bm{p},j,a)\right) 
\]
as a set. 
We equip it with the Gromov topology.
\end{definition2}

\begin{lemma2}
\label{lemma2:m2comp}
For generic choice of perturbation data, $\overline{\mathcal{M}}_2(\bm{p},j,a)$ is regular.
\end{lemma2}
\begin{proof}
Follows from the proof of Proposition \ref{proposition2:nodivintcomp}.

\end{proof}

\subsection{Branched covers}
\label{subsec:2modulibranch}

Let $\phi:(N,E^+) \To (M,D^+)$ be an $\bm{a}$-branched cover of K\"{a}hler$^+$ pairs (see Definition \ref{definition2:2brcov}).

Let $\bm{L}$ be a tuple of anchored Lagrangian branes in $N \setminus E^+$, $\bm{p}$ an associated set of generators, $\bm{E}$ a finite set, and $\ell$ a labelling.
Denote by $\phi(\bm{L})$ the image of these branes in $M \setminus D^+$, and by $\phi(\bm{p})$ the associated set of generators. 
We would like to related the moduli space $\mathcal{M}(\bm{p}, \bm{E}, \ell)$ of disks in $N$ and the moduli space $\mathcal{M}(\phi(\bm{p}),\bm{E},\ell,\bm{a})$ of disks in $M$.

Let us choose perturbation data for the moduli spaces $\mathcal{M}(\phi(\bm{p}),\bm{E},\ell,\bm{a})$ in $M$. 

\begin{condition}
\label{condition:pertham}
On $D^+ \subset M$, the Hamiltonian part of the perturbation datum vanishes. 
\end{condition}

\begin{remark2}
\label{remark2:pertham}
Given a Hamiltonian $H$ on $M$, we can define a Hamiltonian $\phi^*H$ on $N$. 
We denote the associated Hamiltonian vector fields on $M$ and $N$ by $X_H$ and $X_{\phi^* H}$ respectively. 
If $H$ satisfies Condition \ref{condition:pertham}, then $\phi$ respects the the Hamiltonian vector fields, in the sense that
\[ \phi_* X_{\phi^*H} = X_H.\]
That is because, on $\phi: N \setminus E^+ \To M \setminus D^+$, $\phi$ respects the symplectic forms and Hamiltonians by definition; whereas on $\phi: E^+ \To D^+$, the Hamiltonian vector field vanishes.
\end{remark2}

\begin{condition}
\label{condition:pertac}
On $D^+ \subset M$, the almost-complex structure part of the perturbation datum, $J$, is equal to the standard (integrable) almost-complex structure $J_0$. 
It follows that the pullback of $J$ by $\phi$ to $N$ is an almost-complex structure on $N$. 
We require that $J$ be $\omega_M$-tame, and $\phi^*J$ be $\omega_N$-tame. 
\end{condition}

\begin{remark2}
\label{remark2:pertac}
Note the pullback of an arbitrary almost-complex structure by $\phi$ may be singular along the divisors $E$; that is why we need the first part of Condition \ref{condition:pertac}.
\end{remark2}

\begin{remark2}
\label{remark2:pertac2}
Note that the condition that $J$ be $\omega_M$-tame, and $\phi^*J$ be $\omega_N$-tame, is an open condition on $J$ (and the set of such $J$ is non-empty, as it contains $J_0$).
\end{remark2}

\begin{lemma2}
\label{lemma2:pertint}
If our perturbation data in $M$ satisfy Conditions \ref{condition:pertham} and \ref{condition:pertac}, and we use the pulled-back perturbation data to define the moduli space in $N$, then there is an isomorphism of moduli spaces:
\begin{eqnarray*}
\mathcal{M}(\bm{p}, \bm{E}, \ell) & \overset{\cong}{\longrightarrow} & \mathcal{M}(\phi(\bm{p}),\bm{E},\ell,\bm{a}),\\
u & \mapsto & \phi \circ u.
\end{eqnarray*}
\end{lemma2}
\begin{proof}
The pseudo-holomorphic curve equation is defined using the Hamiltonian vector field corresponding to the Hamiltonian part of the perturbation datum, and the almost-complex structure part of the perturbation datum; $\phi$ respects the Hamiltonian vector field by Remark \ref{remark2:pertham} and respects the almost-complex structure by Condition \ref{condition:pertac}.  
It follows that, if $u$ is a solution of the pseudo-holomorphic curve equation with the pulled-back perturbation data, then $\phi \circ u$ is also a solution. 
Hence the map is well-defined (in particular, note that we do not need $\phi$ to respect the symplectic structure for the map to be well-defined).

The map is also injective: given $\phi \circ u$, the anchored brane structures on the boundary labels $\bm{L}$ tell us how to lift the boundary of the disk to $N$, and the rest follows by uniqueness of lifting from $M \setminus D$ to $N \setminus E$.
 
It is also surjective: suppose we are given $u \in \mathcal{M}(\phi(\bm{p}),\bm{E},\ell,\bm{a})$.
It is clear that, locally, $u$ lifts to a pseudo-holomorphic curve in $N \setminus E$, away from the marked points $q_f$. 
At a marked point $q_f$, $u$ intersects the divisor $D_{\ell(f)}$ with multiplicity exactly $a_{\ell(f)}-1$, by Remark \ref{remark2:multinter}, and it follows that a loop around $q_f$ gets mapped to a loop going $a_{\ell(f)}$ times around divisor $D_{\ell(f)}$. 
Therefore, a punctured neighbourhood of $q_f$ lifts to $N \setminus E$. 
By the removable singularity theorem, the point $q_f$ also lifts, so $u$ lifts locally on a neighbourhood of the marked points $q_f$. 
Therefore, since the disk is contractible, $u$ lifts to $N$, and the lift is clearly an element of $\mathcal{M}(\bm{p},\bm{E},\ell)$.
 \end{proof}

We observe that it is possible to achieve regularity of the stable moduli spaces $\mathcal{M}(\bm{p},\bm{E},\ell)$ with perturbation data satisfying Conditions \ref{condition:pertham} and \ref{condition:pertac}, because any element $u \in \mathcal{M}(\bm{p},\bm{E},\ell)$ intersects $M \setminus D^+$ (because the boundary conditions are contained in $M \setminus D^+$), where the perturbation data can be changed essentially arbitrarily.
The transversality argument follows \cite[Section 3.2]{mcduffsalamon}, with modifications as in \cite[Section 9k]{Seidel2008}. 
In the notation of the latter, we choose the neighbourhood $\Omega_r$ to avoid all strip-like ends and marked points, and to lie inside $u^{-1}(M \setminus D^+)$ (compare \cite[Lemma 5.6]{Cieliebak2007}).
 
To achieve regularity of the unstable moduli spaces, i.e., the moduli space of strips which do not intersect the divisors, we must perturb the almost-complex structure part of the Floer datum as  in \cite{Floer1995,oh97}; Condition \ref{condition:pertac} does not affect this argument.

If $N$ is positively  monotone, then the proof of Proposition \ref{proposition2:posmon} carries through, even if our perturbation data satisfy Conditions \ref{condition:pertham} and \ref{condition:pertac}. 
However, the proof of Proposition \ref{proposition2:cyrel} does not, because it is not possible to achieve regularity of the moduli spaces $\mathcal{M}_0(\bm{E},\ell)$ with perturbation data satisfying Condition \ref{condition:pertac}: in particular, on a sphere bubble contained inside one of the divisors, the perturbation datum is required to be equal to the standard integrable complex structure, which may not be regular. 
The proof of Proposition \ref{proposition2:cyrel} used regularity of these moduli spaces in a crucial way, so in general it is not possible to guarantee regularity of the Gromov compactifications $\overline{\mathcal{M}}(\bm{p},\bm{E},\ell)$
of moduli spaces of virtual dimension $\le 1$ using perturbation data satisfying Conditions \ref{condition:pertham} and \ref{condition:pertac}, if $N$ is  Calabi-Yau.

However, the proof of Proposition \ref{proposition2:nodivintcomp} did not require any such regularity, so it continues to hold when we use perturbation data satisfying Conditions \ref{condition:pertham} and \ref{condition:pertac}. 
Therefore, assuming that $D$ has $\ge 2$ irreducible components, we can still achieve regularity of the Gromov compactifications
\[ \overline{\mathcal{M}}_1(\phi(\bm{p}),j,a) \mbox{ and } \overline{\mathcal{M}}_2(\phi(\bm{p}),j,a)\]
with perturbation data satisfying Conditions \ref{condition:pertham} and \ref{condition:pertac}. 
This allows us to prove:

\begin{lemma2}
\label{lemma2:firstordlift}
Suppose that each of $(M,D)$ and $(N,E)$ is either positively  monotone or  Calabi-Yau, and $D$ and $E$ both have $\ge 2$ irreducible components.
Then there exist choices of perturbation data in $M$ and $N$ such that all Gromov compactifications of moduli spaces of pseudo-holomorphic disks of virtual dimension $\le 1$ are regular in both $M$ and $N$, and such that there are furthermore isomorphisms of moduli spaces
\begin{eqnarray*}
\mathcal{M}_1(\bm{p},j,1) &\overset{\cong}{\To} &\mathcal{M}_1(\phi(\bm{p}),j,a_j) \\
u & \mapsto & \phi \circ u
\end{eqnarray*}
for all $\bm{p},j$.
\end{lemma2}
\begin{proof}
First, we choose perturbation data on the moduli spaces $\mathcal{M}_1(\phi(\bm{p}),j,a_j)$ satisfying Conditions \ref{condition:pertham} and \ref{condition:pertac}. 
We define the perturbation data on $\mathcal{M}_1(\bm{p},j,1)$ to be the pullback of these perturbation data under $\phi$.
We can then extend these choices to consistent choices of perturbation data for {\bf all} of the moduli spaces $\mathcal{M}_0$, $\mathcal{M}$, $\mathcal{M}_1$, $\mathcal{M}_2$, separately in $M$ and $N$, such that the moduli spaces and their Gromov compactifications are regular. 
 \end{proof}

\section{The relative Fukaya category}
\label{sec:2relfuks}

In this section, we give our definition of the {\bf relative Fukaya category} of a K\"{a}hler pair $(M,D)$, which we denote $\mathcal{F}(M,D)$. 
It is a (possibly curved) $\bm{G}(M, D)$-graded deformation of the affine Fukaya category, in the sense of Definition \ref{definition2:2graddef}. 
In Section \ref{subsec:2ramcov}, we describe the behaviour of the relative Fukaya category with respect to branched covers.

\subsection{The definition}
\label{subsec:2relfuk}

Suppose that $(M,D)$ is a K\"{a}hler pair, and $\bm{a} = (a_1, \ldots, a_k)$ a tuple of $k$ positive integers (where $k$ is the number of divisors in $D$).

\begin{definition2}
\label{definition2:2relring}
We define the ring 
\[\widetilde{R}(M,D,\bm{a}) := \C [ r_1, \ldots, r_k ]. \]
We equip it with a $\bm{G}(\bm{H}(M, D))$-grading, where $r_j$ has grading 
\[(2(1-a_j),a_j y_j) \in (\Z \oplus Y)/Z\]
(see Definition \ref{definition2:2gradm}).
We then define the $\bm{G}(M,D)$-graded ring $R(M,D,\bm{a})$ to be the completion of $\bm{p}_* \widetilde{R}(M,D,\bm{a})$ with respect to the filtration by order, in the category of $\bm{G}(M,D)$-graded algebras. 
Here
\[ \bm{p}: \bm{G}(\bm{H}(M,D)) \To \bm{G}(M,D)\]
is the morphism of grading data from Definition \ref{definition2:pmorph}. 
We define $R(M,D) := R(M,D,\bm{1})$.
\end{definition2}

\begin{example2}
\label{example2:fermatring}
Suppose that $(M,D) = (M^n_1,D)$ as in Example \ref{example2:fermat}, and $\bm{a} = (a, \ldots, a)$. 
Then we have
\[ R(M^n_1,D,\bm{a}) \cong R^n_a ,\]
where $R^n_a$ is the $\bm{G}(M,D)$-graded ring introduced in Definition \ref{definition2:2ringra}.
\end{example2}

When $M$ is either positively  monotone or  Calabi-Yau, we give a definition of the {\bf relative Fukaya category} $\mathcal{F}(M,D)$, based on \cite{Seidel2002}. 
It is a (possibly curved) $\bm{G}(M, D)$-graded deformation of the affine Fukaya category over $R(M,D)$.

The objects of $\mathcal{F}(M,D)$ are the same as in the affine Fukaya category: anchored Lagrangian branes in $M \setminus D$. 
Generators of the morphism spaces are the same as for the affine Fukaya category. 
The morphism spaces are the free $R(M,D)$-modules generated by the generators. 
Note that the morphism spaces have finite rank, because we are considering compact Lagrangians.

Given a tuple of anchored Lagrangian branes $\bm{L}$ with associated generators $\bm{p}$, and an element $\bm{d} \in \Z^k_{\ge 0}$, we choose a labelling $\ell$ such that $\bm{d}(\ell) = \bm{d}$, and define the coefficient of $r^{\bm{d}}p_0$ in $\mu^s(p_s,\ldots,p_1)$ to be 
\[ \frac{ \# \left(\mathcal{M}(\bm{p},\ell)\right)}{\bm{d}!},\]
where, if $\bm{d} = (d_1, \ldots, d_k)$, we denote
\[ \bm{d}! := d_1! d_2! \ldots d_k!,\]
and $\#$ denotes a signed count of the zero-dimensional part of the moduli space, with signs defined according to the canonical isomorphism of orientation spaces given by Lemma \ref{lemma2:gradmorb}. 
We observe that this is a finite count, by Proposition \ref{proposition2:posmon} (if $M$ is positively  monotone) or Proposition \ref{proposition2:cyrel} (if $M$ is  Calabi-Yau). 
Note that, while the affine Fukaya category is not curved (any holomorphic disk with boundary on an exact Lagrangian is constant, since it has zero energy), the relative Fukaya category may be curved.

It follows from the index computation in Lemma \ref{lemma2:gradmorb} that the structure maps $\mu^s$ define a $\bm{G}(M, D)$-graded $A_{\infty}$ deformation of the affine Fukaya category over $R(M,D)$. 
This is why we can define our coefficient ring $R(M,D)$ to be a completion of the polynomial ring in the category of $\bm{G}(M,D)$-graded algebras, rather than in the category of algebras: the coefficient of $p_0$ in $\mu^s(p_s, \ldots, p_1)$ may involve an infinite sum of monomials $r^{\bm{d}}$, but all of these coefficients will have the same $\bm{G}(M,D)$-degree. 

Observe that the order-$0$ component of $\mu^s$ counts disks that completely avoid the divisors $D$, and therefore coincides with the definition of the structure maps in the affine Fukaya category.

The fact that $\mu \circ \mu = 0$ also follows from Propositions \ref{proposition2:posmon} and \ref{proposition2:cyrel}, since the signed count of boundary points of a compact $1$-dimensional manifold is $0$. 
The sign computation follows directly from that of \cite[Section 12g]{Seidel2008}, essentially because the fibres of the forgetful maps $\mathcal{R}(\bm{L},\bm{E}) \To \mathcal{R}(\bm{L}, \emptyset)$ have a complex structure, hence are canonically oriented.

Observe that we needed the factor $(\bm{d}!)^{-1}$ in the definition of the structure coefficients $\mu^s$. 
This is because, if $\bm{d} = \bm{d}_1 + \bm{d}_2$, then given a labelling $\ell: \bm{F} \To [k]$ with $\bm{d}(\ell) = \bm{d}$, there are $\bm{d}!/(\bm{d}_1! \bm{d}_2!)$ ways of choosing a partition $\bm{F} = \bm{F}_1 \sqcup \bm{F}_2$ such that the restricted labellings $\ell_1$, $\ell_2$ on $\bm{F}_1$ and $\bm{F}_2$ satisfy $\bm{d}(\ell_1) = \bm{d}_1$ and $\bm{d}(\ell_2) = \bm{d}_2$. 
So, in the boundary of the one-dimensional component of the moduli space $\overline{\mathcal{M}}(\bm{p},\ell)$, which consists of nodal disks with two components, there are $\bm{d}!/(\bm{d}_1! \bm{d}_2!)$ ways for the marked points $q_f$ to be distributed between the two components.

It is important to consider in what sense the relative Fukaya category is dependent on the choices (of Floer and perturbation data) involved in its construction.  
We recall the argument of \cite[Section 10a]{Seidel2008}: Let $I$ denote a set of possible choices of Floer and perturbation data. 
For each $i \in I$, we denote by $\mathcal{F}(M,D)^i$ the relative Fukaya category defined using those choices.
We define a new $A_{\infty}$ category, the {\bf total} category $\mathcal{F}(M,D)^{\mathrm{tot}}$, as follows:
\begin{itemize}
\item Objects are pairs $(L,i)$ where $L$ is an object of $\mathcal{F}(M,D)$ and $i \in I$;
\item For each pair of objects we choose a Floer datum, and for each set of labels of objects we choose a perturbation datum;
\item We require that, for a pair $(L_0,i), (L_1,i)$, the Floer datum is that given by $i$;
\item We require that, for a set of labels $(L_0,i),\ldots,(L_k,i)$, the perturbation data are those given by the index $i$;
\item The rest of the Floer and perturbation data we choose arbitrarily.
\end{itemize}
The rest of the construction (of morphism spaces and composition maps) follows that of the relative Fukaya category.
It follows that for each $i \in I$, there is a full embedding 
\[\mathcal{F}(M,D)^i \hookrightarrow \mathcal{F}(M,D)^{\mathrm{tot}}\]
which is given, on the level of objects, by
\[ L \mapsto (L,i).\]
When restricting to the affine Fukaya category ($r_j = 0$), it follows as in \cite[Section 10a]{Seidel2008} that these embeddings are quasi-equivalences, and hence invertible, and therefore that the affine Fukaya category does not depend on the choice of data $i \in I$, up to quasi-equivalence.

This need no longer be the case for the relative Fukaya category (these embeddings need not be quasi-equivalences nor invertible), and in general the relative Fukaya may depend on the data used to define it (compare \cite[Section 8f]{Seidel2003}). 
However, in this paper we are only interested in a certain full subcategory of the relative Fukaya category which is necessarily minimal (for grading reasons). 
Whereas quasi-equivalences are not necessarily invertible over the power series ring $R$ (see Remark \ref{remark2:invertquasi}), quasi-equivalences of minimal categories are (by Lemma \ref{lemma2:mininv}).

Let $\mathscr{L}$ be a set of objects of $\mathcal{F}(M,D)$, and $I$ be some set of possible choices of Floer and perturbation data for the full subcategory $\mathcal{C} \subset \mathcal{F}(M,D)$ with objects $\mathscr{L}$. 
Let us form the total category $\mathcal{C}^{\mathrm{tot}}$, as above. 

\begin{lemma2}
Suppose that we can choose Floer data for $\mathcal{C}^{\mathrm{tot}}$ so that it is minimal.
Then, for any $i, j \in I$, there is a $\bm{G}$-graded quasi-equivalence of minimal $A_{\infty}$ categories over $R$,
\[ \mathcal{C}^i \cong \mathcal{C}^j.\]
\end{lemma2}
\begin{proof}
First, note that we have $A_{\infty}$ embeddings of minimal $A_{\infty}$ categories
\[ \mathcal{C}^i \hookrightarrow \mathcal{C}^{\mathrm{tot}} \hookleftarrow \mathcal{C}^j,\]
as above. 
Now note that, because the objects $(L,i)$ and $(L,j)$ are quasi-isomorphic when restricted to the affine Fukaya category, and $\mathcal{C}^{\mathrm{tot}}$ is minimal, it follows by Lemma \ref{lemma2:quasiorder} that these objects are quasi-isomorphic in $\mathcal{C}^{\mathrm{tot}}$. 

It follows that the above embeddings are quasi-equivalences. 
By Lemma \ref{lemma2:mininv}, quasi-equivalences of minimal $A_{\infty}$ categories can be inverted over $R$. 
It follows that there is a quasi-equivalence $\mathcal{C}^i \cong \mathcal{C}^j$, as required.
 \end{proof}

In other words, if we can choose the category $\mathcal{C}^{\mathrm{tot}}$ to be minimal (e.g., for grading reasons), then $\mathcal{C}$ is independent of the choice of perturbation data $i \in I$ made in its construction, up to quasi-equivalence.

\begin{remark2}
\label{remark2:cygradrel}
If $M$ is Calabi-Yau, then by Remark \ref{remark2:cygrad}, there is a canonical morphism of grading data
\[ \bm{p}: \bm{G}(M, D) \To \bm{G}_{\Z}.\]
Thus, we obtain a canonical $\Z$-grading on the category
\[ \mathcal{F}(M,D) \cong \bm{p}_* \mathcal{F}(M,D),\]
such that the $\Z$-grading of $R$ is zero.
\end{remark2}

More generally, we introduce the {\bf smooth orbifold relative Fukaya category} $\mathcal{F}(M,D,\bm{a})$. 
Morally, it is a (possibly curved) $\bm{G}(M,D)$-graded deformation of $\mathcal{F}(M \setminus D)$ over $R(M,D,\bm{a})$, defined by analogy with $\mathcal{F}(M,D)$, but with structure coefficients counting elements of the moduli spaces $\mathcal{M}(\bm{p},\ell,\bm{a})$ rather than $\mathcal{M}(\bm{p},\ell)$. 
However, we have not proven analogues of Propositions \ref{proposition2:posmon} and \ref{proposition2:cyrel} to ensure that the Gromov compactifications $\overline{\mathcal{M}}(\bm{p},\ell,\bm{a})$ are regular for moduli spaces of virtual dimension $\le 1$, and hence that the $A_\infty$ equations hold; we have only proven Proposition \ref{proposition2:nodivintcomp}, which ensures the result for moduli spaces of disks which do not intersect one of the divisors. 
Assuming that $D$ has $\ge 2$ irreducible components (as is the case for all K\"{a}hler pairs considered in this paper), this allows us to define $\mathcal{F}(M,D,\bm{a})/\mathfrak{m}^2$, which is a $\bm{G}(M,D)$-graded first order deformation of $\mathcal{F}(M \setminus D)$ over $R(M,D,\bm{a})/\mathfrak{m}^2$. 
Recall (Section \ref{subsec:2firstod}) that we denote by $\mathfrak{m} \subset R$ the maximal ideal, and if $\mathcal{F}$ is an $R$-linear $A_{\infty}$ category, then $\mathcal{F}/\mathfrak{m}^2$ is an $R/\mathfrak{m}^2$-linear category, which retains only the information about the first-order part of $\mathcal{F}$.

\begin{remark2}
The smooth orbifold relative Fukaya category is independent of the choices made in its construction, up to first-order quasi-equivalence. 
This follows from the `total category' construction above, together with Lemma \ref{lemma2:similar}.
\end{remark2}

\begin{remark2}
$\mathcal{F}(M,D,\bm{a})$ can certainly be defined in greater generality, but $\mathcal{F}(M,D,\bm{a})/\mathfrak{m}^2$ is all we will need for the purposes of this paper.
\end{remark2}

\subsection{Behaviour with respect to ramified covers}
\label{subsec:2ramcov}

Now we would like to study the relationship between the relative Fukaya category of a K\"{a}hler pair and that of its branched cover. 
As we saw in Section \ref{subsec:2relcase}, we can really only make sense of a branched cover of K\"{a}hler$^+$ pairs. 

\begin{definition2}
If $(M,D^+)$ is a K\"{a}hler$^+$ pair, we define the affine Fukaya category $\mathcal{F}(M \setminus D^+)$, the relative Fukaya category $\mathcal{F}(M,D^+)$, and the orbifold Fukaya category $\mathcal{F}(M,D^+,\bm{a})/\mathfrak{m}^2$, as the full subcategories of $\mathcal{F}(M \setminus D)$, $\mathcal{F}(M,D)$ and $\mathcal{F}(M,D,\bm{a})/\mathfrak{m}^2$ respectively, whose objects are those Lagrangians contained in the interior of $M \setminus D^+$. 
\end{definition2}

\begin{remark2}
The affine Fukaya category $\mathcal{F}(M \setminus D^+)$ is {\bf not} an invariant of the symplectic manifold with boundary $M \setminus D^+$, unless the boundary of $M \setminus D^+$ is convex (for example, if $D^+$ is a super-level set of the K\"{a}hler potential, $\{h > c\}$, where $c$ is a regular value of $h$). 
\end{remark2}

Now let $(M,D^+)$ and $(N,E^+)$ be K\"{a}hler$^+$ pairs, and $\phi: (N,E^+) \To (M,D^+)$ an $\bm{a}$-branched cover of K\"{a}hler pairs, where $\bm{a} = (a_1, \ldots, a_k)$. 
Suppose furthermore that the unbranched cover
\[ \phi: N \setminus E \To M \setminus D\]
satisfies Condition \ref{condition:cov}.
In this section, we will examine the relationship between the following three categories: 
\[ \mathcal{F}(M,D^+) \mbox{, }\mathcal{F}(M,D^+,\bm{a})/\mathfrak{m}^2 \mbox{, and } \mathcal{F}(N,E^+).\]

We first recall the behaviour of the affine Fukaya category under covers.
Because $\phi: N \setminus E \To M \setminus D$ is an (unbranched) cover, and satisfies Condition \ref{condition:cov}, it induces an injective morphism of grading data,
\[ \bm{p}: \bm{G}(N, E) \To \bm{G}(M, D).\]
We recall from Proposition \ref{proposition2:2affcovers} that there is a fully faithful embedding
\[ \bm{p}^*\mathcal{F}(M \setminus D^+) \To \mathcal{F}(N \setminus E^+).\] 

\begin{remark2}
It is necessary for us to restrict to Lagrangians lying in $M \setminus D^+$ and $N \setminus E^+$ for this embedding to exist, even on the level of objects: because the map $\phi: E^+ \To D^+$ need not respect the symplectic structure, the lift of a Lagrangian submanifold which intersects $D^+$ need not be a Lagrangian submanifold in $N$.
\end{remark2}

We now observe that, if it were possible to choose perturbation data in $M$ such that Conditions \ref{condition:pertham} and \ref{condition:pertac} were satisfied, and such that $\mathcal{F}(M,D^+,\bm{a})$ were defined, then Lemma \ref{lemma2:pertint} would imply, by a similar argument to Proposition \ref{proposition2:2affcovers}, that there is a fully faithful embedding
\[ \bm{p}^*\mathcal{F}(M, D^+,\bm{a}) \To \mathcal{F}(N,E^+).\]

\begin{remark2}
\begin{sloppypar}
We observe that $\bm{p}^* R(M,D,\bm{a}) \cong R(N,E)$ is exactly the $\bm{G}(N,E)$-graded coefficient ring over which $\mathcal{F}(N,E^+)$ is defined; this follows immediately from the definition of $R(M,D,\bm{a})$ (Definition \ref{definition2:2relring}) and Lemma \ref{lemma2:branchps}. 
\end{sloppypar}
\end{remark2}

However, we recall (see discussion in Section \ref{subsec:2modulibranch}) that it is {\bf not} possible to guarantee sufficient regularity of all of our moduli spaces under Conditions \ref{condition:pertham} and \ref{condition:pertac}, when $M$ or $N$ is  Calabi-Yau, and nor have we rigorously defined $\mathcal{F}(M,D^+,\bm{a})$, so we can not make this statement.
We circumvent this problem by an ad-hoc method, which we now describe.

We recall, from the discussion at the end of Section \ref{subsec:2modulibranch}, that it {\bf is} possible to obtain regularity for moduli spaces with only a single marked point, under Conditions \ref{condition:pertham} and \ref{condition:pertac} (and assuming that $D$ has $\ge 2$ irreducible components). 
Therefore, we can make sure that the result is true `to first order' (see Section \ref{subsec:2firstod}):

\begin{proposition2}
\label{proposition2:2branchedfirst}
Suppose that:
\begin{itemize}
\item $\phi: (N,E^+) \To (M,D^+)$ is an $\bm{a}$-branched cover of K\"{a}hler$^+$ pairs;
\item $N$ is either positively  monotone or  Calabi-Yau;
\item $E$ and $D$ have $\ge 2$ irreducible components;
\item $\phi|_{N \setminus E}$ satisfies Condition \ref{condition:cov}.
\end{itemize}
Then there exists a $\bm{G}(M, D)$-graded $A_{\infty}$ category $\mathcal{F}(\phi)$ over $R(M,D,\bm{a})$, and valid choices of perturbation data for $\mathcal{F}(N,E^+)$ and $\mathcal{F}(M,D^+,\bm{a})/\mathfrak{m}^2$, such that there exists a strict isomorphism of $R(M,D,\bm{a})/\mathfrak{m}^2$-linear $A_{\infty}$ categories,
\[\mathscr{G}_1: \mathcal{F}(\phi) /\mathfrak{m}^2 \To \mathcal{F}(M,D^+,\bm{a}) /\mathfrak{m}^2,\]
and a full strict embedding of $R(N,E)$-linear $A_{\infty}$ categories,
\[ \mathscr{G}_2: \bm{p}^*\mathcal{F}(\phi) \To \mathcal{F}(N,E^+).\]
\end{proposition2}
\begin{proof}
We define $\mathcal{F}(\phi)$ to be exactly the same as $\mathcal{F}(M,D^+,\bm{a})$ as a $\bm{G}(M,D)$-graded pre-category. 
By Lemma \ref{lemma2:firstordlift} and the preceding argument, we can choose Floer and perturbation data for our first-order moduli spaces in such a way that there is a strict embedding
\[ \bm{p}^*\mathcal{F}(M,D^+,\bm{a})/\mathfrak{m}^2 \To \mathcal{F}(N,E^+)/\mathfrak{m}^2.\]

We now extend these choices of perturbation data for $\mathcal{F}(N,E^+)$ to the moduli spaces of disks with $\ge 2$ marked points. 
We choose our perturbation data to be $\Gamma$-equivariant, where $\Gamma$ is the covering group of $\phi|_{N \setminus E}$, so that the $A_{\infty}$ structure maps in $\mathcal{F}(N,E^+)$ are strictly $\Gamma$-equivariant. 
One might be concerned that the requirements of $\Gamma$-equivariance and regularity of our choices of perturbation data are incompatible. 
In fact they are not: because $\Gamma$ acts freely on $N \setminus E^+$, it acts freely on the set of lifts of a given Lagrangian $L \subset M \setminus D^+$ to $N \setminus E^+$, and hence $\Gamma$ acts freely on the set of choices of perturbation data for moduli spaces of pseudo-holomorphic disks in $N$ (for our purposes, we need only consider Lagrangians in $N \setminus E^+$ which are lifts of Lagrangians in $M \setminus D^+$). 
This means we can choose regular perturbation data for one element of each $\Gamma$ orbit, then extend this to the rest of the orbit via the action of $\Gamma$ (see \cite[Section 8b]{Seidel2008}). 

By counting rigid elements of the resulting moduli spaces of pseudo-ho\-lo\-mor\-phic disks in $N$, we define $A_{\infty}$ structure maps $\mu^*$ on $\bm{p}^*\mathcal{F}(\phi)$ (and they define an $A_\infty$ structure, because $N$ is either positively  monotone or  Calabi-Yau).
Because the perturbation data are chosen to be $\Gamma$-equivariant, the $A_{\infty}$ structure maps lie in the $\Gamma$-equivariant part of the Hochschild cochains:
\[ \mu^* \in CC^2_{\bm{G}(N,E)}(\bm{p}^*\mathcal{F}(\phi))^{\Gamma}.\]
We now recall the isomorphism of Remark \ref{remark2:pullcc}:
\[CC^*_{\bm{G}(N,E)}(\bm{p}^*\mathcal{F}(\phi))^{\Gamma} \cong \bm{p}^*CC^*_{\bm{G}(M,D)}(\mathcal{F}(\phi)).\]
Therefore, the $A_{\infty}$ structure $\mu^*$ on $\bm{p}^*\mathcal{F}(\phi)$, which we constructed above, corresponds to an $A_{\infty}$ structure on $\mathcal{F}(\phi)$.
By Definition \ref{definition2:2pullainf}, and by our choice of perturbation data on the image of $\bm{p}^*\mathcal{F}(\phi)$ in $\mathcal{F}(N,E^+)$, there is a strict $A_{\infty}$ embedding
\[ \mathscr{G}_2: \bm{p}^*\mathcal{F}(\phi) \To \mathcal{F}(N,E^+).\]
This completes the construction of $\mathcal{F}(\phi)$.
 \end{proof}

Now we would like to relate $\mathcal{F}(M,D)/\mathfrak{m}^2$ to $\mathcal{F}(M, D,\bm{a})/\mathfrak{m}^2$. 
Again, we will only relate the first-order part of the deformations. 
We recall from Section \ref{subsec:2firstod} that a first-order deformation of the $A_{\infty}$ category $\mathcal{F} := \mathcal{F}(M \setminus D)$ over $R = R(M,D)$ consists of $A_{\infty}$ structure maps 
\[ \mu^* = \mu_0^* + \mu_1^*,\]
where $\mu_0^*$ gives the structure maps of $\mathcal{F}$, and the $A_{\infty}$ relations say that $\mu_1^*$ defines a class
\[ [\mu_1] \in HH^2_{\bm{G}(M,D)}(\mathcal{F}, \mathcal{F} \otimes R^1).\]

\begin{theorem2}
\label{theorem2:2defclassram}
Let $\mathcal{F} := \mathcal{F}(M \setminus D)$, $\bm{G} := \bm{G}(M ,D)$, and let
\[ [\mu_{1,\bm{1}}] := \sum_{j=1}^k r_j \alpha_j \in HH^2_{\bm{G}}(\mathcal{F},\mathcal{F} \otimes R(M,D)^1)\]
be the first-order deformation class of $\mathcal{F}(M,D,\bm{1}) \cong \mathcal{F}(M,D)$. 
Then the first-order deformation class of $\mathcal{F}(M,D,\bm{a})/\mathfrak{m}^2$ is given by
\[ [\mu_{1,\bm{a}}] = \sum_{j=1}^k \pm r_j \alpha_j^{a_j} \in HH^2_{\bm{G}}(\mathcal{F},\mathcal{F} \otimes R(M,D,\bm{a})^1),\]
where the power is taken with respect to the Yoneda product on $HH^*_{\bm{G}}(\mathcal{F})$ (and we have not determined the signs). 
\end{theorem2}
\begin{proof}
We define elements
\[ \beta_j(b) \in CC^*_{\bm{G}}(\mathcal{F}),\]
for $b\ge 1$, as follows:
Let $\bm{L}$ be a tuple of anchored Lagrangian branes in $M \setminus D$, with associated generators $\bm{p}$. 
Then the coefficient of $p_0$ in $\beta_j(b)^s(p_s, \ldots, p_1)$ is given by the count of rigid elements in the moduli space $\mathcal{M}_1(\bm{p},j,b)$. 
It follows from the fact that the signed count of points in the boundary of the one-dimensional component of the moduli space $\overline{\mathcal{M}}_1(\bm{p},j,b)$ is $0$, that each $\beta_j(b)$ is a Hochschild cocycle.
Furthermore, by definition we have 
\[ [\beta_j(1)] = \alpha_j\]
and
\[ \mu_{1,\bm{a}} = \sum_{j=1}^k r_j \beta_j(a_j) .\]

We also define elements 
\[ H_j(b) \in CC^*_{\bm{G}}(\mathcal{F})\]
by counting rigid elements in $\mathcal{M}_2(\bm{p},j,b)$. 

\begin{lemma2}
\label{lemma2:defyon}
We have 
\[\beta_j(b+1) = \pm \beta_j(b) \smile \beta_j(1) \pm \partial (H_j(b))\]
in $CC^*_{\bm{G}}(\mathcal{F})$, where $\smile$ denotes the Yoneda product and $\partial$ denotes the Hoch\-schild differential.
\end{lemma2}
\begin{proof}
The result follows from the fact that the signed count of points in the boundary of the one-dimensional component of the moduli space $\overline{\mathcal{M}}_2(\bm{p},j,b)$ is $0$. 
See Lemma \ref{lemma2:m2comp} for the description of the boundary components.
The boundary points $\mathcal{M}^{1,\bm{T}}(\bm{p},j,b)$ contribute the term $\partial H_j(b)$ to the sum, the boundary points $\mathcal{M}^{2,\bm{T}}(\bm{p},j,b)$ contribute the term $\beta_j(b+1)$, and the boundary points $\mathcal{M}^{3,\bm{T}}(\bm{p},j,b)$ contribute the term $\beta_j(b) \smile \beta_j(1)$. 
 \end{proof}

It follows that, on the level of Hochschild cohomology,
\[ [\beta_j(b+1)] = \pm [\beta_j(b)] \smile \alpha_j,\]
and hence, by induction, that
\[ [\beta_j(b)] = \pm \alpha_j^{b}.\]
The result follows immediately.
 \end{proof}

\section{Morse-Bott computations in the Fukaya category}
\label{sec:2mb}

In this section, we consider the K\"{a}hler pair $(M,D) = (M^n_1,D)$ of Example \ref{example2:fermat}, and the tuple
\[ \bm{a} := (n, \ldots, n)\]
associated to the branched cover of K\"{a}hler pairs
\[ \phi: (M^n_n,D^+) \To( M^n_1,D^+)\]
from Example \ref{example2:fermatcover}.
We consider the grading datum $\bm{G}:=\bm{G}(M,D)$.
We define the $\bm{G}$-graded rings $R := R^n_1$ and $R_n := R^n_n$ from Definition \ref{definition2:2ringra}, and recall that they are the coefficient rings of $\mathcal{F}(M,D)$ and $\mathcal{F}(M,D,\bm{a})$ respectively (see Example \ref{example2:fermatring}). 
Throughout this section, we denote by
\[ Y := \Z \langle y_1, \ldots, y_n \rangle\]
the abelian group appearing in the pseudo-grading datum $\bm{H}(M,D)$.

We recall that
\[ M := \left\{ \sum_j z_j = 0\right\} \subset \CP{n-1},\]
with the divisors $D_j := \{z_j = 0\}$. 
Thus $M \cong \CP{n-2}$, and $D$ consists of $n$ hyperplanes with normal crossings. 
$M \setminus D$ is called the (generalized) {\bf pair of pants}.
 
We construct an immersed Lagrangian sphere $L^n:S^{n-2} \To M \setminus D$ in the pair of pants (summarising the construction in \cite{Sheridan2011a}). 
The main result of this section (Corollary \ref{corollary2:atypea}) is that the endomorphism algebra $CF^*(L^n,L^n)$, computed in $\mathcal{F}(\phi)$ (see Proposition \ref{proposition2:2branchedfirst}), is of type A (see Definition \ref{definition2:2typea}).

To do this, we first give a Morse-Bott description of the endomorphism algebra $CF^*(L^n,L^n)$ in the relative Fukaya category, $\mathcal{F}(M,D)$, to first order.
Structure coefficients in this description are given by counts of `holomorphic flipping pearly trees' rather than holomorphic disks.  
A holomorphic flipping pearly tree is a Morse-Bott version of a holomorphic disk, made out of holomorphic disks and Morse flowlines. 
We introduce them because it is often possible to explicitly identify moduli spaces of flipping pearly trees, and therefore to make explicit computations of the structure coefficients in the Fukaya category.

The construction is based on \cite{Sheridan2011a} (note that the original idea comes from \cite{Cornea2006}). 
The extra content here is that, whereas \cite{Sheridan2011a} describes the endomorphism algebra in the affine Fukaya category only, we will describe the endomorphism algebra in the first-order relative Fukaya category. 

There are some transversality issues in the definition of moduli spaces of flipping pearly trees, involving the possibility of unstable disk and sphere bubbles. 
In the moduli spaces of holomorphic disks that we used to define the relative Fukaya category, we avoided this problem by introducing extra internal marked points where they intersected the divisors $D$.
This approach is no longer possible for flipping pearly trees: they can intersect the divisors on their boundary. 
For this reason, we can not guarantee transversality and therefore can not give a complete Morse-Bott description of $CF^*(L^n,L^n)$ in the relative Fukaya category. 

However, for moduli spaces with very few intersections with the divisors (in particular, those which intersect only a single divisor), we can rule out any unstable disk or sphere bubbles in an ad hoc way. 
Thus, we are able to give a Morse-Bott description of $CF^*(L^n,L^n)$ to first order, and in particular to identify the first-order deformation class in $\mathcal{F}(M,D)$. 
This allows us to determine the first-order deformation class in $\mathcal{F}(M,D,\bm{a})/\mathfrak{m}^2$ via Theorem \ref{theorem2:2defclassram}, and hence in $\mathcal{F}(\phi)/\mathfrak{m}^2$ by Proposition \ref{proposition2:2branchedfirst}. 
The first-order deformation class is all we need to determine that the algebra is of type A, so the failure of transversality in our Morse-Bott model at higher order does not concern us (to clarify: $CF^*(L^n,L^n)$ is perfectly well-defined at all orders, but our Morse-Bott model for it is well-defined only to first order).

After describing the construction of the Lagrangian $L^n$ in Section \ref{subsec:2ln}, the structure of this section follows that of Section \ref{sec:2moduli}: first we introduce the moduli space of flipping pearly trees (possible domains for a holomorphic flipping pearly tree), then we describe our choice of perturbation data, then we describe the moduli space of holomorphic flipping pearly trees (pseudo-holomorphic maps into $M$), explain why transversality holds, then describe Gromov compactness.

\subsection{The Lagrangian immersion $L^n:S^{n-2} \To M \setminus D$}
\label{subsec:2ln}

We recall \cite{Sheridan2011a} the one-parameter family of Lagrangian immersions 
\[L^n_{\epsilon}: S^{n-2} \To M \setminus D,\]
 for $\epsilon>0$ sufficiently small (actually these were called $L^{n-2}$ in \cite{Sheridan2011a}; apologies for the change in notation, but it makes many formulae cleaner, cf. Remark \ref{remark2:awkward}).
We briefly recall the construction of $L^n_\epsilon$.

We consider the Lagrangian immersion $L':S^{n-2} \To M$ which is the double cover of the real locus $\RP{n-2}$ of $M$. 
If we think of
\[ S^{n-2} := \left\{ \sum_{j=1}^{n} x_j = 0, \sum_{j=1}^{n} x_j^2 = 1\right\} \subset \R^{n},\]
then the immersion is given by
\[ (x_1,\ldots,x_{n}) \mapsto [x_1:\ldots:x_{n}].\]
We construct the immersion $L^n_\epsilon$ by perturbing the immersion $L'$.

Namely, by the Weinstein Lagrangian neighbourhood theorem, $L'$ can be extended to an immersion of the radius-$\eta$ cotangent disk bundle
\[ D^*_{\eta} S^{n-2} \To M,\]
which is $J_0$-holomorphic along the zero section, and such that complex conjugation acts by $-1$ on the covector.

We construct a function $f: S^{n-2} \To \R$ by setting
\[ f(x_1,\ldots,x_{n}) = \sum_{j=1}^{n} g(x_j),\]
where $g: \R \To \R$ has the properties
\begin{enumerate}
\item $g'(x)>0$;
\item $g(-x)=-g(x)$;
\item $g(x) = x$ for $|x| < \delta$;
\item $g'(x)$ is a strictly decreasing function of $|x|$ for $|x|>\delta$;
\item $g'(x) < \delta$ for $|x| > 2\delta$,
\end{enumerate}
for some small $\delta>0$. 

We then define $L^n_{\epsilon} : S^{n-2} \To M$ to be the image of the graph of the exact one-form $\epsilon df$ in $D^*_{\eta} S^{n-2}$, under the immersion into $M$, so that $L' = L^n_0$.
The fact that $\nabla f$ is transverse to the hypersurfaces $\{x_j = 0\}$ implies that the image $L^n_{\epsilon}$ avoids the divisors $\{z_j = 0\}$ for sufficiently small $\epsilon>0$, so we obtain a Lagrangian immersion $L^n_{\epsilon}:S^{n-2} \To M \setminus D$.
It has self-intersections at the critical points of $f$ (where it intersects the other branch of the double cover). 
We define $L^n:=L^n_{\epsilon}$, for some fixed $\epsilon > 0$.
We observe that, for $n \ge 4$, $L^n$ automatically lifts to $\widetilde{\mathcal{G}}(M \setminus D)$, because $\pi_1(S^{n-2}) = 0$. 
We choose such a lift, and hence define an anchored brane structure on the Lagrangian $L^n$.

\begin{figure}
\begin{center}
\includegraphics[width=0.9\textwidth]{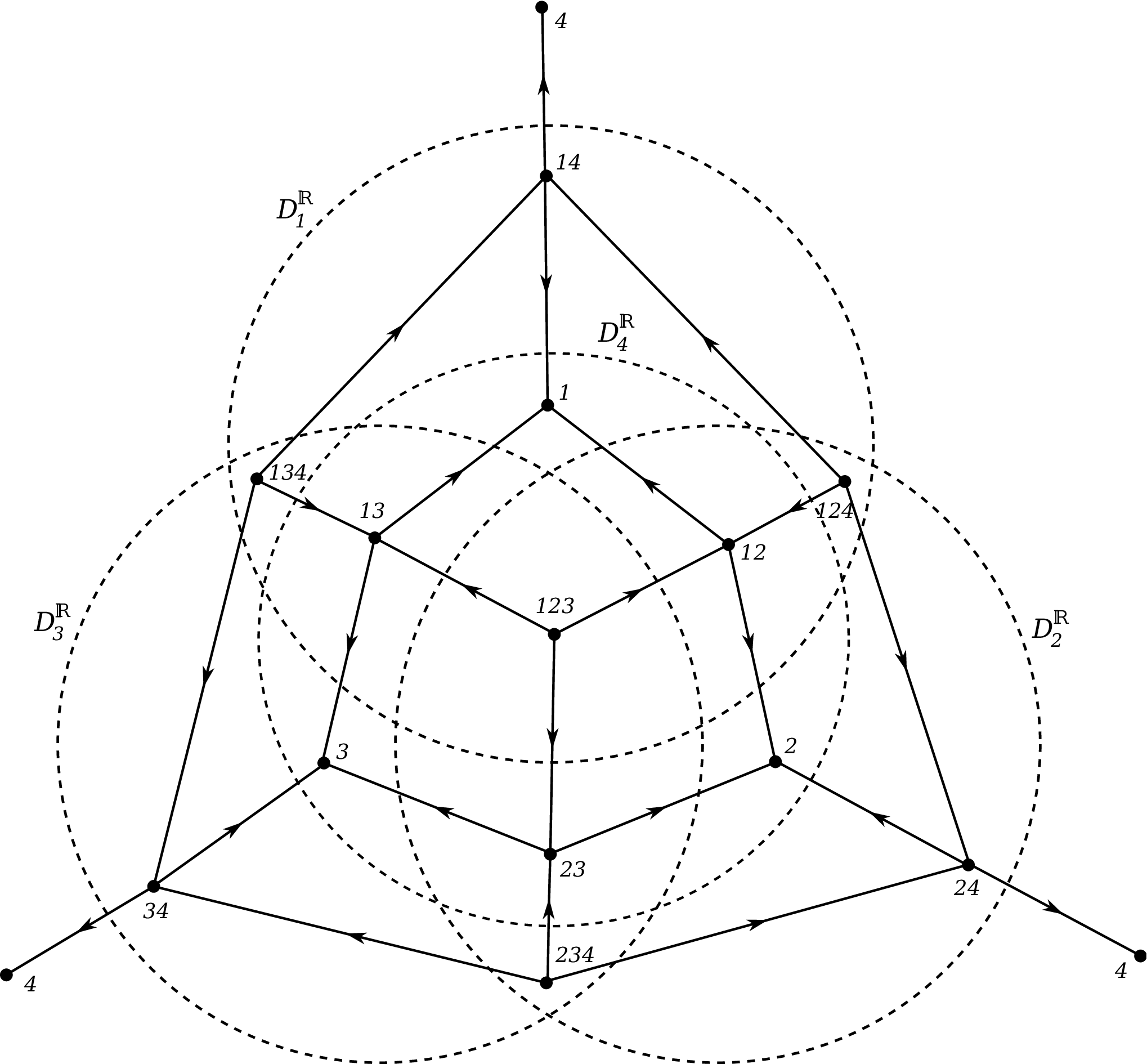}
\end{center}
\caption{The case $n=4$. The dashed circles represent the hypersurfaces $D_j^{\R}$ which are the lifts of $D_j \cap \RP{2}$, as labeled. Each region is labeled with the list of coordinates that are negative in that region (e.g., the label `$124$' means that $x_1<0,x_2<0,x_3>0,x_4<0$ in that region). The arrows represent the index-1 Morse flow lines of $\nabla f$. The dots represent critical points of $f$. The picture really lives on a sphere, and the three points labeled `$4$' should be identified (at infinity).}
\label{fig:2dualcells}
\end{figure}

The flowlines of $\nabla f$ are illustrated in Figure \ref{fig:2dualcells}, in the case $n = 4$. 
The hypersurfaces $\{x_j = 0\}$ split $S^{n-2}$ into $2^{n}-2$ regions, indexed by the proper non-empty sets $K \subset [n]$. 
Namely, $K$ corresponds to the region where coordinates $x_j$ are negative for $j \in K$ and positive for $x_j \notin K$. 
Each region contains a unique critical point $p_K$ of $f$.

The Floer endomorphism algebra $CF^*(L^n,L^n)$ can be defined, despite $L^n$ being immersed (see \cite[Section 3.1]{Sheridan2011a}), and is generated by the self-intersection points of $L^n$ (which are the points $p_K$ indexed by proper non-empty sets $K \subset [n]$), together with the Morse cohomology of $S^n$ (which we choose to have generators $p_{\emptyset}$ and $p_{[n]}$, corresponding to the identity and top class respectively). 
Thus, $CF^*(L^n,L^n)$ has generators $p_K$ indexed by subsets $K \subset [n]$. 

It follows from \cite[Proposition 3.3]{Sheridan2011a} and \cite[Proposition 3.7]{Sheridan2011a} that there is an isomorphism
\[ CF^*(L^n,L^n) \cong A\]
as $\bm{G}$-graded vector spaces, where $A \cong A_n$ is the $\bm{G}$-graded exterior algebra of Definition \ref{definition2:2alga}.

We now observe that we can define the endomorphism algebra $CF^*(L^n,L^n)$ in $\mathcal{F}(M,D)/\mathfrak{m}^2$. 
Most of the rest of this section is concerned with computing $CF^*(L^n,L^n)$ in $\mathcal{F}(M,D)/\mathfrak{m}^2$, using `flipping pearly trees'.

\begin{remark2}
\label{remark2:liftgrom}
Note that, in \cite[Section 3.1]{Sheridan2011a}, well-definedness of $CF^*(L^n,L^n)$ in the affine Fukaya category was proved by passing to the cover $M^n_n \setminus D$ of $M \setminus D$, to bypass proving Gromov compactness for immersed Lagrangians. 
One may worry that this will no longer be valid when we consider disks passing through the divisors $D$, about which the cover $M^n_n \To M$ has some branching. 
However, as we saw in the proof of Proposition \ref{proposition2:2branchedfirst}, for the first-order relative Fukaya category it is possible to choose perturbation data that lift to the branched cover; and therefore we can apply the same trick to rigorously define $CF^*(L^n,L^n)$ in $\mathcal{F}(M,D)/\mathfrak{m}^2$.
\end{remark2}

The main result we will prove is:
 
\begin{proposition2}
\label{proposition2:2fptout}
We have
\[ CF^*_{\mathcal{F}(M,D)}(L^n,L^n)/\mathfrak{m}^2 \cong (A \otimes R/\mathfrak{m}^2,\mu^*),\]
where:
\begin{itemize}
\item $A$ is the $\bm{G}$-graded vector space of Definition \ref{definition2:2alga};
\item $(A,\mu^2_0)$ coincides with the exterior algebra multiplication on $A$;
\item We have
\[ \Phi(\mu^*) = \pm u_1 \ldots u_n + \sum_{j=1}^n \pm r_j u_j \in HH^2_{\bm{G}}(A,A \otimes R / \mathfrak{m}^2)\]
where $\Phi$ is the HKR map (see Definition \ref{definition2:2phihkr}).
\end{itemize}
\end{proposition2}

\begin{remark2}
The main result of \cite{Sheridan2011a} was the zeroth-order part of Proposition \ref{proposition2:2fptout}. It remains to prove that flipping pearly trees can be made to work to first order, and that the first-order deformation classes are as claimed.
\end{remark2}

\subsection{Flipping Pearly trees}
\label{subsec:2pearlyt}

The main step in proving Proposition \ref{proposition2:2fptout} is to define a $\bm{G}$-graded, $R/\mathfrak{m}^2$-linear $A_{\infty}$ category $\mathcal{F}'$. 
$\mathcal{F}'$ has two objects, called $L$ and $L'$. 
 $L$ represents the Lagrangian immersion $L^n:S^{n-2} \To M\setminus D$, and $L'$ represents the Lagrangian immersion $L':S^{n-2} \To M$.

Morally, the endomorphism algebra of $L$ in this category will be 
\[ CF^*_{\mathcal{F}'}(L,L) \cong CF^*_{\mathcal{F}(M,D)/\mathfrak{m}^2}(L^n,L^n),\]
which we have already defined; whereas the endomorphism algebra of $L'$ (denoted $CF^*_{\mathcal{F}'}(L',L')$) will be defined in terms of pearly trees, and explicitly computable. 
This computation will suffice to compute the endomorphism algebra of $L$ up to first-order quasi-equivalence, by Lemma \ref{lemma2:similar}. 

\begin{remark2}
\begin{sloppypar}
Actually, the conventions involved in the definitions of $CF^*_{\mathcal{F}'}(L,L)$ and $CF^*_{\mathcal{F}(M,D)}(L^n,L^n)$ will be slightly different, so we will have to compare them in Lemma \ref{lemma2:fprimef}. 
\end{sloppypar}
\end{remark2} 

For the purposes of this section, $\bm{L}$ will denote a tuple of objects of $\mathcal{F}'$. 
Thus, $\bm{L}$ consists only of two types of entries: $L$ and $L'$.

\begin{definition2}
If $T$ is a semi-stable directed planar tree with labels $\bm{L}$ (see Definition \ref{definition2:2plantree}), we introduce the following notation:
\begin{itemize}
\item $V(T)$ is the set of vertices of $T$;
\item $E(T)$ is the set of edges of $T$;
\item $E'(T) \subset E(T)$ is the subset of edges with both sides labeled $L'$;
\item $F'(T)$ is the set of flags $(v,e)$ of $T$ such that $e \in E'(T)$;
\item $C(T)$ is the set of `segments' between consecutive edges around a vertex (these are indexed by pairs of consecutive flags around a vertex);
\item If $C \in C(T)$, then $L_C \in \bm{L}$ is the label associated to $C$.
\end{itemize}
\end{definition2}

\begin{figure}
\centering
\includegraphics[width=0.7\textwidth]{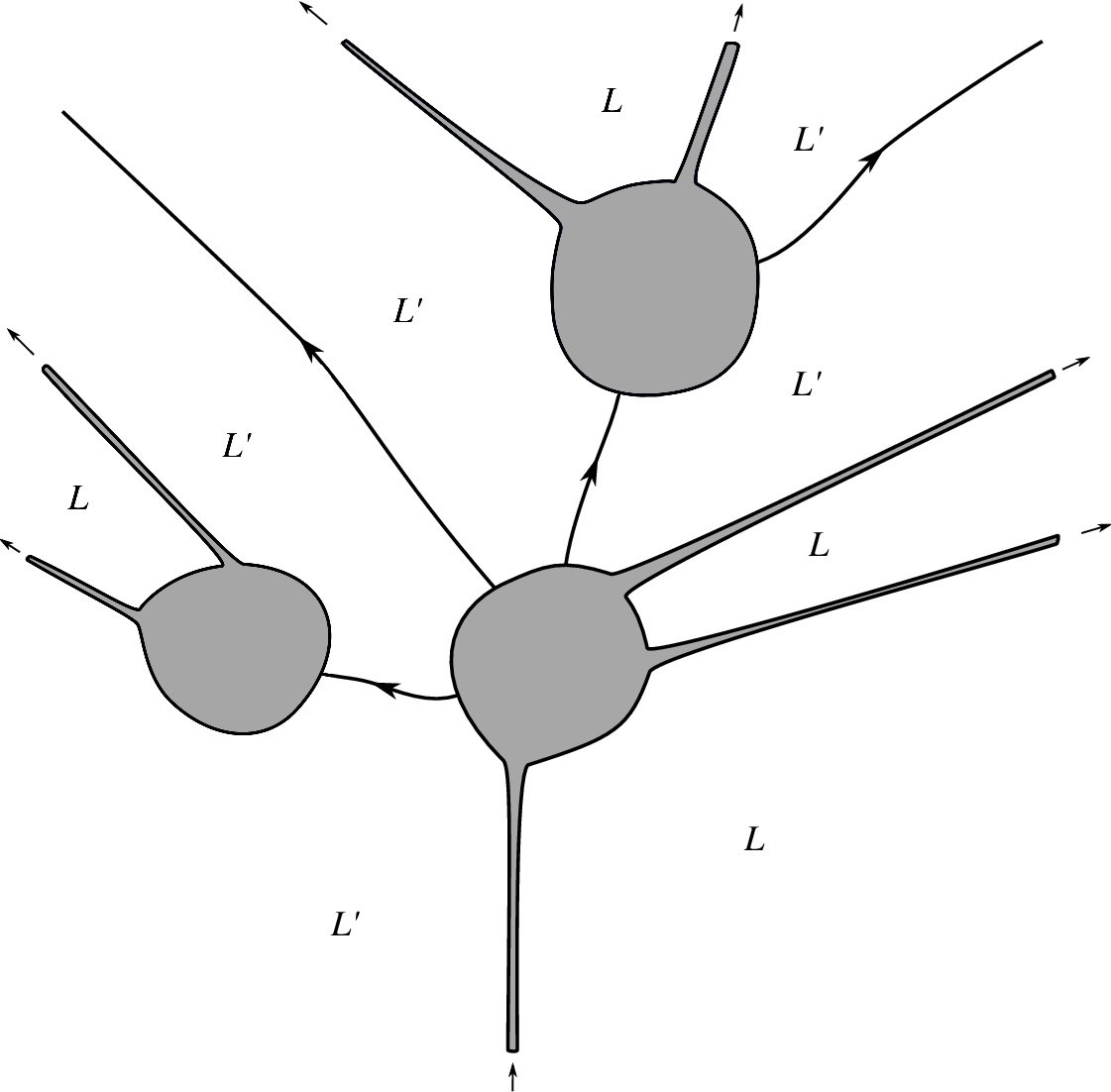}
\caption{A flipping pearly tree $S$. Observe that all internal edges have label $L'$ on either side.
\label{fig:2pearlytree}}
\end{figure}

\begin{definition2}
Let $\bm{L}$ be a tuple. 
We denote by $\mathcal{R}_3(\bm{L})$ the moduli space of {\bf flipping pearly trees}, where a flipping pearly tree $r \in \mathcal{R}_3(\bm{L})$ consists of the following data:
\begin{itemize}
\item A semi-stable directed planar tree $T_r$ with labels $\bm{L}$, such that all internal edges have both sides labeled $L'$ (i.e., all internal edges are contained in $E'(T_r)$); 
\item A designation of each edge $e \in E'(T_r)$ as either {\bf flipping} or {\bf non-flipping};
\item For each stable vertex $v \in V(T_r)$, a point $r_v \in \mathcal{R}(\bm{L}_v,\emptyset)$;
\item For each internal edge $e$, a length parameter $l_e \in [0,\infty)$.
\end{itemize}
See Figure \ref{fig:2pearlytree} for a picture of a flipping pearly tree. 
We also allow one special case: if $\bm{L} = (L',L')$, then we permit $T_r$ to have no vertices, just a single edge with both sides labeled $L'$.
\end{definition2}

Now, given a flipping pearly tree $r \in \mathcal{R}_3(\bm{L})$, we define an associated topological space $S_r$. 
There are a few special cases first:

\begin{definition2}
\label{definition2:2stripsfpt}
If $\bm{L} = (L',L')$, and $T_r$ is the tree with a single edge, then we define $S_r := \R$.
If $|\bm{L}| = 2$ but $\bm{L} \neq (L',L')$, then $S_r := \R \times [0,1]$.
\end{definition2}

Now we define $S_r$ in the remaining cases:

\begin{definition2}
Given a flipping pearly tree $r \in \mathcal{R}_3(\bm{L})$, we define a topological space $S_r$ as follows:
\begin{itemize}
\item For each semi-stable vertex $v \in V(T_r)$ with both sides labeled $L'$, we define $S_v$ to be a disk with two boundary marked points, corresponding to the edges incident to $v$.
\item For each stable vertex $v \in V(T_r)$, we define $S_v$ to be the boundary-marked disk with modulus $r_v$, with all marked points punctured except for those corresponding to edges in $E'(T_r)$ (they remain as marked points). 
These are the `{\bf pearls}'.
\item We define
\begin{eqnarray*}
S^p &:=& \coprod_{v \in V(T_r)} S_v.
\end{eqnarray*}
\item For each internal edge $e$, we define $S_e := [0,l_e]$. 
For each external edge $e$ in $E'(T_r)$, we define $S_e := \R^{\pm}$, with the $+$ or $-$ depending on the orientation of the edge.
\item We define
\begin{eqnarray*}
 S^e &:=& \coprod_{e \in E(r)} S_e.
\end{eqnarray*}
\item For each flag $f = (v,e) \in F'(T_r)$, there is a corresponding marked boundary point $m(f) \in S_v$ and boundary point $b(f) \in S_e$. 
\item We define
\[ S_r:= (S^p \sqcup S^e)/\sim\]
where
\[m(f) \sim b(f) \mbox{ for all $f \in F'(T_r)$}.\]
\end{itemize}
\end{definition2}

\begin{figure}
\centering
\includegraphics[width=0.7\textwidth]{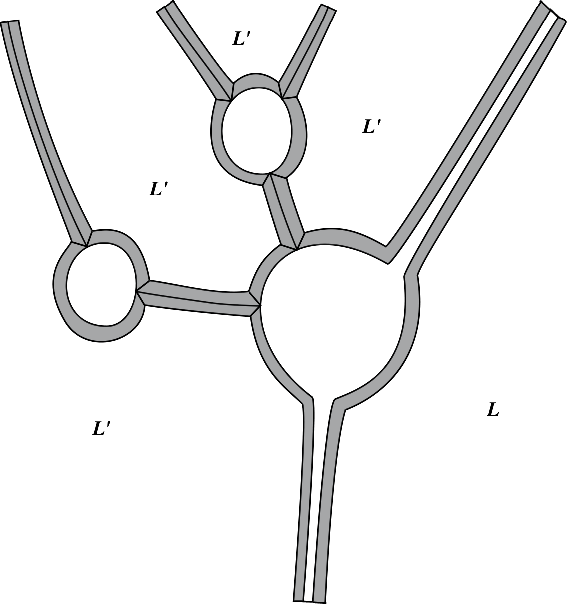}
\caption{Defining the boundary $\partial S \To S$ of a flipping pearly tree $S$, and attaching a strip (shaded in grey) along it.
\label{fig:2addstrip}}
\end{figure}

\begin{definition2}
Given $r \in \mathcal{R}_3(\bm{L})$, we also define a `boundary' $(\partial S)_r$ and a continuous map
\[ (\partial S)_r \To S_r,\]
as follows:
\begin{itemize}
\item For each segment $C \in C(T_r)$ adjacent to vertex $v \in V(T_r)$, we define $(\partial S)_C$ to be the corresponding component of the boundary of $S_v$. 
Thus, $(\partial S)_C$ is an interval, and its two ends correspond to consecutive marked points on the disk with modulus $r_v$. 
If the marked point is punctured in $S_v$ (i.e., if its sides are not both labeled $L'$), then that end of the interval $(\partial S)_C$ is open, and if the marked point remains in $S_v$ (i.e., if its sides are both labeled $L'$), then that end of the interval $(\partial S)_C$ is closed.
\item We define
\begin{eqnarray*}
(\partial S)^p & := & \coprod_{C \in C(T_r)} (\partial S)_C \mbox{, with the obvious map} \\
(\partial S)^p & \To & S^p.
\end{eqnarray*}
\item For each edge $e$, we define $(\partial S)_e := S_e \times \{0,1\}$ (two copies of $S_e$). 
\item We define
\begin{eqnarray*}
(\partial S)^e & := & \coprod_{e \in E(T_r)} (\partial S)^e \mbox{, with the obvious map} \\
(\partial S)^e & \To & S^e.
\end{eqnarray*}
\item For each flag $f = (v,e) \in F'(T_r)$, there are points $\tilde{m}_j(f) \in (\partial S)_v$ for $j = 0,1$, from the boundary components to the right and left of $m(f)$ respectively, and points $\tilde{b}_j(f) = (b(f),j) \in (\partial S)_e$, for $j = 0,1$.
\item We define
\[ (\partial S)_r := ((\partial S)^p \sqcup (\partial S)^e)/\sim,\]
where
\[ \tilde{m}_j(f) \sim \tilde{b}_j(f) \mbox{ for all $f \in F'(r)$, and $j = 0,1$}\]
(see Figure \ref{fig:2addstrip}).
\item It is clear that there is a continuous map
\[ (\partial S)_r \To S_r.\]
\end{itemize}
\end{definition2}

\begin{definition2}
An {\bf automorphism} of a flipping pearly tree is a map $S_r \To S_r$ such that each pearl gets sent to itself by a biholomorphism which preserves the marked points, and each edge gets sent to itself by a translation preserving marked points (in particular, edges are fixed by any automorphism, unless we are in the special case of an edge which is infinite in both directions).
\end{definition2}

In particular, the possible non-trivial automorphisms of a flipping pearly tree are:
\begin{itemize}
\item If $\bm{L} = (L',L')$ and $S_r = \R$, then automorphisms are translations of $\R$;
\item If $|\bm{L}|  =2$ but $\bm{L} \neq (L',L')$, then $S_r = \R \times [0,1]$ and automorphisms are translations in the $\R$-direction;
\item If $v$ is a semi-stable vertex of $T(r)$, then $S_v$ is a disk with two marked boundary points, and there is an $\R$ family of automorphisms (translations) of $S_v$ preserving the marked boundary points.
\end{itemize}

\begin{figure}
\centering
\includegraphics[width=0.7\textwidth]{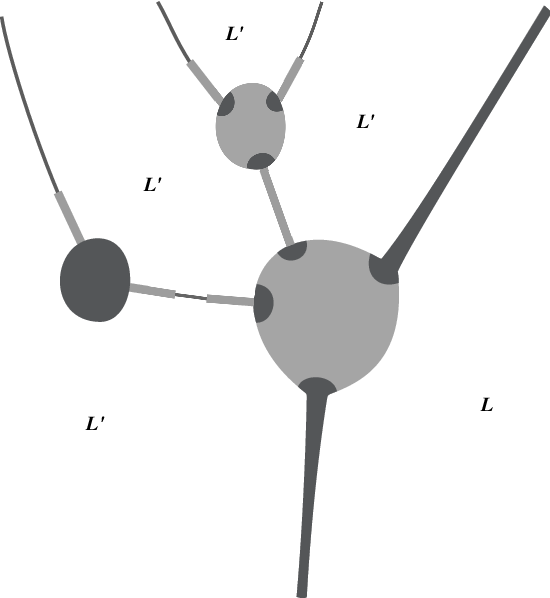}
\caption{The thick and thin regions of the same flipping pearly tree $S$ illustrated in Figure \ref{fig:2addstrip}. 
The thick regions are shown in light grey, and the thin regions in dark grey. 
Note that the unstable disk (with two marked points) is entirely thin, and also that it is possible for an internal edge to be entirely thick (when it has length $\le 2$).
\label{fig:2thickthin}}
\end{figure}

\begin{definition2}
From our universal choice of strip-like ends for the moduli spaces $\mathcal{R}(\bm{L},\emptyset)$, we can define a subset
\[ S^p_{\mathrm{thin}} \subset S^p,\]
called the {\bf thin} region, which consists of the images of strip-like ends under gluing maps (see \cite[Remark 9.1]{Seidel2008}). 
To clarify: the thin region includes a neighbourhood of each boundary marked point of a pearl, and also {\bf all} of any semi-stable pearl $S_v$ (see Figure \ref{fig:2thickthin}).
We define the corresponding {\bf thick} region
\[ S^p_{\mathrm{thick}} := S^p \setminus S^p_{\mathrm{thin}}.\]
\end{definition2}

\begin{definition2}
We define the region
\[ S^e_{\mathrm{thin}} \subset S^e\]
to be the set of points on edges which are distance $>1$ from the boundary of the edge, and
\[ S^e_{\mathrm{thick}} := S^e \setminus S^e_{\mathrm{thin}}\]
(see Figure \ref{fig:2thickthin}).
\end{definition2}

\begin{figure}
\centering
\includegraphics[width=0.9\textwidth]{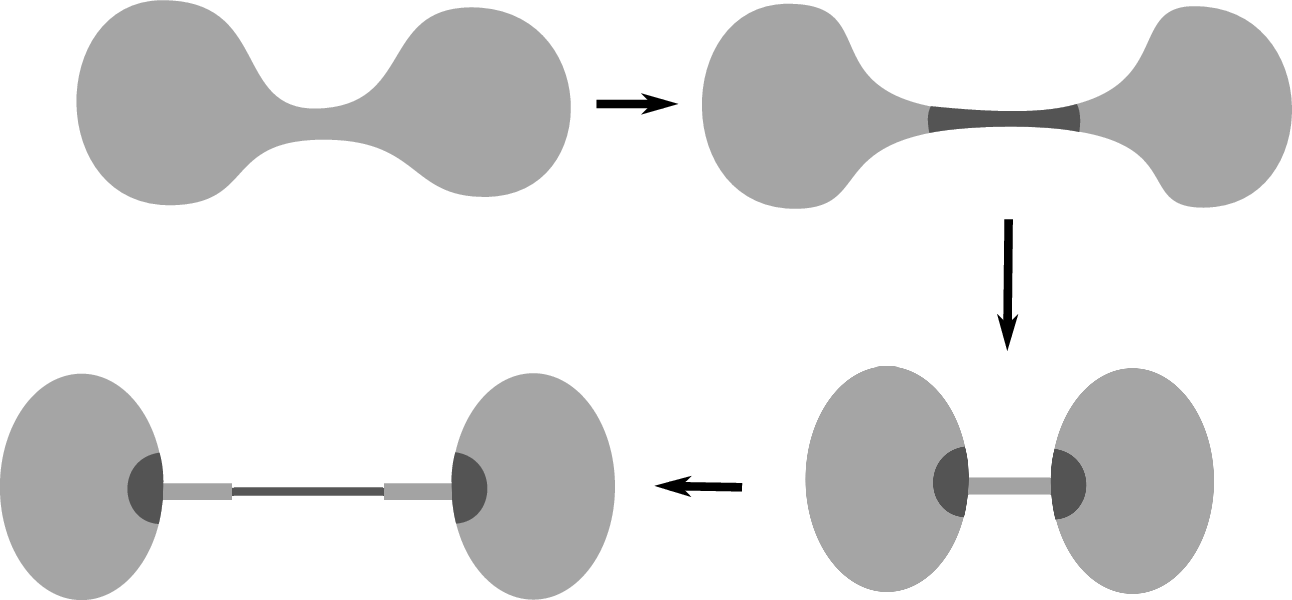}
\caption{A thick region with both sides labelled $L'$ (upper left) can stretch until it becomes a thin region (upper right), then break, becoming a thick internal edge (lower right), which then stretches until it has a thin region in its interior (lower left).
\label{fig:2morseedge}}
\end{figure}

As in \cite[Section 4.1]{Sheridan2011a}, we can define a topology on the moduli space $\mathcal{R}_3(\bm{L})$.
The important point is that thin regions with opposite sides labelled $L'$ can stretch until they `break', then become an internal edge (see Figure \ref{fig:2morseedge}). 
The difference from \cite{Sheridan2011a} is that we now allow semi-stable vertices, and this means that $\mathcal{R}_3(\bm{L})$ no longer has the structure of a manifold with boundary. 
However we will see later that the space of holomorphic maps of pearly trees into our manifold {\bf is} a manifold with boundary, which is what we need to define our algebraic structures. 

We also define a topological space $\overline{\mathcal{R}}_3(\bm{L})$, by analogy with the Deligne-Mumford compactification.
The strata of $\overline{\mathcal{R}}_3(\bm{L})$ are indexed by semi-stable directed planar trees $T$ with labels $\bm{L}$. 
Note there is no requirement that internal edges have opposite sides labeled $L'$ here.
The tree $T$ corresponds to the stratum
\[ \mathcal{R}_3^T(\bm{L}) \cong \prod_{v \in V(T)} \mathcal{R}_3(\bm{L}_v).\]
Points in this stratum correspond to flipping pearly trees, where the pearls are allowed to be nodal and the edges are allowed to have infinite length. 

\begin{remark2}
\label{remark2:semistabcomp}
The topological space $\overline{\mathcal{R}}_3(\bm{L})$ is not compact, because our flipping pearly trees can have arbitrarily many semi-stable vertices. 
In practice (see the proof of Proposition \ref{proposition2:2hfpttransv}), we have an a priori upper bound $N$ on the number of semi-stable vertices that we need consider. 
We consider the subspace 
\[\overline{\mathcal{R}}_3(\bm{L},N) \subset \overline{\mathcal{R}}_3(\bm{L}),\]
consisting of stable flipping pearly trees with $\le N$ semi-stable vertices. 
This subspace {\bf is} compact.
\end{remark2}

\subsection{Floer and perturbation data}

\begin{definition2}
We define 
\[\mathcal{H}^e := C^{\infty}(S^n,\R)\]
(think of this as the space of Morse functions on $S^n$), and
\[\mathcal{H}^p := \left\{H \in C^{\infty}(M,\R): H \mbox{ vanishes, with first derivatives, along each $D_j$}\right\}\]
(think of this as the space of Hamiltonians on $M$), and $\mathcal{J}$, the space of smooth almost-complex structures on $M$ which are compatible with $\omega$, and make the divisors $D_j$ almost-complex submanifolds.
\end{definition2}

\begin{definition2}
\label{definition2:2floerdata}
Let $\bm{L} = (L_0,L_1)$ be a $2$-element tuple. 
For each such tuple, we choose a Floer datum $(H_{\bm{L}},J_{\bm{L}})$ consisting of 
\[H_{\bm{L}} \in C^{\infty}([0,1],\mathcal{H}^p) \mbox{ and }J_{\bm{L}} \in C^{\infty}([0,1],\mathcal{J})\] 
such that:
\begin{itemize}
\item $H_{\bm{L}} = 0 $ unless $\bm{L} = (L,L)$;
\item $H_{(L,L)}$ is small, and its time-$1$ Hamiltonian flow makes $L^n$ transverse to itself;
\item $J_{(L',L')} = J_0$ is constant, equal to the standard integrable complex structure.
\end{itemize}
Now, if $(L_0,L_1) \neq (L',L')$, then we define a {\bf generator} of $CF^*_{\mathcal{F}'}(L_0,L_1)$ to be a path $p:[0,1] \To M$ which is a flowline of the Hamiltonian vector field of $H_{\bm{L}}$, such that $p(0) \in L_0$ and $p(1) \in L_1$.
One defines $CF^*_{\mathcal{F}'}(L_0,L_1)$ to be the $R/\mathfrak{m}^2$-module generated by its generators.
It is $\bm{G}$-graded.
\end{definition2}

\begin{remark2}
\label{remark2:pushing}
The $\bm{G}$-grading is defined exactly as in the affine Fukaya category (Section \ref{subsec:2eqaffuk}). 
One might worry that $L'$ intersects the divisors $D_j$, hence doesn't lie in $M \setminus D$. 
However, we simply push $L'$ off itself using $\nabla f$ (as in the definition of $L^n_{\epsilon}$), and use the pushed-off version of $L'$ for all grading (and index) computations -- see \cite[Proof of Proposition 5.3]{Sheridan2011a}. 
The point is that the grading of the relative Fukaya category arises from index computations of the relevant Fredholm operators, which are computed purely topologically.
\end{remark2}

\begin{definition2}
\label{definition2:2floerdata'}
If $(L_0,L_1) = (L',L')$, then we define the Floer datum to contain additional information, namely:
\begin{itemize}
\item the Morse function $f:S^n \To \R$;
\item another Morse function $h:S^n \To \R$, with exactly two critical points.
\end{itemize}
One defines a {\bf generator} of $CF^*_{\mathcal{F}'}(L',L')$ to be a critical point of one of the Morse functions $f$ or $h$, and $CF^*_{\mathcal{F}'}(L',L')$ to be the $R/\mathfrak{m}^2$-module generated by these critical points. 
We identify the critical points of $f$ as $p_K$ for $K \subset [n]$ proper and non-empty, and the critical points of $h$ as $p_{\emptyset}$ and $p_{[n]}$. 
Then $CF^*_{\mathcal{F}'}(L',L')$ is $\bm{G}$-graded, where $p_K$ has the same grading as the corresponding generator $\theta^K$ of $A$, where $A$ is the $\bm{G}$-graded algebra of Definition \ref{definition2:2alga}.
\end{definition2}

We recall the description of the morphism spaces in $\mathcal{F}'$ from \cite[Proof of Proposition 5.5]{Sheridan2011a}:

\begin{lemma2}
\label{lemma2:morphf}
The $\bm{G}$-graded morphism spaces in $\mathcal{F}'$ are as follows:
\begin{eqnarray*}
CF^*_{\mathcal{F}'}(L,L) &\cong&  CF^*_{\mathcal{F}'}(L',L') \\
& \cong & A \otimes R/\mathfrak{m}^2\\
& \cong & \bigoplus_{K \subset [n]} R/\mathfrak{m}^2 \cdot p_K,
\end{eqnarray*}
and
\begin{eqnarray*}
CF^*_{\mathcal{F}'}(L,L') &\cong&  CF^*_{\mathcal{F}'}(L',L) \\
& \cong & CM^*(f) \otimes R/\mathfrak{m}^2 \oplus CM^*(f) \otimes R/\mathfrak{m}^2 \\
& \cong & \bigoplus_{K \subset [n], K \neq \emptyset, [n]} R/\mathfrak{m}^2 \cdot p_K \oplus \bigoplus_{K \subset [n], K \neq \emptyset, [n]} R/\mathfrak{m}^2 \cdot q_K,
\end{eqnarray*}
where $CM^*(f)$ is the Morse complex of the function $f$ (this equation defines the generators $q_K$). 
The $\bm{G}$-grading of generators labelled by $p_K \in A$ is the same as the grading of the generator $\theta^K$ in Definition \ref{definition2:2alga}. 
The $\bm{G}$-grading of generators labelled by $q_K$ is $(|K|-1,0)$ (where $|K|-1$ is the Morse index of $q_K$). 
\end{lemma2}
\begin{proof}
See \cite[Proof of Proposition 5.5]{Sheridan2011a}.
 \end{proof}

\begin{definition2}
\label{definition2:2flipnonflip}
In each morphism space, we call the generators labelled $p_K$, where $K \subset [n]$ is not equal to $\emptyset$ or $[n]$, the {\bf flipping generators}, and the others (labelled $q_K$, $p_{\emptyset}$ or $p_{[n]}$) the {\bf non-flipping generators}. 
\end{definition2}

\begin{remark2}
The terminology `flipping' and `non-flipping' generators comes from the definition of $L^n$ as a perturbation of the double cover $S^n \To \RP{n} \hookrightarrow \CP{n}$. 
Flipping generators correspond to paths $p$ from one sheet of the cover to the opposite sheet; non-flipping generators correspond to paths from one sheet to the same sheet. 
Recall that $L'$ is the immersion
\[ L': S^{n-2} \To \RP{n-2} \hookrightarrow \CP{n-2},\]
and $L^n$ is the perturbation of $L'$ by the graph of the differential of the Morse function $f$, in a small disk cotangent bundle $D^*_\eta S^{n-2}$. 
The two sheets of $L^n$ intersect at a critical point of $f$; the critical points of $f$ are in one-to-one correspondence with the proper, non-empty subsets $K \subset [k]$, so these give rise to the flipping generators $p_K$ of $CF^*(L,L)$, for $K \neq \emptyset, [n]$. 
We also have a contribution to $CF^*(L,L)$ coming from the Morse complex of $L^n$, which arise by perturbing $L^n$ off itself by a Hamiltonian isotopy; this gives rise to the two non-flipping generators $p_\emptyset$ and $p_{[n]}$, corresponding to the generators of cohomology of $S^{n-2}$. 
$CF^*(L',L')$ is modelled on $CF^*(L,L)$, with the same flipping and non-flipping generators. 
Finally, $CF^*(L,L')$ is generated by intersection points between $L^n$ and $L'$. 
These occur at critical points of $f$, which are (again) indexed by subsets $K \subset [k]$, for $K \neq \emptyset,[n]$.  
Each critical point contributes two generators, corresponding to the two sheets of the double cover $L'$: one is flipping, called $p_K$, and one is non-flipping, called $q_K$.
\end{remark2}

\begin{definition2}
A {\bf perturbation datum} for a fixed flipping pearly tree $r \in \mathcal{R}_3(\bm{L})$ consists of the data $(K^e,K^p,J)$, where:
\begin{itemize}
\item $K^e \in C^{\infty}(S^e,\mathcal{H}^e)$;
\item $K^p \in \Omega^1(S^p,\mathcal{H}^p)$;
\item $J \in C^{\infty}(S^p,\mathcal{J})$,
\end{itemize}
such that, for each boundary component $C$ of a pearl in $S$ with Lagrangian label $L_C$,
\[ K^p(\xi)|_{L_C} = 0 \mbox{ for all $\xi \in TC \subset T(\partial S)$}.\]
\end{definition2}

\begin{definition2}
\label{definition2:2pertdatcomp}
We say that a perturbation datum is {\bf compatible with the Floer data} if, on each component of the thin regions of $S^p$ and $S^e$, the perturbation datum agrees with the corresponding (translation-invariant) Floer datum. 
Explicitly, this means that:
\begin{itemize}
\item On each strip-like end of a pearl, $(K^p,J)$ is given by the translation-invariant extension of the Floer datum $(H_{\bm{L}},J_{\bm{L}})$;
\item In a neighbourhood of each boundary marked point of a pearl, and also on all of each semi-stable pearl with both sides labeled $L'$, we have $(K^p,J) = (0,J_0)$;
\item On each thin region of a flipping edge, $K^e = f$;
\item On each thin region of a non-flipping edge, $K^e = h$.
\end{itemize}
\end{definition2}

\begin{remark2}
\label{remark2:pertaut}
Note that, if our perturbation datum is compatible with the Floer data, then it is preserved by any automorphism of the flipping pearly tree.
\end{remark2}

\begin{definition2}
We define the notion of a {\bf compatible universal choice of perturbation data} for the moduli spaces $\mathcal{R}(\bm{L})$, by analogy with \cite[Section 9i]{Seidel2008}. 
\end{definition2}

\begin{definition2}
\label{definition2:2perthol}
Let $\bm{L}$ be a tuple of objects of $\mathcal{F}'$, and $\bm{p}$ an associated set of generators.
A {\bf holomorphic flipping pearly tree} $\bm{u}$ with ends on $\bm{p}$ consists of the following data:
\begin{itemize}
\item A flipping pearly tree $r \in \mathcal{R}_3(\bm{L})$;
\item For each vertex $v \in V(T(r))$, a smooth map $u_v: S_v \To M$;
\item For each edge $e \in E(r)$ with both sides labeled $L'$, a smooth map $u_e: S_e \To S^n$;
\item A continuous map $\tilde{u}: \partial S_r \To S^n$.
\end{itemize}
We impose the following requirements on these maps:
\begin{itemize}
\item $\bm{u}$ is asymptotic to the generators $\bm{p}$ along the strip-like ends and external edges;
\item For each semi-stable vertex $v$, the map $u_v$ satisfies the pseudo-holomorphic curve equation
\[(Du_v - Y)^{0,1} = 0, \]
where, for $\xi \in TS$, $Y(\xi)$ is the Hamiltonian vector field of the function $K^p(\xi)$;
\item The maps $u_e$ satisfy the Morse flow equation
\begin{eqnarray*}
D u_e - \nabla K^e &=& 0;
\end{eqnarray*}
\item For each boundary component $C$ of a pearl $S_v$, 
\[ L_C \circ \tilde{u}|_C = u_v|_C;\]
\item For each edge $e$ with both sides labeled $L'$,
\[ \tilde{u}|_{S_e \times \{0\}} = u_e,\]
and
\[ \tilde{u}|_{S_e \times \{1\}} = \left\{ \begin{array}{rl} 
					u_e & \mbox{if $e$ is non-flipping} \\
					a \circ u_e & \mbox{if $e$ is flipping}
					\end{array} \right.\]
where $a:S^n \To S^n$ is the antipodal map.
\item If $v\in V(T(r))$ is semi-stable, then the map $u_v$ is non-constant. 
\end{itemize}
\end{definition2}

\begin{definition2}
Two holomorphic flipping pearly trees are {\bf equivalent} if they are related by an automorphism of the domain (recall from Remark \ref{remark2:pertaut} that any automorphism of the domain preserves the perturbation datum and hence acts on the space of holomorphic flipping pearly trees). 
\end{definition2}

\begin{definition2}
\label{definition2:2strip}
Given a flipping holomorphic pearly tree $\bm{u}$ as defined above, one obtains a well-defined homology class $[\bm{u}] \in H_2(M,L^n)$ as follows (see Figure \ref{fig:2addstrip}):
\begin{itemize}
\item Start with the continuous map $\bm{u}: S \To \CP{n}$ associated with the flipping holomorphic pearly tree. 
\item Glue a thin strip along the boundary $\partial S$ of the flipping pearly tree;
\item If the boundary component or edge has label $L$, then it already gets mapped to $L^n$, so we map the strip into $\CP{n}$ by making it constant along its width. 
\item If the boundary component or edge has label $L'$, then by construction, there is a continuous lift $\tilde{u}$ of the boundary of the strip to $S^n$. 
\item Thus, we can map the strip into $\CP{n}$ by letting it interpolate between the zero section and the graph of $\epsilon df$ in the Weinstein neighbourhood $D^*_{\eta} S^n$ used in the construction of $L^n$. 
Thus, boundary components of the strip with label $L'$ now lie on $L^n$.
\end{itemize}
We now define the intersection number $\bm{u} \cdot D_j$ to be the topological intersection number of this class $[\bm{u}] \in H_2(\CP{n},L^n)$ with $D_j \in H_{2n-2}(\CP{n})$, and
\[ \bm{u} \cdot D := \sum_j \left([\bm{u}] \cdot D_j\right) y_j \in Y.\]
\end{definition2}

\begin{definition2}
Let $\bm{u}$ be a holomorphic flipping pearly tree. 
For each $v \in V(T(r))$, the map $u_v$ defines a homology class in $H_2(\CP{n},\RP{n}) \cong \Z$, because its boundary gets mapped to a Weinstein neighbourhood of $\RP{n}$. 
We denote this homology class by $d_v \in \Z_{\ge 0}$ (it is non-negative because holomorphic disks have non-negative area). 
We denote the sum of all homology classes $d_v$ by $d_{\bm{u}} \in \Z_{\ge 0}$.
\end{definition2}

Now we explain how to compute these intersection numbers in a simple way. 
It helps if the holomorphic flipping pearly trees are in general position, in the following sense:

\begin{definition2}
We say that a holomorphic flipping pearly tree $\bm{u}$ is {\bf in general position} if:
\begin{itemize}
\item Each boundary component $C$ with label $L'$ is transverse to the real hypersurfaces $D_j^{\R} \subset S^n$;
\item No flipping marked points lie on the hypersurfaces $D_j^{\R}$.
\end{itemize}
\end{definition2}

\begin{lemma2}
Given a holomorphic flipping pearly tree $\bm{u}$, we can perturb the defining equations of the divisors $D_j$ so that
\begin{itemize}
\item The intersection numbers $\bm{u} \cdot D_j$ do not change;
\item $\bm{u}$ is in general position with respect to the perturbed divisors.
\end{itemize}
\end{lemma2}
\begin{proof}
See the proof of \cite[Proposition 5.1]{Sheridan2011a}.
 \end{proof}

If $\bm{u}$ is in general position, then we can split the surface defining our homology class $[\bm{u}]$ into regions $[u_v]$ corresponding to the pearls $v$, and $[u_e]$ corresponding to the edges $e$ of the pearly tree, in such a way that the boundary of each such region does not intersect the divisors $D_j$. 
We cut the pearls off from the strips in the obvious way -- since they are joined at boundary marked points, which don't lie on the hypersurfaces $D_j^{\R}$, the cuts we introduce do not intersect the divisors $D_j$. 

Thus, each region defines a class in $H_2(M,M \setminus D)$, and $[\bm{u}] \cdot D$ is equal to the sum of $[u_v] \cdot D$ and $[u_e] \cdot D$ over all pearls $v$ and edges $e$ of the pearly tree.

\begin{lemma2}
\label{lemma2:calcint}
Let $\bm{u}$ be a holomorphic flipping pearly tree in general position. 
Then we have:
\begin{itemize}
\item For each non-flipping edge $e$, $[u_e] \cdot D_j = 0$;
\item For each flipping edge $e$, $[u_e] \cdot D_j$ is equal to the topological intersection number of the edge $u_e: [0,l_e] \To S^n$ with the hypersurface $D_j^{\R}$ (this is non-negative because the gradient flow of the function $f$ crosses all hypersurfaces $D_j^{\R}$ positively);
\item For each pearl $v$, $[u_v] \cdot D_j$ is equal to the sum of the number of internal intersections of $u_v$ with $D_j$ (these are counted positively by positivity of intersections), together with $+1$ for each time a boundary lift $\tilde{u}|_C$ with label $L'$ crosses $D_j^{\R}$ in the negative direction (and $0$ if the lift crosses in the positive direction).
\end{itemize}
\end{lemma2}
\begin{proof}
See \cite[Proposition 5.1]{Sheridan2011a}.
 \end{proof}

\begin{corollary2}
\label{corollary2:calcint}
If $\bm{u}$ is a flipping holomorphic pearly tree, then the intersection numbers $\bm{u} \cdot D_j$ are non-negative.
\end{corollary2}

\begin{definition2}
Let $\bm{L}$ be a tuple, $\bm{p}$ an associated set of generators, and $\bm{d} \in Y_{\ge 0}$. 
We define $\mathcal{M}_3(\bm{p},\bm{d})$ to be the moduli space of holomorphic flipping pearly trees $\bm{u}$ with labels $\bm{L}$ and ends on $\bm{p}$, and such that 
\[ \bm{u} \cdot D = \bm{d},\]
modulo equivalence.
\end{definition2}

\begin{proposition2}
\label{proposition2:2pearlenergy}
Let $\bm{L}$ be a tuple, $\bm{p}$ an associated set of generators. 
Suppose that $(p_{K_1}, \ldots, p_{K_k})$ is the tuple obtained from $\bm{p}$ by keeping only the flipping generators (see Definition \ref{definition2:2flipnonflip}).  
Then we have
\[ [\bm{u}] \cdot D = q y_{[n]} + \sum_{j=1}^{k} y_{K_j}\]
for some $q \in \Z$, and the homology class of $[\bm{u}] \in H_2(\CP{n},\RP{n}) \cong \Z$ is given by the formula
\[[\bm{u}] = 2q + k. \]
\end{proposition2}
\begin{proof}
Follows from \cite[Lemma 5.8]{Sheridan2011a}, and a slight modification of \cite[Proposition 5.10]{Sheridan2011a}.
 \end{proof}

\begin{figure}
\centering
\includegraphics[width=0.8\textwidth]{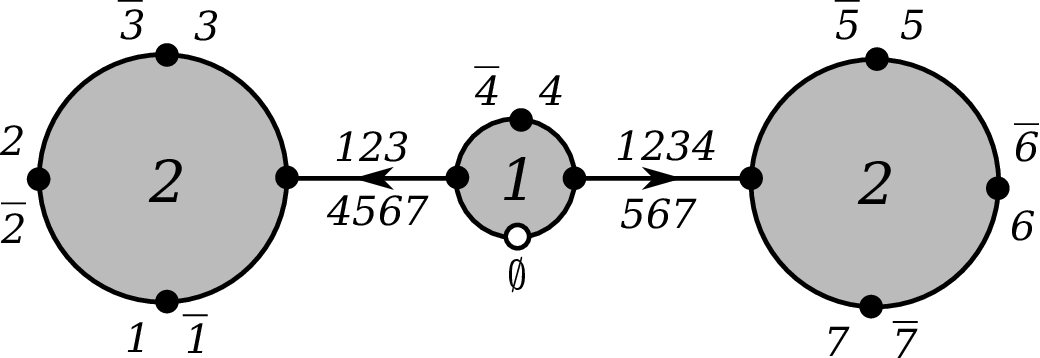}
\caption{An example of a legal labeling of a flipping holomorphic pearly tree, which might contribute to the coefficient of $p_{\emptyset}$ in the $A_{\infty}$ product $\mu^7(p_{\{1\}},\ldots,p_{\{7\}})$. 
We have illustrated a simple case, in which all external flowlines are constant because the points $p_{\{j\}}$ are maxima of the Morse function $f$.
The external label `$1$' means the set $\{1\}$, while `$\overline{1}$' means the complement $\{2,3,4,5,6,7\}$.
The big label `$1$' in the middle of a pearl means that the pearl has degree $1$.
\label{fig:2n5}}
\end{figure}

\begin{definition2}
\label{definition2:2labels}
Let $\bm{u}$ be a holomorphic flipping pearly tree in general position, all of whose boundary components are labelled $L'$. 
For each flag $f = (v,e)$ (recall that a flag consists of a vertex $v$ and an edge $e$ connected to that vertex) such that $e$ is a flipping edge of the tree, we have corresponding marked points in the boundary $(\partial S)_r$, namely
\[\tilde{m}_0(f) \sim \tilde{b}_0(f) \mbox{ and } \tilde{m}_1(f) \sim \tilde{b}_1(f).\]
The boundary map $\tilde{u}$ sends these points to antipodal regions $S^n_K, S^n_{\bar{K}}$ respectively, for some $K \subset [n]$ (recall that $S^n_K$ is defined to be the region where $x_j<0$ for $j \in K$ and $x_j>0$ for $j \notin K$). 
We attach the labels $K$ and $\bar{K}$ to the marked points $\tilde{m}_0(f)$ and $\tilde{m}_1(f)$, respectively. 
Figure \ref{fig:2n5} shows a possible labeling of a flipping holomorphic pearly tree.
\end{definition2}

\begin{lemma2}
\label{lemma2:pearlsen}
Let $\bm{u}$ be a holomorphic flipping pearly tree with all sides labelled $L'$, in general position, and equipped with labels as above. 
We decompose $[\bm{u}]$ into regions $[u_e]$, corresponding to edges, and $[u_v]$, corresponding to pearls, as before.
Then we have:
\begin{itemize}
\item For each flipping edge $e$, the label $K_0$ at the start of the edge contains the label $K_1$ at the end of the edge, and
\[ [u_e] \cdot D = y_{K_1} - y_{K_0};\]
\item For each pearl $v$ such that:
\begin{itemize}
\item The homology class is $d_v \in H_2(\CP{n},\RP{n}) \cong \Z$;
\item The points immediately after the flipping marked points of $u_v$ have labels $K_1,\ldots,K_{k_v}$, if we traverse the boundary in positive direction (in other words, the points $\tilde{m}_1(f)$ for all outgoing flags, and $\tilde{m}_0(f)$ for the incoming flag);
\end{itemize}
then we have:
\begin{itemize}
\item For some $q_v \in \Z$,
\[ [\bm{u}] \cdot D = q_v y_{[n]} + \sum_{j=1}^{k_v} y_{K_j};\]
\item The homology class $d_v$ satisfies
\[ d_v = 2q_v + k_v.\]
\end{itemize}
\end{itemize}
We recall that non-flipping edges $u_e$ do not contribute to the intersection numbers.
\end{lemma2}
\begin{proof}
For edges, the statement follows from the fact that $\nabla f$ only crosses divisors positively, and $[u_e]$ picks up an intersection point with $D_j$ each time $u_e$ crosses $D_j$ positively (compare \cite[Figure 9(b)]{Sheridan2011a}).
For pearls, the statement follows from Proposition \ref{proposition2:2pearlenergy}.
 \end{proof}

\subsection{Transversality and compactness}
\label{subsec:2transcomp}

\begin{proposition2}
\label{proposition2:2hfpttransv}
For generic choice of Floer and perturbation data, the components of $\mathcal{M}_3(\bm{p},\bm{d})$ with $|\bm{d}| \le 1$ are regular, and have the structure of topological manifolds of the expected dimension.
\end{proposition2}
\begin{proof}
The moduli spaces $\mathcal{M}_3(\bm{p},D)$ are constructed by gluing together pieces corresponding to the different possible underlying trees (see \cite[Section 4.4]{Sheridan2011a}). 
If we are to obtain a topological manifold, we need each piece to be cut out transversely, and also for the `seams' along which the pieces are glued (corresponding to holomorphic flipping pearly trees with nodal pearls, where a Morse flowline is about to form as in Figure \ref{fig:2morseedge}) to be regular.

This amounts to checking the following:
\begin{itemize}
\item The Cauchy-Riemann operator $(Du^p - Y)^{0,1}$ on each stable pearl is surjective;
\item For the semi-stable pearls, we require that the moduli space of $J_0$-holo\-mor\-phic disks with two boundary marked points, modulo translation, is regular;
\item The Morse flow operator $(Du^e - \nabla K^e)$ on each edge is surjective;
\item For moduli spaces consisting of a single Morse edge, we require that the moduli space of Morse flowlines of $f$ or $h$ is regular;
\item The restriction that marked points on pearls coincide with ends of edges is cut out transversely;
\item The restriction that marked points on pearls coincide at a node (with both sides labelled $L'$) is cut out transversely (this is the requirement that the `seams' are regular).
\end{itemize}
The Cauchy-Riemann operators on pearls corresponding to stable vertices are generically surjective, by perturbing $K^p$ (as in \cite[Section 9k]{Seidel2008}). 
The moduli spaces of semi-stable pearls with sides labelled anything other than $(L',L')$ are regular for generic choice of Floer data, by the arguments of \cite{Floer1995,oh97}. 
The semi-stable pearls with opposite sides labelled $L'$ are different -- we have required that $K^p = 0$ and $J=J_0$ on these pearls, which is not a generic condition. 
However, the moduli space of such pearls is regular by the `automatic regularity' result of \cite[Proposition 7.4.3]{mcduffsalamon} (the automatic regularity result deals with holomorphic spheres in $\CP{n-2} = M$, and our moduli space of holomorphic disks is the real locus of this moduli space, hence also regular).

The Morse flow operators are surjective for generic choice of Floer and perturbation data. 
 
The intersections of marked points on pearls with the endpoints of edges are generically transverse, by perturbing $K^e$ near the end of the edge. 

Now we deal with intersections of marked points at a node connecting two pearls, which are necessary to make the `seams' along which we glue different parts of our moduli space.
If one of the pearls involved is stable, then the intersection is transverse, because we can perturb $K^p$ on the stable pearl (see \cite[Section 4.6]{Sheridan2011a}) to move the marked point in any direction we please. 
A problem arises if both pearls are semi-stable (with both sides labelled $L'$). 
However, this situation does not arise in the moduli spaces $\mathcal{M}_3(\bm{p},\bm{d})$ with $|\bm{d}| \le 1$: any semi-stable pearl must contribute at least $1$ to one of the intersection numbers $\bm{u} \cdot D_j$. 

This follows from Lemma \ref{lemma2:pearlsen}: if $[\bm{u}] \cdot D = 0$ for a pearl with two boundary marked points with labels $K_0$ and $K_1$, then
\[ q y_{[n]} + y_{K_0} + y_{K_1} = 0.\]
If $K_0 = \emptyset$ or $[n]$, then Lemma \ref{lemma2:pearlsen} shows that $d_v = 0$ so the unstable pearl $u_v$ has zero energy and must be constant, which is not allowed. 
If $K_0 \neq \emptyset$ or $[n]$, then necessarily $q = -1$, and Lemma \ref{lemma2:pearlsen} shows that
\[ d_v = 2q_v + k_v = 0,\]
so again $u_v$ is constant.
 \end{proof}

\begin{proposition2}
\label{proposition2:nobub}
For generic choice of Floer and perturbation data, the components of $\mathcal{M}_3(\bm{p},\bm{d})$ of virtual dimension $\le 1$, and with $|\bm{d}| \le 1$, have the structure of compact topological manifolds with boundary, of the expected dimension.
The boundary strata of the one-dimensional moduli spaces correspond to nodal holomorphic flipping pearly trees (see \cite[Definition 4.33, 4.34]{Sheridan2011a}).
\end{proposition2}
\begin{proof}
This essentially follows from Gromov compactness, as outlined in \cite[Proposition 4.6]{Sheridan2011a}. 
However, we must do a bit more work in this case: we must check that the number of semi-stable pearls is bounded, so that we are only gluing together finitely many pieces to make our moduli space (see Remark \ref{remark2:semistabcomp}). 
This is true because the homology class $d_{\bm{u}}$ of any $\bm{u} \in \mathcal{M}_3(\bm{p},\bm{d})$ is fixed, by Proposition \ref{proposition2:2pearlenergy}, and any semi-stable pearl $v$ contributes at least $1$ to $d_{\bm{u}}$, because it is required to be non-constant.

Furthermore, we must rule out the possibility of sphere and disk bubbling in our moduli spaces. 
Sphere bubbling is easy to rule out: any non-constant sphere bubble must intersect each divisor $D_j$ at least once, hence contribute at least $n$ to $|\bm{d}|$. 
So sphere bubbling does not happen in moduli spaces with $|\bm{d}| \le 1$. 

Disk bubbling needs a little more work.  
Suppose a holomorphic disk bubbles off some pearl in a holomorphic flipping pearly tree. 
Let us denote by $u_1$ the disk, and by $u_2$ the rest of the holomorphic flipping pearly tree.
We will assume that $u_1$ has boundary on $L'$ (the case with boundary on $L^n$ is not very different). 
We can regard this configuration as a holomorphic flipping pearly tree with a vertex $v$ of valence $1$, connected by a Morse edge of length $0$ to the rest of the holomorphic flipping pearly tree.
We will show that, if $u_1$ is non-constant, then its virtual dimension is $\ge n+1$.
This is the same as showing that its Maslov index is at least $3$. 
This will show that disk bubbling generically does not happen, because the rest of the holomorphic flipping pearly tree has virtual dimension $<0$.
 
Firstly, by Lemma \ref{lemma2:pearlsen} (which works exactly the same if there are pearls of valence $1$), if the edge is non-flipping then the disk has zero energy, and hence is constant.

So let us assume that the edge is flipping, and the label attached to it is $K \subset [n]$ (as in Definition \ref{definition2:2labels}). 
By Lemma \ref{lemma2:pearlsen}, we have
\[ [u_1] \cdot D = q y_{[n]}+ y_K.\]
Because $|\bm{d}| \le 1$, it must be that $[u_1] \cdot D = y_j$, $q = 0$, and $K = \{j\}$ for some $j \in [n]$. 
It follows, by Lemma \ref{lemma2:pearlsen}, that the homology class of $u_1$ is $d_1 = 1$. 
Thus, the virtual dimension of $u_1$ is (by the real analogue of Lemma \ref{lemma2:indor0})
\[ \mathrm{v.d.}(u_1) = n + d_1(n+1) - 3 + 1 = 2n-1 \ge n+1\]
(because $n \ge 2$). 

The connect sum formula now says that the virtual dimension of $\bm{u}$ is
\[ \mathrm{v.d.}( \bm{u}) = \mathrm{v.d.}(u_1) + \mathrm{v.d.}(u_2) + 1 - n \le 1\]
(since we are considering moduli spaces of virtual dimension $\le 1$), and hence
\[  \mathrm{v.d.}(u_2) \le n - \mathrm{v.d.}(u_1) \le -1.\]
We saw in Proposition \ref{proposition2:2hfpttransv} that moduli spaces of holomorphic flipping pearly trees are generically regular, hence a moduli space of virtual dimension $<0$ is generically empty. 
So we may conclude that no disk bubbling occurs for generic choices of Floer and perturbation data.
 \end{proof}

\subsection{Morse-Bott model for the first-order Fukaya category}

This section contains the proof of Proposition \ref{proposition2:2fptout}.

We first complete the definition of the $\bm{G}$-graded, $R/\mathfrak{m}^2$-linear $A_{\infty}$ category $\mathcal{F}'$.
The $\bm{G}$-graded morphism spaces are the $R/\mathfrak{m}^2$-module freely generated by their generators, and are described explicitly in Lemma \ref{lemma2:morphf}.
Given a tuple $\bm{L}$ with associated generators $\bm{p} = (p_0,p_1,\ldots,p_s)$, the coefficient of $r^{\bm{d}}p_0$ in the $A_{\infty}$ product $\mu^s(p_s,\ldots,p_1)$ is given by the signed count of holomorphic flipping pearly trees in the zero-dimensional component of the moduli space $\mathcal{M}_3(\bm{p},\bm{d})$ (by analogy with the usual definition of the Fukaya category). 
By the standard argument, the composition maps $\mu^s$ satisfy the $A_{\infty}$ relations. 

Furthermore, they are $\bm{G}$-graded, by the same argument as for the relative Fukaya category (keeping in mind Remark \ref{remark2:pushing}).
To see why, recall that the grading depended on the index theory of a Cauchy-Riemann operator coming from a bundle pair $(D^2,u^*TM,F)$, where $u:D^2 \To M \setminus D$ was some smooth map, and $F$ a lift of $\partial u$ to $\mathcal{G}(M \setminus D)$. 
Given a holomorphic flipping pearly tree $\bm{u}$, we can construct (by a modification of the construction of the smooth surface representing the homology class $[\bm{u}]$) a smooth boundary-punctured disk $\tilde{\bm{u}}$ mapping to $M \setminus D$, with boundary on $L^n$, such that the index of the associated Cauchy-Riemann operator is equal to the index of the Fredholm operator defining the moduli space of holomorphic flipping pearly trees near $\bm{u}$. 
Note that $\tilde{\bm{u}}$ will {\bf not} be holomorphic itself: we are simply using it to compute the index topologically.

\begin{lemma2}
\label{lemma2:llprimeiso}
For small enough $\epsilon>0$, the objects $L^n_{\epsilon}$ and $L'$ are quasi-isomorphic in $\mathcal{F}' /\mathfrak{m}$.
\end{lemma2}
\begin{proof}
See \cite[Proposition 5.5]{Sheridan2011a}.
 \end{proof}

\begin{corollary2}
\label{corollary2:llprimedef}
The $R/\mathfrak{m}^2$-linear $A_{\infty}$ algebras $CF^*_{\mathcal{F}'}(L,L)$ and $CF^*_{\mathcal{F}'}(L',L')$ are first-order quasi-equivalent in the sense of Definition \ref{definition2:similar}.
\end{corollary2}
\begin{proof}
Follows from Lemma \ref{lemma2:llprimeiso} and Lemma \ref{lemma2:similar}.
\end{proof}

\begin{lemma2}
\label{lemma2:fprimef}
For appropriate choices of perturbation data, the $R/\mathfrak{m}^2$-linear $A_{\infty}$ algebras $CF^*_{\mathcal{F}(M,D)/\mathfrak{m}^2}(L^n,L^n)$ and $CF^*_{\mathcal{F}'}(L,L)$ are first-order quasi-equiv\-a\-lent in the sense of Definition \ref{definition2:similar}.
\end{lemma2}
\begin{proof}
Suppose that $\bm{L}$ is a tuple, all of whose entries are $L^n$, $|\bm{L}| \ge 2$, and $\ell$ is a labelling such that $\bm{d}(\ell) = y_j$ for some $j$. 
Then there is a forgetful map $\mathcal{R}(\bm{L},\ell) \To \mathcal{R}_3(\bm{L})$, obtained by simply forgetting the internal marked point. 
We choose Floer and perturbation data for the moduli spaces $\mathcal{R}(\bm{L},\ell)$, with $|\bm{L}| \ge 2$, by pulling back the corresponding perturbation data from $\mathcal{R}_3(\bm{L})$ via the forgetful map. 

We use these pulled-back Floer and perturbation data to define moduli spaces of pseudo-holomorphic disks $\mathcal{M}(\bm{p},\ell)$, for $|\bm{p}| \ge 2$. 
Note that, under these assumptions, there is a forgetful map $\mathcal{M}(\bm{p},\ell) \To \mathcal{M}_3(\bm{p},\bm{d}(\ell))$. 
This map is clearly bijective: given a holomorphic disk with topological intersection number $1$ with the divisor $D_j$, there exists a unique internal point which gets mapped to $D_j$ (uniqueness follows from positivity of intersections). 
The case $|\bm{L}| = 2$ is slightly different: $\mathcal{M}_3(\bm{p},y_j)$ is a quotient of a moduli space of holomorphic strips by translation, whereas $\mathcal{M}(\bm{p},\ell)$ is a moduli space of holomorphic strips with an internal marked point getting mapped to $D_j$. 
However, this map is still bijective: given a strip with topological intersection number $1$ with divisor $D_j$, there is a unique internal point which gets mapped to $D_j$. 
Introducing this as a marked point defines a class in $\mathcal{M}(\bm{p},\ell)$, and the result is independent of translation of the domain in the original strip. 

Therefore, the $R/\mathfrak{m}^2$-linear $A_{\infty}$ algebra whose structure maps count rigid elements of $\mathcal{M}(\bm{p},\ell)$ coincides with the $A_{\infty}$ algebra whose structure maps count rigid elements of $\mathcal{M}_3(\bm{p},\bm{d}(\ell))$, which by definition is $CF^*_{\mathcal{F}'}(L,L)$. 
Note that we never defined moduli spaces $\mathcal{M}(\bm{p},\ell)$ with $|\bm{p}| = 1$; however this doesn't matter, because no disks can bubble off in the moduli spaces we are considering, by the argument in the proof of Proposition \ref{proposition2:nobub}.

We would now like to say that this choice of perturbation data on the moduli spaces $\mathcal{R}(\bm{p},\ell)$ is a valid choice of perturbation data in the definition of the endomorphism algebra $CF^*(L^n,L^n)$, computed in $\mathcal{F}(M,D)/\mathfrak{m}^2$; if this were true, then the result would follow immediately. 
Unfortunately, it is not true: in the definition of $\mathcal{F}(M,D)$, one allows holomorphic disks with one boundary component and a single internal marked point. 
Perturbation data on $\mathcal{R}(\bm{L},\ell)$ are required to be consistent with respect to the Gromov compactification $\overline{\mathcal{R}}(\bm{L},\ell)$, which includes strata in which a disk with one boundary component and an internal marked point bubble off an element of $\mathcal{R}(\bm{L})$. 
If the perturbation data on $\mathcal{R}(\bm{L},\ell)$ are pulled back from $\mathcal{R}_3(\bm{L},\bm{d}(\ell))$, then consistency would require the perturbation data on the disk bubble with one boundary component to be constant. 
However, the perturbation data on this disk bubble is also required to restrict to the Floer data on the strip-like end; these requirements are incompatible. 

Instead, we follow the strategy used in showing the relative Fukaya category is independent of choices made in its construction up to quasi-equivalence: we consider an $R/\mathfrak{m}^2$-linear category with two objects, both of which are identified with the Lagrangian $L^n$. 
If $\bm{L}$ is a tuple consisting entirely of copies of the first object, then we define Floer and perturbation data for the moduli space $\mathcal{R}(\bm{L},\ell)$ as if we were defining the relative Fukaya category $\mathcal{F}(M,D)/\mathfrak{m}^2$. 
If $\bm{L}$ is a tuple consisting entirely of copies of the second object, then we define Floer and perturbation data for the moduli space $\mathcal{R}(\bm{L},\ell)$ by pulling back from the moduli spaces $\mathcal{R}_3(\bm{L},\ell)$. 
We extend these to choices of Floer and perturbation data for the moduli spaces $\mathcal{R}(\bm{L},\ell)$, where $\bm{L}$ is now allowed to be any tuple consisting of copies of either object. 
We use these to define moduli spaces of pseudo-holomorphic disks, whose counts define the coefficients of the structure maps of an $R/\mathfrak{m}^2$-linear $A_{\infty}$ category. 
Note that the argument in the proof of Proposition \ref{proposition2:nobub} shows that no disks with a single boundary component can bubble off in the resulting moduli spaces $\mathcal{M}(\bm{p},\ell)$, so the $A_{\infty}$ relations do hold.

It then follows from Lemma \ref{lemma2:similar} that the endomorphism algebras of the two objects are first-order quasi-equivalent; this completes the proof.
 \end{proof}

Combining Corollary \ref{corollary2:llprimedef} and Lemma \ref{lemma2:fprimef}, we obtain:

\begin{corollary2}
\label{corollary2:llprimesimilar}
The $R/\mathfrak{m}^2$-linear $A_{\infty}$ algebras $CF^*_{\mathcal{F}(M,D)}(L^n,L^n)/\mathfrak{m}^2$ and $CF^*_{\mathcal{F}'}(L',L')$ are first-order quasi-equivalent in the sense of Definition \ref{definition2:similar}.
\end{corollary2}

Now we recall (Lemma \ref{lemma2:morphf}) that the underlying $\bm{G}$-graded vector space of $CF^*_{\mathcal{F}'}(L',L')$ is $A$. 
So $CF^*_{\mathcal{F}'}(L',L')$ has the form $(A,\mu^*)$.

\begin{lemma2}
As a $\bm{G}$-graded algebra, we have
\[ (A,\mu^2_0) \cong (A, \wedge),\]
where `$\wedge$' denotes the usual product on the exterior algebra.
\end{lemma2}
\begin{proof}
See \cite[Theorem 5.12]{Sheridan2011a}.
 \end{proof}

Now recall from Section \ref{subsec:2comp} that there is a map
\[ \Phi: CC^*_{\bm{G}}(A, A \otimes R) \To R \llbracket U \rrbracket  \otimes A .\]

\begin{lemma2}
\label{lemma2:lprimedef}
We have
\[ \Phi \left( \mu^{*} \right) =  \pm u_1 \ldots u_n + \sum_{j=1}^n  \pm r_j u_j + \mathfrak{m}^2.\]
\end{lemma2}
\begin{proof}
See \cite[Proposition 5.15]{Sheridan2011a} for the proof that 
\[ \Phi(\mu_0^{*}) = \pm u_1 \ldots u_n.\]
Our aim now is to calculate the first-order terms in $\Phi(\mu^*)$.
By Lemma \ref{lemma2:defclasshh}, we know that the degree-$2$ part of
\[ HH^*_{\bm{G}}(A, A \otimes R^1) \subset R^1 \llbracket U \rrbracket  \otimes A  \]
is generated by the elements $r_j u_j$.
Thus, the first-order part of $\Phi(\mu^*)$ can be written as
\[\sum_{j=1}^{n} c_j r_j u_j \]
for some numbers $c_j \in \C$.

The number $c_j$ is given by the count of holomorphic flipping pearly trees $\bm{u}$ in the moduli space $\mathcal{M}_3((p_{\emptyset}, p_{\{j\}}), y_j)$.
Such a holomorphic flipping pearly tree must be a chain of semi-stable $J_0$-holomorphic pearls (see Definition \ref{definition2:2pertdatcomp}).

By Proposition \ref{proposition2:2pearlenergy}, the homology class of such $\bm{u}$ is $1$. 
Because semi-stable pearls are not allowed to be constant, this means there can only be a single semi-stable pearl. 
This pearl is a $J_0$-holomorphic disk with boundary on $\RP{n}$,
\[ u_v: (D, \partial D) \To (\CP{n},\RP{n}),\]
together with two marked boundary points, considered up to translation. 
The homology class $[u_v] \in H_2(\CP{n},\RP{n})$ is $1$. 
Thus, $u_v$ is one half of a $J_0$-holomorphic sphere of degree $1$ in $\CP{n}$. 
That is, it is one half of a complex line in $\CP{n}$, and its boundary is a real line in $\RP{n}$.

The Morse index of the input $p_{\{j\}}$ is $n$ (see \cite[Corollary 2.11]{Sheridan2011a}), so the corresponding Morse edge must be constant: this means we just have a point constraint that our real line must pass through the point $L'(p_{\{j\}})$. 
Similarly, the Morse index of $p_{\emptyset}$ is $0$, so the corresponding Morse edge must be constant: this means we have a point constraint that our real line must pass through the point $L'(p_{\emptyset})$.
There is a unique real line through the points $L'(p_{\{j\}})$ and $L'(p_{\emptyset})$, so the pearl must be one half of the corresponding complex line. 

Furthermore, the pearl must admit a lift $\tilde{u}_v$ of the boundary to $S^n$, which changes sheets at $p_{\{j\}}$ but not at $p_{\emptyset}$. 
Since $u_v$ is one half of a complex line, this lift $\tilde{u}_v$ must be one half of a great circle, from  $p_{\{j\}}$ to its antipode $p_{\overline{\{j\}}}$, and passing through the point $p_{\emptyset}$. 
Clearly, there is exactly one such half-great circle. 
Furthermore, the orientation of this half-great circle determines uniquely which half of the complex line we must take. 
Thus, we have uniquely determined our holomorphic flipping pearly tree $\bm{u}$. 

It follows from Proposition \ref{proposition2:2pearlenergy} that $[\bm{u}] \cdot D = y_j$.

We have shown that the moduli space $\mathcal{M}((p_{\{j\}},p_{\emptyset}), y_j)$ contains a unique element $\bm{u}$. 
Thus, each coefficient $c_j$ must be $\pm 1$.
 \end{proof}

\begin{remark2}
We could also have calculated $[\bm{u}] \cdot D$ using Lemma \ref{lemma2:calcint}, and the exercise helps us to get a picture of $\bm{u}$: the edges $u_e$ are constant, hence do not contribute to $[\bm{u}] \cdot D$. 
The pearl $u_v$ does not intersect any divisor $D_i$ in its interior (it is half of a complex line, hence intersects the divisor $D_i$ exactly once, and that intersection is on the real locus $\RP{n}$). 
So we get no contributions to $[\bm{u}] \cdot D$ from interior intersection points.
The boundary lift $\tilde{u}_v$ moves along a great circle from $p_{\{j\}}$ to $p_{\overline{\{j\}}}$, hence crosses the divisors $D_i$ positively for $i \neq j$, and the divisor $D_j$ once negatively. 
Therefore, we have $[\bm{u}] \cdot D = y_j$. 
See Figure \ref{subfig:2l1tear} for the picture in the one-dimensional case.
\end{remark2}

Combining Lemma \ref{lemma2:lprimedef} with Corollary \ref{corollary2:llprimesimilar} completes the proof of Proposition \ref{proposition2:2fptout}. 

\begin{corollary2}
\label{corollary2:atypea}
Let $D^+$ be an open neighbourhood of $D \subset M$ which makes $(M,D^+)$ into a K\"{a}hler pair, such that the image of $L^n$ lies in the interior of $M \setminus D^+$ (and hence is an object of $\mathcal{F}(M,D^+,\bm{a})$).
Let $\phi: (M^n_n,D^+) \To (M,D^+)$ be the branched cover of K\"{a}hler$^+$ pairs of Example \ref{example2:fermatcover}. 
Denote
\[ \mathscr{A} := CF^*_{\mathcal{F}(\phi)}(L^n,L^n),\]
where $\mathcal{F}(\phi)$ is the category defined in Proposition \ref{proposition2:2branchedfirst}.
Then $\mathscr{A}$ satisfies all of the conditions required to be of type A, in the sense of Definition \ref{definition2:2typea}, except it may not be strictly $H$-equivariant.
\end{corollary2}
\begin{proof}
By Proposition \ref{proposition2:2branchedfirst}, there is a quasi-isomorphism of $\C$-linear $A_\infty$ algebras 
\[ \mathcal{A} := CF^*_{\mathcal{F}(\phi)}(L^n,L^n)/\mathfrak{m} \cong CF^*_{\mathcal{F}(M,D^+,\bm{a})}(L^n,L^n)/\mathfrak{m},\]
It follows from Proposition \ref{proposition2:2fptout} that the underlying $\bm{G}$-graded vector space of $\mathscr{A}$ is $A$, and the product $\mu^2_0$ is the exterior product. 
Furthermore, if $\mathscr{A} = (A,\mu^*)$, it follows that
\[ \Phi(\mu^*) = \pm u_1 \ldots u_n +  \mathfrak{m}.\]

We now recall from Lemma \ref{lemma2:ssconv} that the spectral sequence induced by the length filtration on $CC^*(\mathcal{A})$ has  $E_2$ page
\[ E_2^{*,*} \cong HH^*_c(A),\]
and converges to $HH^*(\mathcal{A})$.
It follows from Proposition \ref{proposition2:2fptout} that the first-order deformation class of
\[ CF^*_{\mathcal{F}(M,D^+)}(L^n,L^n)\]
is given by
\[ [\mu_1] = \sum_{j=1}^{n} \pm r_j u_j + \mbox{(higher-order in length filtration)}.\]
It follows from Theorem \ref{theorem2:2defclassram}, and the fact that the spectral sequence induced by the length filtration respects the multiplication given by the Yoneda product, that the first-order deformation class of
\[ CF^*_{\mathcal{F}(M,D^+,\bm{a})}(L^n,L^n)\]
 is given by 
\[ [\mu_{1,\bm{a}}] = \sum_{j=1}^{n} \pm  r_j u_j^{n} + \mbox{(higher-order in length filtration)}.\]
By Proposition \ref{proposition2:2branchedfirst}, it follows that
\[ \Phi(\mu^*) = \pm u_1 \ldots u_{n} + \sum_{j=1}^{n} \pm r_j u_j^{n} + \mathfrak{m}^2,\]
from which the result follows.
 \end{proof}

We now have to deal with the fact that $\mathscr{A}$ may not be strictly $H$-equivariant. 
It is clear that $H$ acts on $(M,D)$, preserving the anchored Lagrangian brane $L^n$, because all of our constructions of $(M,D)$ and $L^n$ have been symmetric with respect to permuting the coordinates. 
However, it may not be possible to choose our perturbation data $H$-equivariantly and still achieve transversality, so we may only have $H$-equivariance `up to homotopy'.
We can fix this using the arguments of Appendix \ref{sec:2strict} (based on \cite[Section 8b]{Seidel2003}), which says that we can replace $\mathscr{A}$ by a quasi-equivalent algebra which is strictly $H$-equivariant: essentially, we just apply the proof of Theorem \ref{theorem2:2typea} to this strictly $H$-equivariant replacement.

\begin{corollary2}
\label{corollary2:atypeaworks}
Suppose that $\mathscr{B}$ is a $\bm{G}$-graded $A_{\infty}$ algebra over $R$ of type A. 
Then there exists $\psi\llbracket T \rrbracket \in \mathrm{Aut}(R)$, and an $A_{\infty}$ quasi-isomorphism
\[ \mathscr{A} \cong \psi \cdot \mathscr{B}.\]
\end{corollary2}
\begin{proof}
First, observe that Corollary \ref{corollary2:aunique} does not require strict $H$-equi\-var\-i\-ance, so there is a quasi-isomorphism
\[ \mathscr{A}/\mathfrak{m} \cong \mathscr{B}/\mathfrak{m} =: \mathcal{A}.\]

Now, by a version of Proposition \ref{proposition2:2streq}, the subcategory of $\mathcal{F}(\phi)$ with object $L^n$ embeds, $H$-equivariantly, into a strictly $H$-equivariant $A_{\infty}$ category, in such a way that the order-$0$ component of the embedding is a quasi-equivalence. 
Now recall that $\mathscr{A}$ is necessarily minimal (by Lemma \ref{lemma2:mindef}), and it follows easily from the proof of Proposition \ref{proposition2:2streq} that this strictly $H$-equivariant $A_{\infty}$ category can be chosen to be minimal too. 
Any $A_{\infty}$ functor between minimal $A_{\infty}$ categories, whose order-$0$ component is a quasi-equivalence, is necessarily a quasi-equivalence by Lemma \ref{lemma2:quasiorder}. 
Using the fact that quasi-equivalences of minimal $A_{\infty}$ categories can be inverted (Lemma \ref{lemma2:mininv}), we can apply \cite[Lemma 4.3]{Seidel2003} to prove that $\mathscr{A}$ is quasi-equivalent to a strictly $H$-equivariant $A_{\infty}$ algebra of type A.  
The result now follows from Theorem \ref{theorem2:2typea}.
 \end{proof}

\begin{definition2}
\label{definition2:2atild}
We denote by
\[ \widetilde{\mathscr{A}} \subset \mathcal{F}(M^n_n,D)\]
the full $\bm{G}^n_n$-graded subcategory whose objects are the lifts of $L^n$. 
\end{definition2}

\begin{corollary2}
\label{corollary2:btildatild}
If $\mathscr{B}$ is any $A_{\infty}$ algebra of type A, then there exists $\psi\llbracket T \rrbracket  \in \mathrm{Aut}(R)$ and an $A_{\infty}$ quasi-isomorphism of $\bm{G}^n_n$-graded $R$-linear $A_{\infty}$ categories,
\[  \widetilde{\mathscr{A}} \cong \psi \cdot \bm{p}_1^* \underline{\mathscr{B}}.\]
\end{corollary2}
\begin{proof}
Consider the branched cover
\[ \phi: (M^n_n,D^+) \To (M^n_1,D^+)\]
of Example \ref{example2:fermatcover}.
By Proposition \ref{proposition2:2branchedfirst}, there is a fully faithful embedding
\[  \bm{p}_1^* \mathcal{F}(\phi) \To \mathcal{F}(M^n_n,D^+) \To \mathcal{F}(M,D),\]
and in particular we have a quasi-equivalence
\[ \bm{p}_1^* \underline{\mathscr{A}} \cong \widetilde{\mathscr{A}}.\] 
The result now follows from Corollary \ref{corollary2:atypeaworks}.
 \end{proof}

\section{The $B$-model}
\label{sec:2bmodel}

The aim of this section is to prove Theorem \ref{theorem2:2main2}.

\subsection{Homological perturbation lemma}
\label{subsec:2hpl}

We will use a version of the homological perturbation lemma which is not quite the usual one (for which see, for example, \cite{Markl2006}, \cite{Markl2001}), but rather the slightly modified version used in \cite{Seidel2008a}, so we feel it is as well to state it. 

Suppose we are given:
\begin{itemize}
\item An $A_{\infty}$ algebra $(B,\mu^*)$ (over a $\C$-algebra $R$);
\item A map $\partial: B \To B$ that is a Maurer-Cartan element for $(B,\mu^*)$, in the sense that $\tilde{\mu}^* := (\mu^1+\partial,\mu^2,\ldots)$ is an $A_{\infty}$ structure on $B$;
\item A chain complex $(C,d_C)$;
\item Chain maps
\[(C,d_C) \overset{i}{\underset{p}{\rightleftarrows}} (B,\mu^1);\]
\item A map
\[h: B  \To B,\]
\end{itemize}
such that 
\begin{itemize}
\item $pi = \mathrm{id}$;
\item $h$ defines a homotopy between $ip$ and the identity, which just means that
\[ ip = \mathrm{id} - [\mu^1,h];\]
\item the {\bf side conditions}
\begin{eqnarray*}
h^2 &=& 0, \\
hi &=& 0 \mbox{, and}\\
ph &=& 0.
\end{eqnarray*}
are satisfied;
\item there exists some integer $N$ such that $(\partial h)^N = 0$.
\end{itemize}

Then we construct:
\begin{itemize}
\item An $A_{\infty}$ structure $\nu^*$ on $C$;
\item An $A_{\infty}$ morphism $I^*$ from $(B,\tilde{\mu}^*)$ to $(C,\nu^*)$;
\item An $A_{\infty}$ morphism $P^*$ from $(C,\nu^*)$ to $(B, \tilde{\mu}^*)$,
\end{itemize}
such that $I^*$ and $P^*$ are mutually inverse $A_{\infty}$ quasi-isomorphisms. 
In fact, we can show that $P^1 \circ I^1 = \mathrm{id}$, and we can construct an $A_{\infty}$ homotopy $H^*$ such that
\[ I^* \circ P^* = \mathrm{id} - [\tilde{\mu}^*,H^*].\]

This result is proved in the case $\partial = 0$ in \cite{Markl2006}. 
The operations $\nu^*,I^*,P^*,H^*$ are defined by certain counts over stable directed planar trees (we use the opposite orientation convention from Definition \ref{definition2:2plantree}, so trees have $s$ incoming edges and a single outgoing edge). 
We attach the operation $\mu^k$ to each vertex of arity $k$ (arity $=$ number of incoming edges).
If $\partial \neq 0$, then we make exactly the same construction, but sum instead over {\bf semistable} directed planar trees, and attach the operation $\partial$ to each vertex of arity $1$. 
The assumption that $(\partial h)^N = 0$ ensures that we need only sum over a finite number of trees, because a tree with a sufficiently long chain of vertices of arity $1$ does not contribute to the sum.

For example, to define $\nu^s$, we sum over semistable directed planar trees with $s$ incoming edges. 
We attach operations to each vertex and edge of such a tree, as follows (omitting signs):
\begin{itemize}
\item to each vertex of arity $1$, attach $\partial$;
\item to each vertex of arity $k \ge 2$, attach $\mu^k$;
\item to each internal edge, attach $h$;
\item to each incoming edge, attach $i$;
\item to each outgoing edge, attach $p$.
\end{itemize}
Composing the operations as prescribed by the tree  determines a map $C^{\otimes s} \To C$.
Summing these maps, over all such trees, defines $\nu^s$.

The modifications in the definitions of $I^*,P^*,H^*$, and the proofs that $\nu^*$ is an $A_{\infty}$ structure, that $I^*$ and $P^*$ are $A_{\infty}$ morphisms, and that $H^*$ defines an $A_{\infty}$ homotopy from $I^* \circ P^*$ to $\mathrm{id}$, should all be clear from \cite{Markl2006}. 
The fact that $P^1 \circ I^1 = \mathrm{id}$ follows easily from the side conditions.

\subsection{Matrix factorization computations}
\label{subsec:2mfdefclass}

Let $k$ be a field of characteristic $0$, and $R$ a commutative $k$-algebra. 
Consider the polynomial $R$-algebra $S := R [x_1,\ldots,x_n]$. 
Suppose we are given $w \in S$.
We consider the differential $\Z_2$-graded category of matrix factorizations $MF(S,w)$. 
Objects are finitely-generated free $\Z_2$-graded $S$-modules $K$, equipped with a `differential' $\delta_K : K \To K$ of odd degree such that $\delta_K^2 = w \cdot \mathrm{id}$. 
Morphisms are $S$-module homomorphisms, with the standard differential and compositions.

By \cite[Theorem 3.9]{orlov2004}, there is an exact equivalence between $\mathrm{Ho}(MF(S,w))$ (where `$\mathrm{Ho}$' denotes the homotopy category) and Orlov's `derived category of singularities' $D^b_{\mathrm{Sing}}(w^{-1}(0))$.  
We consider a matrix factorization $(B,\delta_B)$ which corresponds to the ideal 
\[ \mathcal{O}_0 := ( x_1, \ldots x_n )\]
under this equivalence. 
Following the method described in \cite[Section 2.3]{Dyckerhoff2009}, we take $B$ to be the free finitely-generated $\Z_2$-graded algebra generated by odd supercommuting variables $\theta_1,\ldots,\theta_n$ (with $S$ in even degree). 
That is,
\[ B := S [ \theta_1, \ldots, \theta_n ].\]
We define the differential on $B$ to be
\[ \delta := \delta_0 + \delta_1,\]
where
\begin{eqnarray*}
\delta_0 &=& \sum_{j=1}^{n} x_j \del{}{\theta_j},\\
\delta_1 &=& \sum_{j=1}^{n} w_j \theta_j,
\end{eqnarray*}
where $w_j \in S$ are elements chosen such that
\[ w = \sum_{j=1}^{n} w_j x_j.\]
Observe that
\begin{eqnarray*}
\delta_0^2 &=& 0, \\
\delta_1^2 &=& 0\mbox{, and} \\
{[}\delta_0,\delta_1{]} &=& w.
\end{eqnarray*}
It follows that $\delta^2 = w \cdot \mathrm{id}$ as required. 

Now consider the differential $\Z_2$-graded algebra
\[ \mathcal{B} := hom^*_{MF(S,w)} (B,B).\]
Again following \cite[Section 2.3]{Dyckerhoff2009}, we take the underlying vector space to be the algebra of differential operators
\[ \mathcal{B} := S \left[ \theta_1, \ldots, \theta_n, \del{}{\theta_1},\ldots,\del{}{\theta_n} \right],\]
acting on $B$ in the obvious way, with the natural multiplication and the differential $d := d_0 + d_1$, where $d_j := [\delta_j, -]$. 
$\mathcal{B}$ is freely generated, as an $R$-module, by generators $x^{\bm{b}} \theta^J \partial^K$, where $\bm{b} = (b_1,\ldots,b_n)$ is a multi-index, $J \subset [n]$, $K \subset [n]$. 
We will use the shorthand $\partial_j$ for $\partial / \partial \theta_j$.

Now we use the homological perturbation lemma to construct a minimal $A_{\infty}$ model for $\mathcal{B}$. 

To put ourselves in the situation of Section \ref{subsec:2hpl}, let us consider $\mathcal{B}$ with the $A_{\infty}$ (in fact, differential graded) structure given by the differential $ \mu^1 = d_0$ and standard multiplication $\mu^2$, and let $\partial := d_1$. 
Then $(\mathcal{B},\tilde{\mu}^*) = \mathcal{B}$.
Furthermore, let
\begin{eqnarray*}
C &:= & R [ \partial_1 , \ldots, \partial_n ],\\
d_C &=& 0,\\
i:(C,d_C) & \To & (\mathcal{B},d_0) \mbox{ the obvious inclusion,} \\
p:(\mathcal{B},d_0) & \To & (C,d_C) \mbox{ the projection defined by} \\
p(x^{\bm{b}} \theta^J \partial^K) &=& \left\{ \begin{array}{rl}
							\partial^K & \mbox{if $\bm{b} = 0$ and $J = \emptyset$,} \\
							0 & \mbox{otherwise,}
							\end{array} \right.
\end{eqnarray*}
so that $pi = \mathrm{id}$.
We define 
\begin{eqnarray*}
\tilde{h} : B &\To &B,\\
\tilde{h}(f \theta^J \partial^K) &:=& \left(\sum_{j=1}^n \del{f}{x_j} \theta_j \right)\theta^J \partial^K.
\end{eqnarray*}
One can check that
\[ [d_0,\tilde{h}]\left(x^{\bm{b}} \theta^J \partial^K\right) = (|\bm{b}| + |J|)\left(x^{\bm{b}} \theta^J \partial^K\right).\]
Therefore, if we define
\[ h\left(x^{\bm{b}} \theta^J \partial^K\right) :=  \left\{ \begin{array}{rl}
							0 & \mbox{if $\bm{b} = 0$ and $J = \emptyset$,} \\
							\frac{1}{|\bm{b}| + |J|}\tilde{h}\left(x^{\bm{b}} \theta^J \partial^K\right) & \mbox{otherwise,}
							\end{array} \right.
\]
then we have
\[ ip = \mathrm{id} - [d_0,h].\]
Furthermore, we can check that the side conditions are satisfied:
\begin{eqnarray*}
h^2 &=& 0 \mbox{ (for the same reason the exterior derivative squares to $0$),} \\
hi &=& 0 \mbox{, and}\\
ph &=& 0.
\end{eqnarray*}
Finally, observe that $h \partial$ decreases the grading $|x^{\bm{b}} \theta^J \partial^K| := |K|$ by $1$, so $(h \partial)^{n+1} = 0$.

Thus, we can apply Section \ref{subsec:2hpl} to construct an $A_{\infty}$ structure $\nu^*$ on $C$, which is quasi-isomorphic to $\mathcal{B}$. 
If $w$ has degree $\ge 3$, then the differential $\nu^1 = 0$ and the product $\nu^2$ is the standard (exterior algebra) product on $C$. 
Thus, $\nu^{\ge 3}$ defines a Maurer-Cartan element in the Hochschild cochain complex $CC^*(C)$.

Now we recall, from Definition \ref{definition2:2phihkr}, the Hochschild-Kostant-Rosenberg map from the Hochschild cohain complex to the space of polyvector fields,
\[ \Phi: CC^*(C) \To S [ \partial_1, \ldots, \partial_{n} ],\]
given by
\[ \Phi(\alpha) := \sum_{s \ge 0} \alpha^s(\bm{x}, \ldots, \bm{x}),\]
where we denote
\[ \bm{x} := \sum_{j=1}^{n} x_j \partial_j.\]

\begin{proposition2}
\label{proposition2:2mfdef} 
(see \cite[Proposition 7.1]{Efimov2009})
The image of the Maurer-Car\-tan element $\nu^{\ge 3}$ under $\Phi$ is exactly the superpotential $w$.
\end{proposition2}
\begin{proof}
We recall the construction of the maps $\nu^k$ from Section \ref{subsec:2hpl}, by summing over trees.
The only trees that give a non-zero contribution to a product $\nu^k(\partial_{i_1}, \ldots, \partial_{i_k})$ are those depicted in Figure \ref{fig:2perttree}.

\begin{figure}
\centering
\includegraphics[width=0.7\textwidth]{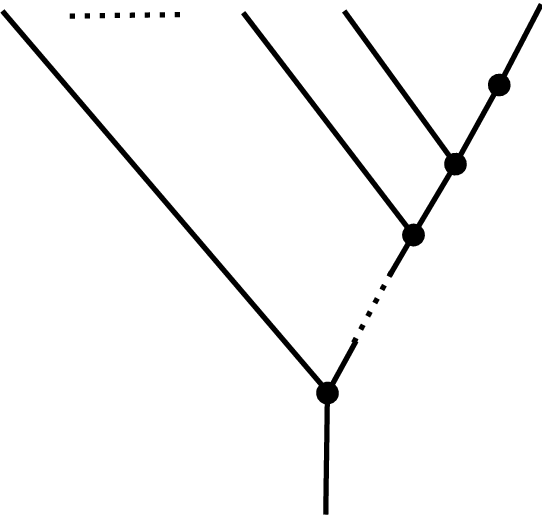}
\caption{The only trees contributing to the homological perturbation lemma computation.
\label{fig:2perttree}}
\end{figure}

Given such a tree, with inputs $\partial_{i_1}, \ldots, \partial_{i_k}$, one easily determines the output: it is the constant term of
\[ \partial_{i_k} \ldots \partial_{i_2} w_{i_1},\]
divided by $(k-1)!$ (coming from the terms in the denominator of $h$). 
Given a monomial $x^{\bm{b}}$ in $w_j$, there are $|\bm{b}|!/\bm{b}!$ ways of choosing the order of the inputs $\partial_{i_m}$ we take from the term
\[ \nu^{|\bm{b}| + 1}(\bm{x}, \ldots, \bm{x})\] 
of $\Phi(\nu^{\ge 3})$, in such a way that $j = i_1$ comes first, and the constant term of 
\[\partial_{i_k} \ldots \partial_{i_2} w_j\]
is non-zero.
Thus, each monomial $x^{\bm{b}}$ of $w_j$ contributes a term
\[ x_j\frac{1}{|\bm{b}|!} \frac{|\bm{b}|!}{\bm{b}!} \left(\partial^{\bm{b}} x^{\bm{b}}\right) x^{\bm{b}} = x_j x^{\bm{b}}\]
to $\Phi(\nu^{\ge 3})$, and the result follows.
 \end{proof}

\subsection{$\bm{G}$-Graded matrix factorizations}
\label{subsec:2gradmf}

For the purposes of this section, let $\bm{G}$ be a grading datum:
\[ \Z \overset{f}{\To} Y \To X \To 0,\]
with sign morphism $\bm{\sigma}$. 
 
Let $S$ be a $\bm{G}$-graded algebra such that $\bm{\sigma}_*S$ is concentrated in degree $0 \in \Z_2$, and let $w \in S$ be an element of degree $f(2) \in Y$. 

\begin{definition2}
A {\bf $\bm{G}$-graded matrix factorization} of $w \in S$ is a $\bm{G}$-graded finitely-generated free $S$-module $K$, together with a homomorphism
\[ \delta_K \in \mathrm{Hom}_S(K,K)\]
of degree $f(1) \in Y$, such that
\[ \delta_K^2 = w \cdot \mathrm{id}.\]
\end{definition2}

\begin{definition2}
\label{definition2:ggradmat}
We define the {\bf differential $\bm{G}$-graded category of matrix factorizations}, $MF^{\bm{G}}(S,w)$:
\begin{itemize}
\item Objects are $\bm{G}$-graded matrix factorizations of $w$;
\item Morphisms are $S$-module homomorphisms:
\[ hom^*((K,\delta_K),(L,\delta_L)) := \mathrm{Hom}_S(K,L);\]
\item Differential on morphism spaces is as usual:
\[ \partial(F) :=  \delta_L \circ F - (-1)^{\sigma(F)} F \circ \delta_K;\]
\item Composition is composition of $S$-module homomorphisms.
\end{itemize}
We note that the morphism spaces are naturally $\bm{G}$-graded $S$-modules, the differential and composition maps have degrees $f(1)$ and $f(0) \in Y$ respectively, and they satisfy the Leibniz rule.
It follows that $MF^{\bm{G}}(S,w)$ is a $\bm{G}$-graded $A_{\infty}$ category over $S$ (see Remark \ref{remark2:ainfgrad}). In fact it is a differential $\bm{G}$-graded category, since all $\mu^{\ge 3}$ are zero.
\end{definition2}

We observe that ordinary matrix factorizations are nothing more than $\bm{G}_{\sigma}$-graded matrix factorizations.
It follows that there is a fully faithful embedding
\[ \bm{\sigma}_* MF^{\bm{G}}(S,w) \To MF(S,w).\]

Now let us introduce our main example. 
We will introduce graded matrix factorizations mirror to the smooth orbifold relative Fukaya category (compare Section \ref{subsec:2ln}). 
Let $\bm{G}:= \bm{G}^n_1$ be the grading datum introduced in Example \ref{example2:ses1}.
Let 
\[ R:= R^n_n \]
be the $\bm{G}$-graded power series ring introduced in Definition \ref{definition2:2ringra}.

Let $U$ be the $\bm{G}$-graded vector space of Example \ref{example2:vecu}. 
We consider the $\bm{G}$-graded algebra
\[ S := R[U] \cong R[u_1, \ldots, u_{n}].\]
We consider the element
\[ w = u_1 \ldots u_{n} + \sum_{j=1}^{n} r_j u_j^n \in S,\]
which has degree $f(2) \in Y$, hence we can define the $\bm{G}$-graded category of matrix factorizations $MF^{\bm{G}}(S,w)$.

We consider the $\bm{G}$-graded $S$-module
\[ K := R[U] \otimes \Lambda \left( U^{\vee} \right) \cong R [u_1, \ldots,u_{n}][ \theta_1, \ldots, \theta_{n}],\]
where the variables $\theta_j$ anti-commute. 
We introduce the differential 
\begin{eqnarray*}
\delta_K : K & \To & K \\
\delta_K &=& \sum_{j=1}^{n} u_j \del{}{\theta_j} + w_j \theta_j,
\end{eqnarray*}
where
\[ w_j = \frac{u_1 \ldots u_{n}}{nu_j} + r_j u_j^{n-1}.\]
We observe that $\delta_K$ has degree $f(1) \in Y$, and that
\[ \delta_K^2 = w \cdot \mathrm{id}.\]
Thus, $(K,\delta_K)$ is a $\bm{G}$-graded matrix factorization of $w$. 
We denote it by $\mathcal{O}_0$.
Finally, we observe that $\delta_K$ is $H$-invariant, where $H$ is the symmetric group acting in the obvious way.

\begin{corollary2}
\label{corollary2:btypea}
Let us define the $\bm{G}$-graded $A_{\infty}$ algebra over $R$:
\[ \mathscr{B} := hom^*_{MF^{\bm{G}}(S,w)}(\mathcal{O}_0,\mathcal{O}_0).\]
Then $\mathscr{B}$ is quasi-isomorphic to an $A_\infty$ algebra of type A, in the sense of Definition \ref{definition2:2typea}.
\end{corollary2}
\begin{proof}
Follows from Proposition \ref{proposition2:2mfdef}. 
We observe that it is necessary to check that the $\bm{G}$-grading and $H$-equivariance interact appropriately with the homological perturbation lemma construction, but this is clear.
 \end{proof}

\begin{corollary2}
\label{corollary2:alltypea}
Let
\[ \mathscr{A} := CF^*_{\mathcal{F}(\phi)}(L^n,L^n),\]
where $\mathcal{F}(\phi)$ is the category of Proposition \ref{proposition2:2branchedfirst} and $\phi: (M^n_n,D^+) \To (M,D^+)$ is the branched cover of Example \ref{example2:fermatcover}. 
Let
\[ \mathscr{B} := hom^*_{MF^{\bm{G}}(S,w)}(\mathcal{O}_0,\mathcal{O}_0)\]
as above. 
Both $\mathscr{A}$ and $\mathscr{B}$ are $\bm{G}$-graded $A_{\infty}$ algebras over $R$, and there exists a power series $\psi \in \C\llbracket T \rrbracket $, with $\psi(0) = 1$ (recalling $T = r_1 \ldots r_{n}$), and an $A_{\infty}$ quasi-isomorphism 
\[\mathscr{A} \cong \psi \cdot \mathscr{B}.\]
\end{corollary2}
\begin{proof}
Follows from Corollary \ref{corollary2:btypea} and Corollary \ref{corollary2:atypeaworks}.
 \end{proof}

\subsection{Equivariant matrix factorizations}
\label{subsec:2equivmf}

Suppose that $\bm{p}: \bm{G}' \To \bm{G}$ is an injective morphism of grading data. 
Then we consider the $\bm{G}'$-graded category 
\[ \bm{p}^*MF^{\bm{G}}(S,w).\]

\begin{remark2}
Objects of $\bm{p}^*MF^{\bm{G}}(S,w)$ are again $\bm{G}$-graded matrix factorizations, but the morphism spaces are just the parts whose $Y$-grading lies in the image of $p : Y' \To Y$. 
Thus $\bm{p}^*MF^{\bm{G}}(S,w)$ embeds, fully faithfully, in the category of $\mathrm{coker}(p)^*$-equivariant matrix factorizations.
\end{remark2}

For example, let $\bm{p}_1: \bm{G}^n_n \To \bm{G}^n_1$ be the morphism defined in Lemma \ref{lemma2:squaregrad}. 

\begin{definition2}
\label{definition2:2btild}
We denote 
\[\widetilde{\mathscr{B}} := \bm{p}_1^* \underline{\mathscr{B}} \subset \bm{p}_1^*MF^{\bm{G}}(S,w).\]
It is the full subcategory whose objects are $\mathcal{O}_0$ and its shifts by elements $y \in Y$. 
\end{definition2}

We obtain the following:

\begin{corollary2}
\label{corollary2:relfukeqmat}
There exists a power series $\psi \in \C\llbracket T \rrbracket $, $\psi(0) = 1$, and an $A_{\infty}$ quasi-isomorphism of $\bm{G}^n_n$-graded $R$-linear $A_{\infty}$ categories,
\[ \widetilde{\mathscr{A}} \cong \psi \cdot \widetilde{\mathscr{B}}.\]
\end{corollary2}
\begin{proof}
Follows from Corollary \ref{corollary2:btypea} and Corollary \ref{corollary2:btildatild}.
  \end{proof}

\subsection{Coherent sheaves}
\label{subsec:2coh}
 
In this section, we will explain how to relate equivariant categories of coherent sheaves on a projective variety to equivariant categories of graded matrix factorizations. 
We state a result closely related to \cite[Proposition 1.2.2]{Polishchuk2011} (see also \cite[Section 2]{Caldararu2010}), in the language of $\bm{G}$-graded matrix factorizations. 
First, we recall Orlov's {\bf category of graded matrix factorizations} (see \cite[Section 3.1]{Orlov2009}). 

Let $S$ be a $\Z$-graded ring, and $w \in S$ homogeneous of degree $d$.
Recall from Example \ref{example2:grmfd} the grading datum $\bm{G}_{MF(d)}$. 
We define
\begin{eqnarray*}
 q: \Z & \To & \Z \oplus \Z/(2,-d), \\
q(j) &=& (0,j),
\end{eqnarray*}
so we can equip $S \cong q_* S$ with a $\bm{G}_{MF(d)}$-grading.
Then $w$ has degree $(0,d) \sim (2,0) =f(2)\in Y$, so we can define the category of $\bm{G}_{MF(d)}$-graded matrix factorizations. 
Now there is a unique injective morphism of grading data, $\bm{p}: \bm{G}_{\Z} \To \bm{G}_{MF(d)}$. 
We define
\[ \mathrm{GrMF}(S,w) := \bm{p}^* MF^{\bm{G}_{MF(d)}}(S,w).\]
It is an $S$-linear, DG category.
It is not hard to see that this definition coincides with Orlov's category of graded matrix factorizations of $w$ (actually Orlov defines graded matrix factorizations to be the homotopy category of this DG category). 
Namely, if we denote 
\[ \pi: \Z \oplus \Z \To \Z \oplus \Z/ (2,-d),\]
then given a $\bm{G}_{MF(d)}$-graded $S$-module $K$, $\pi^*K$ is a quasi-periodic complex of graded $S$-modules, as in the usual definition of the category of graded matrix factorizations.

\begin{sloppypar}
We recall the relationship of $\mathrm{GrMF}$ to coherent sheaves. 
Suppose that $k$ is a field, $S = k[u_1, \ldots, u_n]$ is the $\Z$-graded polynomial ring, and $w \in S$ is homogeneous of degree $n$. 
Suppose that the variety
\[ X := \{ w = 0\} \subset \mathbb{P}^{n-1}_k \]
is smooth. 
Then, because $X$ is Calabi-Yau, \cite[Theorem 3.11]{Orlov2009} says that there is an equivalence of triangulated categories,
\[ D^b Coh\left(X\right) \cong \mathrm{Ho}(\mathrm{GrMF}(S,w)).\]
In fact, by the uniqueness of DG enhancements proved in \cite{Lunts2010}, there is a quasi-equivalence of $A_{\infty}$ categories
\[ D^bCoh (X) \cong \mathrm{GrMF}(S,w),\]
where the left-hand side denotes some appropriate DG enhancement of $D^bCoh(X)$ (see Remark \ref{remark2:lunts}).
\end{sloppypar}

\begin{remark2}
\label{remark2:beilinson}
We observe that the matrix factorization $\mathcal{O}_0$ is always $\bm{G}_{MF(n)}$-graded. 
If $\mathcal{O}_0[j]$ denotes the shift of $\mathcal{O}_0$ by the integer $j$ (i.e., by $f(j)$), then we can arrange that the object
\[ \bigoplus_{j=0}^{n-1} \mathcal{O}_0[j] \in Ob(\mathrm{GrMF}(S,w))\]
corresponds, under Orlov's equivalence, to the restriction of the Beilinson exceptional collection
\[ \bigoplus_{j=0}^{n-1} i^* \Omega^j(j),\]
where $i: X \hookrightarrow \mathbb{P}^{n-1}_k$ denotes the inclusion. 
See \cite[Remark 5.20]{Ballard2010}, also \cite[Section IV.A]{Aspinwall2007}.
\end{remark2}

Now let us consider the situation of Section \ref{subsec:2equivmf}. 
Recall the commutative square of grading data of Lemma \ref{lemma2:squaregrad}. 
We have the $\bm{G}^n_1$-graded ring
\[ S := R[u_1, \ldots, u_{n}],\]
and $w \in S$ has degree $f(2) \in Y$. 
Recall that the Novikov field $\Lambda$ is an $R$-algebra, via the map
\begin{eqnarray*}
R & \To & \Lambda,\\
r_j & \mapsto & r.
\end{eqnarray*}
We define
\[S_{nov}  :=  S \otimes_R \Lambda \cong \Lambda[u_1, \ldots, u_n],\]
with the object
\[w_{nov}  :=  w \otimes 1 \in S_{nov}.\]
Equip $S_{nov}$ with the standard $\Z$-grading (where $\Lambda$ is concentrated in degree $0$, and each $u_j$ has degree $1$). 
Consider the variety
\[ \widetilde{N}^n_{nov} := \left\{ w_{nov} = 0 \right\} \subset \mathrm{Proj}(S_{nov}) \cong \mathbb{P}^{n-1}_{\Lambda}.\]
Define the group $\Gamma_n$ to be the kernel of the map
\[ (\Z_n)^n/y_{[n]} \To \Z_n,\]
and equip $\widetilde{N}^n_{nov}$ with the action of the character group $\Gamma_n^*$, acting by multiplying the coordinates $u_j$ by $n$th roots of unity. 
We denote
\[ N^n_{nov} := \widetilde{N}^n_{nov}/\Gamma_n^*.\]
We observe that these definitions coincide with those given in the Introduction.

\begin{lemma2}
\label{lemma2:orlovemb}
There is a fully faithful embedding of $\Lambda$-linear DG categories,
\[\bm{q}_{1*} \bm{p}_1^* MF^{\bm{G}}(S,w) \otimes_R \Lambda \To D^b Coh\left(N^n_{nov}\right).\]
\end{lemma2}
\begin{proof}
First, observe that $\bm{q}_{2*} S$ is a $\bm{G}_{MF(n)}$-graded ring. 
Furthermore, it is easy to check that $\bm{q}_{2*}R$ is concentrated in degree $0$, and $u_i$ has degree $(0,1)$. 
This coincides with the $\bm{G}_{MF(n)}$-grading of $S$ induced by the standard $\Z$-grading.

There is a fully faithful embedding
\[ \bm{q}_{2*} MF^{\bm{G}^n_1}(S,w) \To MF^{\bm{G}_{MF(n)}}(\bm{q}_{2*}S,w),\]
which sends
\[ (K,\delta_K) \mapsto (\bm{q}_{2*}K, \delta_K).\]
Therefore, there is a fully faithful embedding
\[ \bm{p}_2^* \bm{q}_{2*} MF^{\bm{G}}(S,w) \To \mathrm{GrMF}(S,w).\]
Now note that $\bm{p}_2^* \bm{q}_{2*} R \cong R$ is concentrated in degree $0 \in \Z$, so the morphism $R \To \Lambda$ respects the $\Z$-grading, and we obtain a fully faithful embedding 
\[ \left(\bm{p}_2^* \bm{q}_{2*} MF^{\bm{G}}(S,w)\right) \otimes_R \Lambda \To \mathrm{GrMF}(S,w) \otimes_R \Lambda.\]

There is a fully faithful embedding
\[ \mathrm{GrMF}(S,w) \otimes_R \Lambda \To \mathrm{GrMF}(S_{nov},w_{nov}),\]
which sends
\[ K \mapsto K \otimes_R \Lambda\]
on the level of objects (recall that $K$ is by definition a {\bf free} $S$-module). 
By Orlov's theorem \cite[Theorem 3.11]{Orlov2009} together with \cite{Lunts2010}, there is an $A_\infty$ quasi-equivalence of $\Lambda$-linear DG categories,
\[ \mathrm{GrMF}(S_{nov},w_{nov}) \cong D^bCoh \left( \widetilde{N}^n_{nov} \right)\]
(cf. Remark \ref{remark2:lunts}).

It follows that there is an $A_\infty$ quasi-equivalence of $\Lambda$-linear DG categories
\[ \mathrm{GrMF}(S_{nov},w_{nov})^{\Gamma_n^*} \cong D^bCoh^{\Gamma_n^*}\left(\widetilde{N}^n_{nov}\right).\]
Hence, by the argument above, there is a fully faithful embedding
\[ \bm{p}_2^* \bm{q}_{2*} MF^{\bm{G}}(S,w)^{\Gamma_n^*} \otimes_R \Lambda \To D^bCoh^{\Gamma_n^*}\left(\widetilde{N}^n_{nov}\right) \cong D^bCoh(N^n_{nov}).\]

Now we recall Lemma \ref{lemma2:gradsquare}. 
It shows that there is an isomorphism
\[ \bm{p}_2^* \bm{q}_{2*} MF^{\bm{G}}(S,w)^{\Gamma_n^*} \cong \bm{q}_{1*} \bm{p}_{1}^* MF^{\bm{G}}(S,w),\]
where
\[ \Gamma_n:= \mathrm{ker}(q_{2,X})/\mathrm{im}(p_{1,X}).\]
In this case, an examination of Lemma \ref{lemma2:squaregrad} shows that $\Gamma_n$ is indeed the kernel of the map
\[ (\Z_n)^n/y_{[n]} \To \Z_n\]
given by summing the coordinates, and the $\Gamma_n$-gradings of morphism spaces correspond under Orlov's equivalence.
The result follows.
 \end{proof}

\begin{corollary2}
\label{corollary2:mfcoh}
There is an $A_\infty$ quasi-equivalence of $\Lambda$-linear DG categories,
\[ D^{\pi} \left( \bm{q}_{1*} \widetilde{\mathscr{B}}  \otimes_R \Lambda \right) \cong D^bCoh\left(N^n_{nov}\right)\]
(where we use `$D^{\pi}$' to denote the split-closure of the category of twisted complexes of a DG category).
\end{corollary2}
\begin{proof}
It follows from Lemma \ref{lemma2:orlovemb} that there is a fully faithful embedding of $\bm{q}_{1*} \widetilde{\mathscr{B}} \otimes_R \Lambda$ into the right-hand side. 
Furthermore, the image of this embedding split-generates, by combining Remark \ref{remark2:beilinson}, together with (the $\Gamma_n^*$-equivariant version of) \cite[Lemma 5.4]{Seidel2003}, together with (the $\Gamma_n^*$-equivariant version of) Beilinson's generation result \cite{Beilinson1978}. 
Applying \cite[Lemma 2.5]{Seidel2003} shows that there is a quasi-equivalence
\[ D^{\pi} \left( \bm{q}_{1*} \widetilde{\mathscr{B}}  \otimes_R \Lambda \right) \cong D^{\pi}Coh\left(N^n_{nov}\right).\]
Now \cite[Lemma 5.3]{Seidel2003} shows that the superscript `$\pi$' on the right-hand side is unnecessary, and can be replaced by a `$b$', so the proof is complete. 
 \end{proof}

\begin{definition2}
\label{definition2:2atildnov}
We denote 
\[ \widetilde{\mathscr{A}}_{nov} := \bm{q}_{1*} \widetilde{\mathscr{A}}\otimes_R \Lambda.\]
It is a full subcategory of $\bm{q}_{1*} \mathcal{F}(M^n_n,D) \otimes_R \Lambda$.
\end{definition2}

\begin{corollary2}
\label{corollary2:pfmingen}
There exists $\psi \in \C\llbracket r^n \rrbracket $, and a quasi-equivalence of triangulated $\Lambda$-linear $\Z$-graded $A_{\infty}$ categories,
\[ D^{\pi} \left(\widetilde{\mathscr{A}}_{nov} \right) \cong \psi \cdot D^bCoh(N^n_{nov})\]
\end{corollary2}
\begin{proof}
Follows from Corollary \ref{corollary2:mfcoh} and Corollary \ref{corollary2:relfukeqmat}.
 \end{proof}

\section{The full Fukaya category}
\label{sec:2splitgen}

In this Section, we consider the full Fukaya category $\mathcal{F}(M^n)$, where $M^n:=M^n_n$. 
We explain why there is an embedding
\[ \mathcal{F}(M^n,D) \otimes_R \Lambda \To \mathcal{F}(M^n),\]
and prove that the full subcategory 
\[\widetilde{\mathscr{A}} \otimes_R \Lambda \subset \mathcal{F}(M^n)\]
 (see Definition \ref{definition2:2atild}) split-generates, using the criterion of \cite{Abouzaid2012}. 
This allows us to complete the proof of Theorem \ref{theorem2:2main}.

\subsection{Relating the full and relative Fukaya categories}

Let $(M,\omega)$ be a compact symplectic manifold, satisfying $c_1(M) = 0$.
The Fukaya category $\mathcal{F}(M)$ is not defined in full generality in the literature; a definition that is sufficiently general for our purposes will be given in \cite{Abouzaid2012}. 
It depends on the choice of a bulk class $\mathfrak{b} \in H^{even}(M;\Lambda_0)$ (where $\Lambda_0$ is the {\bf universal Novikov ring}: the subring of $\Lambda$ containing all formal sums of non-negative powers of the generator $r$), and a background class $st \in H^2(M;\Z_2)$. 
We choose both of these to be $0$. 
The objects of $\mathcal{F}(M)$ are $(L,b)$, where $L \subset M$ is a graded spin Lagrangian submanifold, and $b \in H^1(L;\Lambda_0)$ is a weak bounding cochain.
Composition maps are defined by counting holomorphic disks.
$\mathcal{F}(M)$ is a $\Z$-graded $A_{\infty}$ category.

We would like to relate $\mathcal{F}(M,D)$ to $\mathcal{F}(M)$. 
First, we observe that there is a canonical morphism of grading data,
\[ \bm{q}: \bm{G}(M,D) \To \bm{G}_{\Z},\]
by Remark \ref{remark2:cygrad}. 
Thus we can define a $\bm{G}_{\Z}$-graded (i.e., $\Z$-graded) $A_{\infty}$ category
\[ \bm{q}_* \mathcal{F}(M,D).\]

$\mathcal{F}(M,D)$ is defined over the coefficient ring $R$, which is a certain completion of the polynomial ring $\C [r_1, \ldots, r_k ]$. 
It is a simple task to show that the coefficient ring $\bm{q}_* R$ has degree $0 \in \Z$. 
Therefore, the ring homomorphism
\begin{eqnarray*}
R & \To & \Lambda, \\
r_j & \mapsto & r^{l_j},
\end{eqnarray*}
respects the $\Z$-grading, because both have degree $0 \in \Z$ (here $l_j$ are the linking numbers, as defined in Definition \ref{definition2:2kahlp}).
This makes $\Lambda$ into an $R$-algebra, and means that we can define the $\Z$-graded, $\Lambda$-linear $A_{\infty}$ category
\[ \bm{q}_*\mathcal{F}(M,D) \otimes_R \Lambda.\]

We remark that, up until this point, we have given a complete definition of the relative Fukaya category $\mathcal{F}(M,D)$ of a K\"{a}hler pair (satisfying assumptions as in Definition \ref{definition2:2kahlp}), using explicit domain-dependent perturbations of the holomorphic curve equation. 
To define the full Fukaya category $\mathcal{F}(M)$ however, we need some virtual perturbation scheme as in \cite{fooo,Abouzaid2012}. 
These two categories should be related as follows:

\begin{assumption}
\label{assumption:relfull}
There is a fully faithful embedding
\[ \bm{q}_*\mathcal{F}(M,D) \otimes_R \Lambda \To \mathcal{F}(M)\]
of $\Z$-graded, $\Lambda$-linear $A_{\infty}$ categories. 
\end{assumption}

To prove this, we would have to relate the two different perturbation schemes (explicit domain-dependent perturbations versus perturbations of Kuranishi structures), which would take us beyond the scope of this paper.
So we leave it as an assumption. 
We do, however, provide the following justification:

\begin{remark2}
\label{remark2:relfull}
We observe that there is an obvious map on the level of unobstructed objects (Lagrangians with $\mu^0 = 0$):
\begin{eqnarray*}
 Ob(\bm{q}_*\mathcal{F}(M,D))_{unob} & \To & Ob(\mathcal{F}(M)),\\
L & \mapsto & (L,0),
\end{eqnarray*}
and anchored Lagrangian branes automatically come with a grading (recall that taking $\bm{q}_*$ of a category involves identifying certain objects; in this case, this exactly means that we identify all anchored Lagrangian branes which have the same grading). 
However, there is no map in the other direction: objects of $\mathcal{F}(M)$ may intersect the divisors $D$.
Suppose now that:
\begin{itemize}
\item $\bm{L}$ is a tuple of exact, transversely-intersecting anchored Lagrangian branes in $M \setminus D$;
\item $\bm{p}$ is an associated set of generators (intersection points) of $\bm{L}$;
\item the moduli space of rigid, boundary-punctured holomorphic disks in $M$ with boundary on $\bm{L}$, asymptotic to $\bm{p}$, is regular.
\end{itemize}
We show that a rigid holomorphic disk in this moduli space contributes the same term to an $A_{\infty}$ structure map $\mu^s$ in $\mathcal{F}(M,D) \otimes_R \Lambda$ and in $\mathcal{F}(M)$. 
We have
\[ \omega(u) = -\alpha(p_0) + \sum_{j=1}^s \alpha(p_j) + \sum_{j=1}^k l_j (u \cdot D_j)\]
by Lemma \ref{lemma2:stokes}, where $\alpha(p)$ denotes the symplectic action of generator $p$. 
In particular, if we define the map
\begin{eqnarray*}
 CF^*_{\mathcal{F}(M,D)\otimes_R \Lambda}(L_0,L_1) &\To& CF^*_{\mathcal{F}(M)}(L_0,L_1), \\
 p & \mapsto & r^{\alpha(p)} p,
\end{eqnarray*}
then the holomorphic disk $u$ contributes the same term $\pm r^{\omega(u)}$ to $\mu^s$ in both categories. 
Note that, in $\mathcal{F}(M,D)$, we have $(u \cdot D)!$ choices for the labelling of the marked points mapping to the divisors $D$, so this disk in fact contributes $(u \cdot D)!$ identical terms to $\mu^s$, each of which is $\pm r^{\omega(u)}/(u \cdot D)!$ (see the definition in Section \ref{subsec:2relfuk}). 
Thus its total contribution is exactly $\pm r^{\omega(u)}$, which is the same as its contribution to $\mathcal{F}(M)$.
\end{remark2}

\subsection{Split-generation}
\label{subsec:splitgens}

We recall the subcategory 
\[ \widetilde{\mathscr{A}}_{nov} \subset \bm{q}_{1*} \mathcal{F}(M^n,D) \otimes_R \Lambda\]
from Definition \ref{definition2:2atildnov}. 
By abuse of notation, we will identify this category with its image under the embedding of Assumption \ref{assumption:relfull}. 
Our aim in this section is to prove that $\widetilde{\mathscr{A}}_{nov}$ split-generates $D^{\pi}\mathcal{F}(M^n)$. 

\begin{definition2}
We summarize the definition the {\bf closed-open string map} from the (small) quantum cohomology ring of $M$ (with coefficients in $\Lambda$) to the Hochschild cohomology of $\mathcal{F}(M)$, following \cite{Abouzaid2012}. 
The map
\[ \mathcal{CO}: QH^*(M) \To HH^*(\mathcal{F}(M)),\]
is defined as follows: let $\alpha \in H^j(M;\C)$ be Poincar\'{e} dual to a smooth cycle $A \subset M$. 
Let $\bm{L}$ be a tuple of objects with associated generators $\bm{p}$.
We consider the moduli space $\mathcal{M}_4(\bm{p},A)$, whose objects consist of pairs $(r,u)$, where $r \in \mathcal{R}_1(\bm{L})$ and $u: S_r \To M$ is a smooth map, such that
\begin{itemize}
\item $u$ satisfies the (perturbed) holomorphic curve equation, with Lagrangian boundary conditions given by the labels $\bm{L}$;
\item $u$ is asymptotic to the generators $\bm{p}$ at the boundary punctures;
\item $u(q) \in A$, where $q \in S_r$ is the internal marked point.
\end{itemize}
Then each rigid disk $u \in \mathcal{M}_4(\bm{p},A)$ contributes a term $\pm r^{\omega(u)}p_0$ to
\[ \mathcal{CO}(\alpha)(p_s, \ldots, p_1). \]
We remark that
\begin{itemize}
\item \begin{sloppypar} When we say `perturbed' holomorphic curve equation, it really means we must define a Kuranishi structure on $\mathcal{M}_4(\bm{p},A)$ and introduce virtual perturbations thereof (see \cite{Abouzaid2012});
\item The map $\mathcal{CO}$ is a homomorphism of $\Z$-graded $\Lambda$-algebras, where the product on $QH^*(M)$ is quantum cup product $*$, and the product on $HH^*(\mathcal{F}(M))$ is the Yoneda product, and the $\Z$-gradings are the standard ones. \end{sloppypar}
\end{itemize}
\end{definition2}

We now aim to apply the following result, which is due to \cite{Abouzaid2012}:

\begin{theorem2} 
\label{theorem2:2afooo} 
If $(M,\omega)$ is a compact $2d$-dimensional Calabi-Yau symplectic manifold, $\mathscr{L}$ a full subcategory of $\mathcal{F}(M)$ with some finite set of objects, and the map
\[ \mathcal{CO}^{2d}: QH^{2d}(M) \To HH^{2d}(\mathscr{L})\]
is non-zero, then $\mathscr{L}$ split-generates $\mathcal{F}(M)$. 
\end{theorem2}

We consider the subcategory $\widetilde{\mathscr{A}}_{nov} \subset \mathcal{F}(M^n)$ (actually this is not a finite collection of Lagrangians, because we include all shifts, but it will suffice to choose one representative of each geometric lift of $L^n$ to $M^n$).
We aim to understand the map 
\[\mathcal{CO}^{2(n-2)}:QH^{2(n-2)}(M^n) \To HH^{2(n-2)}\left( \widetilde{\mathscr{A}}_{nov}\right)\]
by first understanding the degree-$2$ part of the map, $\mathcal{CO}^2$, then using the fact that $\mathcal{CO}$ is a $\Lambda$-algebra homomorphism. 
It is expected that the image of the class of the symplectic form, 
\[ \mathcal{CO}([\omega]) \in HH^2(\mathcal{F}(M)),\]
should be the class corresponding to the deformation of $\mathcal{F}(M)$ given by scaling the Novikov parameter $r$. 
In fact, in our relative setting, we can make a statement with a cleaner proof:

\begin{lemma2}
\label{lemma2:ocrel}
Let $(M,D)$ be a K\"{a}hler pair. 
Consider the full subcategory 
\[ \bm{q}_*\mathcal{F}(M,D) \otimes_R \Lambda \subset \mathcal{F}(M)\]
of Assumption \ref{assumption:relfull}. 
Then $\mathcal{CO}(P.D.([D_j]))$ is the image of the class
\[ \left( r_j \del{\mu^*}{r_j} \right) \otimes 1 \in HH^*( \mathcal{F}(M,D)) \otimes_R \Lambda\]
in $HH^*(\mathcal{F}(M,D) \otimes_R \Lambda)$, for any $j$.
\end{lemma2}
To clarify: $\mu^* \in CC^*(\mathcal{F}(M,D))$ is the $A_{\infty}$ structure map, and
\[ r_j \del{\mu^*}{r_j} \in CC^*( \mathcal{F}(M,D))\]
is a Hochschild cochain, as can be seen by applying $r_j \partial/\partial r_j$ to the $A_{\infty}$ associativity equation $\mu^* \circ \mu^* = 0$. 
Thus, it defines a class in $HH^*(\mathcal{F}(M,D))$, and we consider the image of this class under the map
\[ HH^*(\mathcal{F}(M,D)) \otimes_R \Lambda \To HH^*(\mathcal{F}(M,D) \otimes_R \Lambda).\]
\begin{proof}
We ignore issues of transversality of our moduli spaces, as in Remark \ref{remark2:relfull}.
Given Lagrangian branes $\bm{L}$ with associated generators $\bm{p}$, each rigid holomorphic disk $u$ with boundary on $\bm{L}$, asymptotic to $\bm{p}$, contributes a term
\[ \pm r_1^{u \cdot D_1} \ldots r_k^{u \cdot D_k}\]
to the coefficient of $p_0$ in $\mu^s(p_s, \ldots, p_1)$, and hence a term
\[ \pm (u \cdot D_j)  r_1^{u \cdot D_1} \ldots r_k^{u \cdot D_k}\]
to the corresponding coefficient of $r_j \partial \mu^s/ \partial r_j$, and hence a term
\[ \pm (u \cdot D_j) r^{\omega(u)}\]
to the corresponding coefficient of $(r_j \partial \mu^s/ \partial r_j) \otimes_R 1$ (see Remark \ref{remark2:relfull}). 

So, by definition of the map $\mathcal{CO}$, each such holomorphic disk $u$ together with an internal marked point $q$ mapping to $D_j$, contributes a term $\pm r^{\omega(u)}$ to the corresponding coefficient of $\mathcal{CO}(P.D.([D_j]))$. 
There are $u \cdot D_j$ choices for the internal marked point $q$, so the total contribution of each such holomorphic disk $u$ is $\pm (u \cdot D_j) r^{\omega(u)}$. 

Therefore, 
\[ \mathcal{CO}(P.D.([D_j])) = \left( r_j \del{\mu^*}{r_j} \right) \otimes_R 1 \]
as required.
 \end{proof}

\begin{remark2}
\label{remark2:ocrel}
Again, Lemma \ref{lemma2:ocrel} should perhaps be thought of as an assumption, for the same reason that we make Assumption \ref{assumption:relfull}.
\end{remark2}

\begin{proposition2}
\label{proposition2:2ocmn}
Let $M^n$ be the Calabi-Yau Fermat hypersurface of Example \ref{example2:fermat}, and recall the full subcategory $\widetilde{\mathscr{A}}_{nov} \subset \mathcal{F}(M^n)$.
The map
\[ \mathcal{CO}^{2(n-2)}: QH^{2(n-2)}(M^n) \To HH^{2(n-2)}\left(\widetilde{\mathscr{A}}_{nov}\right)\]
is non-zero.
\end{proposition2}
\begin{proof}
\begin{sloppypar}
We recall from Lemma \ref{lemma2:hhtypea} that there is an action of $\tilde{\Gamma}_n^*$ on $HH^*(\widetilde{\mathscr{A}}_{nov})$, and the $\tilde{\Gamma}_n^*$-invariant part is
\[ HH^*_{\bm{G}_{\Z}} \left( \widetilde{\mathscr{A}}_{nov} \right) ^{\tilde{\Gamma}_n^*} \cong \Lambda [\alpha]/\alpha^{n-1}\]
as a $\Z$-graded $\Lambda$-algebra, where $\alpha$ has degree $2$.

It follows from Lemma \ref{lemma2:ocrel} and Lemma \ref{lemma2:defom} that the image of $\mathcal{CO}(P.D.([D_j]))$ under this isomorphism is $g \cdot \alpha$, for some invertible $g \in \Lambda^*$.  
Therefore, because $\mathcal{CO}$ is a $\Lambda$-algebra homomorphism, the image of $\mathcal{CO}(P.D.([D_j])^{n-2})$ under this isomorphism is $g^{n-2} \alpha^{n-2}$, which does not vanish in $\Lambda[\alpha]/\alpha^{n-1}$. 
Because it has degree $2(n-2)$, this completes the proof.
\end{sloppypar}
 \end{proof}

\begin{corollary2}
\label{corollary2:splgen}
The full subcategory 
\[\widetilde{\mathscr{A}}_{nov} \subset \mathcal{F}(M^n)\]
split-generates the Fukaya category.
\end{corollary2}
\begin{proof}
Follows from Proposition \ref{proposition2:2ocmn} and Theorem \ref{theorem2:2afooo}.
 \end{proof}

Theorem \ref{theorem2:2main} now follows from Corollary \ref{corollary2:pfmingen}, Assumption \ref{assumption:relfull}, and Corollary \ref{corollary2:splgen}.

\appendix

\section{Strict group actions on Fukaya categories}
\label{sec:2strict}

This section essentially reproduces the argument of \cite[Section 8b]{Seidel2003}, in the context of the relative Fukaya category.
Suppose that we have a finite group $\Gamma$, which acts on a K\"{a}hler pair $(M,D)$, permuting the divisors $D$, and preserving the Liouville one-form $\alpha$. 
We consider the relative Fukaya category $\mathcal{F}(M,D)$ defined in Section \ref{subsec:2relfuk}. 
There is an obvious action of $\Gamma$ on the objects of the relative Fukaya category.
We would like to say that this action extends to an action on the relative Fukaya category, in a suitable sense. 

Recall the notion of a strictly $\Gamma$-equivariant $A_{\infty}$ structure from Definition \ref{definition2:2stricteq}.
Na\-\"{i}vely, one might try to argue that $\Gamma$ acts on the moduli spaces of pseudo-holomorphic disks used to define the structure maps of the Fukaya category, and hence the structure maps are strictly $\Gamma$-equivariant. 
However, this does not work: for $\Gamma$ to act on the moduli spaces of pseudo-holomorphic disks, we would have to make a $\Gamma$-equivariant choice of perturbation data, which would destroy our chances of achieving transversality. 
Instead, we have the following:

\begin{proposition2}
\label{proposition2:2streq}
In the situation described above, there is a fully faithful embedding of $\mathcal{F}(M,D)$ into a strictly $\Gamma$-equivariant $A_{\infty}$ category. 
The order-$0$ part of this embedding is a quasi-equivalence, and respects the action of $\Gamma$ on objects, in the sense that $F(\gamma \cdot L)$ is quasi-isomorphic to $\gamma \cdot F(L)$.
\end{proposition2}
\begin{proof}
We consider a category $\widetilde{\mathcal{F}}(M,D)$ with objects $(\gamma,L)$, where $L$ is an object of $\mathcal{F}(M,D)$ and $\gamma \in \Gamma$. 
Think of $(\gamma,L)$ as representing the object $\gamma \cdot L$ of $\mathcal{F}(M,D)$, but we now have $|\Gamma|$ copies of each object.
 
We define an action of $\Gamma$ on these objects, via
\[ \gamma_1 \cdot (\gamma_2,L) :=  (\gamma_1 \cdot \gamma_2,L).\]

Now for each pair of objects $((1,L_0),(\gamma,L_1))$ of $\widetilde{\mathcal{F}}(M,D)$, we choose a regular Floer datum for the objects $(L_0,\gamma \cdot L_1)$ of $\mathcal{F}(M,D)$. 
We then define Floer data for pairs of objects $((\gamma_0,L_0),(\gamma_1,L_1))$ by acting with $\gamma_0$ on the Floer data for $((1,L_0),(\gamma_0^{-1}\cdot \gamma_1,L_1))$. 
We thus define morphism spaces
\[ CF^*((\gamma_0,L_0),(\gamma_1,L_1))\]
for all pairs of objects. 
We define the Floer differential $\mu^1_0$ as before, and note that it is now strictly $\Gamma$-equivariant.

Now for each tuple of objects $\bm{L} := ((1,L_0),(\gamma_1,L_1),\ldots,(\gamma_k,L_k))$, with associated generators $\bm{y}$, we choose regular, consistent perturbation data on the moduli spaces $\mathcal{R}(\bm{p},\ell)$. 
We then define perturbation data for tuples $\bm{L} = ((\gamma_0,L_0),\ldots,(\gamma_k,L_k))$ by acting with $\gamma_0$ on the perturbation data chosen for $((1,L_0),(\gamma_0^{-1}\cdot \gamma_1,L_1),\ldots,(\gamma_0^{-1}\cdot \gamma_k,L_k))$. 
This allows us to define the rest of the Floer products $\mu^k$. 
Note that they are strictly $\Gamma$-equivariant.

Observe now that the full subcategory with objects $(1,L)$ is equivalent to $\mathcal{F}(M,D)$ (making the corresponding choice of perturbation data). 
Thus, we have an inclusion of $\mathcal{F}(M,D)$ as a full subcategory of $\widetilde{\mathcal{F}}(M,D)$. 
Furthermore, if we restrict to the affine Fukaya category, then this inclusion is a quasi-equivalence, because each object $(\gamma,L)$ is quasi-isomorphic to the element $(1,\gamma \cdot L)$ of the subcategory.
This concludes the proof.
 \end{proof}

\section{List of notation}

\begin{longtable}{rp{3.8in}}
$[k]$ & The set $\{1,2,\ldots,k\}$ \\
$\Z_k$ & The abelian group $\Z/(k)$ \\
$(M,D)$ & K\"{a}hler pair, $M$ is a K\"{a}hler manifold, $D \subset M$ a normal-crossings divisor (Definition \ref{definition2:2kahlp}). In Section \ref{sec:2mb}, this notation is reserved for the specific K\"{a}hler pair $(M^n_1,D)$.\\
$D_K$ & Intersection of all divisors $D_j \subset M$ such that $j \in K$.\\
$(M,D^+)$ & K\"{a}hler$^+$ pair, a K\"{a}hler pair with a neighbourhood $D^+$ of $D$ (Definition \ref{definition2:kahlerplus}). \\
$(M^n_a,D)$ & The K\"{a}hler pair where $M^n_a \subset \CP{n-2}$ is the degree-$a$ Fermat hypersurface, $D \subset M^n_a$ is the union of coordinate divisors $\{z_j = 0\}$ (Example \ref{example2:fermat}). \\
$\phi_a$ & Branched cover of K\"{a}hler pairs $\phi_a: (M^n_a,D^+) \To (M^n_1,D^+)$ (Example \ref{example2:fermatcover}).\\
$M^n$ & Same as $M^n_n$ (notation used in Sections \ref{sec:intro} and \ref{sec:2splitgen}). \\
$\bm{G}$ & Grading datum: an abelian group $Y$ with a morphism $\Z \To Y$ (Definition \ref{definition2:gradingdat}). In Section \ref{sec:2mb}, this notation is reserved for the specific grading datum $\bm{G}^n_1$. \\
$\bm{H}$ & Pseudo-grading datum (Definition \ref{definition2:2pseudograd}), which induces a grading datum $\bm{G}(\bm{H})$ (Definition \ref{definition2:gradfromps}). \\
$\bm{G}(M)$ & Grading datum coming from Lagrangian Grassmannian of $M$ (Definition \ref{definition2:2hm}). \\
$\bm{G}(M,D)$ & Same as $\bm{G}(M \setminus D)$ (Definition \ref{definition2:gmd}).\\
$\bm{H}(M,D)$ & Pseudo-grading datum coming from $(M,D)$ (Definition \ref{definition2:2gradm}), which comes with a map $\bm{G}(\bm{H}(M,D)) \To \bm{G}(M,D)$ (Definition \ref{definition2:pmorph}) which is an isomorphism under certain conditions (Lemma \ref{lemma2:psgrel}). \\
$\bm{G}^n_a$ & Same as $\bm{G}(M^n_a,D)$ (Example \ref{example2:ses1} and Lemma \ref{lemma2:fermatgrad}). \\
$R(M,D,\bm{a})$ & Coefficient ring of the orbifold relative Fukaya category $\mathcal{F}(M,D,\bm{a})$ (Definition \ref{definition2:2relring}): a $\bm{G}(M,D)$-graded formal power series ring with one generator for each component of $D$. \\
$R(M,D)$ & Coefficient ring of the relative Fukaya category $\mathcal{F}(M,D)$ (Definition \ref{definition2:2relring}). \\
$R^n_a$ & The $\bm{G}^n_1$-graded coefficient ring of $\mathcal{F}(M^n_1,D,(a,\ldots,a))$ (Definition \ref{definition2:2ringra} and Example \ref{example2:fermatring}). \\
$R_a$ & Same as $R^n_a$ (notation only used in Section \ref{subsec:2comp}). \\
$R^j_a$ & Order-$j$ component of $R_a$ (notation only used in Section \ref{subsec:2comp}); apologies for the clash with notation $R^n_a$. \\
$R$ & Notation reserved for the ring $R^n_n$ in Section \ref{subsec:2comp}; also notation reserved for the ring $R^n_1$ in Section \ref{sec:2mb}); apologies for the clash.\\
$\Lambda$ & Universal Novikov field (Definition \ref{definition2:2novikov}).\\
$U_n$ & A certain $n$-dimensional $\bm{G}^n_1$-graded vector space (Example \ref{example2:vecu}).\\
$A_n$ & The $\bm{G}^n_1$-graded exterior algebra of $U_n$ (Definition \ref{definition2:2alga}).\\
$S$ & The $\bm{G}^n_1$-graded algebra $R_n[U_n]$ (Section \ref{subsec:2gradmf}).\\
$w \in S$ & The element $w = u_1 \ldots u_n + \sum r_j u_j^n$ (Section \ref{subsec:2gradmf}). \\
$S_{nov}$ & The $\Z$-graded $\Lambda$-algebra $S \otimes_{R_n} \Lambda$ (Section \ref{subsec:2coh}). \\
$w_{nov}$ & The element $w \otimes 1 \in S_{nov}$ (Section \ref{subsec:2coh}). \\
$\tilde{\Gamma}^*_n$ & The group $(\Z_n)^n/(1,\ldots,1)$, where $\Z_n := \Z/(n)$ (Definition \ref{definition2:2mirror}), which acts on $\mathrm{Proj}(S)$ and $\mathrm{Proj}(S_{nov})$. \\
$\Gamma^*_n$ & The kernel of the sum map $\tilde{\Gamma}^*_n \To \Z_n$ (Definition \ref{definition2:2mirror}). \\
$\widetilde{N}^n$ & The subscheme of $\mathrm{Proj}(S)$ defined by $\{w = 0\}$ (Section \ref{subsec:introb}). \\
$N^n$ & The quotient $\widetilde{N}^n/\Gamma^*_n$ (Section \ref{subsec:introb}). \\
$\widetilde{N}^n_{nov}$ & The subscheme of $\mathrm{Proj}(S_{nov})$ defined by $\{w_{nov} = 0\}$ (Definition \ref{definition2:2mirror}). \\
$N^n_{nov}$ & The quotient $\widetilde{N}^n_{nov} / \Gamma^*_n$ (Definition \ref{definition2:2mirror}). \\
$\mathcal{F}'(M)$ & In Section \ref{sec:2affeqi}, the exact Fukaya category of an exact symplect manifold $M$.\\
$\mathcal{F}'$ & In Section \ref{sec:2mb}, an auxiliary category used to compute the Fukaya category via pearly trees. \\
$\mathcal{F}(M \setminus D)$ & The affine Fukaya category; closely related to $\mathcal{F}'(M \setminus D)$ (Section \ref{subsec:2eqaffuk}). \\
$\mathcal{F}(M,D)$ & The relative Fukaya category (Section \ref{subsec:2relfuk}). \\
$\mathcal{F}(M,D,\bm{a})$ & The orbifold relative Fukaya category: $\bm{a} = (a_1,\ldots,a_k)$, where $a_j$ is the degree of orbifolding about divisor $D_j$ (Section \ref{subsec:2relfuk}). \\
$\mathcal{F}(M)$ & The full Fukaya category (Section \ref{sec:2splitgen}).\\
$L^n$ & An immersed Lagrangian sphere, $L^n: S^{n-2} \To M^n_1 \setminus D$ (Section \ref{subsec:2ln}). \\
$\mathcal{A}$ & In Section \ref{sec:2deftheory}, a $\C$-linear $A_\infty$ algebra; In Section \ref{sec:2mb} and subsequently, the endomorphism algebra of $L^n$ in $\mathcal{F}(M^n_1 \setminus D)$ (Corollary \ref{corollary2:atypea}). \\
$\mathscr{A}$ & In Section \ref{sec:2deftheory}, a deformation of a $\C$-linear $A_\infty$ algebra to one over a formal power series ring; In Section \ref{sec:2mb} and subsequently, the endomorphism algebra of $L^n$ in $\mathcal{F}(M^n_1, D,(n,\ldots,n))$ (Corollary \ref{corollary2:atypea}).\\
$\widetilde{\mathscr{A}}$ & The full subcategory of $\mathcal{F}(M^n_n,D)$ whose objects are lifts of $L^n$ under the branched cover $\phi_n: M^n_n \To M^n_1$ (Definition \ref{definition2:2atild}). \\
$\widetilde{\mathscr{A}}_{nov}$ & The full subcategory of $\mathcal{F}(M^n)$ whose objects are lifts of $L^n$ (Section \ref{subsec:splitgens}). \\
$MF^{\bm{G}}(S,w)$ & Differential $\bm{G}$-graded category of matrix factorizations of $w \in S$ (Definition \ref{definition2:ggradmat}). \\
$\mathcal{O}_0$ & The matrix factorization corresponding to the skyscraper sheaf at the origin (Section \ref{subsec:2mfdefclass}). \\
$\mathscr{B}$ & The endomorphism algebra of $\mathcal{O}_0$ in $MF^{\bm{G}}(S,w)$ \\
$\widetilde{\mathscr{B}}$ & Full subcategory of $\bm{p}_1^*MF^{\bm{G}}(S,w)$ (`$\tilde{\Gamma}^*_n$-equivariant matrix factorizations') whose objects are equivariant lifts of $\mathcal{O}_0$. \\
$\widetilde{\mathscr{B}}_{nov}$ & Full subcategory of (a DG enhancement of) $D^bCoh(N^n_{nov})$ whose objects are $\Gamma^*_n$-equivariant twists of the Beilinson exceptional collection $\Omega^j(j)$ for $j = 0, \ldots, n-1$.
\end{longtable}

\bibliographystyle{plain}
\bibliography{library2}

\end{document}